\documentclass[reqno,twoside,openright,12pt]{amsbook}
\usepackage{amssymb,latexsym,epsfig,verbatim,xypic,amsmath,mathrsfs}
\usepackage[french]{babel}
\usepackage[T1]{fontenc}
\usepackage{indentfirst}

\def\C{\mathbf{C}}

\def\vs{\vspace*{4mm}}

\def\P{\mathbf{P}}

\theoremstyle{plain}
\newtheorem{thm}{Th\'eor\`eme}[chapter]
\newtheorem{pro}[thm]{Proposition}
\newtheorem{cor}[thm]{Corollaire}
\newtheorem{lem}[thm]{Lemme}

\theoremstyle{definition}
\newtheorem*{defi}{D\'efinition}

\newtheorem*{rem}{Remarque}
\newtheorem*{rems}{Remarques}
\newtheorem*{defis}{D\'efinitions}

\begin{document}

\thispagestyle{empty}

\vspace{6mm}

\begin{center}
{\large \textbf{FEUILLETAGES ET ACTIONS DE GROUPES}}
\end{center}
\vspace{4mm}

\begin{center}
{\large \textbf{SUR LES ESPACES PROJECTIFS} }
\vspace{10mm}
\end{center}

\begin{center}
JULIE DESERTI  ET DOMINIQUE CERVEAU
\end{center}

\vspace*{20mm}

ABSTRACT. A holomorphic foliation $\mathscr{F}$ on 
a compact complex
manifold $M$ is said to be an $\mathscr{L}$-foliation if there
exists an action of a complex Lie group $G$ such that the generic
leaf of $\mathscr{F}$ coincides with the generic orbit of $G$. We
study $\mathscr{L}$-foliations of codimension one, in particular
in projective space, in the spirit of classical invariant theory,
but here the invariants are sometimes transcendantal ones. We give
a bestiary of examples and general properties. Some
classification results are obtained in low dimensions.

\footnotetext{IRMAR, UMR 6625 du CNRS,
Universit\'e de Rennes I, 35042 Rennes, 
France.}

\tableofcontents
\addcontentsline{toc}{chapter}{Introduction}

\chapter*{Introduction}

\vspace*{8mm}

La th\'eorie classique des invariants, au sens o\`u l'on peut
l'entendre au XIX$^{eme}$ si\`ecle, propose \'etant donn\'ee une
action alg\'ebrique d'un groupe alg\'ebrique $G$ sur une
vari\'et\'e alg\'ebrique compacte $M$ de d\'ecrire le corps des
fonctions rationnelles $R$ sur $M$ invariantes sous l'action de
$G$.

Parmi les pionners de la th\'eorie des invariants, on retiendra en
particulier Sylvester et Hesse, puis Gordan et N\oe ther qui font
une approche que l'on qualifierait aujourd'hui d'effective. Ces
travaux accompagnent le d\'eveloppement de l'alg\`ebre lin\'eaire,
de la th\'eorie des matrices et de la g\'eom\'etrie projective.
Cette approche change radicalement, non sans pol\'emique, avec les
travaux de Hilbert et Hurwitz.

La g\'eom\'etrie classique produit tout un folklore d'exemples.
Nous en rappelons certains, en liaison avec la classification des
objets alg\'ebriques ; l'un des plus populaires est sans doute
l'invariant $j$ des courbes elliptiques que l'on peut voir de
diff\'erentes mani\`eres : action de $\P\mbox{GL}(2,\C)$ sur
$$\C \P(4) \simeq \{ \mbox{quadruplets de points
sur la sph\`ere de Riemann} \}$$ ou action de $\P\mbox{GL}
(3,\C)$ sur $$\C \P(9)\simeq
\{\mbox{courbes de degr\'e $3$ dans } \C \P(2)\}$$ etc. Nous
mentionnerons aussi l'exemple de Gordan-N\oe ther en r\'eponse \`a une
affirmation de Hesse. L'ubiquit\'e de cet exemple est tout \`a fait \'etonnante.

Dans tout le discours qui suit les objets sont d\'efinis sur le
corps $\C$ des nombres complexes. Si Aut $M$ d\'esigne le groupe
des automorphismes de $M$, on dispose d'un morphisme $\varphi
\hspace{1mm} \colon \hspace{1mm} G \to \mbox{Aut }M$ et $R$ est
invariante si $$R(\varphi(g).m)=R(m)$$ pour tout $m \in M$ et $g
\in G$. Si $\chi(M)$ d\'esigne l'espace des champs de vecteurs
holomorphes sur $M$, on sait, lorsque $M$ est compacte, que
$\chi(M)$ est une $\C$-alg\`ebre de Lie de dimension finie qui
s'identifie naturellement \`a $\mathscr{L}ie(\mbox{Aut }M)$. Par
suite, $\varphi$ induit un morphisme d'alg\`ebres de Lie
$$D\varphi \hspace{1mm} \colon \hspace{1mm} \mathscr{L}ie \hspace{1mm}
G \to \mathscr{L}ie (\mbox{Aut }M) \simeq \chi(M)$$ et produit
une sous-alg\`ebre de Lie $\mathscr{L}$ de $\chi(M)$ : l'image de $D\varphi$.
Notons que l'action de $G$ sur $M$ produit un \guillemotleft
\hspace{1mm} feuilletage \guillemotright \hspace{1mm} singulier $\mathscr{F}$ de
la vari\'et\'e $M$. Nous dirons que $\mathscr{F}$ est un
$\mathscr{L}$-feuilletage. Evidemment ce feuilletage peut \^etre trivial,
par exemple lorsque l'action est transitive ou \`a l'inverse lorsque $G$ est un
groupe discret. On s'int\'eressera au cas o\`u ce feuilletage est de codimension
 $1$ dans une situation g\'en\'erale, en particulier lorsque l'action n'est plus
alg\'ebrique, c'est-\`a-dire lorsque les feuilles du feuilletage associ\'e sont
transcendantes. Sur une vari\'et\'e alg\'ebrique $M$ compacte, un
$\mathscr{L}$-feuilletage de codimension $1$ est associ\'e \`a une $1$-forme
ferm\'ee rationnelle $\Omega$. C'est une cons\'equence \`a peu pr\`es directe du
fait que le groupe d'automorphismes de $M$ est alg\'ebrique. Comme nous l'a
signal\'e Serge Cantat, ce r\'esultat persiste si $M$ est compacte
k\"{a}hlerienne. Par contre ce n'est plus vrai dans le cas non k\"{a}hler.
L'exemple le plus standard est le suivant : soit $\Gamma \subset \mbox{SL}(2,
\C)$ un sous-groupe discret cocompact. La vari\'et\'e
$M=\mbox{SL}(2,\C)/\Gamma$ est munie d'une $\mathscr{L}$-feuilletage
correspondant \`a l'action du groupe triangulaire sur $M$. Ce feuilletage n'est
pas d\'efini par une $1$-forme ferm\'ee mais est \guillemotleft \hspace{1mm}
transversalement projectif \guillemotright. Il est probable que ce soit un fait
g\'en\'eral, ceci \'etant confort\'e par un r\'esultat r\'ecent de
\cite{[C-LN-L-P-T]}.

Donnons maintenant quelques propri\'et\'es des
$\mathscr{L}$-feuilletages de codimension $1$ sur les espaces
projectifs $\C\P(n)$ (chapitre 1). Le degr\'e d'un
$\mathscr{L}$-feuilletage sur $\C\P(n)$ est major\'e par $n-1$ ;
lorsqu'il est maximal, on a alors $\dim \mathscr{L}=n-1$ et on
peut profiter pour $n$ petit de la classification des alg\`ebres
de Lie pour donner une description des \guillemotleft \hspace{1mm}
invariants g\'en\'eralis\'es \guillemotright \hspace{1mm} de
$\mathscr{F}$. Par invariant g\'en\'eralis\'e on entend une
primitive de la forme ferm\'ee rationnelle $\Omega$ d\'efinissant
$\mathscr{F}$ ou toute fonction \guillemotleft \hspace{1mm}
\'el\'ementaire \guillemotright \hspace{1mm} int\'egrale
premi\`ere de $\mathscr{F}$. Remarquons que les feuilles d'un
$\mathscr{L}$-feuilletage de codimension $1$ sur $\C\P(n)$ sont
domin\'ees par $\C^{n-1}$ ; il suffit de choisir $n-1$
\'el\'ements g\'en\'eriques $X_1$, $\ldots$, $X_{n-1}$ dans
$\mathscr{L}$ et de consid\'erer l'application
$$(t_1,\ldots,t_{n-1})\mapsto \exp(t_1X_1)\circ\ldots\circ
\exp(t_{n-1}X_{n-1})(z)$$ qui est d'image dense dans la feuille
passant par $z$. Ceci explique pourquoi les formes ferm\'ees
$\Omega$ intervenant ici sont tr\`es sp\'eciales ; par exemple le
compl\'ement des p\^oles de $\Omega$ (les p\^oles sont invariants
par $\mathscr{F}$) ne peut \^etre Kobayashi hyperbolique. Les
feuilles d'un $\mathscr{L}$-feuilletage peuvent \^etre denses ou
d'adh\'erence des vari\'et\'es r\'eelles Levi-plates.

Nous montrons qu'un $\mathscr{L}$-feuilletage sur $\C\P(n)$, $n
\geq 3$, poss\'edant un point singulier isol\'e est n\'ecessairement de degr\'e
$1$ et que dans une carte affine ad-hoc les feuilles sont les niveaux d'une
forme quadratique de rang maximum. C'est une cons\'equence plus ou moins directe
 du th\'eor\`eme de Frobenius singulier de B. Malgrange. Ce r\'esultat devrait
avoir un analogue en codimension sup\'erieure.

Nous d\'ecrivons ensuite quelques propri\'et\'es des
$\mathscr{L}$-feuilletages li\'ees \`a la nature alg\'ebrique de
$\mathscr{L}$. Par exemple si tous les \'el\'ements de
$\mathscr{L}$ sont nilpotents, le feuilletage associ\'e poss\`ede
une int\'egrale premi\`ere rationnelle. Il faut cette hypoth\`ese
forte car m\^eme dans le cas ab\'elien ce r\'esultat ne persiste
pas. Modulo des hypoth\`eses de r\'egularit\'e naturelles et
n\'ecessaires, si $\mathscr{L}$ est r\'esoluble, $\mathscr{F}$
poss\`ede un hyperplan invariant. Enfin si
$[\mathscr{L},\mathscr{L}]=\mathscr{L}$, ce qui est le cas si
$\mathscr{L}$ est semi-simple, alors $\mathscr{F}$ poss\`ede
encore un invariant rationnel.

Le chapitre 2 est consacr\'e \`a la description d'exemples, issus d'actions de
groupes, pour la plupart classiques.

Dans le chapitre 3 nous montrons que les feuilletages de degr\'es $0$ et $1$ sur
$\C\P(n)$ sont des $\mathscr{L}$-feuilletages. Nous donnons la
liste des invariants g\'en\'eralis\'es correspondants.

Le chapitre 4 est consacr\'e \`a la description compl\`ete des
$\mathscr{L}$-feuilletages sur $\C\P(3)$. De cette \'etude on
peut d\'eduire le th\'eor\`eme suivant : 

\vs

\begin{thm} Soit $G$ un groupe de Lie complexe
connexe agissant sur $\C \P(3)$. On est dans l'un
des cas suivants :
\begin{itemize}
\item[(i)] l'action est triviale,

\item[(ii)] il existe une orbite dense,

\item[(iii)] l'orbite g\'en\'erique est de dimension $2$ et
l'action de $G$ poss\`ede un invariant g\'en\'eralis\'e (fonction
constante sur les feuilles, \'eventuellement multivalu\'ee)
appartenant \`a conjugaison pr\`es \`a la liste qui suit :
\begin{eqnarray}
& & z_0^{\lambda_0} z_1^{\lambda_1} z_2^{\lambda_2} \mbox{ o\`u }
\displaystyle \sum_{i=0}^2 \lambda_i=0,\hspace{2mm}
z_0^{\lambda_0}z_1^{\lambda_1} z_2^{\lambda_2}z_3^{\lambda_3}
\mbox{o\`u } \displaystyle \sum_{i=0}^3 \lambda_i=0,\hspace{2mm} \frac{z_0}{z_1}
,\nonumber \\
& & \nonumber \\
& & \frac{z_1}{z_3}\exp\left(\frac{z_0z_3+z_2z_1}{z_1
z_3}\right),\hspace{2mm}
\frac{z_0}{z_3}\exp\left(\frac{z_2^2-z_1z_3}{z_3^2}\right),
\hspace{2mm}
\frac{z_0}{z_3}\exp\left(\frac{z_2^2-z_1z_3}{z_0z_3}\right), \nonumber \\
& & \nonumber \\
& &\frac{z_1z_3^{\kappa-1}}{z_0^{\kappa}}\exp\left(\frac{z_2}
{z_3}\right),\hspace{2mm}\frac{z_0}{z_1}\exp\left(\frac{z_2}{z_1}
\right),\hspace{2mm}\frac{z_3^2}{z_1z_3-z_2^2}
\exp\left(\frac{z_0}{z_3}\right),\nonumber \\
& & \nonumber \\
& &\frac{z_2}{z_3}\exp\left(\frac{z_0z_3-z_1z_2}{z_3^2}\right),
\hspace{2mm} \frac{z_0z_3-z_1z_2}{z_3^2}\exp\left(\frac{z_2}{z_3}
\right), \nonumber \\
& & \nonumber \\
& & \frac{Q}{z_0^2} \mbox{ o\`u $Q$ est une forme quadratique}, \nonumber \\
& & \nonumber \\
& &\frac{z_0z_3^2+z_1z_2z_3+z_2^3}{z_3^3}, \hspace{2mm}
\frac{z_0^2z_3^{2\kappa}}{z_3^2(z_1z_3-z_2^2)^\kappa},
 \nonumber \\
& & \nonumber \\
& &\frac{z_0z_3+z_2z_1}{z_1z_3}, \hspace{2mm}
\frac{{(z_0z_3-z_1z_2)}^\kappa}{z_2^{\kappa+1}z_3^{\kappa-1}}
\mbox{ avec } \kappa \not= \pm 1 \mbox{ et } \nonumber  \\
& &  \nonumber \\
&
&\frac{(z_0z_3^2-z_1z_2z_3+\frac{z_2^3}{3})^2}{(z_1z_3-\frac{z_2^2}{2})^3}
\mbox{ o\`u }\kappa \mbox{ et } \lambda_i \in \C. \nonumber
\end{eqnarray}

\item[(iv)] l'orbite g\'en\'erique est de dimension $1$ ; le
feuilletage associ\'e est celui d'un champ de vecteurs lin\'eaire.
\vs
\end{itemize}
\end{thm}

Ce th\'eor\`eme permet d'envisager une description topologique syst\'ematique
des feuilles des $\mathscr{L}$-feuilletages sur
$\C\P(3)$.\vs

On note l'\'etrange fait suivant : toutes les surfaces invariantes qui
apparaissent ci-dessus (ce sont les p\^oles de $\Omega$) sont donn\'ees par les
z\'eros de polyn\^omes \`a coefficients entiers.\vs

Rappelons qu'un polyn\^ome $P \in \C[z_1,\ldots,z_n]$ est dit minimal
s'il est \`a fibre g\'en\'erique connexe. Tout polyn\^ome $P$ se factorise
comme suit : $P=p(Q)$ o\`u $p\in \C[t]$ et $Q$ est un polyn\^ome
minimal.

Du th\'eor\`eme pr\'ec\'edent on d\'eduit le : \vs

\begin{cor} Soit $P \in \C[z_1,z_2,z_3]$
un polyn\^ome minimal ; notons $$G^P:=\{g \hspace{1mm} \colon
\hspace{1mm} \C^3 \to \C^3 \mbox{ affine},
\hspace{1mm} P \circ g=P\}.$$ On suppose que l'orbite
g\'en\'erique de $G^P$ est de dimension $2$, i.e. $P$ est
pr\'ecis\'ement l'invariant de $G^P$. Alors \`a conjugaison pr\`es
$P$ appartient \`a la liste qui suit :
\begin{eqnarray}
& & z_0,\hspace{2mm} z_0^{p_0}z_1^{p_1}z_2^{p_2},\hspace{2mm} z_0^{p_0}z_1^{p_1}
\nonumber \\
& & \nonumber \\
& & Q,\hspace{2mm} z_0+z_1z_2+z_2^3  \nonumber
\vs
\end{eqnarray}
o\`u $Q$ d\'esigne un polyn\^ome de degr\'e $2$, $p_i \in
\mathbb{N}$ et $\gcd(p_0,\ldots)=1$.
\end{cor}

\vs

Le chapitre 5 est consacr\'e aux $\mathscr{L}$-feuilletages de
degr\'e $2$. Actuellement nous n'en avons pas la description
g\'en\'erale pour $n \geq 4$. On pr\'esente quelques exemples et
on traite certains cas sp\'eciaux.

Dans le chapitre 6 on propose ce qui devrait \^etre la classification des
$\mathscr{L}$-feuilletages de degr\'e $3$ sur $\C\P(4)$. Comme
on l'a dit $\mathscr{L}$ est ici de dimension trois. Nous avons utilis\'e le
logiciel Maple pour d\'ecrire certaines sous-alg\`ebres r\'esolubles de
dimension $3$ de l'alg\`ebres des matrices complexes $5 \times 5$. Une fois
cette description, c'est-\`a-dire les diff\'erentes possibilit\'es
d'alg\`ebres $\mathscr{L}$, obtenue, les calculs des invariants sont localement
triviaux mais tr\`es lourds. Nous en avons pr\'esent\'e quelques-uns, m\^eme au
risque de rebuter le lecteur... La taille des calculs explose de la dimension
trois \`a la dimension quatre et certains \guillemotleft \hspace{1mm}
r\'esultats \guillemotright \hspace{1mm} ne seront finalement que des
observations. Toutefois m\^eme si l'on ne peut pr\'etendre obtenir des
\'enonc\'es au sens classique, ceci permet de pr\'esenter une grande liste
d'exemples que l'on peut esp\'erer \^etre exhaustive. Bien s\^ur, chaque fois
qu'un \'enonc\'e d\'epend de l'utilisation de Maple, nous le signalerons.

Pr\'ecisons tout de m\^eme quelques r\'esultats sur $\C\P(4)$.
Si $\mathscr{F}$ est un $\mathscr{L}$-feuilletage de degr\'e
$3$ sur $\C\P(4)$ alors l'alg\`ebre de Lie $\mathscr{L}$ fait
partie de la liste suivante : \begin{itemize}
\item[1. ] $s\ell(2,\C)$
\item[2. ] $\C^3$
\item[3. ] $\mathscr{L}_\alpha$ dont la pr\'esentation est
$\{[X,Y]=Y, \hspace{1mm} [X,Z]=\alpha Z, \hspace{1mm} [Y,Z]=0\}$
\item[4. ] $\mathscr{L}_0=\C \oplus A$ o\`u $A$ est l'alg\`ebre
de Lie du groupe affine.
\end{itemize}

On note qu'il manque ici quelques alg\`ebres de dimension $3$ ; en
fait ce sont des d\'eg\'en\'erescences des pr\'ec\'edentes qui ne
produisent pas de $\mathscr{L}$-feuilletages de degr\'e $3$. C'est
le traitement des cas $3$ et $4$ qui a n\'ecessit\'e l'usage de
Maple. Nous avons trouv\'e $17$ mod\`eles dans l'\'eventualit\'e
$3$ et $9$ dans l'\'eventualit\'e $4$ L'\'etude des situations $1$
et $2$ se fait de fa\c{c}on classique. On obtient les deux
r\'esultats suivants : \vs

\begin{thm} Soit $\mathscr{F}$ un
$\mathscr{L}$-feuilletage de degr\'e $3$ sur
$\C\P(4)$. Si $\mathscr{L}$ est ab\'elienne, alors
le feuilletage $\mathscr{F}$ admet \`a conjugaison pr\`es l'une
des int\'egrales premi\`eres suivantes :
\begin{eqnarray}
& & z_0^{\kappa_0}
z_1^{\kappa_1} z_2^{\kappa_2} z_3^{\kappa_3} z_4^{\kappa_4}
\hspace{8mm} \sum_{j=0}^4 \kappa_j=0 \nonumber \\
& &\nonumber \\
& &\frac{z_1^{\kappa+\xi+1}}{z_2^\kappa z_3^\xi z_4}\exp\left(
\frac{z_0}{z_1}\right) \nonumber \\
& &\nonumber \\
& &\frac{z_3^\kappa z_4}{z_1^{1+\kappa}}\exp\left(\frac{z_1z_2-z_0z_3}
 {z_1z_3}\right) \nonumber \\
& &\nonumber \\
& &\frac{z_0z_4^2-z_1z_2z_4+z_1z_3^2}{z_1z_4^2} \nonumber \\
& &\nonumber \\
& &\left(\frac{z_1}{z_4}\right)
\exp\left(\frac{z_0z_4^2-z_1z_2z_4+z_1z_3^2}{z_1z_4^2}\right)
\nonumber \\
& &\nonumber \\
& & \frac{z_1^\kappa z_4^{1-\kappa}}{z_0}
\exp\left(\frac{z_3^2-z_2z_4}{z_4^2}\right)\nonumber \\
& &\nonumber \\
& & \frac{z_0}{z_4}\exp\left(\frac{z_1z_4^2+z_2z_3z_4-
z_3^3}{z_4^3}\right)\nonumber \\
& &\nonumber \\
& & \frac{z_0z_4^3+z_1z_3z_4^2-\frac{z_2^2z_4^2}{2}-\frac{z_3^4}
{4}+z_2z_3^2z_4}{z_4^4} \nonumber
\end{eqnarray}
Les $\kappa_i$, $\kappa$ et $\xi$ sont des nombres complexes.
\end{thm}

Dans le cas o\`u $\mathscr{L}$ est isomorphe \`a $s\ell(2,\C)$
 on montrera le : \vs

\begin{thm} Soit $\mathscr{F}$ un
$\mathscr{L}$-feuilletage de degr\'e $3$ sur
$\C\P(4)$. Si $\mathscr{L}$ est isomorphe \` a
$s\ell(2,\C)$, alors le feuilletage $\mathscr{F}$ admet
\`a conjugaison pr\`es l'une des int\'egrales premi\`eres
suivantes : $$\frac{z_3^2z_4^2-4z_2z_4^3-4
z_3^3z_5+18z_2z_3z_4z_5}{z_1^4}$$
$$\frac{{(z_3z_2^2+z_5z_1^2-z_1z_2
z_4)}^2}{{\left(z_3z_5-\frac{z_4^2}{4}\right)}^3}$$
\end{thm}

Nous terminons par un court chapitre d'exemples de $\mathscr{L}$-feuilletages de
codimension sup\'erieure ou \'egale \`a $2$.

Notre classification est exprim\'ee en termes de feuilletages,
plus pr\'ecis\'ement en termes d'int\'egrales premi\`eres. Nous
n'avons pas list\'e les actions de groupes ; mais dans le texte
on trouvera les diff\'erentes sous-alg\`ebres de Lie qui
produisent des $\mathscr{L}$-feuilletages. On peut en d\'eduire
ais\'ement les actions de groupes correspondantes.

Les $\mathscr{L}$-feuilletages jouent aussi un r\^ole important
dans l'\'etude locale des singularit\'es de feuilletages. En
dimension trois, il existe une th\'eorie de r\'eduction des
singularit\'es (\cite{[Ca-Ce]}, \cite{[Ca]}) ; les mod\`eles
locaux apr\`es r\'eduction sont les $\mathscr{L}$-feuilletages
associ\'es \`a des alg\`ebres de Lie ab\'eliennes. Ils
interviennent aussi de fa\c{c}on naturelle sur les vari\'et\'es
compactes complexes sans fonction m\'eromorphe non constante mais
poss\'edant suffisamment d'automorphismes (tores g\'en\'eraux et
plus g\'en\'eralement certains quotients de groupes de Lie par des
r\'eseaux cocompacts etc). Sur ces vari\'et\'es $M$ tout
feuilletage est associ\'e \`a la donn\'ee d'une sous-alg\`ebre de
Lie de l'alg\`ebre $\chi(M)$. Initialement nous avons entrepris
cette \'etude pour la raison suivante. Il existe un feuilletage
$\mathscr{F}_\Gamma$ dit exceptionnel sur $\C \P(3)$ associ\'e \`a
une action du groupe des transformations affines de la droite. Ce
feuilletage est \guillemotleft \hspace{1mm} stable \guillemotright
\hspace{1mm} et par cons\'equent l'adh\'erence de l'orbite de
$\mathscr{F}_\Gamma$ sous l'action naturelle de $\mbox{Aut
}\C\P(3) \simeq \P\mbox{GL}(4,\C)$ est une composante de l'espace
des feuilletages de codimension $1$ sur $\C\P(3)$. Les
adh\'erences des feuilles de $\mathscr{F}_\Gamma$ sont des
surfaces de degr\'e six, degr\'e \guillemotleft \hspace{1mm}
anormalement \guillemotright \hspace{1mm} grand. La question est
de savoir s'il existe d'autres feuilletages exceptionnels
associ\'es \`a des actions de groupe. Nous avons trouv\'e parmi
les $\mathscr{L}$-feuilletages sur $\C\P(4)$ un candidat \`a
\^etre exceptionnel ; il est de degr\'e trois et ses feuilles sont
des hypersurfaces de degr\'e douze. Elles sont donn\'ees par les
niveaux de la fonction rationnelle
$$\frac{{\left(z_0z_4^3-\frac{z_2^2z_4^2}{2}-z_3z_4^2+z_2z_3^2z_4-
\frac{z_3^4}{4}\right)}^3}{{\left(z_1z_4^2-z_2z_3z_4+\frac{z_3^3}{3}
\right)}^4}$$ On conjecture que ce feuilletage est stable.

L'essentiel du m\'emoire est \guillemotleft \hspace{1mm}
self-contained \guillemotright \hspace{1mm} ; nous avons laiss\'e
dans le texte des calculs explicites qui, bien que tous
\'el\'ementaires, illustrent diverses techniques d'approche.
\vs

Nous remercions J. Pereira pour sa gentillesse et sa disponibilit\'e
permanentes.

%%%%%%%%%%%%%%%%%%%%%%%%%%%%%%%%%%%%%%%%%%%%%%%%%%%%%%%%%%%%%%%%%%%%%%%%%%%%
%%%%%%%%%%%%%%%%%%%%%%%%%%%%%%%%%%%%%%%%%%%%%%%%%%%%%%%%%%%%%%%%%%%%%%%%%%%%

\chapter{$\mathscr{L}$-feuilletages.}

\vs

\section[Premi\`eres propri\'et\'es.] 
{Feuilletages de codimension $1$ de $\C
\P(n)$ associ\'es \`a une alg\`ebre de Lie de champs de vecteurs.
Premi\`eres propri\'et\'es.}

\vs

On note $\pi \hspace{1mm} \colon \hspace{1mm} \C^{n+1} \setminus
\{0\} \to \C\P(n)$ la projection canonique. Si $\omega=
\displaystyle \sum_{k=0}^n A_kdz_k$ est
une $1$-forme int\'egrable sur $\C^{n+1}$ avec $A_k$
polyn\^ome homog\`ene de degr\'e $\nu$ sur $\C^{n+1}$, on d\'esigne par
$\mathscr{S}ing \hspace{1mm} \omega$ le lieu singulier de
$\omega$ $$\mathscr{S}ing \hspace{1mm} \omega:=\{A_0=\ldots=A_n=0\}$$ Dans toute
la suite on suppose, quitte \`a diviser les $A_i$ par $\gcd(A_0,\ldots,A_n)$,
que codim$_\C \mathscr{S}ing \hspace{1mm} \omega \geq 2$. On dira que
$\nu$ est le degr\'e de $\omega$.
Le feuilletage associ\'e \`a la $1$-forme $\omega$ descend
\`a $\C\P(n)$ si et seulement si $\displaystyle
\sum_{k=0}^n z_kA_k=0$ et tout feuilletage $\mathscr{F}$ de dimension $1$ sur
$\C\P(n)$ est ainsi d\'efini. Le feuilletage d\'efini par
$\omega$ sur $\C^{n+1}$ est not\'e $\tilde{\mathscr{F}}$ mais on dira
indiff\'eremment que $\omega$ d\'efinit $\mathscr{F}$ ou $\tilde{\mathscr{F}}$.
L'\'egalit\'e Aut $\C \P(n)=\P\mbox{GL}(n+1,\C)$
conduit \`a
$$\chi(\C \P(n))=\mathscr{L}\mbox{ie Aut }
\C \P(n) \simeq \mathscr{M}(n+1,\C)/\C.Id$$ o\`u
 $\chi(\C \P(n))$
d\'esigne l'ensemble des champs de vecteurs holomorphes sur
$\C \P(n)$ et $\mathscr{M}(n+1,\C)$ l'espace des
matrices complexes $(n+1) \times (n+1)$. Ainsi un champ de vecteurs
holomorphes sur $\C\P(n)$ peut \^etre vu comme un champ
de vecteurs lin\'eaire sur $\C^{n+1}$ modulo $\C.R$ o\`u
$R$ d\'esigne le champ radial $R=\displaystyle \sum_{i=0}^n z_i
\frac{\partial}{\partial z_i}$. Notons encore que l'identification
champs-matrices est naturelle puisque le flot d'un champ lin\'eaire $X$
est du type $e^{tA}x$ o\`u $A$ est \guillemotleft \hspace{1mm} la
matrice \guillemotright \hspace{1mm} de $X$. On prendra garde au fait
que si les matrices $A$ et $B$ correspondent aux champs $X$ et $Y$, alors
le crochet de Lie $[X,Y]$ correspond \`a $BA-AB$.  \vs

Soit $\mathscr{F}$ un feuilletage de codimension $1$ sur l'espace
projectif $\C \P(n)$. On introduit la
\vs

\begin{defi}
On dit que $\mathscr{F}$ est un $\mathscr{L}$-feuilletage de
codimension $1$ s'il existe une sous-alg\`ebre de Lie
$\mathscr{L}$ de $\chi(\C \P(n))$ telle qu'en tout point
g\'en\'erique $z \in \C \P(n)$ on ait la propri\'et\'e :
$$\mathscr{L}(z) \mbox{ est l'espace tangent \`a la feuille de }
\mathscr{F} \mbox{ en } z. \hspace{6mm} (*)$$ En particulier, on a
$\dim_{\C} \mathscr{L}(z)=\dim_{\C} \{ X(z),
\hspace{1mm} X \in \mathscr{L}\}=n-1$ pour $z$ g\'en\'erique.
\end{defi}

\vs

\hspace{-6.5mm} Dans la suite on suppose $\mathscr{L}$ maximale
pour la propri\'et\'e ($*$).
\vs

On peut de m\^eme d\'efinir les $\mathscr{L}$-feuilletages sur
$\C^n$ : il suffit de demander que les \'el\'ements de
$\mathscr{L}$ soient des champs de vecteurs affines (on
prend une d\'efinition rigide compatible avec la prise de carte
affine).\vs

On peut \'etendre la notion de $\mathscr{L}$-feuilletage sur
$\C\P(n)$ ou $\C^n$ \`a la codimension $q$
 ; dans ce cas $\dim \mathscr{L}(z)=n-q$. Sauf mention expresse du contraire
 tous les $\mathscr{L}$-feuilletages seront de codimension $1$ et l'on
 dira simplement $\mathscr{L}$-feuilletage.
\vs

\begin{rem}
Soit $Q(z_1,z_2,z_3)=z_1^2+z_2^2+z_3^2$ la forme quadratique
standard sur $\C^3$. Le feuilletage sur $\C \P(3)$ d\'ecrit dans
la carte affine $\C^3 \subset \C \P(3)$ par la forme $dQ$ est un
$\mathscr{L}$-feuilletage dont l'alg\`ebre de Lie $\mathscr{L}$
est isomorphe \`a $so(3;\C)$ ; une base de $\mathscr{L}$ est
$$X_{k,j}=z_k\frac{\partial}{\partial z_j}-z_j\frac{\partial}{\partial
z_k},\hspace{2mm} 1 \leq k<j \leq 3$$ en particulier l'alg\`ebre
$\mathscr{L}$ est de dimension $3$. Soit $m \in \C^3
\setminus \{0\}$, alors $\mathscr{L}(m)=<X_{k,j}(m)>$ est de
dimension $2$ car co\"{\i}ncide avec le plan tangent \`a la fibre
$Q^{-1}(m)$. Ainsi l'\'egalit\'e $\dim \mathscr{L}(m)=2$ en tout
point $m$ g\'en\'erique n'entra\^ine pas $\dim \mathscr{L}=2$. Par
ailleurs, $so(3;\C)$ est isomorphe \`a
$s\ell(2;\C)$ et la sous-alg\`ebre de dimension
$2$ : $$\mathscr{L}_0
:=\{\left(%
\begin{array}{cc}
  a & b \\
  0 & -a \\
\end{array}%
\right), \hspace{1mm} a, \hspace{1mm} b \in \C\} \subset
s\ell(2;\C)$$ peut \^etre vue dans $so(3;\C)$. Soit $\{X_1,X_2\}$
une base de $\mathscr{L}_0$. En un point g\'en\'erique $m$, les
vecteurs $X_1(m)$ et $X_2(m)$ sont ind\'ependants ; en particulier
$\mathscr{L}_0 \subsetneq so(3;\C)$ et $so(3;\C)$ produisent le
m\^eme feuilletage. C'est pour cette raison que l'on introduit une
notion de maximalit\'e pour $\mathscr{L}$.
\end{rem}

\vs

Supposons que
$\mathscr{F}$ soit un $\mathscr{L}$-feuilletage de $\C
\P(n)$ et que $X_1,\ldots,X_s$ soit une base de
$\mathscr{L} \subset \chi(\C\P(n))$. On note
$\tilde{X_1},\ldots,\tilde{X_s}$ des relev\'es de $X_1,\ldots,X_s$
\`a $\C^{n+1}$. Les $\tilde{X_k}$ sont des champs de
vecteurs lin\'eaires et le $\C$-espace vectoriel
engendr\'e par $R$ et les $\tilde{X_k}$ est une alg\`ebre de Lie
not\'ee $\tilde{\mathscr{L}}$ de dimension $s+1$ ; l'application
$\pi$ induit un morphisme d'alg\`ebres $$\pi_* \hspace{1mm} \colon
\hspace{1mm} \tilde{\mathscr{L}} \to \mathscr{L}$$ satisfaisant
$$\pi_*\tilde{X_k}=X_k, \hspace{2mm} \pi_*R=0 \mbox{ et }\ker
\pi_* =\C R$$

\hspace{-6.5mm} On remarque que si $\tilde{z} \in
\C^{n+1}$ est un point g\'en\'erique alors
$\tilde{\mathscr{L}}(\tilde{z})$ est l'espace tangent \`a la
feuille du feuilletage $\tilde{\mathscr{F}}$, relev\'e de
$\mathscr{F}$, passant par $\tilde{z}$.
\vs

Introduisons la sous-alg\`ebre de $\tilde{\mathscr{L}}$ d\'efinie
par
$$\mathscr{L}':=\{\tilde{X} \mbox{ champs lin\'eaires sur }
\C^{n+1}, \hspace{1mm} i_{\tilde{X}}d\omega=0\}.$$ o\`u
$\omega$ d\'efinit le feuilletage $\tilde{\mathscr{F}}$.

\hspace{-5.8mm}On note que $R \not \in \mathscr{L}'$ ; en
effet, l'identit\'e d'Euler assure que $i_R d\omega=(\nu+2) \omega
\not =0$ o\`u l'entier $\nu+1$ d\'esigne le degr\'e des composantes
de $\omega$.\vs

\begin{rem}
Soit $\omega$ une $1$-forme int\'egrable sur $\C^n$ d\'efinissant
un feuilletage $\tilde{\mathscr{F}}$ ; alors $d\omega$, si elle
est non nulle, d\'efinit un feuilletage de codimension $2$
\guillemotleft \hspace{1mm} contenu dans le feuilletage
$\mathscr{F}$ \guillemotright \hspace{1mm}. En effet le
th\'eor\`eme de Frobenius assure qu'au voisinage d'un point
r\'egulier~$m$, il existe des coordonn\'ees locales $x_1$,
$\ldots$, $x_n$ dans lesquelles la $1$-forme $\omega$ s'\'ecrit
$gdx_1$ o\`u $g \in \mathscr{O}^*(\C^n,m)$. Si $d\omega(m) \not
=0$, alors $dg \wedge dx_1 \not =0$ ; on peut donc supposer que
$g=1-x_2$, i.e. $\omega=(1-x_2)dx_1$ et $d\omega=dx_1\wedge dx_2$.
Ainsi localement l'espace des champs tangents au feuilletage
d\'efini par $\omega$ co\"{\i}ncide avec l'espace engendr\'e par
les champs $$\frac{\partial}{\partial
x_2},\ldots,\frac{\partial}{\partial x_n}$$ tandis que l'espace
des champs tangents au feuilletage d\'efini par $d\omega$ est
engendr\'e par
$$\frac{\partial }{\partial x_3},\ldots,\frac{\partial }{\partial x_n}$$ On
remarque au passage que si $i_Xd\omega=0$, alors $i_X\omega=0$ ; ainsi
$\mathscr{L}'$ est associ\'ee au $\mathscr{L}$-feuilletage de
codimension $2$ de $\C^{n+1}$ d\'ecrit par $d\omega$ et $\mathscr{L}'
\subset \tilde{\mathscr{L}}$.
\end{rem}

\vs

On d\'esigne par $\Omega^q(X)$ l'ensemble des $q$-formes holomorphes sur
l'espace $X$. Les deux \'enonc\'es suivants nous seront utiles :
\vs

\begin{lem} (\cite{[Sa]}) Soient $\alpha
\in \Omega^1(\C^n,0)$ tel que codim$_{\C}
\hspace{1mm} \mathscr{S}ing \hspace{1mm} \alpha \geq p+1$ et
$\beta \in \Omega^q(\C^n,0)$ avec $q \leq p$. Les
assertions suivantes sont \'equivalentes :

(i) $\alpha \wedge \beta=0$

(ii) Il existe $\gamma \in \Omega^{q-1}(\C^n,0)$ tel que
$\beta=\alpha \wedge \gamma$.

\hspace{-6mm} En particulier si $\alpha$ et $\beta$ sont deux
germes de $1$-formes holomorphes en $0$ tels que
codim$_\C$ $\mathscr{S}$ing $\alpha \geq 2$ et $\beta
\wedge \alpha=0$. Alors $\beta$ s'\'ecrit $h \alpha$ o\`u $h$ est
une fonction holomorphe.
\vs
\end{lem}

\begin{proof}[D\'emonstration]
Faisons la pour $q=1$. Ecrivons $\alpha=\displaystyle \sum_{k=1}^n
a_kdz_k$ et supposons d'abord que $$\mathscr{S}ing \hspace{1mm}
\alpha=\{a_1=\ldots=a_n=0\}=\emptyset$$ i.e. par exemple
$a_1 \not = 0$. Soit $X:=${\large
$\frac{1}{a_1}\frac{\partial}{\partial z_1}$} ; alors $i_X
\alpha=1$ et l'\'egalit\'e $\alpha \wedge \beta=0$ conduit \`a
$i_X(\alpha \wedge \beta)=0$ et $\gamma=-i_X \beta$ convient.

Si maintenant codim$_\C\mathscr{S}ing$ $\alpha \geq~2$
alors pour tout $m \not \in \mathscr{S}ing\mbox{ } \alpha$, on
a $$\beta(m)=h(m)\alpha(m)$$ o\`u $h$ est une fonction holomorphe sur
$(\C^n \setminus\mathscr{S}ing\mbox{ } \alpha,0)$. Comme codim
$\mathscr{S}ing\mbox{ }\alpha \geq 2$, le th\'eor\`eme de prolongement
d'Hartogs assure que $h$ se prolonge en une fonction holomorphe sur
$(\C^n,0)$ encore not\'ee $h$. R\'eciproquement, si $\beta=h\alpha$
alors $\alpha \wedge \beta=0$.

Pour le cas g\'en\'eral on renvoie le lecteur \`a \cite{[Sa]}.
\end{proof}

\begin{cor}\label{deRham} Si $\alpha$ et $\beta$ sont deux
$1$-formes homog\`enes de m\^eme degr\'e telles que $\alpha \wedge
\beta=0$. Alors $\beta=h\alpha$ o\`u $h$ est une fonction
rationnelle de degr\'e $0$. Si codim $\mathscr{S}$ing $\alpha \geq
2$, la fonction $h$ est constante.
\vs
\end{cor}

\begin{pro} Le morphisme $\pi_*$ induit un
isomorphisme de $\mathscr{L}'$ sur $\mathscr{L}$. \vs
\end{pro}

\begin{proof}[D\'emonstration]
Par la remarque qui pr\'ec\`ede, tout \'el\'ement de
$\mathscr{L}'$ descend en un \'el\'ement de $\mathscr{L}$. Soient
$X_1, \ldots, X_s$ une base de $\mathscr{L}$ et $\tilde{X}_1,
\ldots, \tilde{X}_s$ des relev\'es de $X_1, \ldots, X_s$ \`a
$\C^{n+1}$ par $\pi$. Les $X_k$ forment une base de $\mathscr{L}$
donc sont tangents au feuilletage ; par cons\'equent $\tilde{X}_k$
est tangent \`a $\omega$, i.e. $i_{\tilde{X}_k}\omega=0$.
Puisque $\omega \wedge d\omega=0$, on a
$$i_{\tilde{X}_k}d\omega \wedge \omega=0$$  et le corollaire \ref{deRham}
assure que
$$i_{\tilde{X}_k}d\omega=\mu(\tilde{X}_k)\omega$$ o\`u
$\mu(\tilde{X}_k) \in \C.$ Comme $L_R\omega=(\deg\omega +1)
\omega$, quitte \`a changer $\tilde{X}_k$ en $\tilde{X}_k-\frac{\mu
(\tilde{X}_k)}{deg\hspace{1mm} \omega+1}R$, on peut choisir
$\tilde{X}_k \in \mathscr{L}'$ et la proposition est \'evidente.
\end{proof}

On rappelle la
\vs

\begin{defi}
Un feuilletage $\mathscr{F}$ est de degr\'e $\nu$ sur $\C \P(n)$
s'il est d\'ecrit en coordonn\'ees homog\`enes par une $1$-forme
int\'egrable homog\`ene $\omega$ de degr\'e $\nu+1$ satisfaisant
codim$_{\C} \mathscr{S}ing \hspace{1mm} \omega \geq 2$.
\end{defi}

\vs

Donnons une interpr\'etation g\'eom\'etrique du degr\'e. Soient
$\mathscr{F}$ un feuilletage sur $\C \P(n)$,
$\mathscr{D}$ une droite g\'en\'erique de $\C\P
(n)$, $z$ un point de $\mathscr{D}$ et $\mathscr{F}_z$ la feuille
de $\mathscr{F}$ passant par $z$ ; le degr\'e de $\mathscr{F}$ est
\'egal au nombre de points $z \in \mathscr{D}$ tels que
l'hyperplan tangent \`a $\mathscr{F}_z$ en $z$ contienne $\mathscr{D}$. En
effet, soit $\omega$ une $1$-forme de degr\'e $\nu+1$ sur
$\C^{n+1}$ d\'efinissant $\mathscr{F}$ ; $\omega$
s'\'ecrit
$$\sum_{k=0}^n A_kdz_k$$ o\`u les $A_k$ sont des polyn\^omes homog\`enes
de degr\'e $\nu+1$. Dans la carte affine $z_0=1$, on note
$$\omega_0:=\omega_{|z_0=1}=\sum_{k=1}^nA_k(1,z_1,\ldots,z_n)dz_k$$
 Supposons que la droite $\mathscr{D}:=
\{z_2=\ldots=z_n=0\}$ soit une droite g\'en\'erique ; dans la
carte $z_0=1$ le fait que $\omega$ annule le champ radial se traduit sur
$\mathscr{D}$ par
$$A_0(1,z_1,0,\ldots,0)+z_1A_1(1,z_1,0,\ldots,0)=0.$$ Ainsi comme
g\'en\'eriquement (sur le choix de $\mathscr{D}$) le polyn\^ome
$A_0(1,z_1,0,\ldots,0)$ est de degr\'e $\nu+1$, on trouve que
$A_1(1,z_1,0,\ldots,0)$ est de degr\'e $\nu$. Or
$\omega_{0_{|\mathscr{D}}}=A_1(1,z_1,0,\ldots,0)dz_1$, donc
$\omega_{0_{|\mathscr{D}}}$ s'annule en $\nu$ points, i.e. le
nombre de points de tangence entre le feuilletage et la droite
$\{z_2=\ldots=z_n=0\}$ est $\nu$. \vs

La proposition suivante contr\^ole le degr\'e d'un
$\mathscr{L}$-feuilletage $\mathscr{F}$ associ\'e \`a une
sous-alg\`ebre $\mathscr{L} \subset \chi(\C\P(n))$
ainsi que la dimension de $\mathscr{L}$ en degr\'e maximal :
\vs

\begin{pro}\label{degre} 
Soit $\mathscr{F}$ un
$\mathscr{L}$-feuilletage de degr\'e $\nu$ sur $\C
\P(n)$ d\'efini par une $1$-forme $\omega$. Alors $\nu
\leq n-1$. De plus, si $\nu=n-1$, alors $\dim \mathscr{L}=n-1$.
\vs
\end{pro}

\begin{proof}[D\'emonstration] Soient $X_1$, $\ldots$, $X_s$ une base de $\mathscr{L}$ et
$\tilde{X}_k \in \chi(\C^{n+1})$ des relev\'es des $X_k$.
Soit $m$ un point g\'en\'erique, i.e. $\dim \mathscr{L}(m)=n-1$. On
peut supposer, quitte \`a r\'eindicer, que l'espace vectoriel
engendr\'e par
$$R(\tilde{m}),\tilde{X}_1(\tilde{m}),\ldots,\tilde{X}_{n-1}(\tilde{m})$$
(o\`u $\pi(\tilde{m})=m$) est de dimension $n$, de sorte
que la $1$-forme $\bar{\omega}$ d\'efinie par
$$i_Ri_{\tilde{X}_1}\ldots i_{\tilde{X}_{n-1}}dz_0
\wedge \ldots \wedge dz_n$$ est homog\`ene de degr\'e $n$, non
identiquement nulle et $\ker
\bar{\omega}(\tilde{m})=<R(\tilde{m})$, $\tilde{X}_1(\tilde{m})$,
$\ldots$, $\tilde{X}_{n-1}(\tilde{m})>$. Par ailleurs, $\ker
\omega(\tilde{m})$ est l'hyperplan tangent \`a la feuille passant
par $m$.
Ainsi si $\tilde{m}$ est un point g\'en\'erique, on a $\ker
\omega(\tilde{m})=\ker \bar{\omega}(\tilde{m})$ : les deux
$1$-formes $\bar{\omega}$ et $\omega$ sont colin\'eaires. Le lemme
de de Rham-Saito assure l'existence d'un polyn\^ome homog\`ene
$f_{n-\nu-1}$ tel que $\bar{\omega}=f_{n-\nu-1} \omega$ ; comme
$\omega$ est de degr\'e $\nu+1$ et $\bar{\omega}$ de degr\'e $n$,
le polyn\^ome homog\`ene $f_{n-\nu-1}$ est de degr\'e $n-\nu-1$ et
$n-1 \geq \nu$.

Supposons maintenant que $\nu=n-1$ ; avec les notations
pr\'ec\'edentes, on a $f_{n-\nu-1}=f_0\in\C^*$. Par suite
les champs lin\'eaires $R, \tilde{X}_1,\ldots,\tilde{X}_{n-1}$
sont lin\'eairement ind\'ependants en tout point $\tilde{m}$
r\'egulier pour $\tilde{\mathscr{F}}$ et les champs holomorphes
$X_1(m),\ldots,X_{n-1}(m)$ sont lin\'eairement ind\'ependants en
tout point $m \in \C\P(n) \setminus \mathscr{S}ing
\hspace{1mm} \mathscr{F}$. Soient $X \in \mathscr{L}$ et $m \in
\C \P(n) \setminus \mathscr{S}ing \hspace{1mm}
\mathscr{F}$, on a $$X(m)=\sum_{k=1}^{n-1} \alpha_k(m)X_k(m)$$
o\`u $\alpha_k(m) \in \C$ et les applications $\alpha_k$ sont
holomorphes sur $\C \P(n)
\setminus \mathscr{S}ing \hspace{1mm} \mathscr{F}$. Comme
codim$_{\C}$ $\mathscr{S}ing \hspace{1mm} \mathscr{F} \geq
2$, le th\'eor\`eme de prolongement de Hartogs assure que les
$\alpha_k$ se prolongent \`a $\C \P(n)$ tout
entier et sont donc constants. \end{proof}

On en d\'eduit le
\vs

\begin{cor} Soit $\mathscr{F}$ un
$\mathscr{L}$-feuilletage de degr\'e $2$ sur $\C
\P(3)$ associ\'e \`a $\mathscr{L}\subset\chi(\C
\P(3))$. Alors $\mathscr{L}$ est ab\'elienne de dimension
$2$ ou isomorphe \`a l'alg\`ebre du groupe des transformations
affines $\{z \mapsto az+b\}$.
\vs
\end{cor}

Ainsi on d\'eduira la classification des
$\mathscr{L}$-feuilletages de degr\'e $2$ sur $\C
\P(3)$ de celle de certaines sous-alg\`ebres de Lie de dimension
$2$ de $s\ell(4;\C)$.
\vs

\begin{defi}
Un feuilletage $\mathscr{F}$ d\'ecrit en coordonn\'ees homog\`enes
par la $1$-forme $\omega$ poss\`ede un facteur int\'egrant s'il
existe un polyn\^ome homog\`ene $P$ non trivial tel que la
$1$-forme {\Large $\frac{\omega}{P}$} soit ferm\'ee.
\end{defi}

\vs

\begin{pro} (\cite{[Ce-Ma]}) Soit $\mathscr{F}$ un feuilletage sur $\C
\P(n)$ donn\'e en coordonn\'ees homog\`enes par la $1$-forme
$\omega$. Supposons que $\mathscr{F}$ poss\`ede un facteur
int\'egrant $P$ ; si $P=P_1^{n_1+1}\ldots P_s^{n_s+1}$ est la
d\'ecomposition de $P$ en facteurs irr\'eductibles, alors
\begin{eqnarray}
\frac{\omega}{P_1^{n_1+1}\ldots P_s^{n_s+1}}&=&\sum_{k=1}^s
\lambda_k \frac{dP_k}{P_k} + d(\frac{H}{P_1^{n_1}\ldots
P_s^{n_s}}) \nonumber \\
&=&d(\sum_{k=1}^s \lambda_k \log P_k + \frac{H}{P_1^{n_1}\ldots
P_s^{n_s}}) \nonumber
\end{eqnarray} o\`u $\lambda_k \in \C$ et $H$ d\'esigne un
polyn\^ome homog\`ene. \vs
\end{pro}

Bien que techniquement difficile, l'\'enonc\'e est \guillemotleft
\hspace{1mm} essentiellement \guillemotright \hspace{1mm} le
th\'eor\`eme de d\'ecomposition en \'el\'ements simples des
fractions rationnelles \`a une variable. Il y a un \'enonc\'e
analogue pour les germes holomorphes (\cite{[Ce-Ma]}).
\vs

\begin{rems}
(i) Le feuilletage associ\'e \`a la forme {\Large
$\frac{\omega}{P}$} a en g\'en\'eral des singularit\'es, les
croisements $P_j=P_k=0$ notamment. Les feuilles de ce feuilletage
sont les composantes connexes des niveaux de la fonction
multivalu\'ee $$\displaystyle \sum_{k=1}^s \lambda_k \log P_k +
\frac{H}{P_1^{n_1}\ldots P_s^{n_s}}$$ \vs

(ii) Si les $\lambda_k$ sont tous nuls, le feuilletage
$\mathscr{F}$ admet une int\'egrale premi\`ere rationnelle
$H/P_1^{n_1} \ldots P_s^{n_s}$. Il suit de l'identit\'e d'Euler
$i_R \omega=0$ que $H/P_1^{n_1} \ldots P_s^{n_s}$
est homog\`ene de degr\'e $0$. Finalement, $\deg(P)=\deg \omega+1=
\nu+2$.\vs

(iii) Si l'un des $\lambda_k$ est non nul, l'\'egalit\'e
$$\frac{\omega}{P}=\displaystyle \sum_{k=1}^s \lambda_k
\frac{dP_k}{P_k} + d(\frac{H}{P_1^{n_1}\ldots P_s^{n_s}})$$ conduit
aussi \`a $\deg(P)=\deg \omega+1=\nu+2$.\vs

(iv) Il suit de (ii) et (iii) que $\frac{\omega}{P}$ d\'efinit une $1$-forme
ferm\'ee rationnelle sur $\C\P(n)$.\vs

(v) Si le polyn\^ome $P$ est r\'eduit, i.e. $n_k=0$ pour tout $k
\in \{1,\ldots,s\}$, alors toujours pour des raisons
d'homog\'en\'eit\'e le polyn\^ome $H$ est constant et on peut
prendre $H=0$. On dit alors que le feuilletage $\mathscr{F}$
d\'efini par {\Large $\frac{\omega}{P}$}$=d(\displaystyle
\sum_{k=1}^s \log P_k)$ est logarithmique ; l'identit\'e d'Euler
implique que $\displaystyle \sum_{k=1}^s \lambda_k \deg(P_k)=0$.\vs

(vi) A l'int\'egrale premi\`ere $\displaystyle \sum_{k=1}^s
\lambda_k \log P_k + \frac{H}{P_1^{n_1}\ldots P_s^{n_s}}$ on
pr\'ef\`ere parfois son exponentielle $\displaystyle \prod_{k=1}^s
P_k^{\lambda_k} \exp(\frac{H}{P_1^{n_1}\ldots P_s^{n_s}})$.
\vs
\vs
\end{rems}

La proposition qui suit est bien classique \cite{[Ce-Ma]} : \vs

\begin{pro} Soient $\mathscr{F}$ un
feuilletage sur $\C \P(n)$ d\'ecrit par une
$1$-forme $\omega$, $X \in \chi(\C\P (n))$ et
$\exp(tX)$ le flot de $X$. Supposons que $X$ ne soit pas tangent
\`a $\mathscr{F}$ et que $\mathscr{F}$ soit invariant par $X$,
i.e. $(\exp tX)^*\mathscr{F}=\mathscr{F}$. Alors $i_{\tilde{X}}
\omega$ est un facteur int\'egrant de $\mathscr{F}$, o\`u
$\tilde{X}$
 est un relev\'e de $X$. \vs
\end{pro}

\begin{proof}[D\'emonstration] Les hypoth\`eses se traduisent par :
\begin{itemize}
\item[(i)] $i_{\tilde{X}} \omega =P$ est un polyn\^ome homog\`ene non trivial

\item[(ii)] $(\exp t\tilde{X})^* \omega \wedge \omega=0$
\end{itemize}
En d\'erivant (ii), on obtient $$(L_{\tilde{X}} \omega)\wedge \omega=0$$
ou encore $$(i_{\tilde{X}} d\omega+ d(i_{\tilde{X}} \omega))
\wedge \omega=0$$ Comme $\omega$ est int\'egrable, on a
$$i_{\tilde{X}} (\omega \wedge d\omega)=0$$ soit
$$(i_{\tilde{X}} \omega)d\omega+i_{\tilde{X}}d\omega \wedge
\omega=0 $$
qui \'equivaut \`a $$d(\frac{\omega}{P})=0$$
\end{proof}

\begin{defi}
Un champ de vecteurs $X$ comme ci-dessus est appel\'e sym\'etrie
de $\mathscr{F}$.
\end{defi}

\vs

Via la proposition pr\'ec\'edente, on peut \'etablir qu'un
$\mathscr{L}$-feuilletage de $\C \P(n)$ poss\`ede
un facteur int\'egrant, c'est-\`a-dire est d\'efini par une $1$-forme
ferm\'ee rationnelle : \vs

\begin{thm} Soit $\mathscr{F}$ un
$\mathscr{L}$-feuilletage de $\C \P(n)$. Alors
$\mathscr{F}$ poss\`ede un facteur int\'egrant. Plus
pr\'ecis\'ement, on a l'alternative :
\begin{itemize}
\item[(i)] le feuilletage $\mathscr{F}$ poss\`ede une sym\'etrie

\item[(ii)] le feuilletage $\mathscr{F}$ poss\`ede une int\'egrale
premi\`ere rationnelle non triviale.
\end{itemize}
En particulier $\mathscr{F}$ est d\'efini par une $1$-forme
ferm\'ee rationnelle sur $\C\P(n)$.\vs
\end{thm}

\begin{proof}[D\'emonstration] On introduit les deux sous-groupes de Aut $\C
\P(n)$ suivants : $$G_1:=<\exp tY,\hspace{1mm} Y\in\mathscr{L}>$$ le
groupe analytique engendr\'e par les flots des champs tangents \`a
$\mathscr{F}$ et
$$G_2:=\{\phi \in \mbox{Aut } \C \P(n),
\hspace{1mm} \phi^* \mathscr{F}=\mathscr{F}\}=\{\phi \in \mbox{Aut
} \C \P(n), \hspace{1mm} \phi^* \omega \wedge
\omega=0\}$$ o\`u $\omega$ d\'efinit le feuilletage $\mathscr{F}$.

Alors que $G_2$ est un sous-groupe alg\'ebrique de Aut $\C
\P(n)$, il se peut que $G_1$ ne le soit pas. On a
l'inclusion $G_1 \subset G_2$. Notons $G_1^Z$ l'adh\'erence de
Zariski de $G_1$ et $G_2^0$ la composante connexe
neutre de $G_2$ ; comme $G_1$ est connexe, $G_1^Z$ aussi et
$G_1^Z \subset G_2^0$. On a l'alternative

- $\mathscr{L}ie \hspace{1mm} G_1 \subsetneq \mathscr{L}ie
\hspace{1mm} G_1^Z \subset \mathscr{L}ie \hspace{1mm} G_2$

- $G_1=G_1^Z$.

Dans le premier cas, si $X \in \mathscr{L}ie \hspace{1mm} G_1^Z
\setminus \mathscr{L}ie \hspace{1mm} G_1$, alors $X$ est non
tangent \`a $\mathscr{F}$ (maximalit\'e de $\mathscr{L}$) et $\exp
tX \in G_2^0$. Par suite $X$ est une sym\'etrie.

Dans le deuxi\`eme cas, le groupe $G_1$ est alg\'ebrique ; comme
il agit sur $\C \P(n)$ par action alg\'ebrique,
l'adh\'erence ordinaire de toute orbite $G_1.m$ est alg\'ebrique
comme adh\'erence de l'image d'un morphisme alg\'ebrique. Par suite
toute feuille de $\mathscr{F}$ est d'adh\'erence alg\'ebrique. Il
r\'esulte du th\'eor\`eme de Darboux-Jouanolou (\cite{[Ce-Ma]}, \cite{[Jo]}) que
$\mathscr{F}$ poss\`ede une int\'egrale premi\`ere rationnelle.
\end{proof}

Ainsi \'enonc\'e ce th\'eor\`eme est semblable \`a un r\'esultat
de J. Pereira et P. S\'{a}nchez (\cite{[Pe-Sa]}). Il se
g\'en\'eralise comme suit (la preuve est analogue) :
\vs

\begin{thm}{\it Un $\mathscr{L}$-feuilletage de
codimension un d'une vari\'et\'e alg\'ebrique est d\'efini par
une forme ferm\'ee rationnelle.}\vs
\end{thm}

\begin{lem}\label{loc} Soient $\mathscr{F}$ un
$\mathscr{L}$-feuilletage sur $\C \P(n)$ d\'ecrit
par une $1$-forme $\omega$, $X$ une sym\'etrie et $Y$ un champ
tangent au feuilletage. Alors $i_{[X,Y]} \omega=0$, i.e. $[X,Y]$
est tangent au feuilletage $\mathscr{F}$. \vs
\end{lem}

\begin{proof}[D\'emonstration] En un point r\'egulier, le th\'eor\`eme de Frobenius assure
que $\omega$ s'\'ecrit $udx_0$ o\`u $u$ est une unit\'e et $x_0$
une submersion. Si on compl\`ete les coordonn\'ees \`a l'aide de
$x_1,\ldots,x_n$ alors les champs locaux tangents \`a $\omega$ sont
engendr\'es par $\frac{\partial}{\partial x_1},\ldots,
\frac{\partial}{\partial x_n}$. Comme $X=\displaystyle
\sum_{k=0}^n X_k \frac{\partial}{\partial x_k}$ est une
sym\'etrie, la $1$-forme $\frac{\omega}{i_X\omega}$ est ferm\'ee. Or
$i_X\omega=uX_0$ ; donc $\frac{dx_0}{X_0}$ est ferm\'ee, i.e :
$$X=\ell(x_0)\frac{\partial} {\partial x_0}+\displaystyle
\sum_{k=1}^n X_k \frac{\partial}{\partial x_k}$$ Puisque le champ $Y$
\'etant tangent au feuilletage, il s'\'ecrit $\displaystyle
\sum_{k=1}^n Y_k\frac{\partial}{\partial x_k}$ ; alors
$$[X,Y]=[\ell(x_0)\frac{\partial} {\partial x_0}+\displaystyle
\sum_{k=1}^n X_k \frac{\partial}{\partial x_k},\displaystyle
\sum_{k=1}^n Y_k\frac{\partial}{\partial x_k}]$$ et $[X,Y]$ annule
$\omega$. \end{proof}

\begin{rem}
Le lemme pr\'ec\'edent est \'evidemment de nature locale.
\end{rem}

\vs

\begin{cor} Soient $\mathscr{F}$ un
$\mathscr{L}$-feuilletage sur $\C \P(n)$ et $X$ un
champ de vecteurs non tangent au feuilletage $\mathscr{F}$. Le
champ $X$ est une sym\'etrie de $\mathscr{F}$ si et seulement si
$[X,\mathscr{L}] \subset \mathscr{L}$.\vs
\end{cor}

\begin{proof}[D\'emonstration] Soit $Y \in \mathscr{L}$, le lemme 
\ref{loc} assure que si $X$
est une sym\'etrie alors le champ $[X,Y]$ est tangent au
feuilletage $\mathscr{F}$. Par maximalit\'e de l'alg\`ebre
$\mathscr{L}$, le crochet $[X,Y] \in \mathscr{L}$.

R\'eciproquement, le th\'eor\`eme de Frobenius assure que
localement au voisinage d'un point r\'egulier $m$ la $1$-forme
$\omega$ d\'ecrivant le feuilletage s'\'ecrit $udx_0$. L'alg\`ebre
$\tilde{\mathscr{L}}=\mathscr{L}' \oplus \C R$ est de
dimension ponctuelle $n$. Le module des champs tangents en $m$ au
feuilletage co\"{\i}ncide avec le module engendr\'e par
$\frac{\partial}{\partial x_1}, \ldots, \frac{\partial} {\partial
x_n}$ et avec le module engendr\'e par $\tilde{\mathscr{L}}$. Soit $\tilde{X}$
un relev\'e de $X$ ; la
condition $[\tilde{X},\tilde{\mathscr{L}}] \subset \tilde{\mathscr{L}}$
implique que pour tout champ $Y$ local tangent \`a
$\tilde{\mathscr{F}}$, le champ $[\tilde{X},Y]$ est tangent au feuilletage. En
particulier, c'est le cas pour chaque $[\tilde{X},\frac{\partial}
{\partial x_i}]$, $i \geq 1$. Supposons que le
champ $\tilde{X}$ s'\'ecrive $$X_0 \frac{\partial}{\partial x_0} +\ldots+
X_n \frac{\partial}{\partial x_n}$$ Comme la $1$-forme $\omega$
annule $$[X,\frac{\partial} {\partial x_k}]=-\frac{\partial
X_0}{\partial x_k}\frac{\partial} {\partial x_0}+\ldots$$ on a
$-\frac{\partial X_0}{\partial x_k}=0$ pour tout $k \geq 1$ ; on note puisque
$X$ n'est pas tangent \`a $\tilde{\mathscr{F}}$ que
$X_0=\ell(x_0) \not =0$. Ainsi {\Large $\frac{\omega}{i_X \omega}=
\frac{dx_0}{\ell(x_0)}$} est ferm\'ee, ce qui implique que $X$ est
une sym\'etrie.
\end{proof}

Les th\'eor\`emes qui suivent proposent un proc\'ed\'e de construction
 en cascade : \vs

\begin{thm} Soit $\mathscr{F}$ un
$\mathscr{L}$-feuilletage de degr\'e $\nu$ sur $\C
\P(n)$ d\'efini en coordonn\'ees homog\`enes sur
$\C^{n+1}$ par la $1$-forme $\omega$ de degr\'e $\nu+1$.
Supposons que le $\mathscr{L}$-feuilletage $\mathscr{F}$ n'ait pas
d'int\'egrale premi\`ere rationnelle, alors $\tilde{\mathscr{F}}$
a une sym\'etrie lin\'eaire $Z$ sur $\C^{n+1}$. Cette
sym\'etrie induit un facteur int\'egrant $P$ de degr\'e $\nu+2$
qui d\'efinit un $\mathscr{L}$-feuilletage de degr\'e $\leq \nu+1$
sur $\C^{n+1}$ que l'on peut prolonger \`a $\C
\P(n+1)$.\vs
\end{thm}

\begin{proof}[D\'emonstration] Soit $X \in\mathscr{L}'$ ; puisque $P$ est un facteur int\'egrant de
$\omega$ on a $$i_X(d\omega-\frac{dP}{P} \wedge \omega)=0$$ et comme $X \in
\mathscr{L}'$, on a $i_Xd\omega=i_X\omega=0$. Par cons\'equent, on obtient
$i_XdP=0$.

Soit $G \subset \mbox{GL}(n+1;\C)$ le groupe analytique de
$\mathscr{L}'$. Ce groupe n'est pas alg\'ebrique sinon les
feuilles de $\mathscr{F}$ seraient alg\'ebriques et par le
th\'eor\`eme de Darboux-Jouanolou, le feuilletage $\mathscr{F}$
admettrait une int\'egrale premi\`ere rationnelle ce qui est
exclu. On consid\`ere alors $G^Z$ l'adh\'erence de Zariski de $G$
; on a $\dim \hspace{1mm} \mathscr{L}ie \hspace{1mm} G<\dim
\hspace{1mm} \mathscr{L}ie \hspace{1mm} G^Z$. Soit $G^P$ le groupe
alg\'ebrique d\'efini par
$$G^P=\{h \in \mbox{GL}(n+1;\C), \hspace{1mm} P \circ h=P\}
\hspace{1mm};$$ pour tout $X \in \mathscr{L}'$, le champ $X$
annule $P$ donc $P(\exp tX)=P$ autrement dit $G \subset G^P$.
Ainsi $G^Z \subset G^P$ et $$\mathscr{L}'=\mathscr{L}ie \hspace{1mm} G
\varsubsetneq \mathscr{L}ie \hspace{1mm} G^Z \subset \mathscr{L}ie \hspace{1mm}
G^P$$ Par maximalit\'e de $\mathscr{L}$ la dimension ponctuelle de
$\mathscr{L}ie\hspace{1mm} G^P$ est $n$ ; en effet un champ lin\'eaire qui
appartient \`a $\mathscr{L}'$ en tout point g\'en\'erique appartient \`a
$\mathscr{L}'$. Le facteur int\'egrant $P$, de degr\'e $\nu+2$, apparait donc
comme invariant d'un sous-groupe alg\'ebrique de $\mbox{GL}(n+1;\C)$ \`a
orbites g\'en\'eriques de dimension $n$ ; le $\mathscr{L}$-feuilletage associ\'e
\`a $P$ est un $\mathscr{L}$-feuilletage de codimension $1$ et de degr\'e $\leq
\nu+1$ sur $\C^{n+1}$ que l'on peut prolonger \`a
$\C \P(n+1)$.
\end{proof}

\begin{rem}
L'in\'egalit\'e $\leq \nu+1$ est due au fait que l'on ne suppose
pas le polyn\^ome $P$ r\'eduit ; si $P=P_1^{n_1+1}\ldots
P_s^{n_s+1}$ est de degr\'e $\nu+2$, la $1$-forme $$P_1\ldots P_s
\displaystyle \sum_{k=1}^s (n_k+1) \frac{dP_k}{P_k}$$
d\'efinissant le feuilletage associ\'e \`a $P$ est de degr\'e
strictement inf\'erieur au degr\'e de $P$ d\`es que l'un des $n_k$
est strictement positif.
\end{rem}

\vs

On a une r\'eciproque lorsque $P$ est r\'eduit : \vs

\begin{thm} Soit $P$ un polyn\^ome
homog\`ene r\'eduit de degr\'e $k$ sur $\C^{n+1}$ dont la
d\'ecomposition en polyn\^omes irr\'eductibles s'\'ecrit $P_1
\ldots P_s$, $s \geq 2$. Si $P$ d\'efinit un
$\mathscr{L}$-feuilletage sur $\C^{n+1}$, il existe des
$\mathscr{L}$-feuilletages sur $\C \P(n)$ de
degr\'e $k-1$ ayant pour facteur int\'egrant $P$.
\vs
\end{thm}

\begin{proof}[D\'emonstration] On d\'efinit l'alg\`ebre de Lie $g$ de dimension ponctuelle
$n$ par
$$g:=\{X \mbox{ champs lin\'eaires, } i_X dP=0\}.$$ Si $X\in g$, le flot de $X$
est tangent aux niveaux de $P$ et en particulier tangent \`a
$P^{-1}(0)$ d'o\`u $$X(P_k)=\mu_k(X) P_k$$ avec $\mu_k \hspace{1mm} \colon
\hspace{1mm} g \to \C$ lin\'eaire. Soit
$(\lambda_k)_{k=1,\ldots,s}$ satisfaisant $\displaystyle
\sum_{k=1}^s \lambda_k \deg (P_k)=0$, alors la $1$-forme
$\Omega=\displaystyle \sum_{k=1}^s \lambda_i \frac{dP_k}{P_k}$
d\'efinit un $\mathscr{L}$-feuilletage sur $\C
\P(n)$. En effet, soit
$$g'=\{X \in g, \hspace{1mm} \sum_{k=1}^s
\lambda_k \mu_k(X)=0\}\subset g$$ Cette alg\`ebre de Lie
tangente au feuilletage logarithmique $\mathscr{F}_\Omega$ est de
codimension $1$ dans $g$ et de dimension ponctuelle $n-1$. Comme
$P$ est un polyn\^ome, ses niveaux (except\'e $P^{-1}(0)$)
n'adh\`erent pas \`a l'origine donc $R \not \in g$. L'alg\`ebre
$\C R+g'$ est de dimension ponctuelle $n$ et tangente au
feuilletage $\mathscr{F}_\Omega$. Pour
$(\lambda_k)_{k=1,\ldots,s}$ g\'en\'erique, la $1$-forme
$\omega=P\Omega$ satisfait codim$_\C \mathscr{S}ing
\hspace{1mm} \omega \geq 2$ et $\mathscr{F}_\Omega$ est de degr\'e
$k-1$.
\end{proof}

\begin{rem}
Si $s=1$, alors $\mu_1(X) \equiv 0$ pour tout $X \in g$.

Si $s \geq 2$, alors $\mu_i(X) \not \equiv 0$ $\forall
\hspace{1mm} i \in \{1,\ldots,s\}$ ; en effet, si $\mu_i(X) \equiv 0$
on a $X(P_i)=0$ $\forall \hspace{1mm} X \in g$. Ainsi les feuilletages
d\'ecrits par $P_i$ et $P$ co\"{\i}ncident et $s=i=1$.
\end{rem}

\vs

\section[$\mathscr{L}$-feuilletages avec une singularit\'e
isol\'ee.]{$\mathscr{L}$-feuilletages avec une singularit\'e
isol\'ee en dimension $n \geq 3$.}

Rappelons la d\'efinition d'un feuilletage \`a singularit\'e
isol\'ee.
\vs

\begin{defi} Soit $\mathscr{F}$ un feuilletage sur
un ouvert de $\C^n$. On dit que $\mathscr{F}$ a une singularit\'e
isol\'ee en $m$ si $m$ est isol\'e dans $\mathscr{S}ing$
$\mathscr{F}$.
\end{defi}

\vs

On d\'emontre de fa\c{c}on \'el\'ementaire que si $\mathscr{F}$
est d\'efini par la $1$-forme
$\omega=\displaystyle \sum_{k=1}^n a_k dz_k$ alors $\mathscr{F}$
est \`a singularit\'e isol\'ee en $m$ si et seulement si le germe
en $m$ de $\{a_1=\ldots=a_n=0\}$ est r\'eduit \`a $\{m\}$.\vs

On montre ici qu'un $\mathscr{L}$-feuilletage sur $\C
\P(n)$ ayant une singularit\'e isol\'ee est associ\'e \`a
une forme quadratique de rang maximum. Les outils n\'ecessaires
sont le th\'eor\`eme de Frobenius de Malgrange (\cite{[Ma]}) ainsi que
 le r\'esultat suivant d'alg\`ebre commutative :

\vs

\begin{pro} (\cite{[To]}) Soient $a_1,\ldots,a_n$
des \'el\'ements de $\mathscr{O}(\C^n,0)$. Les assertions
suivantes sont \'equivalentes
\begin{itemize}
\item[(i)] $\{a_1=\ldots=a_n=0\}=\{0\}$

\item[(ii)] le module des relations $$M_a:=\{(b_1,\ldots,b_n)
\hspace{1mm} | \hspace{1mm} \displaystyle \sum_{k=1}^n a_kb_k=0, b_i \in
\mathscr{O}(\C^n,0)\}$$ est
engendr\'e par les relations triviales
$$b_{kj}=(0,\ldots,0,a_j,0,\ldots,0,-a_k,0,\ldots,0).$$
\end{itemize}
\end{pro}

\vs

Le th\'eor\`eme de Frobenius singulier de Malgrange donne
 une condition g\'eom\'etrique d'existence d'int\'egrale premi\`ere
  holomorphe pour un feuilletage singulier :

\vs

\begin{thm} (\cite{[Ma]}) Soit $\omega$ un germe
de $1$-forme holomorphe int\'egrable en $0$. Si
codim$_{\C} \hspace{1mm} \mathscr{S}ing \hspace{1mm}
\omega \geq 3$, il existe $f \in \mathscr{O}(\C^n,0)$ et
$g \in \mathscr{O}^*(\C^n,0)$ tels que $\omega=gdf$, i.e.
$f$ est int\'egrale premi\`ere du feuilletage associ\'e \`a
$\omega$.
\vs
\end{thm}

Le th\'eor\`eme de Malgrange est inop\'erant en dimension $2$ : la
condition codim$_{\C} \mathscr{S}ing$ $\omega \geq 3$ ne
peut \^etre satisfaite sauf lorsque $\omega$ n'a pas de
singularit\'e et l\`a il s'agit du th\'eor\`eme de Frobenius
classique. \vs

Enon\c{c}ons le th\'eor\`eme principal de cette partie :
\vs

\begin{thm} \label{cormal}
Soit $\mathscr{F}$ un
$\mathscr{L}$-feuilletage sur $\C \P(n)$, $n \geq
3$, ayant un point singulier isol\'e. Alors $\mathscr{F}$ est de
degr\'e $1$ et il existe une carte affine $\C^n \subset
\C\P(n)$ telle que $\mathscr{F}_{|\C^n}$
soit associ\'e \`a une forme quadratique de rang maximum.
\vs
\end{thm}

\begin{rem}
L'\'enonc\'e ne dit pas que le lieu singulier est r\'eduit \`a un
point mais qu'il existe un point singulier isol\'e.
\end{rem}

\vs

\begin{proof}[D\'emonstration] Supposons que le point singulier de $\mathscr{F}$ soit
l'origine $0$ d'une carte affine $\C^n \subset\C
\P(n)$. Si $X \in \mathscr{L}$, on sait que dans la carte
affine $\C^n \subset \C \P(n)$ le champ
$X$ s'\'ecrit
$$\sum_{k=1}^n a_k(z)\frac{\partial}{\partial z_k}+\tilde{\ell}(z)R$$
o\`u $\tilde{\ell}$ est lin\'eaire, les $a_k$ affines et $R$ le champ radial de
$\C^n$.

Le champ $X$ s'annule en $0$ sinon toute la trajectoire de $X$
passant par $0$ serait singuli\`ere pour $\mathscr{F}$, ce qui
contredirait l'hypoth\`ese $\mathscr{F}$ admet une singularit\'e
isol\'ee. Ainsi les $a_k$ sont lin\'eaires, i.e. $X$
s'\'ecrit
$$X'+\ell(X)R$$ o\`u $\ell \hspace{1mm} \colon \hspace{1mm}
\mathscr{L} \to (\C^n)^*$ est lin\'eaire et $X'$
 un champ lin\'eaire.

Soient $X_1, \ldots, X_{n-1} \in \mathscr{L}$ tels que pour
$z$ g\'en\'erique, on ait $$\dim_{\C} \hspace{1mm}
<X_1(z),\ldots,X_{n-1}(z)>=n-1$$ La $1$-forme $\Omega$ d\'efinie
par $$i_{X_1}\ldots i_{X_{n-1}} dz_1\wedge \ldots \wedge dz_n$$
est int\'egrable et s'\'ecrit $$\Omega_{n-1}+\Omega_n$$ o\`u les
$\Omega_k$ sont homog\`enes de degr\'e $k$,
$\Omega_{n-1}=i_{X'_1}\ldots i_{X'_{n-1}} dz_1\wedge \ldots \wedge
dz_n$ et $i_R \Omega_n=0$. Il se peut que $\Omega$ s'annule sur une
hypersurface et se laisse diviser par un polyn\^ome ; c'est en
fait le cas. Le feuilletage $\mathscr{F}$ est d\'efini par une
$1$-forme \`a singularit\'e isol\'ee en $0$ ; le
th\'eor\`eme de Malgrange assure l'existence de $f \in
\mathscr{O}(\C^n,0)$ telle que $df$ d\'efinisse localement $\mathscr{F}$
en l'origine. Si $X=\displaystyle \sum_{k=1}^n b_k \frac{\partial}{\partial z_k}
\in \mathscr{L}$, alors $X.f=i_Xdf=0$, autrement dit
$\displaystyle \sum_{k=1}^n b_k \frac{\partial f} {\partial z_k}=0$. Comme
n\'ecessairement $df$ est \`a singularit\'e isol\'ee, la proposition
1.2.3. assure que le module des relations des $\frac{\partial
f}{\partial z_k}$ est engendr\'e par les relations triviales ; on
a $$X=\sum \mu_{kj}\left(\frac{\partial f}{\partial
x_j}\frac{\partial }{\partial x_k}-\frac{\partial f}{\partial
x_k}\frac{\partial }{\partial x_j}\right)$$ o\`u $\mu_{kj} \in
\mathscr{O}(\C^n,0)$. Comme $\dim_{\C}
\hspace{1mm} <X_1(z),\ldots,X_{n-1}(z)>=n-1$ pour $z$
g\'en\'erique, les $X'$ ne
sont pas tous identiquement nuls et puisque $f$ a une
singularit\'e isol\'ee \`a l'origine, $f$ est non submersive. De sorte
qu'en calculant le $1$-jet d'un \'el\'ement $X$ g\'en\'erique on
constate que
$$f=q+\mbox{ termes d'ordres sup\'erieurs}$$ o\`u $q$ est une
forme quadratique non triviale. Comme $\Omega \wedge df=0$, le
lemme de de Rham-Saito assure l'existence de $G \in
\mathscr{O}(\C^n,0)$ tel que
$$\Omega=\Omega_{n-1}+\Omega_n=Gdf=G(dq+\mbox{ termes d'ordre
sup\'erieur}).$$ Remarquons que $\Omega_{n-1} \not = 0$ sinon on
aurait $$i_R \Omega=i_R \Omega_{n-1}+i_R \Omega_n=0$$ qui implique
$i_Rdf=0$, ce qui est absurde. Comme $\Omega$ est d'ordre $n-1$
et $dq$ d'ordre $1$, l'\'egalit\'e
$$\Omega=\Omega_{n-1}+\Omega_n=G(dq + \mbox{termes d'ordre
sup\'erieur})$$ implique que $G$ est d'ordre $n-2 \geq 1$. On en
d\'eduit que $G$ s'annule sur une hypersurface et par cons\'equent
$\Omega$ aussi ; or $\Omega$ est polyn\^omiale donc cette
hypersurface est alg\'ebrique. Ainsi $G$ s'\'ecrit $UP$ avec $U
\in \mathscr{O}^*(\C^n,0)$ et $P$ un polyn\^ome :
$$P_{n-2}+\ldots+P_N$$ o\`u les $P_i$ sont homog\`enes de degr\'e $i$ et
$P_{n-2} \not =0$. Finalement $\Omega$ s'\'ecrit
$$\Omega=\Omega_{n-1}+\Omega_n=Gdf=PUdf=(P_{n-2}+\ldots+P_N)(\omega_1
+ \ldots+ \omega_{N'})=P\omega$$ o\`u les $\omega_k$ sont des $1$-formes
homog\`enes de degr\'e $k$, $\omega_1=dq$ avec $P_N$ et $\omega_{N'}$
non nuls.

\hspace{-6.5mm} Par identification, on a
$$\left\{
\begin{array}{c}
\Omega_{n-1}=P_{n-2}dq \hspace{24mm} \\
\Omega_n=P_{n-2}\omega_2+P_{n-1}dq \hspace{11mm} \\
0=P_{n-2}\omega_3+P_{n-1}\omega_2+P_ndq \\
\vdots \\
0=P_N\omega_{N'} \hspace{31mm}
\end{array} \right.$$

\hspace{-6.5mm} Il en r\'esulte que $$\left\{
\begin{array}{l} N+N'=n \\ 1 \leq N' \\ n-2 \leq N \end{array} \right.$$
 autrement dit on a l'alternative :
\begin{itemize}
\item[- ] $ N'=1 \mbox{ et } N=n-1$

\item[- ] $ N'=2 \mbox{ et } N=n-2.$
\end{itemize}

Dans le premier cas, $P=P_{n-2}+P_{n-1}$, $\omega=\omega_1$ et
$\mathscr{F}$ est un feuilletage de degr\'e $1$. Dans le second cas,
$P=P_{n-2}$ et $\omega=\omega_1+\omega_2$ avec $i_R\omega_2=0$ car $i_R
\Omega_n=0$. Ce qui signifie ici encore que $\mathscr{F}$
est un feuilletage de degr\'e $1$. La classification des
feuilletages de degr\'e $1$ (que nous rappelons au chapitre 3) nous permet de
trouver une carte affine dans laquelle $\omega_2=0$. Ainsi $\mathscr{F}$ est
d\'efini par $dq$ o\`u $q$ est une forme quadratique n\'ecessairement de rang
maximum.
\end{proof}

\vs

\begin{rem}
Ce th\'eor\`eme n'est pas correct en dimension $2$ ; sur $\C
\P(2)$ tous les feuilletages de degr\'e $1$ sont des
$\mathscr{L}$-feuilletages et poss\`edent une singularit\'e
isol\'ee (trois compt\'ees avec multiplicit\'e) sans toutefois
\^etre associ\'es \`a une forme quadratique.
\end{rem}

\vs

\section{Th\'eor\`emes g\'en\'eraux.}

\vs

On rappelle quelques r\'esultats \'el\'ementaires classiques sur
les alg\`ebres de Lie et on pr\'esente leur traduction en terme de
$\mathscr{L}$-feuilletages. On renvoie aux livres classiques pour
ces r\'esultats.
\vs

\begin{pro}\label{semisimple} 
Soit $g$ une sous-alg\`ebre de
Lie de $\mathscr{M}(n,\C)$ dont tous les \'el\'ements sont
semi-simples. Alors $g$ est isomorphe \`a une alg\`ebre diagonale.
\vs
\end{pro}

On en d\'eduit ais\'ement la \vs

\begin{pro} Soit $g$ une alg\`ebre de Lie
de champs de vecteurs lin\'eaires sur $\C^{n+1}$. Si tous
les \'el\'ements de $g$ sont semi-simples, alors $g$ est
ab\'elienne et $\dim_{\C} g= \dim_{\C} g(m)$ pour
$m \in \C^{n+1}$ g\'en\'erique.
\vs
\end{pro}

\begin{proof}[D\'emonstration] Par la proposition \ref{semisimple},
$g$ est isomorphe \`a une
alg\`ebre de matrices diagonales ; ainsi si $X_1,\ldots,X_s$ est
une base de $g$, on peut trouver un syst\`eme de coordonn\'ees
$z_0,\ldots,z_n$ tel que
$$X_k=\sum_{j=0}^n \mu_{j,k} z_j \frac{\partial}{\partial z_j}$$
Alors la matrice $$(\mu_{j,k})_{\stackrel{0\leq j \leq n}{1\leq k
\leq s}}$$ a pour rang $s$ la dimension de $g$. Evaluons
l'\'el\'ement $\displaystyle \sum_{k=1}^s \lambda_k X_k$ de
l'alg\`ebre $g$ au point g\'en\'erique $(1,\ldots,1)$ (chaque
$m=(m_1,\ldots,m_n)$ avec $m_j \not =0$ se ram\`ene \`a
$(1,\ldots,1)$ par un isomorphisme
 diagonal). Sous forme matricielle on obtient : $$\left(%
\begin{array}{ccc}
  \displaystyle \sum_{k=1}^s \lambda_k \mu_{0,k} &  & \mbox{{\LARGE $0$}} \\
   & \ddots &  \\
  \mbox{{\LARGE $0$}} &  & \displaystyle \sum_{k=1}^s \lambda_k \mu_{n,k} \\
\end{array}%
\right) \left(%
\begin{array}{c}
  1 \\
  \vdots \\
  1 \\
\end{array}%
\right)%= \left(%
%\begin{array}{c}
%  \displaystyle \sum_{k=1}^s \lambda_k \mu_{0,k} \\
%  \vdots \\
%  \displaystyle \sum_{k=1}^s \lambda_k \mu_{n,k} \\
%\end{array}%
%\right)
=(\mu_{j,k})_{\stackrel{0 \leq j \leq n}{1\leq k \leq s}}\left(%
\begin{array}{c}
  \lambda_1 \\
  \vdots \\
  \lambda_s \\
\end{array}%
\right)$$ Par suite l'espace engendr\'e par les $\displaystyle
\sum_{k=1}^s \lambda_k X_k$ \'evalu\'e en $(1,\ldots,1)$ est de
dimension $s$.
\end{proof}

On obtient comme cons\'equence directe : \vs

\begin{cor} Soit $\mathscr{F}$ un
$\mathscr{L}$-feuilletage (de codimension $1$) sur
$\C\P(n)$ associ\'e \`a $\mathscr{L} \subset
\chi(\C\P(n))$. Si tous les \'el\'ements de
$\mathscr{L}$ sont semi-simples, alors $\mathscr{F}$ est
lin\'eairement conjugu\'e \`a un feuilletage logarithmique du type
$$z_0\ldots z_n \sum_{j=0}^n \lambda_j \frac{dz_j}{z_j},
\hspace{4mm} \lambda_j \in \C$$ o\`u les $z_j$ sont des
coordonn\'ees.
\vs
\end{cor}

\begin{proof}[D\'emonstration] Soit $\omega$ une $1$-forme d\'efinissant le feuilletage
$\mathscr{F}$. On identifie $\mathscr{L}$ \`a $\mathscr{L}'$. La proposition
1.3.1. assure l'existence d'un
syst\`eme de coordonn\'ees $z_0,\ldots,z_n$ dans lequel
$\mathscr{L}'$ est diagonale. Si $X \in \mathscr{L}'$, $X$ s'\'ecrit
$$\sum_{j=0}^n \mu_j(X)z_j\frac{\partial }{\partial z_j}$$ o\`u
$\mu_j(X) \in \C$. Alors l'alg\`ebre $\tilde{\mathscr{L}}$ est engendr\'ee par
les champs $X_1,\ldots,X_{n-1}$ base
de $\mathscr{L}'$ et $X_0=R$. Si $(\lambda_0,\ldots, \lambda_n)$
satisfait
$$\left\{
\begin{array}{l} \displaystyle \sum_{j=0}^n \lambda_j=0 \\
\displaystyle \sum_{j=0}^n \mu_j(X_k) \lambda_j=0 \mbox{ pour }
k=1,\ldots,n-1\end{array} \right.$$ alors $\bar{\omega}=z_0\ldots
z_n \displaystyle \sum_{j=0}^n \lambda_j \frac{dz_j}{z_j}$
d\'efinit $\mathscr{F}$.
\end{proof}

\begin{rems}
(i) Tous les $\mathscr{L}$-feuilletages sur $\C \P(n)$ associ\'es
\`a une alg\`ebre ab\'elienne ne sont pas conjugu\'es \`a un
feuilletage du type $z_0\ldots z_n \displaystyle \sum_{j=0}^n
\lambda_j \frac{dz_j}{z_j}$. Par exemple,
 les champs $z_j \frac{\partial }{\partial z_0}$, $j \not = 1$,
engendrent une alg\`ebre $g$ ab\'elienne satisfaisant $\dim g(X)
\leq~1$ pour tout $X \in g$. Ainsi l'hypoth\`ese tous les
\'el\'ements sont semi-simples est importante pour le dernier
point, exiger que l'alg\`ebre soit ab\'elienne ne suffit pas.\vs

(ii) L'adh\'erence des feuilletages conjugu\'es aux feuilletages
logarithmiques du type $$z_0\ldots z_n \sum_{k=0}^n \lambda_k
\frac{dz_k}{z_k} \hspace{6mm}
(**)$$ forme une composante irr\'eductible de l'adh\'erence de
l'ensemble des feuilletages de degr\'e $n-1$ sur $\C
\P(n)$ (\cite{[Ce-LN]}).\vs

(iii) On note aussi que les limites de feuilletages d\'efinis par des
$1$-formes du type $(**)$ ne sont pas n\'ecessairement des
$\mathscr{L}$-feuilletages ; soit $\mathscr{F}$ le feuilletage sur
$\C^4=\{z_0,z_1,z_2,z_3\}$ d\'ecrit par la $1$-forme
$$\omega=z_0z_1(z_1-z_0)z_3 \left(\kappa_0 \frac{dz_0}{z_0}+\kappa_1
\frac{dz_1}{z_1} +\kappa_2\frac{d(z_1-z_0)}{z_1-z_0} + \kappa_3
\frac{ dz_3} {z_3}\right)$$ avec $\displaystyle \sum_{i=0}^3
\kappa_i=0$. On remarque que $\omega$ appara\^it comme
$\displaystyle \lim_{\varepsilon \to 0} \omega_\varepsilon$ o\`u
$$\omega_\varepsilon=z_0z_1(\varepsilon z_2+ z_1-z_0)z_3 \left(\kappa_0
 \frac{dz_0}{z_0}+\kappa_1 \frac{dz_1}{z_1} +\kappa_2
 \frac{d(\varepsilon z_2+z_1-z_0)}{z_1-z_0} + \kappa_3
\frac{ dz_3} {z_3}\right)$$ avec $\displaystyle \sum_{i=0}^3
\kappa_i=0$ et donc $\mathscr{F}$ est limite des $\mathscr{L}$-feuilletages
$\mathscr{F}_\varepsilon$. On peut montrer que $\omega$ ne d\'efinit pas
un $\mathscr{L}$-feuilletage de la fa\c{c}on suivante ; si
\begin{eqnarray}
\varphi \hspace{1mm} \colon \hspace{1mm} \C &\to& \C^4 \nonumber
\\
z &\mapsto& (\varphi_0(z),\varphi_1(z),\varphi_2(z),\varphi_3(z))\nonumber
\end{eqnarray}
est une application enti\`ere
qui \'evite les hyperplans $z_0=0$, $z_1=0$ et $z_0=z_1$, il r\'esulte du
th\'eor\`eme de Picard qu'il existe une relation non triviale
$$a_0\varphi_0+a_1\varphi_1=0, \hspace{2mm} a_i \in \C$$ En appliquant
cette remarque aux flots des champs lin\'eaires tangents \`a $\omega$, on
constate que la condition de dimension ponctuelle trois ne peut \^etre
r\'ealis\'ee.
\vs
\end{rems}

On rappelle la \vs

\begin{defi}
Soit $g$ une alg\`ebre de Lie, on note $(g_i)_{i\in\mathbb{N}}$ la
suite d\'efinie par $$g_0=g, \hspace{1mm} g_{i+1}=[g_i,g_i]$$ On
dit que $g$ est r\'esoluble si $g_p=\{0\}$ pour un certain $p$.
\end{defi}

\vs

Rappelons aussi le th\'eor\`eme de triangulation des alg\`ebres de Lie
r\'esolubles : \vs

\begin{thm}\label{resoluble}
Toutes les alg\`ebres de Lie complexes r\'esolubles de
matrices $g$ sont triangulables, i.e. il existe une matrice
inversible $P$ telle que $PgP^{-1}$ soit triangulaire. \vs
\end{thm}

\begin{defi}
Soit $\mathscr{F}$ un $\mathscr{L}$-feuilletage sur $\C\P(n)$. On
dit que $\mathscr{L}$ est r\'eguli\`ere si pour tout point $m \in
\C\P(n) \setminus \mathscr{S}ing \hspace{1mm} \mathscr{F}$ la
dimension ponctuelle de $\mathscr{L}$ en $m$ est $n-1$.
\end{defi}

\vs

On remarque qu'un $\mathscr{L}$-feuilletage de degr\'e maximal $n-1$ est
n\'ecessairement r\'egulier. \vs

On d\'eduit de \ref{resoluble} le :
\vs

\begin{thm} Soit $\mathscr{F}$ un
feuilletage sur $\C \P(n)$ associ\'e \`a
$\mathscr{L} \subset \chi(\C\P(n))$ r\'eguli\`ere.
Si $\mathscr{L}$ est r\'esoluble, alors $\mathscr{F}$ poss\`ede un
hyperplan invariant $\mathscr{H}$. En particulier, il existe une
carte affine $\C^n \simeq \C \P(n)
\setminus \mathscr{H}$ dans laquelle les \'el\'ements de
$\mathscr{L}$ sont affines.
\vs
\end{thm}

\begin{proof}[D\'emonstration] L'alg\`ebre $\mathscr{L}'$ est une alg\`ebre de Lie r\'esoluble isomorphe
\`a une sous-alg\`ebre de $\mathscr{M}(n+1,\C)$. D'apr\`es
le th\'eor\`eme \ref{resoluble} cette alg\`ebre est triangulable et admet
donc pour hyperplan invariant $\tilde{\mathscr{H}}:=\{z_n=0\}$,
o\`u $\{z_0,\ldots,z_n\}$ est un syst\`eme de coordonn\'ees
triangularisant de $\C^{n+1}$. Puisque $\mathscr{L}$ est r\'eguli\`ere,
en un point g\'en\'erique $m$ de $\tilde{\mathscr{H}}$ la dimension ponctuelle
de $\mathscr{L}'$ est $n-1$, et $\tilde{\mathscr{H}}$ est invariant par
 le feuilletage ; d'o\`u
l'existence via la projection canonique $\pi \hspace{1mm} \colon \hspace{1mm}
\C^{n+1}\setminus \{0\} \to \C \P(n)$ d'un
hyperplan invariant pour $\mathscr{F}$.

\hspace{-6mm} Enfin, on montre facilement que si le feuilletage
$\mathscr{F}$ admet un hyperplan invariant $\mathscr{H}$, alors dans
la carte affine $\C\P(n)\setminus \mathscr{H}$ les
 champs tangents \`a $\mathscr{F}$ sont affines.
\end{proof}

L'hypoth\`ese de r\'egularit\'e est n\'ecessaire. Nous d\'etaillerons en 
5.1.4 l'exemple d'un feuilletage de degr\'e deux sur $\C \P(4)$ ayant
une int\'egrale premi\` ere de type $\frac{Q_1}{Q_2}$, $Q_i$ quadratique,
associ\'e \`a une alg\`ebre de Lie r\'esoluble non r\'eguli\`ere maximale ;
l'alg\`ebre poss\'ede un hyperplan invariant, le feuilletage non.\vs

On rappelle aussi les \vs

\begin{defis} 1) Soit $g$ une alg\`ebre de Lie ; il
existe un unique id\'eal r\'esoluble not\'e $\mathscr{R}(g)$
appel\'e radical de $g$ contenant tous les id\'eaux r\'esolubles
de $g$.

2) Une alg\`ebre $g$ est semi-simple si elle ne contient pas
d'id\'eal r\'esoluble non trivial ce qui est \'equivalent \`a
$\mathscr{R}(g) =\{0\}$.
\end{defis}

Il est assez naturel de regarder les
$\mathscr{L}$-feuilletages associ\'es aux alg\`ebres semi-simples : \vs

\begin{thm} Soit $\mathscr{F}$ un
$\mathscr{L}$-feuilletage associ\'e \`a la sous-alg\`ebre
$\mathscr{L} \subset \chi(\C\P(n))$. Si
$\mathscr{L}$ est semi-simple, ou plus g\'en\'eralement si
$[\mathscr{L},\mathscr{L}]=\mathscr{L}$, alors $\mathscr{F}$
poss\`ede une int\'egrale premi\`ere rationnelle.
\vs
\end{thm}

\begin{proof}[D\'emonstration] Comme d'habitude on rel\`eve $\mathscr{L}$ \`a
$\C^{n+1}$ en $\mathscr{L}'$. Le feuilletage $\mathscr{F}$ poss\`ede un
 facteur int\'egrant $P=P_1^{n_1+1}\ldots P_s^{n_s+1}$, avec
$\gcd(P_1,\ldots,P_s)=1$, i.e.
$$\frac{\omega}{P}=\sum_{k=1}^s \lambda_k \frac{dP_k}{P_k}+d\left(
\frac{H}{P_1^{n_1}\ldots P_s^{n_s}}\right)$$ o\`u $\omega$ d\'efinit
$\mathscr{F}$. Si tous les
$\lambda_i$ sont nuls, le feuilletage $\mathscr{F}$ poss\`ede une
int\'egrale premi\`ere rationnelle ; si ce n'est pas le cas
supposons par exemple que
$\lambda_1 \not =0$. Soit $X \in \mathscr{L}'$ ; le champ $X$ est
tangent \`a $\tilde{\mathscr{F}}$ et par suite tangent aux hypersurfaces
invariantes par $\tilde{\mathscr{F}}$. Or le lieu des z\'eros de $P_k$ en est
une pour $k=1,\ldots,s$ ; ceci se traduit par :
$$X(P_k)=i_XdP_k=\mu_k(X)P_k.$$ A priori $\mu_k(X) \in
\mathscr{O}(\C^{n+1},0)$ mais $P_k$ et $X(P_k)$ sont
homog\`enes de m\^eme degr\'e donc $\mu_k(X) \in\C$. On
constate que si $X$ est un commutateur, i.e. si $X$ s'\'ecrit
$[Y,Z]$ avec $Y$, $Z \in \mathscr{L}'$, alors $\mu_k(X)=0$. En fait chaque
$\mu_k \hspace{1mm} \colon \hspace{1mm} \mathscr{L}' \to
\C$ est un morphisme d'alg\`ebres de Lie donc trivial puisque
$\hspace{1mm} [\mathscr{L}',\mathscr{L}']=\mathscr{L}'$. Il se
trouve que $P$ admet au moins deux facteurs premiers entre eux\
dans sa d\'ecomposition, i.e. $s \geq 2$. En effet supposons que
$P=P_1^{n_1+1}$, alors
$$\frac{\omega}{P}=\frac{\omega}{P_1^{n_1+1}}=\lambda_1
\frac{dP_1}{P_1}+d\left(\frac{H}{P_1^{n_1}}\right)$$ En particulier
$\deg d\left(\frac{H}{P_1^{n_1}}\right)=-1$ et $\deg
\frac{H}{P_1^{n_1}}=0$. L'\'egalit\'e $\deg \frac{H}{P_1^{n_1}}=0$
et l'identit\'e d'Euler impliquent que $\lambda_1 \deg P_1=0$,
i.e. $\lambda_1=0$ ce qui est exclu, ou $\deg P_1=0$ ce qui est
encore exclu.

Finalement pour chaque $X \in \mathscr{L}$ on a $X(P_1)=X(P_2)=0$
et si $\nu_i=\deg P_i$, la fonction
rationnelle $\frac{P_1^{\nu_2}}{P_2^{\nu_1}}$ est homog\`ene
non constante de degr\'e $0$ donc passe \`a $\C
\P(n)$ et est int\'egrale premi\`ere de tout \'el\'ement
de $\mathscr{L}$ donc de $\mathscr{F}$.
\end{proof}

\begin{rems}
(i) La construction pr\'ec\'edente donne une \guillemotleft
\hspace{1mm} construction
\`a deux branches \'eventuellement multiples \guillemotright
\hspace{1mm} : un
$\mathscr{L}$-feuilletage associ\'e \`a une alg\`ebre
$\mathscr{L}$ satisfaisant $[\mathscr{L},\mathscr{L}]=\mathscr{L}$
admet une int\'egrale premi\`ere du type
$\frac{P_1^{\nu_2}}{P_2^{\nu_1}}$ avec $P_1$ et $P_2$ deux
polyn\^omes irr\'eductibles et $\nu_i=\deg P_i$. On peut en effet
 lorsque les $\lambda_i$ sont nuls adapter le raisonnement
pr\'ec\'edent.\vs

(ii) Dans le cas o\`u $\mathscr{L}$ est semi-simple, on sait que
$\mathscr{L}$ s'int\`egre en un groupe alg\'ebrique ce qui donne
imm\'ediatement le r\'esultat. En fait l'\'enonc\'e donne plus que
le cas semi-simple et donne le (i) de la remarque qui a un
int\'er\^et en soi. 
\end{rems}

\vs

\begin{thm} Soit $\mathscr{F}$ un
$\mathscr{L}$-feuilletage sur $\C\P(n)$ d\'ecrit
par la $1$-forme $\omega$. Si tous les \'el\'ements de
$\mathscr{L}$ sont nilpotents, alors $\mathscr{F}$ poss\`ede une
int\'egrale premi\`ere rationnelle.
\vs
\end{thm}

\begin{proof}[D\'emonstration] Comme toujours on rel\`eve l'alg\`ebre $\mathscr{L}$ \`a $\C^{n+1}$
en l'alg\`ebre $\mathscr{L}'$. Puisque tous les \'el\'ements
de $\mathscr{L}$ sont nilpotents, $\mathscr{L}'$ est nilpotente et apr\`es
triangulation les hyperplans
$z_n=$ cte sont invariants par les \'el\'ements de $\mathscr{L}'$. Par
cons\'equent, les feuilles du feuilletage $\mathscr{F}'$ associ\'e \`a
$\mathscr{L}'$ sont contenues dans les hyperplans $z_n=$ cte et
sont, puisque $\mathscr{L}'$ produit un feuilletage de codimension
$2$, de codimension $1$ dans ces hyperplans. La forme
$d\omega_{|z_n=cte}$ est nulle sinon le feuilletage restreint
serait de codimension $2$. On a ainsi
l'\'egalit\'e $d\omega \wedge~dz_n=~0$
qui conduit \`a $d(\omega \wedge dz_n)=0$ ; autrement dit $\omega$
est ferm\'ee dans les fibres $z_n=$ cte donc par le th\'eor\`eme
de Poincar\'e exacte dans ces fibres. Par suite, la $1$-forme
$\omega$ s'\'ecrit $\alpha dz_n+dH$. Par ailleurs, l'\'egalit\'e
$\deg H=~\deg \alpha+~1=\deg \omega +1$ et l'identit\'e d'Euler
$\alpha z_n+H\deg
H=0$ impliquent que $H=-\frac{\alpha z_n}{\deg \omega +1}$ ; donc la
$1$-forme $\omega$ s'\'ecrit $$(\deg \omega) \alpha dz_n-z_n
d\alpha$$ Ainsi le feuilletage $\mathscr{F}$ admet {\Large
$\frac{\alpha} {z_n^{\deg \omega}}$} comme int\'egrale premi\`ere.
\end{proof}

\begin{rem}
Dans la carte affine $z_n=1$ le feuilletage admet pour
int\'egrale premi\`ere polynomiale $\alpha$.
\end{rem}

%%%%%%%%%%%%%%%%%%%%%%%%%%%%%%%%%%%%%%%%%%%%%%%%%%%%%%%%%%%%%%%%%%%%%%%%
%%%%%%%%%%%%%%%%%%%%%%%%%%%%%%%%%%%%%%%%%%%%%%%%%%%%%%%%%%%%%%%%%%%%%%%%

\chapter{Exemples.}

\vs

\section{Principe de construction d'exemples.}

\vs

On va d\'ecrire un principe qui permet de construire des
$\mathscr{L}$-feuilletages \`a partir d'exemples classiques venant
en particulier de la th\'eorie des invariants et des actions de
$\C^{n-1}$ sur $\C^n$.

On note $E_k:=\C^{n_k+1}$, $n_k \in \mathbb{N}^*$, et
$\P(E_k) \simeq \C \P(n_k)$,
$k=1,\ldots,q$. Sur $\P(E_k)$, on se donne un
$\mathscr{L}$-feuilletage $\mathscr{F}_k$ associ\'e \`a une
alg\`ebre de Lie $\mathscr{L}_k$ de champs de vecteurs lin\'eaires
; pour $m_k \in \P(E_k)$, g\'en\'erique, on a
$$\dim \mathscr{L}_k(m_k)=n_k-1$$ On a vu que
chaque $\mathscr{F}_k$ est d\'efini par des formes rationnelles
ferm\'ees $$\frac{\omega_k}{P_k}=\sum_{j=1}^{s_k} \lambda_{k,j}
\frac{dP_{k,j}}{P_{k,j}}+d\left(\frac{H_k}{P_{k,1}^{n_{k,1}}
\ldots P_{k,s_k}^{n_{k,s_k}}}\right)$$ avec les conditions
$\displaystyle \sum_{j=1}^{s_k} \lambda_{k,j} \deg P_{k,j}=0$ et
$\deg(H_k)=\displaystyle \sum_{\ell=1}^{s_k} n_{k,\ell} \deg
P_{k,\ell}$.

Posons $N=\displaystyle \sum_{k=1}^q n_k$ ; on consid\`ere sur
$E_1\oplus \ldots \oplus E_q \simeq \C^{N+q}$ le
feuilletage $\mathscr{F}_\mu$ associ\'e \`a la $1$-forme ferm\'ee
$$\frac{\Omega}{P}=\sum_{k=1}^q \mu_k \frac{\omega_k}{P_k}$$ o\`u
 $P=\displaystyle \prod_{k=1}^q P_k$ et les $\mu_k \in \C$
 sont non tous nuls. On remarque que $\mathscr{F}_\mu$ descend sur
$\C\P(N+q-1)$. En effet, si $R_k$ d\'esigne le
champ radial sur $E_k$, alors $R=\displaystyle \sum_{k=1}^q R_k$
est le champ radial sur $E_1\oplus \ldots \oplus E_q$ ; comme
chaque $\omega_k$ annule $R_j$ pour $j=1,\ldots,q$, la $1$-forme
$\Omega$ annule $R$.
\vs

\begin{pro}\label{cons} Si chaque $\mathscr{F}_k$
poss\`ede une sym\'etrie $X_k$, alors le
feuilletage $\mathscr{F}_\mu$ descend en un
$\mathscr{L}$-feuilletage sur $\P(E_1\oplus \ldots \oplus
E_q)$.
\vs
\end{pro}

\begin{proof}[D\'emonstration] Soit $\tilde{m}=(\tilde{m}_1,
\ldots,\tilde{m}_q) \in E_1
\oplus \ldots \oplus E_q$ un point g\'en\'erique. La
sous-alg\`ebre
$\mathscr{L}^0~=~\tilde{\mathscr{L}}_1~\oplus~\ldots~\oplus
~\tilde{\mathscr{L}}_q$ est form\'ee de champs tangents \`a
$\mathscr{F}_\mu$. On a
$$\dim \mathscr{L}^0(\tilde{m})=\sum_{k=1}^q \dim \tilde
{\mathscr{L}}_k(\tilde{m_k})=\sum_{k=1}^q n_k=N$$
Soit $\tilde{X}_k$ un relev\'e de $X_k$ et $P_k=i_{\tilde{X}_k}\omega_k$.
Consid\'erons maintenant un champ $\tilde{X}_\alpha$ de la forme
$$\tilde{X}_\alpha=\sum_{k=1}^q \alpha_k \tilde{X}_k$$ avec $\alpha=(\alpha_1,
\ldots,\alpha_q) \in \C^q$ et $<\alpha,\mu>=0$. On a
$$i_{\tilde{X}_\alpha}\frac{\Omega}{P}= \sum_{k=1}^q \mu_k
\alpha_k \frac{i_{\tilde{X}_\alpha} \omega_k}{P_k}= \sum_{k=1}^q
\mu_k \alpha_k$$ On remarque que $\tilde{X}_\alpha$ est
tangent \`a $\mathscr{F}_\mu$. Soit $\tilde{\mathscr{L}}$
l'alg\`ebre de Lie engendr\'ee par $\mathscr{L}^0$ et par les
champs $\tilde{X}_\alpha$ ;
l'alg\`ebre $\tilde{\mathscr{L}}$ est form\'ee de champs
lin\'eaires tangents \`a $\mathscr{F}_\mu$. De plus, en un point
$\tilde{m} \in E_1 \oplus \ldots \oplus E_q$ g\'en\'erique, on a
$$\dim \tilde{\mathscr{L}}(m)=\dim \mathscr{L}^0(\tilde{m})+q-1
\nonumber =\sum_{k=1}^q n_k+q-1= N+q-1$$ donc $\mathscr{F}_\mu$
induit un $\mathscr{L}$-feuilletage.
\end{proof}

\vs

Il se peut que $\tilde{\mathscr{L}}$ ne soit pas maximale alors
que tous les $\tilde{\mathscr{L}}_k$ le sont ; nous allons voir
que modulo une hypoth\`ese d'irr\'eductibilit\'e elle l'est.
\vs

\begin{defi}
Soit $\mathscr{F}$ un feuilletage sur $\C \P(n)$ d\'efini en
coordonn\'ees homog\`enes par $\omega=\displaystyle \sum_{k=0}^n
a_k dz_k$. On dit que $\mathscr{F}$ est irr\'eductible si $\omega$
n'annule pas de champ constant non trivial.
\end{defi}

\vs

Remarquons que si $\mathscr{F}$ n'est pas irr\'eductible, alors
$\mathscr{F}$ est le pull-back lin\'eaire d'un feuilletage
$\mathscr{F}'$ sur un $\C \P(\ell)$, $\ell < n$.
\vs

Soient $\mathscr{F}_\mu$, $\mathscr{F}_k$, $\mathscr{L}^0$ et
$\tilde{\mathscr{L}}$ comme dans la proposition \ref{cons} ; on a la :
\vs

\begin{pro}\label{consbis}
Si chaque $\mathscr{F}_k$ est
irr\'eductible, alors $\tilde{\mathscr{L}}$ est maximale.

De plus, si on note $g=\{\displaystyle \sum_{k=1}^q \alpha_k
\tilde{X}_k, \hspace{1mm} <\alpha,\mu>=0\}$, on a
$\tilde{\mathscr{L}}=g \oplus \mathscr{L}^0$ et $[g,\mathscr{L}^0]
\subset \mathscr{L}^0$.
\vs
\end{pro}

\begin{proof}[D\'emonstration] Soit $Y$ un champ de vecteurs lin\'eaire tangent \`a
$\tilde{\mathscr{F}}_\mu$. Montrons que $Y \in \tilde{\mathscr{L}}$ ; le
champ $Y$ s'\'ecrit $$Y_1+\ldots+Y_q$$ o\`u chaque $Y_i$ est un
champ lin\'eaire sur $E_1 \oplus \ldots \oplus E_q$ et parall\`ele
au facteur $E_i$. L'\'egalit\'e
$$i_{Y} \frac{\Omega}{P}=0$$ conduit \`a $$\sum_{k=1}^q
\frac{\mu_k}{P_k} i_Y \omega_k =0$$ et \`a
$$\sum_{k=1}^q \frac{\mu_k}{P_k} i_{Y_k} \omega_k
=0$$

\vs

Par suite, $i_{Y_k} \omega_k$ est divisible par $P_k$; comme $Y_k$
est lin\'eaire, on a $i_{Y_k} \omega_k=c_kP_k$ o\`u $c_k \in
\C$ et $\displaystyle \sum_{k=1}^q \mu_k c_k=0$. Puisque
$i_{\tilde{X}_k} \omega_k=P_k$, le champ $\tilde{Y}_k:=Y_k - c_k
\tilde{X}_k$ est annul\'e par $\omega_k$. Les composantes de
$\tilde{Y}_k$ pourraient d\'ependre de variables autres que celles
de $E_k$ ; montrons qu'en fait ce n'est pas le cas. Soient
$u=(u_1,\ldots,u_s)$ des coordonn\'ees de $E_k$ et
$v=(v_{s+1},\ldots,v_r)$ des coordonn\'ees de $E_1 \oplus \ldots
\oplus E_{k-1} \oplus E_{k+1} \oplus \ldots \oplus E_q$. Le champ
$\tilde{Y}_k$ s'\'ecrit
$$\sum_{j=1}^s a_j(u) \frac{\partial}{\partial u_j}+\sum_{j=1}^s
b_j(v) \frac{\partial} {\partial u_j}$$ o\`u les $a_j$ et $b_j$
sont lin\'eaires respectivement en $u$ et $v$. L'\'egalit\'e
$i_{\tilde{Y}_k} \omega_k=0$ conduit \`a $$i_Z \omega_i=0$$ o\`u
$Z=\displaystyle \sum_{j=1}^s b_j \frac{\partial}{\partial u_j}$.
Pour $v$ g\'en\'eral fix\'e ceci produit un champ de vecteurs
constant non nul, d\`es que l'un des $b_j$ est non trivial ; ce
champ est tangent \`a $\mathscr{F}_k$ ce qui est exclu puisque
le feuilletage $\mathscr{F}_k$ est suppos\'e irr\'eductible.
Finalement, $\tilde{Y}_k$ s'\'ecrit
$$\sum_{j=1}^s a_j(u) \frac{\partial}{\partial u_j}$$
Or l'alg\` ebre $\tilde{\mathscr{L}}_k$ est par hypoth\`ese
maximale donc $\tilde{Y}_k \in \tilde{\mathscr{L}_k}$. De sorte
que
$$Y=\sum_{k=1}^q Y_k%=\sum_{k=1}^q c_k\tilde{X}_k+\tilde{Y}_k
=\sum_{k=1}^q c_k \tilde{X}_k + \tilde{Y}$$ avec
$\tilde{Y}=\displaystyle \sum_{k=1}^q \tilde{Y}_k
\in\mathscr{L}^0$.

Pour terminer il suffit de montrer que
$[\tilde{X}_k,\mathscr{L}^0] \subset \mathscr{L}^0$ ; pour cela il
suffit de v\'erifier que les champs lin\'eaires
$[\tilde{X}_k,Z_k]$, o\`u $Z_k \in \tilde{\mathscr{L}}_k$,
annulent $\omega_k$ et sont donc dans $\tilde{\mathscr{L}}_k$.
Mais comme $\tilde{X}_k$ est une sym\'etrie lin\'eaire de
$\mathscr{F}_k$ et $i_{Z_k}\omega_k=0$, on a
$i_{[\tilde{X}_k,Z_k]}\omega_k=0$.
\end{proof}

\vs

Plus g\'en\'eralement on montre de la m\^eme mani\`ere le :
\vs

\begin{pro} Soit $\omega_k$ une $1$-forme
int\'egrable homog\`ene d\'efinissant un $\mathscr{L}$-feuilletage
sur $\C^{n_k}$ admettant une sym\'etrie $X_k$ et $P_k=i_{X_k}
\omega_k$. Si les $\mu_k \in \C \setminus \{0\}$ satisfont
$\displaystyle \sum_{k=1}^q \mu_k \frac{i_{R_k}
\omega_k}{P_k}~=~0$, alors la $1$-forme $\Omega=P_1\ldots P_q
\displaystyle \sum_{k=1}^q \mu_k \frac{\omega_k}{P_k}$ d\'efinit
un $\mathscr{L}$-feuilletage sur $\P(\C^{n_1} \oplus \ldots \oplus
\C^{n_q})$.
\vs
\end{pro}

\begin{rems}
(i) Si $i_{R_k} \omega_k \not =0$, on a une sym\'etrie automatique
$X_k=R_k$ et dans ce cas $i_{R_k} \frac{\omega_k}{P_k}=1$.\vs

(ii) Si $i_{R_k} \omega_k \not =0$, alors $Q_k=i_{R_k} \omega_k$
et $P_k$ sont des facteurs int\'egrants de $\omega_k$. Si la
$1$-forme $\omega_k$ n'a pas d'int\'egrale premi\`ere rationnelle,
alors $Q_k=\mbox{ cte }P_k$.
\end{rems}

\vs

\section{Exemples logarithmiques.}

\vs

Les exemples logarithmiques vus en 1.3. s'obtiennent trivialement
 via le principe pr\'ec\'edent. \vs

On consid\`ere sur $E_k=\C$ la $1$-forme $\omega_k=dz_k$,
le polyn\^ome $P_k=z_k$ et $\mathscr{L}_k=\{0\}$ ; le
\guillemotleft \hspace{1mm} feuilletage \guillemotright
\hspace{1mm} $\mathscr{F}_k$ en points est d\'ecrit par la
$1$-forme
$$\frac{\omega_k}{P_k}=\frac{dz_k}{z_k}$$

On construit sur $E_1 \oplus \ldots \oplus E_{n+1}$ le feuilletage
logarithmique d\'ecrit par la $1$-forme
$$z_1 \ldots z_{n+1} \sum_{k=1}^{n+1} \mu_k
\frac{dz_k}{z_k}$$ $\mu_k \not =0$, qui descend sur
$\P\C^{n+1} =\C \P(n)$ d\`es que
$\displaystyle \sum_{k=1}^{n+1} \mu_k=0$. On a ici $\mathscr{L}
\simeq \C^{n-1}$.

Les feuilles, qui sont les composantes connexes des niveaux de la
fonction multivalu\'ee $\displaystyle \sum_{k=1}^{n+1} \mu_k \log
z_k$, sont en g\'en\'eral denses pour $n \geq 3$. On remarque que
tout champ $\displaystyle \sum_{k=1}^{n+1} \lambda_k z_k
\frac{\partial}{\partial z_k} $ tel que $\displaystyle
\sum_{k=1}^{n+1} \lambda_k \mu_k \not = 0$ est une sym\'etrie.
\vs

\section{Exemples issus de formes quadratiques.}

\vs

Sur $\C^{n_j}$ on consid\`ere une forme quadratique $Q_j$
de rang maximal, i.e. quitte \`a composer par un isomorphisme on
peut supposer que $Q_j(z)=\displaystyle \sum_{k=1}^{n_j} z_{j,k}^2$.

Soit $\mathscr{F}_j$ le $\mathscr{L}$-feuilletage sur
$\C^{n_j}$ d\'efini par la $1$-forme $dQ_j$ et associ\'e
\`a l'alg\`ebre $so(Q_j) \simeq~so(n_j;\C)$. Le champ
radial $R_j$ de $\C^{n_j}$ est une sym\'etrie de
$\mathscr{F}_j$ : c'est l'identit\'e d'Euler. %$$i_{R_j}
%dQ_j=2i_{R_j}\sum_{\ell=1}^{n_j}
%z_\ell dz_\ell=2 \sum_{k=1}^{n_j} z_k \frac{\partial}{\partial
%z_k} \sum_{\ell=1}^{n_j} z_\ell dz_\ell=2Q_j.$$

On construit sur $\C^N=\C^{n_1} \times \ldots
\C^{n_q}$ le feuilletage d\'ecrit par $$Q_1 \ldots Q_q
\sum_{j=1}^q \mu_j \frac{dQ_j}{Q_j}$$ (qui descend en un
$\mathscr{L}$-feuilletage sur $\C \P(N-1)$ d\`es
que $\displaystyle \sum_{j=1}^q \mu_j=0$). La proposition
\ref{consbis}
assure que ce dernier est associ\'e \`a l'alg\`ebre
\begin{eqnarray}
\mathscr{L} &=& \{ \sum_{j=1}^q \kappa_j R_j, \hspace{1mm}
<\kappa,\mu>=0\} \oplus so(Q_1) \oplus \ldots \oplus
so(Q_q) \nonumber \\
&\simeq& \C^{q-1} \oplus so(n_1;\C) \oplus \ldots
\oplus so(n_q;\C) \hspace{8mm} \mbox{(somme directe
d'alg\`ebres de Lie)}\nonumber
\end{eqnarray}

En consid\'erant $\C^N$ comme une carte affine de $\C\P(N)$ on
constate que ce feuilletage s'\'etend \`a $\C \P(N)$ et que ce
prolongement est associ\'e \`a l'alg\`ebre $\mathscr{L}$.

Si $q \geq 3$, les feuilles sont en g\'en\'eral denses. Si $\mu_j
\in \mathbb{Z}$, on obtient des exemples avec int\'egrales premi\`eres
rationnelles.
\vs

On peut aussi coller des feuilletages logarithmiques avec des
feuilletages associ\'es \`a des formes quadratiques ;
consid\'erons par exemple la $1$-forme
$$\omega=\frac{d(z_1^2+z_2^2+z_3^2)}{z_1^2+z_2^2+z_3^2}-\kappa
\frac{dz_0}{z_0}$$ On note $R_1=z_1\frac{\partial}{\partial z_1}+
z_2\frac{\partial}{\partial z_2}+z_3\frac{\partial}{\partial
z_3}$. Soit $Z$ un champ annulant $\omega$ ; visiblement $Z$
s'\'ecrit $\mu z_0\frac{\partial}{\partial z_0}+Z'$ avec $Z' \in
<\frac{\partial}{\partial z_1}, \hspace{1mm} \frac{\partial}
{\partial z_2}, \hspace{1mm} \frac{\partial} {\partial z_3}>$.
Alors $Z-\mu(z_0\frac{\partial} {\partial z_0}+\frac{\kappa}{2}
R_1)$ annule $\omega$ et n'a pas de composante en $\frac{\partial}
{\partial z_0}$ ; en particulier il annule $$\frac{d(z_1^2+z_2^2+
z_3^2)}{z_1^2+z_2^2+ z_3^2}$$ La $1$-forme $\omega$ d\'ecrit donc
un $\mathscr{L}$-feuilletage de degr\'e $2$ et de codimension $1$
sur $\C^4$ associ\'e \`a $\mathscr{L}=so(3;\C)
\oplus \C$. Ce feuilletage s'\'etend en un
$\mathscr{L}$-feuilletage sur $\C \P(4)$ d\'ecrit
par la $1$-forme
$$\frac{d(z_1^2+z_2^2+z_3^2)}{z_1^2+z_2^2+z_3^2}-\kappa
\frac{dz_0}{z_0}+(\kappa -2)\frac{dz_4}{z_4} $$

\vs

\section{Exemples issus d'actions classiques.}

\vs

On pr\'esente quelques exemples d'actions classiques qui peuvent
\^etre des bases pour le principe de construction d'exemples ; la
plupart de ces exemples se trouvent dans les livres (\cite{[Ha]}, \cite{[Do]},
\cite{[Po-Vi]}). Il est plus commode ici de pr\'esenter les groupes
eux-m\^emes que les sous-alg\`ebres de Lie correspondantes.
\vs

\subsection{Action de $\P\mbox{GL}(2,\C)$ sur
$\C \P(3)$ et construction d'un
$\mathscr{L}$-feuilletage sur $\C\P(4)$.}

\vs

On note V$_4$ l'ensemble des polyn\^omes homog\`enes de degr\'e
trois en deux variables et de mani\`ere plus g\'en\'erale V$_k$
l'ensemble des polyn\^omes homog\`enes de degr\'e $k-1$ en deux
variables. Le groupe $\mbox{GL}(2,\C)$ agit sur V$_k$
$$g.P=P \circ g^{-1}$$ o\`u $P \in \mbox{V}_k$
 et $g \in \mbox{GL}(2,\C)$. Soient $\lambda$, $\mu
\in \C^*$, on a $$(\lambda g).(\mu P)=\mu
P(\frac{1}{\lambda}g^{-1})=\frac{\mu}{\lambda^k}Pg^{-1}=\frac{\mu}
{\lambda^k}g.P$$ L'action de $\mbox{GL}(2,\C)$ sur $V_k$
 induit donc une action de $\P
\mbox{GL}(2,\C)$ sur $\P(V_k)\simeq\C \P(k-1)$. Or un \'el\'ement
de $\C \P(3) ~\simeq~\P(V_4)$ est la donn\'ee d'un polyn\^ome $P$
de degr\'e $3$ \`a multiplication par un scalaire pr\`es,
autrement dit la donn\'ee de trois droites vectorielles
$\Delta_1$, $\Delta_2$ et $\Delta_3$ dans $\C^2$ non
num\'erot\'ees (les z\'eros de $P$), soit encore la donn\'ee de
trois points de $\C \P(1)$. Ainsi l'action de $\P\mbox{GL}(2,\C)$
sur $\C \P(3)$ correspond \`a l'action de $\P\mbox{GL}(2,\C)\simeq
\mbox{Aut }\C \P(1)$ sur les triplets de points de $\C \P(1)$.

D\'ecrivons les orbites de cette action.

On compte trois orbites : l'orbite des triplets de
points distincts, l'orbite des triplets de points dont deux
seulement sont confondus et l'orbite des points triples.

On commence tout d'abord par s'int\'eresser au sous-ensemble de
$\C \P(3)$ correspondant aux points triples ; un
point triple c'est la donn\'ee d'un polyn\^ome homog\`ene $P$ de
degr\'e $3$ qui est un cube : $$P(z_0,z_1)=(\beta z_0-\alpha
z_1)^3$$ c'est encore la donn\'ee du point $[\alpha:\beta] \in \C
\P(1)$. On peut donc param\'etrer les points triples par
l'application $$(\alpha, \beta) \mapsto (\beta z_0-\alpha
z_1)^3=\beta^3 z_0^3-3\beta^2\alpha z_0^2 z_1+ 3\beta\alpha^2 z_0
z_1^2-\alpha^3z_1^3$$ qui induit l'application
$$[\alpha:\beta] \mapsto [\beta^3:-3\beta^2 \alpha:3\beta
\alpha^2:-\alpha^3]$$ dont l'image est \guillemotleft \hspace{1mm}
la \guillemotright \hspace{1mm} cubique gauche de $\C \P(3)$
not\'ee $\Gamma$.

D\'ecrivons maintenant l'orbite des triplets de points form\'es
d'un point double et d'un point simple. On consid\`ere le
polyn\^ome $z_0^3$, c'est le point $m=[1:0:0:~0] \in \Gamma$. La
tangente \`a $\Gamma$ en $m$ est donn\'ee par $$t \mapsto
(1,0,0,0) + t\frac{\partial}{\partial \alpha}(1,-3\alpha,
3\alpha^2,-\alpha^3)_{|\alpha=0}=(1,-3t,0,0)$$ qui correspond au
polyn\^ome $z_0^2(z_0-3tz_1)$ ; pour $t\not = 0$, ce polyn\^ome
est le mod\`ele \`a conjugaison pr\`es de tout polyn\^ome de
degr\'e trois ayant un terme carr\'e.
\vs

On remarque que si $\varphi$ est un \'el\'ement de Aut $\C \P(3)$
laissant $\Gamma$ invariante et si $m \in \Gamma$, l'image par
$\varphi$ de la tangente \`a $\Gamma$ en $m$ est la droite
tangente \`a $\Gamma$ en $\varphi(m)$.
\vs

\hspace{-6.5mm} Ainsi l'orbite des polyn\^omes, dont le lieu des
z\'eros est une droite double et une droite simple, est l'ensemble
des tangentes $\Sigma$ \`a la cubique $\Gamma$ priv\'e de
$\Gamma$. La surface $\Sigma \setminus \Gamma$ est invariante sous
l'action de $\P\mbox{GL}(2,\C)$, elle est param\'etr\'ee par
l'application
$$(\alpha,t) \mapsto
(-3\alpha,3\alpha^2,-\alpha^3)+t(-3,6\alpha,-3\alpha^2)$$ et est
de degr\'e $4$ (\cite{[Ha]}). On v\'erifie facilement que dans une
carte affine ad-hoc, $\Sigma$ est isomorphe \`a un cylindre
cuspidal $(T,s) \mapsto (T,s^2,s^3)$.
\vs

L'action de $\P\mbox{GL}(2,\C)$ sur $\C \P(3)$ produit une unique
orbite g\'en\'erique $\C \P(3) \setminus \Sigma$ et par
cons\'equent ne produit pas de feuilletage. Toutefois elle nous
sera utile plus loin. Par contre l'action de $\mbox{SL}(2,\C)$ sur
$V_4$ admet un invariant $\Delta$ (le discriminant). On obtient
ainsi un $\mathscr{L}$-feuilletage sur $\C^4$ dont une int\'egrale
premi\`ere de ce feuilletage est
\begin{eqnarray}
\Delta(\alpha_0,\alpha_1,\alpha_2,\alpha_3)&:=&\Delta(\alpha_0
z_0^3+\alpha_1z_0^2z_1+\alpha_2z_0z_1^2+\alpha_3z_1^3) \nonumber
\\
&=& \alpha_1^2\alpha_2^2-4\alpha_0
\alpha_2^3-4\alpha_1^3\alpha_3-27\alpha_0^2\alpha_3^2+18 \alpha_0
\alpha_1 \alpha_2 \alpha_3. \nonumber
\end{eqnarray}
On le prolonge de fa\c{c}on naturelle \`a $\C\P(4)$. Ce
feuilletage est associ\'e \`a une sous-alg\`ebre $\mathscr{L}
\subset \chi(\C \P(4))$ isomorphe \`a $s\ell(2,\C)$ ; il est de
degr\'e $3$.
\vs

\subsection{Action de $\P\mbox{GL}(3,\C)$ sur
$\C \P(5)$ et construction d'un
$\mathscr{L}$-feuilletage sur $\C\P(6)$.}

Le groupe $\mbox{GL}(3,\C)$ agit sur l'espace des matrices
sym\'etriques $\mathscr{S}ym(3,\C) \simeq \C^6$
comme suit
$$g.P=(^tg^{-1}) P g^{-1}$$ o\`u $P \in \mathscr{S}ym(3,
\C)$ et $g \in \mbox{GL}(3,\C)$. Cette action correspond \`a
 l'action naturelle de $\mbox{GL}(3,\C)$ sur les formes
 quadratiques. Soient $\lambda$, $\mu
\in \C^*$ ; on a $$(\lambda g).(\mu P)=\frac{\mu}{\lambda^2}
(^tg^{-1}) P g^{-1}$$ Ainsi $\P \mbox{GL}(3,\C)$ agit sur
$\P\mathscr{S}ym(3,\C)~ \simeq~\C
\P(5)$.

D\'ecrivons les orbites de cette action. Comme $\C \P(5)
\simeq \{\mbox{coniques de } \C\P(2) \}$, l'\'etude des
orbites se ram\`ene \`a celles des droites doubles, des paires de droites
distinctes et des coniques lisses. L'orbite des
droites doubles est l'image de l'application de V\'eron\`ese
\begin{eqnarray}
\C \P(2) \hspace{3mm} &\to& \hspace{17mm}
\C \P(5) \nonumber
\\
\hspace{1mm} [z_0:z_1:z_2] &\mapsto&
[z_0^2:z_1^2:z_2^2:z_0z_1:z_0z_2:z_1z_2] \nonumber
\end{eqnarray}
appel\'ee surface de V\'eron\`ese ; l'orbite des paires de droites
distinctes est la vari\'et\'e s\'ecante \`a la surface de
V\'eron\`ese (c'est l'hypersurface constitu\'ee des droites qui
coupent cette surface) priv\'ee de celle-ci. Enfin l'orbite des
coniques lisses est l'ouvert dense de $\C \P(5)$ :
$$\C \P(5) \setminus \{\mbox{vari\'et\'e s\'ecante
\`a la surface de V\'eron\`ese}\}.$$ Pour obtenir un feuilletage
de codimension $1$ on consid\`ere comme pr\'ec\'edemment l'action de
$\mbox{SL}(3,\C)$ sur $\mathscr{S}ym(3,\C)$ : deux
matrices $A$, $B \in \mathscr{S}ym(3,\C)$ sont
conjugu\'ees si et seulement si il existe $P \in
\mbox{SL}(3,\C)$ telle que $^t PAP =B$ c'est-\`a-dire si et
seulement si $\det A=\det B$. Donc
$$\delta(z)
 = \det \left(%
\begin{array}{ccc}
  z_1 & z_2 & z_3 \\
  z_2 & z_4 & z_5 \\
  z_3 & z_5 & z_6 \\
\end{array}%
\right)=z_1z_4z_6-z_1z_5^2-z_2^2z_6+2z_2z_3z_5-z_3^2z_4
$$ est un invariant de degr\'e $3$ de
l'action de $\mbox{SL}(3,\C)$ sur
$\mathscr{S}ym(3,\C) \simeq \C^6$. En prolongeant
$\delta$ \`a $\C \P(6)$, on d\'efinit un
$\mathscr{L}$-feuilletage de degr\'e $2$ et de codimension $1$ sur
$\C \P(6)$ dont une int\'egrale premi\`ere est
$$\frac{z_1z_4z_6-z_1z_5^2-z_2^2z_6+2z_2z_3z_5-z_3^2z_4}{z_0^3}$$
Ce $\mathscr{L}$-feuilletage est associ\'e \`a l'alg\`ebre de Lie
$s\ell(3,\C)$ de dimension $8$ ; en effet, l'alg\`ebre
$\mathscr{L}ie \hspace{1mm}
\P\mbox{GL}(3,\C)=s\ell(3,\C)$ est simple
donc le morphisme $$\mathscr{L}ie \hspace{1mm}
\P\mbox{GL}(3,\C)=s\ell(3,\C) \to
\mathscr{L}ie \hspace{1mm} \mbox{Aut } \C
\P(5)=\chi(\C \P(5))$$ est injectif et
$\chi(\C \P(5)) \supset \mathscr{L} \simeq
s\ell(3,\C)$.

Plus g\'en\'eralement l'action de $\mbox{SL}(n,\C)$ sur
$\mathscr{S}ym(n;\C)$ produit un $\mathscr{L}$-feuilletage
de degr\'e $n-1$ et de codimension $1$ sur $\C
\P(\frac{n(n+1)}{2})$ associ\'e \`a l'alg\`ebre
$s\ell(n,\C)$ dont une int\'egrale premi\`ere est
$$\frac{\det A}{z_0^n}$$ avec $A \in \mathscr{S}ym(n,\C)$.
\vs

\subsection{Action de $\P\mbox{GL}(3,\C)$ sur
$\C \P(9)$ et construction d'un
$\mathscr{L}$-feuilletage sur $\C\P(9)$.}

Soit W$_4$ l'ensemble des polyn\^omes homog\`enes de degr\'e trois
en trois variables. Le groupe $\mbox{GL}(3,\C)$ agit sur
W$_4$
$$g.P=P \circ g^{-1}$$ o\`u $P \in \mbox{W}_4$
 et  $g \in \mbox{GL}(3,\C)$. Comme pr\'ec\'edemment
cette action de $\mbox{GL}(3,\C)$ sur W$_4$ se laisse
projectiviser en une action de $\P\mbox{GL}(3,\C)$
sur $\P(\mbox{W}_4) \simeq \C \P(9)$.

Il y a huit orbites sp\'eciales : celle des droites triples, de l'union d'une
droite double et d'une droite, de l'union de trois droites en position
g\'en\'erale, de l'union de trois droites concourrantes, d'une conique et d'une
droite en position g\'en\'erale, d'une conique et d'une droite tangente, des
cubiques cuspidales et des cubiques \`a point double.

Les orbites g\'en\'eriques sont celles des courbes elliptiques. Une courbe
elliptique $E_\lambda$ s'\'ecrit \`a composition par un automorphisme pr\`es
$$z_1^2z_2=z_0(z_0-z_2)(z_0-\lambda z_2).$$ On constate que deux
courbes elliptiques
$E_\lambda$ et $E_{\lambda'}$ sont \'equivalentes si et seulement
si les ensembles $$\Lambda=\{\lambda, \frac{1}{\lambda},
1-\lambda, \frac{1}{1-\lambda}, \frac{\lambda}{1-\lambda},
\frac{\lambda}{\lambda-1}\}$$ (orbite de $\lambda$ sous l'action
de $z \mapsto \frac{1}{z}$ et $z \mapsto 1-z$) et
$$\Lambda'=\{\lambda', \frac{1}{\lambda'}, 1-\lambda',
\frac{1}{1-\lambda'}, \frac{\lambda'}{1-\lambda'},
\frac{\lambda'}{\lambda'-1}\}$$ co\"{\i}ncident ; autrement dit si
et seulement si $j(\lambda)=j(\lambda')$  o\`u $j$ d\'esigne la
fonction $$\lambda \mapsto 256 \frac{(\lambda^2-\lambda+1)^3}
{\lambda^2(\lambda-1)^2}$$ En effet, $j(\frac{1}{\lambda})=j(\lambda)$ et
 $j(1-\lambda)=j(\lambda)$.  La fonction $j$ est
donc bien d\'efinie sur les quadruplets de points de $\C
\P(1)$ et induit une fonction rationnelle $J \hspace{1mm}
\colon \hspace{1mm} \C \P(9) \dashrightarrow \C
\P(1)$ qui est la composition de $j$ avec $\lambda$. Ceci
est reli\'e au fait que l'on peut voir la donn\'ee d'une
courbe elliptique comme celle d'un rev\^etement au dessus de $\C
\P(1)$ ramifi\'e en quatre points, ce que nous
pr\'ecisons dans le paragraphe suivant. L'adh\'erence de l'orbite
d'une courbe elliptique est une fibre de $J$. C'est une hypersurface de degr\'e
$6$ sauf pour les deux valeurs critiques de $j$ : $0$ et $1728$ (les points
critiques de $j$, qui sont $2$, $\frac{1}{2}$, $1$ de multiplicit\'e $1$ et
$-j$, $-j^2$ de multiplicit\'e $2$, correspondent aux r\'eseaux
sp\'eciaux $\mathbb{Z} \oplus \mathrm{i}\mathbb{Z}$ et $\mathbb{Z}
\oplus j \mathbb{Z}$).

Les orbites des cubiques singuli\`eres satisfont des conditions
d'incidence ; par exemple l'adh\'erence de l'orbite d'une courbe elliptique
contient le lieu des cubiques cuspidales.
\vs

L'action de $\P\mbox{GL}(3,\C)$ sur $\C \P(9)$ produit donc un
$\mathscr{L}$-feuilletage de degr\'e $8$ et de codimension $1$ sur
$\C \P(9)$ associ\'e \`a l'alg\`ebre $\mathscr{L} \simeq
s\ell(3,\C)$. Le calcul du degr\'e n'est pas si facile et a
\'et\'e effectu\'e par J. V. Pereira (\cite{[Pe]}) en utilisant une
m\'ethode due \`a Darboux. Les descriptions des feuilles et du
lieu singulier s'obtiennent \`a partir des consid\'erations
pr\'ec\'edentes sur les orbites.
\vs

On trouvera la formule explicite de la fonction $J$ dans le livre
de Dolgachev \cite{[Do]} page 160. \vs

\subsection{Action de $\P\mbox{GL}(2,\C)$ sur
$\C \P(4)$.}

Comme on l'a dit pr\'ec\'edemment $\P\mbox{GL}(2,\C)$ agit sur
$\P(V_5)\simeq \C\P(4)$. On
retrouve et on pr\'ecise les r\'esultats du paragraphe
ci-dessus. \vs

On compte cinq cas de figures diff\'erents : les quadruplets
de points tous distincts, les paires de points doubles, les
quadruplets form\'es d'un point triple et d'un point double, les
quadruplets form\'es d'une paire de points doubles et d'une paire
de points simples et enfin les points quadruples.

Les quatre premi\`eres configurations produisent quatre orbites.
En effet, une homographie permet d'\'echanger deux triplets de
points donn\'es mais envoie quatre points distincts $z_i$, $i=1,
\ldots, 4$, sur $0$, $1$, $\infty$ et le birapport $\lambda$ des
$z_i$. Toutefois le calcul du birapport tient compte de l'ordre
dans lequel on s'est donn\'e les $z_i$ ; permuter ces quatre
points a pour effet de changer le birapport $\lambda$ en
$\frac{1}{\lambda}$, $1-\lambda$, $\frac{1}{1-\lambda}$,
$\frac{\lambda}{1-\lambda}$, $\frac{\lambda}{\lambda-1}$ ou
$\frac{\lambda}{\lambda-1}$. Ainsi deux quadruplets de points
distincts non num\'erot\'es peuvent \^etre envoy\'es l'un sur
l'autre par un \'el\'ement de $\P\mbox{GL}(2;\C)$ si et seulement
si les ensembles
$$\Lambda=\{\lambda, \frac{1}{\lambda}, 1-\lambda,
\frac{1}{1-\lambda}, \frac{\lambda}{1-\lambda},
\frac{\lambda}{\lambda-1}\}$$ et
$$\Lambda'\{\lambda', \frac{1}{\lambda'},1-\lambda', \frac{1}{1-\lambda'},
\frac{\lambda'}{1-\lambda'},\frac{\lambda'}{\lambda'-1}\}$$ co\"{\i}ncident.
Or $\Lambda$ et $\Lambda'$ co\"{\i}ncident si et seulement si
$j(\lambda)=j(\lambda')$. On retrouve ainsi la situation du 2.4.3.

D\'ecrivons par exemple l'orbite d'un point quadruple : se donner
un point quadruple revient \`a se donner un polyn\^ome homog\`ene
de degr\'e $4$ de la forme $$P(z_0,z_1)=(\beta z_0-\alpha z_1)^4$$
ce qui d\'etermine un point $[\alpha:\beta] \in
\C \P(1)$. On peut donc param\'etrer les points
quadruples par l'application
$$(\alpha, \beta) \mapsto (\beta z_0-\alpha z_1)^4=\beta^4
z_0^4-4\beta^3\alpha z_0^3z_1+6\beta^2 \alpha^2 z_0^2 z_1^2-4\beta
\alpha^3 z_0z_1^3+\alpha^4 \beta^4$$ qui induit l'application de
$\C\P(1)$ dans $\C\P(4)$ :
$$[\alpha : \beta] \mapsto [\beta^4:-4\beta^3\alpha:6\beta^2
\alpha^2:-4\beta \alpha^3 :\alpha^4 \beta^4]$$ Son image est une
courbe de degr\'e $4$ appel\'ee courbe de V\'eron\`ese.

L'action de $\P\mbox{GL}(2,\C)$ sur $\C
\P(4)$ induit un $\mathscr{L}$-feuilletage associ\'e \`a
$\mathscr{L}~\simeq~s\ell(2,\C)$. On cherche une
expression de l'int\'egrale premi\`ere de ce
$\mathscr{L}$-feuilletage ; on profite de l'existence
d'invariants de l'action de $\mbox{GL}(2,\C)$ sur $V_5$.
Il existe deux tels invariants homog\`enes de degr\'e
respectivement $2$ et $3$ (\cite{[Ha]}):
$$P(\alpha_0 z_0^4+\alpha_1 z_0^3z_1+\alpha_2 z_0^2z_1^2+\alpha_3
z_0z_1^3+\alpha_4z_1^4) =\alpha_0 \alpha_4-\frac{1}{4} \alpha_1
\alpha_3 +\frac{1}{12} \alpha_2^2 \hspace{1mm}$$ et
\begin{eqnarray}
& & H(\alpha_0 z_0^4+\alpha_1 z_0^3z_1+\alpha_2
z_0^2z_1^2+\alpha_3 z_0z_1^3+\alpha_4z_1^4) =\det \left(%
\begin{array}{ccc}
  \alpha_0 & \frac{\alpha_1}{4} & \frac{\alpha_2}{6} \\
\frac{\alpha_1}{4} & \frac{\alpha_2}{6} & \frac{\alpha_3}{4} \\
 \frac{\alpha_2}{6} & \frac{\alpha_3}{4} & \alpha_4 \\
\end{array}%
\right) \nonumber \\
& & \hspace{22mm} = \frac{\alpha_0 \alpha_2
\alpha_4}{6}-\frac{\alpha_0
\alpha_3^2}{16}-\frac{\alpha_1^2\alpha_2}{16}+\frac{\alpha_1
\alpha_2\alpha_3}{96}+\frac{\alpha_1 \alpha_2 \alpha_3}{96}
-\frac{\alpha_2^3}{216} \nonumber
\end{eqnarray}

\vs

\hspace{-6.5mm} Ainsi on trouve que $\frac{P^3}{\Delta}$ o\`u
$\Delta=2^8(P^3-27H^2)$ est une int\'egrale premi\`ere du
$\mathscr{L}$-feuilletage induit par l'action de
$\P\mbox{GL}(2,\C)$ sur $\C
\P(4)$.
 On en d\'eduit que $\frac{P^3}{\Delta}$ est une fonction
lin\'eaire de l'invariant $j$ : $\frac{P^3}{\Delta}=\frac{1}{1728}j$.
\vs

\subsection{Action de $\P\mbox{GL}(3,\C)$ sur
$\C \P(7)$.}

On consid\`ere l'espace V des formes diff\'erentielles
$$\omega=\sum_{k=0}^2 A_k dz_k$$ o\`u les $A_k$ sont des
polyn\^omes homog\`enes de degr\'e $2$ tels que $\displaystyle
\sum_{k=0}^2 A_k z_k=0$. Cet espace est de dimension $8$ ; un
\'el\'ement g\'en\'erique de $\mbox{V}$, i.e. un
\'el\'ement tel que $\gcd(A_0,A_1,A_2)=1$, repr\'esente un
feuilletage de degr\'e $1$ sur $\C \P(2)$ (cf.
chapitre 3)

Le groupe $\mbox{GL}(3,\C)$ agit sur V de fa\c{c}on
naturelle
$$g.\omega=g^* \omega$$ o\`u $\omega
\in V$ et $g \in \mbox{GL}(3,\C)$. Soient
$\lambda$, $\mu \in \C^*$, on a $$(\lambda g).(\mu
\omega)=\frac{\mu^2}{\lambda}g.\omega \hspace{1mm} $$ et l'action
de $\mbox{GL}(3,\C)$ sur V induit donc une action de
$\P \mbox{GL}(3,\C)$ sur $\P(\mbox{V})
\simeq \C \P(7)$.

Donnons une br\`eve description de cette action.

\hspace{-7.4mm} Commen\c{c}ons par nous int\'eresser aux \'el\'ements
$\omega_{n,g}$ non g\'en\'eriques de $V$, i.e. aux \'el\'ements
tels que $\gcd(A_0,A_1,A_2) \not = 1$ ; la $1$-forme
$\omega_{n,g}$ s'\'ecrit $$L(\sum_{k=0}^2 L_k dz_k)$$ avec
$\displaystyle \sum_{k=0}^2 L_k z_k=0$. On remarque que
$\displaystyle \sum_{k=0}^2 L_k dz_k$ est lin\'eairement conjugu\'e \`a
$$z_1dz_2-z_2dz_1$$ Par suite, les $1$-formes $\omega_{n,g}$ sont
lin\'eairement conjugu\'ees \`a des formes du type
$$L(z_1dz_2-z_2dz_1)$$ o\`u $L$ est lin\'eaire. Si $L$ d\'epend seulement de
$(z_1,z_2)$, on se ram\`ene au mod\`ele
$$\omega_1=z_1(z_1dz_2-z_2dz_1)$$ et sinon \`a
$$\omega_0=z_0(z_1dz_2-z_2dz_1)$$

\vs

Se donner un \'el\'ement de $\mathscr{O}rb [\omega_1]$, o\`u
$\mathscr{O}rb [\omega_1]$ d\'esigne l'orbite de $[\omega_1]$,
revient \`a se donner une droite (la droite $z_1=0$) du pinceau de
droites associ\'e au feuilletage de degr\'e $0$ d\'ecrit par
$z_1dz_2-z_2dz_1$ ; ainsi $\dim \mathscr{O}rb [\omega_1]=3$.

Un \'el\'ement de $\mathscr{O}rb [\omega_0]$ correspond au point
de base du pinceau de droites associ\'e au feuilletage de degr\'e
$0$ d\'ecrit par $z_1dz_2-z_2dz_1$ et \`a une droite g\'en\'erique ; ainsi
$\dim \mathscr{O}rb [\omega_0]=4$.

On constate que $\mathscr{O}rb [\omega_1] \subset
\overline{\mathscr{O}rb [\omega_0]}$ et que
$\overline{\mathscr{O}rb [\omega_0]}$ est biholomorphe \`a
$\C \P(2) \times \check{\C \P(2)}
\simeq~\C \P(2) \times~\C \P(2)$,
o\`u $\check{\C \P(2)}$ d\'esigne l'espace
projectif dual. Dans $\C \P(2) \times
\check{\C \P(2)}$, l'orbite de $[\omega_1]$
s'identifie \`a la vari\'et\'e d'incidence
$$\{(x,\mathscr{D}) \in \C \P(2) \times
\check{\C \P(2)}, \hspace{1mm} x \in
\mathscr{D}\}$$ i.e. \`a
$$\{([z_0:z_1:z_2],[\alpha_0:\alpha_1:\alpha_2]) \in \C \P(2)
\times \check{\C \P(2)}, \hspace{1mm}
z_0\alpha_0+z_1\alpha_1+z_2\alpha_2=0\}$$ \vs Soit $\omega$
un \'el\'ement g\'en\'erique de V ; le feuilletage $\mathscr{F}$
d\'efini par $\omega$ a alors $3$ points singuliers. Il laisse
donc trois droites distinctes invariantes (les droites joignant
les trois points singuliers) que l'on peut supposer \^etre, quitte
\`a faire une transformation lin\'eaire, les droites d'\'equation
$z_0=0$, $z_1=0$, $z_2=0$. Dans la carte affine $z_0=1$, la
$1$-forme $\omega$ s'\'ecrit
$$B_1(1,z_1,z_2)dz_1+B_2(1,z_1,z_2)dz_2+B_3(1,z_1,z_2)(z_1dz_2-z_2dz_1)$$
o\`u $B_1$ et $B_2$ d\'esignent des fonctions affines et $B_3$ une
fonction lin\'eaire. La droite d'\'equation $z_1=0$ \'etant
invariante par le feuilletage, la fonction affine $B_2$ est
divisible par $z_1$ donc s'\'ecrit cte.$z_1$ ; de m\^eme, on
obtient que $B_1$ s'\'ecrit cte.$z_2$. Enfin l'invariance de la
droite \`a l'infini par le feuilletage implique que $B_3=0$.
Finalement $\omega$ s'\'ecrit $$\omega_\lambda=z_2dz_1-\lambda
z_1dz_2$$ En consid\'erant
\begin{eqnarray}
\phi \hspace{1mm} \colon \hspace{3mm} \C^2 \hspace{3mm}
&\to&
\C^2 \nonumber \\
(z_1,z_2) &\mapsto& (z_2,z_1) \nonumber
\end{eqnarray}
on constate que $\mathscr{O}rb [\omega_\lambda]=\mathscr{O}rb
[\omega_{\frac{1}{\lambda}}]$ ; dit autrement les feuilletages
$\mathscr{F}_\lambda$ et $\mathscr{F}_{\frac{1}{\lambda}}$ d\'ecrits
respectivement par $\omega_\lambda$ et $\omega_{\frac{1}{\lambda}}$ sont
conjugu\'es. En permutant les droites invariantes
entre elles, on obtient : les deux feuilletages
$\mathscr{F}_\lambda$ et $\mathscr{F}_{\lambda'}$ sont conjugu\'es
si et seulement si $\Lambda=\Lambda'$ autrement dit si et
seulement si $j(\lambda)=j(\lambda')$ avec les notations
habituelles.

Comme les orbites g\'en\'eriques sont de dimension $6$, l'action
de $\P \mbox{GL}(3,\C)$ sur $\C \P(7)$ produit un
$\mathscr{L}$-feuilletage de degr\'e $3$ et de codimension $1$ sur
$\C \P(7)$ associ\'e \`a $s\ell(3,\C)$. Le calcul du degr\'e a
aussi \'et\'e effectu\'e par J. V. Pereira en utilisant encore la
m\'ethode de Darboux. On constate, mais ce n'est pas nouveau,
l'ubiquit\'e de la fonction $j$.
\vs

\subsection{Action de $\mbox{SL}(2n,\C)$ sur les matrices
antisym\'etriques en dimension paire
$\mathscr{A}sym(2n,\C)$.}

Soit $M \in \mathscr{A}sym(2n,\C)$, alors $M$ s'\'ecrit $$\left(%
\begin{array}{ccc}
  0 & & m_{i,j} \\
   & \ddots &  \\
  -m_{i,j} &  & 0 \\
\end{array}%
\right)$$ et $\det M$ est un polyn\^ome de degr\'e $2n$ en
$m_{i,j}$ qui s'av\`ere \^etre un carr\'e. Il se trouve que
$\sqrt{\det}$ est l'unique invariant de l'action de
$\mbox{SL}(2n,\C)$ sur $\mathscr{A}sym (2n,\C)$ ce qui produit un
$\mathscr{L}$-feuilletage de degr\'e $n-1$ sur $\C^{n(2n-1)}$ dont
les feuilles sont les niveaux g\'en\'eriques de $\sqrt{det}$. On
peut prolonger ce feuilletage \`a $\C \P(n(2n-1))$. Pour $n=3$ on
construit ainsi un $\mathscr{L}$-feuilletage de degr\'e $2$ sur
$\C \P(15)$.
\vs

\subsection{Action de $\mbox{GL}(2,\C)$ sur
$\mbox{V}_5$.}

L'action naturelle de $\mbox{GL}(2,\C)$ sur $\mbox{V}_5$ poss\`ede
l'invariant homog\`ene de degr\'e $2$ suivant (\cite{[Ha]})
\begin{eqnarray}
& & H(\alpha_0 z_0^4+\alpha_1 z_0^3z_1+\alpha_2
z_0^2z_1^2+\alpha_3
z_0z_1^3+\alpha_4z_1^4)=\det \left(%
\begin{array}{ccc}
  \alpha_0 & \frac{\alpha_1}{4} & \frac{\alpha_2}{6} \\
\frac{\alpha_1}{4} & \frac{\alpha_2}{6} & \frac{\alpha_3}{4} \\
 \frac{\alpha_2}{6} & \frac{\alpha_3}{4} & \alpha_4 \\
\end{array}%
\right) \nonumber \\
& & \hspace{26mm}=\frac{\alpha_0 \alpha_2
\alpha_4}{6}-\frac{\alpha_0 \alpha_3^2}
{16}-\frac{\alpha_1^2\alpha_4}{16}+\frac{\alpha_1\alpha_2^2}{144}
+\frac{\alpha_1\alpha_2\alpha_3}{216}-\frac{\alpha_2^3}{216}
\nonumber
\end{eqnarray} Cette action induit un $\mathscr{L}$-feuilletage de degr\'e $2$
et de codimension $1$ sur $\C^5$ qui se prolonge \`a
$\C \P(5)$.
\vs

\subsection{Feuilletage exceptionnel.}

Soit $\Gamma$ la cubique gauche de $\C^3$ param\'etr\'ee par $$t
\mapsto (t,\frac{t^2}{2},\frac{t^3}{6})$$ On note
$\mbox{Aut}_\Gamma \C^3$ les automorphismes affines de $\C^3$
laissant $\Gamma$ invariante. En fait $\mbox{Aut}_\Gamma \C^3$ est
isomorphe au groupe affine de la droite.

L'orbite d'un point $m$ sous l'action de $\mbox{Aut}_\Gamma \C^3$
est une vari\'et\'e d'adh\'erence alg\'ebrique. Cette action
induit donc un $\mathscr{L}$-feuilletage, not\'e
$\mathscr{F}_\Gamma$, sur $\C^3$ dont le lieu singulier est
$\Gamma$. Pour mieux comprendre $\mathscr{F}_\Gamma$ on redresse
$\Gamma$ par l'automorphisme polynomial $$\phi \hspace{1mm} \colon
\hspace{1mm} (x,y,z) \mapsto (x,y-\frac{x^2}{2},z-\frac{x^3}{3})$$
On constate que $\phi^*\Gamma$ est une droite et $\phi^*
\mathscr{F}_\Gamma$ est un feuilletage dont les feuilles sont des
cylindres au dessus de cusp. On peut prolonger ce feuilletage en
un $\mathscr{L}$-feuilletage de degr\'e $2$ et de codimension $1$
sur $\C \P(3)$ dont le lieu singulier est l'union de $\Gamma$,
d'une droite et d'une conique. L'alg\`ebre $\mathscr{L}$ est
isomorphe \`a l'alg\`ebre de Lie du groupe des transformations
affines. Le th\'eor\`eme de Darboux-Jouanolou assure que le
feuilletage $\mathscr{F}_\Gamma$ a une int\'egrale premi\`ere
rationnelle qui se calcule en fait explicitement (\cite{[Ce-LN]})
:
$$\frac{{\left(z_0z_3^2-z_1z_2z_3+\frac{z_2^3}{3}\right)}^2}
{{\left(z_1z_3-\frac{z_2^2}{2}\right)}^3}$$

On note $\omega_\Gamma$ une $1$-forme d\'efinissant le feuilletage
$\mathscr{F}_\Gamma$. On dit que le feuilletage $\mathscr{F}'$ de
degr\'e $2$ d\'efini par $[\omega']$ (classe de $\omega$ dans
l'espace projectif des formes de degr\'e $3$) est voisin de
$\mathscr{F}_\Gamma$ si $[\omega']$ est proche de
$[\omega_\Gamma]$. Si $\mathscr{F}'$ est un feuilletage de degr\'e
$2$ voisin de $\mathscr{F}_\Gamma$ alors $\mathscr{F}'$ s'\'ecrit
$\sigma^* \mathscr{F}_\Gamma$ o\`u $\sigma$ est un automorphisme
de $\C \P(3)$ (\cite{[Ce-LN]}). On dit que le feuilletage
$\mathscr{F}_\Gamma$ est stable et on l'appelle feuilletage
exceptionnel. De la stabilit\'e r\'esulte que la vari\'et\'e alg\'ebrique
$$\overline{\{\sigma^* \mathscr{F}_\Gamma, \hspace{1mm} \sigma
\in \mbox{Aut } \C \P(3)\}}$$ est une composante
irr\'eductible de l'espace des feuilletages de degr\'e $2$ sur
$\C\P(3)$ (\cite{[Ce-LN]}).
\vs

Le feuilletage exceptionnel apparait aussi lorsqu'on consid\`ere
l'action naturelle du groupe triangulaire $$\mathscr{T}_2=\left\{\left(%
\begin{array}{cc}
  \alpha & \beta \\
  0 & \frac{1}{\alpha} \\
\end{array}%
\right) \hspace{1mm} | \hspace{1mm} \alpha \in \C^*, \hspace{1mm}
\beta \in \C \right\} \subset SL(2;\C)$$ sur $V_4$. Elle induit
une action de $\P\mathscr{T}_2$ sur $\C \P(3)$ dont les orbites
sont de dimension $2$, autrement dit un $\mathscr{L}$-feuilletage
de degr\'e $2$ et de codimension $1$ sur $\C \P(3)$ associ\'e \`a
l'alg\`ebre du groupe des transformations affines. On peut montrer
que c'est le feuilletage exceptionnel.
\vs

Il y a encore une autre fa\c{c}on de voir le feuilletage
exceptionnel. On consid\`ere l'action de $\P\mbox{GL}(2,\C)$ sur
$\C \P(4)=\{ \mbox{quadruplets de points sur } \C \P(1)\}$. On a
vu que l'invariant de cette action est l'invariant $j$ des courbes
elliptiques ; fixer la position d'un des quatre points revient \`a
se donner un hyperplan $\C \P(3) \subset \C \P(4)$. La restriction
de l'invariant $j$ \`a ce $\C \P(3)$ donne encore le feuilletage
exceptionnel.

\vs

\section{Exemple de Gordan-N\oe ther.}

\vs

Hesse affirme que si $P$ est un polyn\^ome homog\`ene en les
variables $z_1,\ldots,z_n$, alors les assertions suivantes sont
\'equivalentes

\vs

\begin{itemize}
\item[(i) ] $\det \mathscr{H}\mbox{ess } P=0$

\item[(ii) ] le polyn\^ome $P$ \guillemotleft \hspace{1mm}
d\'epend de moins de $n$ variables \guillemotright.
\vs
\end{itemize}

Bien s\^ur $\mathscr{H}\mbox{ess } P$ d\'esigne la matrice
hessienne de $P$. \vs

Sylvester commence par corriger l'affirmation de Hesse ; si $P$
est un polyn\^ome homog\`ene en $z_1,\ldots,z_n$, les affirmations
suivantes sont \'equivalentes :
\begin{itemize}
\item[(a) ] Il existe $\kappa_1,\ldots,\kappa_n \in \C$
non tous nuls tels que $\displaystyle \sum_{j=1}^n \kappa_j
\frac{\partial P}{\partial z_j}=0$

\item[(b) ] le polyn\^ome $P$ \guillemotleft \hspace{1mm} d\'epend
de moins de $n$ variables \guillemotright.
\vs
\end{itemize}

Puis Gordan et N\oe ther donnent un contre-exemple \`a
l'affirmation de Hesse : le polyn\^ome homog\`ene de degr\'e $3$
sur $\C^5$
$$P(z_1,z_2,z_3,z_4,z_5)=z_1^2z_3+z_1z_2z_4+z_2^2z_5$$ v\'erifie
(i) mais pas (ii). Il se trouve que $P$ d\'efinit un
$\mathscr{L}$-feuilletage $\mathscr{F}_{GN}$ sur $\C^5$.
Les champs
$$X_1=z_2\frac{\partial}{\partial z_3}-z_1\frac{\partial}{\partial z_4},
 \hspace{2mm} X_2=z_2\frac{\partial}{\partial z_4}-z_1\frac{\partial}
 {\partial z_5},$$ $$X_3=z_1\frac{\partial}{\partial z_1}-2
 z_3\frac{\partial}{\partial z_3}-z_4\frac{\partial}
 {\partial z_4}, \hspace{2mm} X_4=z_2\frac{\partial}
 {\partial z_2}-z_4\frac{\partial}{\partial z_4}-2z_5\frac{\partial}
 {\partial z_5}$$ \vs

\hspace{-7mm}annulent $P$, sont lin\'eairement ind\'ependants et v\'erifient les relations
\vs

 $$[X_1,X_2]=0, \hspace{1mm} [X_1,
 X_3]=-2X_1,\hspace{1mm} [X_1,X_4]=-X_1,$$ $$[X_2,X_3]=-X_2,
 \hspace{1mm} [X_2,X_4]=-2X_2 \mbox{ et } [X_3,X_4]=0 $$ \vs

\hspace{-7mm} Ils engendrent une alg\`ebre de Lie maximale $\mathscr{L}$ qui
est ponctuellement de dimension $4$. Le
 $\mathscr{L}$-feuilletage $\mathscr{F}_{GN}$
  se prolonge en un $\mathscr{L}$-feuilletage de degr\'e $2$ sur $\C
\P(5)$ dont une int\'egrale premi\`ere est
$$\frac{z_1^2z_3+z_1z_2z_4+z_2^2z_5}{z_0^3}$$

Dans ce feuilletage il y a un sous-feuilletage en $2$-plans qui admet pour
int\'egrale premi\`ere $$z_1=\mbox{cte},\hspace{1mm} z_2=\mbox{cte},\hspace{1mm}
 z_1^2z_3+z_1z_2z_4+z_2^2 z_5=\mbox{cte}$$ Ce feuilletage en
 $2$-plans produit un $\mathscr{L}$-feuilletage de codimension $3$
 sur $\C^5$ et $\C\P(5)$. L'alg\`ebre $\mathcal{L}$ est
 engendr\'e par les champs $X_1$ et $X_2$.

Il produit aussi un $\mathscr{L}$-feuilletage de codimension $2$ sur
 $\C\P(4)$ d'alg\`ebre $\tilde{\mathscr{L}}=\hspace{1mm}<X_1,
 \hspace{1mm} X_2, \hspace{1mm} R>$ avec pour int\'egrales
 premi\`eres $$\frac{z_1}{z_2} \mbox{ et } \frac{P}{z_1^3}$$ et
 donc d\'efini en coordonn\'ees homog\`enes par la $2$-forme de degr\'e
$3$ $$z_1dz_2\wedge dP-z_2 dz_1\wedge dP+3Pdz_1\wedge dz_2$$ Il
s'obtient comme \guillemotleft \hspace{1mm} intersection
\guillemotright \hspace{1mm} de deux feuilletages ; cependant la
$2$-forme $\Omega$ n'est pas du type $\Omega=\alpha \wedge
\beta$ avec $\alpha$ et $\beta$ deux $1$-formes homog\`enes.

Revenons au polyn\^ome de Gordan-N\oe ther ; le lieu singulier de
$P$ est l'intersection des deux hyperplans $\{z_1=0\}$ et
$\{z_2=0\}$, i.e. un espace lin\'eaire de codimension $2$ dans
$\C^5$. \vs

Soit $f$ un polyn\^ome homog\`ene de degr\'e $3$ sur $\C^n$ dont le
lieu singulier est l'espace lin\'eaire de codimension $2$ d\'efini par
$\{z_1=z_2=0\}$. Alors $f=z_1 A+z_2 B$ ; comme le lieu des z\'eros
de
$$\frac{\partial f}{\partial z_1}=A+z_1\frac{\partial A}{\partial
z_1}+z_2 \frac{\partial B}{\partial z_1} \hspace{2mm} \mbox{ et }
\hspace{2mm} \frac{\partial f}{\partial z_2}=z_1\frac{\partial A}
{\partial z_2}+B+z_2 \frac{\partial B}{\partial z_2}$$ est
$\{z_1=z_2=0\}$, $f$ s'\'ecrit $$\alpha z_1^2+\beta z_1 z_2
+\gamma z_2^2$$ o\`u $\alpha$, $\beta$ et $\gamma$ sont des formes
lin\'eaires. On constate alors que l'exemple de Gordan-N\oe ther
n'est jamais que l'exemple g\'en\'erique : $z_1$, $z_2$, $\alpha$,
$\beta$ et $\gamma$ sont ind\'ependantes ; ceci n'arrive qu'\`a
partir de la dimension $5$. D'o\`u la
\vs

\begin{pro} Soit $f$ un polyn\^ome homog\`ene de
degr\'e $3$ sur $\C^n$, $n\geq 5$. Supposons que
$\mathscr{S}ing$ $f$ soit un espace lin\'eaire de codimension $2$.
Alors g\'en\'eriquement $f$ est lin\'eairement conjugu\'e \`a
l'exemple de Gordan et N\oe ther $z_1^2z_2+z_1z_2z_4+z_2^2
z_5$.

En particulier le feuilletage associ\'e \`a $f$ est un
$\mathscr{L}$-feuilletage.
\vs
\end{pro}

%%%%%%%%%%%%%%%%%%%%%%%%%%%%%%%%%%%%%%%%%%%%%%%%%%%%%%%%%%%%%%%%%%%%%%%%
%%%%%%%%%%%%%%%%%%%%%%%%%%%%%%%%%%%%%%%%%%%%%%%%%%%%%%%%%%%%%%%%%%%%%%%%

\chapter{$\mathscr{L}$-feuilletages de degr\'es $0$ et
$1$ sur $\C \P(n)$.}

\vs

\section{Feuilletages de degr\'e $0$ sur $\C
\P(n)$.}

\vs

La proposition suivante donne une description des feuilletages de
degr\'e $0$ sur $\C \P(n)$ : \vs

\begin{pro} Un feuilletage $\mathscr{F}$
de degr\'e $0$ sur $\C\P(n)$ est associ\'e \`a un
pinceau d'hyperplans, i.e. poss\`ede une int\'egrale premi\`ere
$\frac{\ell_1}{\ell_2}$ o\`u $\ell_1$ et $\ell_2$ sont
des formes lin\'eaires ind\'ependantes. En particulier,
$\mathscr{F}$ est un $\mathscr{L}$-feuilletage.
\vs
\end{pro}

\begin{proof}[D\'emonstration] Soient $\omega$ une $1$-forme d\'efinissant le feuilletage
$\mathscr{F}$ et $m \not \in \mathscr{S}ing \hspace{1mm} \omega$.
Le th\'eor\`eme de Frobenius assure l'existence au voisinage du
point $m$ d'une submersion $z_1$ et d'une unit\'e $u$ que l'on
peut \'ecrire $\mbox{cte} - z_2$ telles que
$$\omega=(\mbox{cte} - z_2)dz_1$$ d'o\`u l'\'egalit\'e
$$d\omega=dz_1\wedge dz_2$$ Mais comme $d\omega$ est une $2$-forme
\guillemotleft \hspace{1mm} constante \guillemotright \hspace{1mm},
un d\'eveloppement de Taylor au voisinage
de $m$ g\'en\'erique assure l'existence de deux formes lin\'eaires
$\ell_1$ et $\ell_2$ telles que $d\omega=d\ell_1 \wedge d\ell_2$.

Il en r\'esulte que
$\omega=\frac{1}{2}(\ell_1d\ell_2-\ell_2d\ell_1)+dQ$ o\`u $Q$ est
une forme quadratique. Comme $\omega$ annule le champ
radial, la forme quadratique $Q$ est nulle.

Autrement dit \`a conjugaison pr\`es par un automorphisme de
$\C\P(n)$, il y a un unique feuilletage
$\mathscr{F}_0$ de degr\'e $0$ dans $\C\P(n)$, le
feuilletage associ\'e \`a $$\omega_0=z_0dz_1-z_1dz_0$$ i.e. au
pinceau d'hyperplans $\frac{z_1}{z_0}=\mbox{cte}$.

Soit $\tilde{\mathscr{L}}$ l'alg\`ebre de Lie des champs
lin\'eaires sur $\C^{n+1}$ annul\'es par la $1$-forme
$$z_0dz_1-z_1dz_0$$ On constate que $$R,\hspace{1mm} z_j \frac{\partial
}{\partial z_i} \hspace{2mm} i\geq 2, \hspace{2mm} j\geq 0$$
forment une base de $\tilde{\mathscr{L}}$ qui est visiblement
maximale. On peut \'evidemment pr\'eciser la nature alg\'ebrique
de $\tilde{\mathscr{L}}$ et donc de $\mathscr{L}$ : l'alg\`ebre
$\mathscr{L}$ s'identifie \`a celle des matrices $(n+1)\times
(n+1)$ dont les deux premi\`eres lignes sont nulles.
\end{proof}

\vs

\section{ Feuilletages de degr\'e $1$ sur $\C
\P(n)$.}

\vs

On commence par rappeler qu'un feuilletage de degr\'e $1$ sur
$\C \P(2)$ poss\`ede une int\'egrale premi\`ere de
l'un des $3$ types suivants :

$$z_0^{\lambda_0}z_1^{\lambda_1} z_2^{\lambda_2} \mbox{ o\`u }
\lambda_i \in \C \mbox{ et }
\displaystyle \sum_{i=0}^2 \lambda_i=0$$
$$\frac{z_0}{z_1} \exp\left(\frac{z_2}{z_1}\right)$$
$$\frac{Q}{z_0^2}$$ o\`u $Q$ est une forme quadratique de rang
maximum.
\vs

En effet, soit $\mathscr{F}$ un feuilletage de degr\'e $1$ sur
$\C \P(2)$ d\'ecrit par la $1$-forme $\omega$. Le
feuilletage $\mathscr{F}$ admet trois points singuliers donc une
droite invariante. En effet supposons que ces trois points ne soient pas
confondus, alors la droite passant par deux points singuliers
distincts est invariante car a deux points de contacts avec le
feuilletage qui est de degr\'e $1$. Si les trois points
singuliers sont confondus, alors chaque droite passant par le
point singulier triple a un contact d'ordre $3$ donc est
invariante. Le feuilletage est donc un pinceau d'hyperplans, i.e.
 de degr\'e $0$ ce qui est exclu.

Si les trois points singuliers du feuilletage $\mathscr{F}$
sont distincts et non align\'es, notons $\mathscr{D}$ une droite passant par
deux de ces points ; dans la carte affine $\C
\P(2) \setminus \mathscr{D}$, le feuilletage $\mathscr{F}$
est d\'ecrit par la $1$-forme
$$\omega_{|\C \P(2) \setminus
\mathscr{D}}=f(z_1,z_2) dz_1+g(z_1,z_2) dz_2$$ avec $f$ et $g$
affines. Si on prend comme origine le troisi\`eme point singulier
de $\mathscr{F}$, le feuilletage
$\mathscr{F}$ est alors d\'ecrit dans la carte affine par la
$1$-forme
$$\tilde{\omega}=(\alpha z_1+\beta z_2)dz_1+(\gamma z_1+ \delta
z_2)dz_2$$ avec $\alpha$, $\beta$, $\gamma$, $\delta \in
\C$. On constate alors que le feuilletage est invariant
sous l'action du flot du champ radial $$(z_1,z_2) \mapsto (e^t
z_1,e^t z_2)$$ autrement  dit $R$ est une sym\'etrie et $P=i_R
\tilde{\omega}$ est un facteur int\'egrant. On a alors
\`a conjugaison pr\`es l'alternative $P=z_0z_1$ ou
$P=z_0^2$ ; si $P=z_0z_1$ alors le
feuilletage $\mathscr{F}$ admet une int\'egrale premi\`ere de type
$z_0^{\lambda_0} z_1^{\lambda_1} z_2^{\lambda_2}$ ; dans le second
cas le
feuilletage admet une int\'egrale premi\`ere de type $\frac{z_1}{z_0}
\exp(\frac{z_2}{z_1})$. En fait, dans cette situation, on v\'erifie que
$\mathscr{F}$ n'a que deux points singuliers.

Si les trois points singuliers du feuilletage sont distincts et
 align\'es, on se donne, pour tout $t \in \C$, une $1$-forme
 $\omega_t$ d\'efinissant un feuilletage $\mathscr{F}_t$ dont les trois
points singuliers ne sont pas align\'es telle que $\displaystyle
\lim_{t\to 0} \omega_t=\omega$. Par ce qui pr\'ec\`ede le
feuilletage $\mathscr{F}_t$ poss\`ede un facteur int\'egrant $P_t$
; par compacit\'e, $\displaystyle \lim_{t \to
0} [P_t]$ a un sens et produit un facteur int\'egrant $P$ de
$\omega$. L'identit\'e d'Euler assure que le polyn\^ome $P$ est r\'eductible,
donc appartient \`a la liste suivante (\`a conjugaison pr\`es) :
$$ z_0^3,\hspace{2mm} z_0z_1z_2,
\hspace{2mm} z_0z_1^2,\hspace{2mm} z_0Q$$ o\`u $Q$ est une forme quadratique de
rang maximum. Si $P=z_0Q$ o\`u $Q$ est une forme
quadratique de rang maximum, alors :
$$\frac{\omega}{P}=\lambda_0 \frac{dz_0}{z_0}+\lambda_1
\frac{dQ}{Q}$$ avec $\lambda_0+2\lambda_1=0$. Quitte \`a prendre
$\lambda_0=2$ et $\lambda_1=-1$, on obtient
$$\frac{\omega}{P}=d\left(\frac{Q}{z_0^2}\right)$$ Consid\'erons
 maintenant le cas o\`u $P=z_0^3$ ; alors
$$\frac{\omega}{P}=d\left(\frac{Q}{z_0^2}\right)$$ Si
$P=z_0z_1z_2$, alors $\mathscr{F}$ est un feuilletage
logarithmique rencontr\'e pr\'ec\'edemment ; enfin si
 $P=z_0z_1^2$, alors $\mathscr{F}$
admet pour int\'egrale premi\`ere
$$\frac{z_0}{z_1}\exp\left(\frac{z_2}{z_1}\right)$$ \vs

Rappelons le r\'esultat suivant classifiant les feuilletages de
degr\'e $1$ sur $\C \P(n)$, $n \geq 3$. Il est
probablement du \`a G. Reeb en tout cas dans le contexte affine ; on
le trouve dans \cite{[Jo]}.
\vs

\begin{thm} Soit $\mathscr{F}$ un
feuilletage de degr\'e $1$ sur $\C \P(n)$ ; on a
l'alternative suivante :

(i) Il existe un feuilletage $\mathscr{F}_1$ de degr\'e $1$ sur
$\C \P(2)$ et une projection $\tau \hspace{1mm}
\colon \hspace{1mm} \C \P(n) \dashrightarrow
\C \P(2)$ telle que $\mathscr{F}=\tau^*
\mathscr{F}_1$.

(ii) Le feuilletage $\mathscr{F}$ poss\`ede une int\'egrale
premi\`ere rationnelle du type $\frac{Q}{L^2}$ avec $\deg(Q)=2$ et
$\deg(L)=1$, i.e. le feuilletage $\mathscr{F}$ est donn\'e en
coordonn\'ees homog\`enes par $$\omega=LdQ-2QdL$$ avec codim $\mathscr{S}ing$
$\omega\geq 2$.
\vs
\end{thm}

Cet \'enonc\'e permet de montrer que l'espace des feuilletages de
degr\'e $1$ sur $\C \P(n)$ a deux composantes
irr\'eductibles correspondant aux deux possibilit\'es de
l'alternative. \vs

\begin{cor} Tout feuilletage $\mathscr{F}$
de degr\'e $1$ poss\`ede un hyperplan $\mathscr{H}$ invariant par
$\mathscr{F}$, autrement dit un hyperplan $\mathscr{H}$ tel que
$\mathscr{H} \setminus \mathscr{S}ing \hspace{1mm} \mathscr{F}$
soit une feuille de $\mathscr{F}$.\vs
\end{cor}

%\begin{proof}[D\'emonstration] Supposons que $\mathscr{F}$ poss\`ede une int\'egrale
%premi\`ere rationnelle du type $\frac{Q}{L^2}$ avec $\deg(Q)=2$ et
%$\deg(L)=1$, i.e. le feuilletage $\mathscr{F}$ est donn\'e en
%coordonn\'ees homog\`enes par $$LdQ-2QdL$$ $L$ ne divisant pas $Q$
%et $Q$ n'\'etant pas un carr\'e. Alors $\mathscr{H}=(L=0)$
%convient.
%
%Sinon, il existe un feuilletage $\mathscr{F}_1$ de degr\'e $1$ sur
%$\C \P(2)$ et une projection $\tau \hspace{1mm}
%\colon \hspace{1mm} \C \P(n) \dashrightarrow
%\C \P(2)$ tels que $\mathscr{F}=\tau^*
%\mathscr{F}_1$. Or un feuilletage de degr\'e $1$ sur
%$\C\P(2)$ a au moins une droite invariante,
%d'o\`u l'existence via $\tau$ d'un hyperplan invariant par
%$\mathscr{F}$.\end{proof}

%\textbf{Remarque. } Si $X$ est un champ de vecteurs sur $\C
%\P(n)$ laissant l'hyperplan $\mathscr{H}$ invariant, alors
%dans des coordonn\'ees affines $z_1,\ldots,z_n$ de la carte
%$\C \P(n) \setminus \mathscr{H} \simeq
%\C^n$, le champ $X$ s'\'ecrit
%$$\sum_{j=1}^n a_j(z) \frac{\partial}{\partial z_j}$$ o\`u les
%$a_j$ sont affines. En particulier si $\mathscr{F}$ est un
%$\mathscr{L}$-feuilletage poss\'edant un hyperplan invariant. Alors
%$\mathscr{L}$ s'identifie \`a une sous-alg\`ebre de champs de
%vecteurs affines sur $\C^n \simeq \C \P(n)
%\setminus \mathscr{H}$.
%\vs

On en d\'eduit la
\vs

\begin{pro} Soit $\mathscr{F}$ un feuilletage de degr\'e $1$ sur $\C
\P(n)$ poss\'edant une int\'egrale premi\`ere de type
$\frac{Q}{L^2}$. Alors $\mathscr{F}$ est un
$\mathscr{L}$-feuilletage. Si $Q$ est g\'en\'erique, l'alg\`ebre
$\mathscr{L}$ associ\'ee est isomorphe \`a $so(n,\C)$.
\vs
\end{pro}

\begin{proof}[D\'emonstration] On prend l'hyperplan invariant $\mathscr{H}=(L=0)$ comme
hyperplan \`a l'infini. Sur $\C \P(n) \setminus \mathscr{H} \simeq
\C^n$, le feuilletage $\mathscr{F}$ a pour int\'egrale premi\`ere
le polyn\^ome $$Q=q_0+q_1(z)+q_2(z)$$ o\`u les $q_i$ sont
homog\`enes de degr\'e $i$. Si $q_2$ est de rang maximum, alors
quitte \`a composer par une translation \`a la source et une au
but, on peut supposer que $q_0=q_1=0$. Ainsi $Q=q_2$ et
l'alg\`ebre des champs affines qui annulent $Q$ est en fait une
alg\`ebre de champs lin\'eaires car $0$ est singularit\'e isol\'ee
de $Q$ ; cette alg\`ebre est $so(Q)\simeq so(n,\C)$.

Dans les cas o\`u $q_2$ n'est pas de rang maximum,
alors $Q$ est conjugu\'e \`a $z_0^2+\ldots+z_k^2+q_1(z)$, avec
$q_1$ affine. Si $q_1$ est fonction de $z_0,\ldots,z_{k}$, alors
le polyn\^ome $Q$ est conjugu\'e \`a $z_0^2+\ldots+
z_k^2$ ; le feuilletage $\mathscr{F}$ d\'ecrit par $Q$ est associ\'e \`a
l'alg\`ebre engendr\'ee par les champs
$$z_i\frac{\partial}{\partial z_j}-z_j\frac{\partial} {\partial
z_i}, \hspace{2mm} 0 \leq i<j \leq k,$$ $$z_\ell
\frac{\partial}{\partial z_{p}}, \hspace{1mm} \ell \geq 0,
\hspace{1mm} p \geq k+1, \mbox{ et } \frac{\partial}{\partial
z_m}, \hspace{1mm} m \geq k+1.$$ Si $q_1$ n'est pas une fonction
de $z_0,\ldots,z_{k}$, alors le polyn\^ome $Q$ s'\'ecrit \`a
conjugaison pr\`es $z_0^2+\ldots+z_{k}^2+z_{k+1}$ ; le feuilletage
est alors associ\'e \`a l'alg\`ebre engendr\'ee par les champs
$$z_i\frac{\partial}{\partial z_j}-z_j\frac{\partial} {\partial
z_i}, \hspace{1mm} 1 \leq i<j \leq k, \hspace{1mm} z_\ell
\frac{\partial}{\partial z_p}, \hspace{1mm} \ell \geq 1,
\hspace{1mm} p \geq k+2,$$ $$\frac{\partial}{\partial z_m},
\hspace{1mm} m \geq k+2, \mbox{ et } 2z_r\frac{\partial}{\partial
z_{k+1}}-\frac{\partial}{\partial z_r}, \hspace{1mm} r \not =
k+1$$\end{proof}

\vs

Lorsque $\mathscr{F}$ est de type \guillemotleft \hspace{1mm}
pull-back \guillemotright \hspace{1mm} on a la : \vs

\begin{pro}\label{feuil}
Soit $\mathscr{F}$ un
feuilletage de degr\'e $1$ sur $\C \P(n)$ du type
$\mathscr{F}=\tau^*\mathscr{F}_1$, o\`u $\tau \hspace{1mm} \colon
\hspace{1mm} \C\P(n) \dashrightarrow \C
\P(2)$ est lin\'eaire et $\mathscr{F}_1$ un feuilletage de
degr\'e $1$ sur $\C \P(2)$. Alors $\mathscr{F}$
est un $\mathscr{L}$-feuilletage. \vs
\end{pro}

\begin{proof}[D\'emonstration] Soit $\mathscr{D}$ une droite invariante par $\mathscr{F}_1$
et $\mathscr{H}=\tau^*(\mathscr{D})$ ; l'hyperplan
$\mathscr{H}$ est invariant par $\mathscr{F}$. Dans les cartes affines
$\C \P(n) \setminus
\mathscr{H}=~\{z_1,\ldots,z_n\}$ et $\C
\P(2)\setminus \mathscr{D}=\{z_1,z_2\}$, l'application
$\tau$ est de la forme $$(z_1,\ldots,z_n) \mapsto (z_1,z_2)$$
Le feuilletage $\mathscr{F}_1$ est d\'ecrit dans la carte affine
$\{z_1,z_2\}$ par $$\omega_1=a(z_1,z_2)dz_1+b(z_1,z_2)dz_2$$ avec
$a$ et $b$ affines. Alors $\mathscr{F}=\tau^*(\mathscr{F}_1)$
est donn\'e par : $$\omega=\tau^*\omega_1=a(z_1,z_2)dz_1+b(z_1,
z_2)dz_2$$ Dans cette m\^eme carte, on constate que la $1$-forme
$\omega$ annule les champs affines $$b(z_1,z_2)\frac{\partial}
{\partial z_1}-a(z_1,z_2)\frac{\partial}{\partial z_2},
\hspace{1mm} \frac{\partial}{\partial z_3},
  \hspace{1mm} \ldots, \hspace{1mm} \frac{\partial}{\partial z_n},$$
  $$ z_j\frac{\partial}{\partial z_k} \mbox{ pour } 3\leq k \leq n,
 \hspace{2mm} 1 \leq j\leq n.$$
\end{proof}

On en d\'eduit le  \vs

\begin{thm} Tout feuilletage de degr\'e $1$ sur
$\C\P(n)$ est un
$\mathscr{L}$-feuilletage.\vs
\end{thm}

La proposition \ref{feuil} ne se g\'en\'eralise pas en degr\'e
sup\'erieur
 : \vs

\begin{pro}\label{obstruction}
Soit $\mathscr{F}$ du type
$\mathscr{F}=\tau^*\mathscr{F}_1$, o\`u $\tau \hspace{1mm} \colon
\hspace{1mm} \C\P(n) \dashrightarrow~\C
\P(2)$ est lin\'eaire et $\mathscr{F}_1$ est un
feuilletage de degr\'e plus grand que $2$ sur $\C
\P(2)$. Alors $\mathscr{F}$ n'est pas un
$\mathscr{L}$-feuilletage. \vs
\end{pro}

\begin{proof}[D\'emonstration] Comme dans la preuve de la proposition
\ref{feuil}, on constate
que $\mathscr{F}$ est d\'ecrit par
$$\omega=A_0(z_0,z_1,z_2)dz_0+A_1(z_0,z_1,z_2)dz_1+A_2(z_0,z_1,z_2)dz_2$$
o\`u les $A_i$ sont de degr\'e $\geq 3$. On remarque que $\omega$
annule le champ radial et les champs $z_j\frac{\partial }{\partial
z_k}$ pour $k \geq 3$ et $j=0,\ldots,n$. Supposons que pour $z$
g\'en\'erique $$\dim_{\C} \mathscr{L}(z)=
 \dim_{\C}\{ X(z), \hspace{1mm} X \in \mathscr{L}\}=n-1$$ Alors
il existe $X$ champ lin\'eaire qui annule $\omega$ et tel que la
$1$-forme $\bar{\omega}$ de degr\'e $2$ d\'efinie par
$$i_Ri_Xdz_0\wedge dz_1 \wedge dz_2$$ soit non identiquement nulle.
En un point $m$ g\'en\'erique, on a $\ker \omega(m)=\ker
\bar{\omega}(m)$, autrement dit $\omega$ et $\bar{\omega}$ sont
colin\'eaires. La $1$-forme $\bar{\omega}$ s'\'ecrit
$Q\bar{\omega}'$ o\`u $\mathscr{S}ing \hspace{1mm} \bar{\omega}'$
est de codimension $2$ et $\deg \bar{\omega}'\leq 2$. La
colin\'earit\'e de $\omega$ et $\bar{\omega}$ entra\^ine celle de
$\omega$ et $\bar{\omega}'$ ; ainsi $\omega$ s'\'ecrit
$P\bar{\omega}'$ avec $\deg P \geq 1$ et $\omega$ s'annule sur
l'hypersurface d\'efinie par le lieu des z\'eros de $P$ ce qui
viole la condition codim$_{\C} \mathscr{S}ing \hspace{1mm} \omega
\geq 2$.
\end{proof}

%%%%%%%%%%%%%%%%%%%%%%%%%%%%%%%%%%%%%%%%%%%%%%%%%%%%%%%%%%%%%%%%%%%
%%%%%%%%%%%%%%%%%%%%%%%%%%%%%%%%%%%%%%%%%%%%%%%%%%%%%%%%%%%%%%%%%%%

\chapter{$\mathscr{L}$-feuilletages sur $\C
\P(3)$.}

\section{Feuilletages de degr\'e $2$ sur $\C\P(3)$.}

On rappelle d'abord la description de l'espace des feuilletages de
degr\'e $2$ sur $\C\P(3)$ (\cite{[Ce-LN]}). Puis
on d\'egage quelques obstructions \`a ce qu'un tel feuilletage
soit un $\mathscr{L}$-feuilletage. On introduit les sous-ensembles
$\Sigma(.;3)$ des formes homog\`enes int\'egrables
 de degr\'e $3$ sur $\C^4$ :
\vs

(i) On note
$$\Sigma(1,1,1,1;3)=\{\ell_0 \ell_1 \ell_2 \ell_3 \displaystyle
\sum_{k=0}^3 \lambda_k \frac{d\ell_k}{\ell_k} \hspace{1mm} | \hspace{1mm}
\lambda_k \in \C^*,\hspace{1mm} \displaystyle \sum_{k=0}^3
\lambda_k=0\}$$ o\`u les $\ell_i$ sont des formes lin\'eaires non nulles. On
d\'efinit alors la sous-vari\'et\'e
irr\'eductible $\overline{\P\Sigma(1,1,1,1;3)}$ de
$\P\Omega^1_3(\C^4)$, o\`u
$\Omega^1_p(\C^n)$ d\'esigne l'espace des $1$-formes
homog\`enes de degr\'e $p$ sur $\C^n$. Un \'el\'ement
$[\omega] \in \overline{\P\Sigma(1,1,1,1;3)}$ d\'efinit un
feuilletage de degr\'e $2$ sur $\C\P(3)$ si
codim$_{\C} \mathscr{S}ing \hspace{1mm} \omega \geq 2$.
Ces \'el\'ements constituent un ouvert de Zariski
$\overline{\P\Sigma(1,1,1,1;3)}\hspace{1mm}^*$ de la
sous-vari\'et\'e unirationnelle
$\overline{\P\Sigma(1,1,1,1;3)}$.
\vs

Plus g\'en\'eralement si $\Theta$ est un sous-ensemble de
$\P\Omega^1_3(\C^4)$, on note $$\Theta^*:=\{[\omega] \in \Theta
\hspace{1mm} | \hspace{1mm} \mbox{codim}_{\C}
\mathscr{S}\mbox{ing} \hspace{1mm} \omega \geq 2\}$$

(ii) On introduit l'ouvert de Zariski
$\overline{\P\Sigma(1,1,2;3)}\hspace{1mm}^*$ de la vari\'et\'e
unirationnelle $\overline{\P\Sigma(1,1,2;3)}$ o\`u
$$\Sigma(1,1,2;3)=\{\ell_0 \ell_1 Q \left(\lambda_0 \frac{d\ell_0}{\ell_0}
+\lambda_1
\frac{d\ell_1}{\ell_1}+\lambda_2\frac{dQ}{Q}\right)\hspace{1mm} | \hspace{1mm}
\lambda_i \in \C^*, \hspace{1mm} \lambda_0 +\lambda_1
+2\lambda_2=0\}$$ o\`u les $\ell_i$ sont des formes lin\'eaires non nulles et
$Q$ une forme quadratique non triviale.
\vs

(iii) On d\'efinit l'ouvert de Zariski $\overline{\P\Sigma(1,3;3)}
\hspace{1mm}^*$ de la vari\'et\'e unirationnelle
$\overline{\P\Sigma(1,3;3)}$ adh\'erence de :
$$\Sigma(1,3;3)=\{\ell C \left(\lambda_0 \frac{d\ell}{\ell} +\lambda_1
\frac{dC}{C}\right) \hspace{1mm} | \hspace{1mm} \lambda_i \in \C^*,
\hspace{1mm} \lambda_0 +3\lambda_1=0\}$$ avec $\ell$ forme lin\'eaire
non nulle et $C$ de degr\'e $3$.
\vs

(iv) On d\'efinit aussi l'ouvert de Zariski
$\overline{\P\Sigma(2,2;3)} \hspace{1mm}^*$ de la vari\'et\'e
unirationnelle $\overline{\P\Sigma(2,2;3)}$ adh\'erence de :
$$\Sigma(2,2;3)=\{Q_0 Q_1 \left(\lambda_0 \frac{dQ_0}{Q_0} +\lambda_1
\frac{dQ_1}{Q_1}\right),\hspace{1mm} | \hspace{1mm} \lambda_i \in
\C^*, \hspace{1mm}\lambda_0+\lambda_1=0\}$$ o\`u les $Q_i$ sont
des formes quadratiques.\vs

On d\'efinit ensuite $\mbox{PB}(2;3)$ l'ensemble des pull back
lin\'eaires comme suit : $$\mbox{PB}(2;3)=\{\tau^*
\tilde{\omega}\hspace{1mm}|\hspace{1mm} \tau \hspace{1mm} \colon
\hspace{1mm} \C^4 \to \C^3 \mbox{ lin\'eaire}, \hspace{1mm}
\tilde{\omega} \in \Omega^1_3(\C^3), \hspace{1mm} i_R
\tilde{\omega}=0\}$$ et l'ouvert de Zariski
$\overline{\P\mbox{PB(2;3)}} \hspace{1mm}^*$ de la
sous-vari\'et\'e unirationnelle $\overline{\P\mbox{PB}(2;3)}$.

On termine par l'adh\'erence de l'orbite du feuilletage exceptionnel :
$$\mbox{Excep}(2,3)=\overline{\P\mbox{GL}
(4,\C).\omega_\Gamma} \cap \{[\omega], \hspace{1mm}
\mbox{codim}_{\C} \mathscr{S}ing \hspace{1mm} \omega \geq 2\}$$
o\`u $\omega_\Gamma$ d\'esigne la $1$-forme d\'efinissant le
feuilletage exceptionnel. Par construction, les \'el\'ements
g\'en\'eriques de Excep$(2,3)$ sont conjugu\'es et correspondent
\`a un $\mathscr{L}$-feuilletage o\`u $\mathscr{L}$ est isomorphe
\`a l'alg\`ebre du groupe des transformations affines comme nous
l'avons vu dans 2.4.8. \vs

On peut maintenant \'enoncer le : \vs

\begin{thm} (\cite{[Ce-LN]}) Les ensembles $\overline{\P\Sigma(.;3)}
\hspace{1mm}^*$, $\overline{\P\mbox{PB}(2;3)}\hspace{1mm}^*$ et
$\mbox{Excep}(2;3)$ sont les composantes irr\'eductibles de
l'espace des feuilletages de degr\'e $2$ sur $\C\P(3)$.
\vs
\end{thm}

Ce th\'eor\`eme se g\'en\'eralise en toute dimension (\cite{[Ce-LN]}).
\vs

Pour se convaincre que tous les feuilletages de degr\'e $2$ ne
sont pas des $\mathscr{L}$-feuilletages nous allons voir des
obstructions plus ou moins \'evidentes \`a ce qu'ils le soient ; ces
obstructions sont utilis\'ees pour v\'erifier certains calculs.
L'\'enonc\'e en dimension $n \geq 4$ pourrait s'av\'erer pertinent
pour la classification des feuilletages quadratiques de
$\C \P(n)$.
\vs

\begin{lem} Un \'el\'ement g\'en\'erique de
$\overline{\P\Sigma(1,3;3)} \hspace{1mm}^*$ ne d\'efinit
pas
 un $\mathscr{L}$-feuilletage de degr\'e $2$ sur $\C
\P(3)$.
\vs
\end{lem}

\begin{proof}[D\'emonstration] Soit $\mathscr{F}$ le feuilletage d\'ecrit par $$\omega=\ell
C\left(\lambda_0 \frac{d\ell}{\ell}+\lambda_1
\frac{dC}{C}\right)$$ avec $\lambda_0+3\lambda_1=0$. Le
feuilletage $\mathscr{F}$ a pour int\'egrale premi\`ere {\large
$\frac{C}{\ell^3}$} ; ainsi dans la carte affine $\C^3=\C \P(3)
\setminus (\ell=0)$, les feuilles de $\mathscr{F}$ sont les
niveaux du polyn\^ome $C$ de degr\'e $3$ priv\'es des points
singuliers. On suppose que $C$ est g\'en\'eral autrement dit que
les niveaux de $C$ sont des surfaces cubiques g\'en\'erales ; une
telle surface contient $27$ droites. Soit $X \in \chi(\C \P(3))$
un champ tangent \`a $\mathscr{F}$, alors $X$ laisse l'hyperplan
$\ell=0$ invariant et par suite est affine dans la carte $\C^3$.
Un tel $X$ est n\'ecessairement tangent aux $27$ droites de chaque
niveau g\'en\'erique. Si $\mathscr{F}$ \'etait un
$\mathscr{L}$-feuilletage, il existerait deux champs affines $X$
et $Y$ tels que
$$\omega_{|\C^3}=\mbox{cte }i_Xi_Y dz_0 \wedge dz_1 \wedge
dz_2$$ en particulier les champs $X$ et $Y$ seraient
colin\'eaires le long d'un nombre fini de courbes. Or on vient de
voir que dans chaque fibre g\'en\'erique de $C$, les champs $X$ et
$Y$ sont colin\'eaires le long de 27 droites.
\end{proof}

\vs

De m\^eme on a le : \vs

\begin{lem}\label{courbeell}
Un \'el\'ement g\'en\'erique de
$\overline{\P\Sigma(2,2;3)} \hspace{1mm}^*$ ne d\'efinit
pas
 un $\mathscr{L}$-feuilletage de degr\'e $2$ sur $\C
\P(3)$.\vs
\end{lem}

\begin{proof}[D\'emonstration] Soit $\mathscr{F}$ un feuilletage de degr\'e $2$ sur
$\C \P(3)$ d\'ecrit par
$$\omega=Q_0Q_1\left(\frac{dQ_0}{Q_0}-\frac{dQ_1}{Q_1}\right).$$
Si $Q_0$ et $Q_1$ sont des formes quadratiques g\'en\'eriques,
alors $\P(Q_0=Q_1=0)$ est une courbe elliptique $\mathscr{E}$ non
plane (quartique). Soit $X \subset \chi(\C \P(3))$ un champ
tangent au feuilletage $\mathscr{F}$ ; comme $\mathscr{E} \subset
\mathscr{S}ing \hspace{1mm} \mathscr{F}$, le champ $X$ est tangent
\`a $\mathscr{E}$. Si $X$ est non identiquement nul sur
$\mathscr{E}$, il induit un champ de vecteurs $X_{\mathscr{E}}$
sur $\mathscr{E}$. Donc $X_\mathscr{E}$ est constant et en
particulier n'admet pas de singularit\'e. Par suite $\mathscr{E}$
est une trajectoire de $X$. La classification des champs
lin\'eaires dit que leurs seules trajectoires compactes sont les
points singuliers ; ainsi $X_{\mathscr{E}}=0$. Le lieu des z\'eros
d'un champ lin\'eaire \'etant une sous-vari\'et\'e lin\'eaire, on
obtient finalement que $X=0$.
\end{proof}

\vs

Le lemme suivant a \'et\'e prouv\'e au chapitre 3 : \vs

\begin{lem} Un \'el\'ement g\'en\'erique de
$\overline{\P\mbox{PB}(2;3)}$ n'est pas un
$\mathscr{L}$-feuilletage de degr\'e $2$ sur $\C
\P(3)$.\vs
\end{lem}

\section{$\mathscr{L}$-feuilletages de degr\'e $2$ sur $\C\P
(3)$.}

On s'int\'eresse \`a la classification des
$\mathscr{L}$-feuilletages sur $\C \P(3)$. Dans le chapitre
pr\'ec\'edent, on a vu que les feuilletages de degr\'es $0$ et $1$
sont des $\mathscr{L}$-feuilletages.

La liste des $\mathscr{L}$-feuilletages de degr\'e $2$ sur $\C
\P(3)$ est donn\'ee par le th\'eor\`eme suivant ; dans
l'\'enonc\'e $(z_0,z_1,z_2,z_3)$ d\'esigne un syst\`eme de
coordonn\'ees.
\vs

\begin{thm}\label{clas}
Un $\mathscr{L}$-feuilletage de
degr\'e $2$ sur $\C \P(3)$ est \`a conjugaison
pr\`es d\'ecrit par l'une des formes ferm\'ees rationnelles
suivantes
\begin{eqnarray}
&Ab (i)& z_0^{\lambda_0}z_1^{\lambda_1}z_2^{\lambda_2}
z_3^{\lambda_3} \nonumber \\
 & & \nonumber \\
&Ab (ii)&\frac{z_1z_2+z_0z_3}{z_1z_3} \nonumber \\
& & \nonumber \\
&Ab (ii)& \frac{z_1}{z_3} \exp\left(\frac{z_0z_3+z_2z_1}{z_1
z_3}\right) \nonumber \\
& & \nonumber \\
&Ab (iii)& \frac{z_0}{z_3}\exp \left(\frac{z_2^2
-z_1z_3}{z_3^2}\right)\nonumber \\
& & \nonumber \\
&Ab (iv)& \frac{z_1z_3^{\kappa-1}}{z_0^{\kappa}} \exp\left(
\frac{z_2}{z_3}\right) \nonumber \\
& & \nonumber \\
&Ab (v)& \frac{z_0z_3^2+z_1z_2z_3+z_2^3}{z_3^3}\nonumber \\
& & \nonumber \\
&Af (i)& \mbox{le feuilletage exceptionnel donn\'e par
}\nonumber \\
& & \nonumber \\
& & \frac{{\left(z_0z_3^2-z_1z_2z_3+\frac{z_2^3}{3}
\right)}^2}{{\left(z_1z_3-\frac{z_2^2}{2}\right)}^3}
\nonumber \\
& & \nonumber \\
&Af (ii)& \frac{z_0^2z_3^{2(\kappa-1)}}{{\left(z_1z_3-z_2^2
\right)}^\kappa}\nonumber \\
& & \nonumber \\
&Af (iii)& \frac{z_3^2}{z_1z_3-z_2^2}\exp\left(\frac{z_0}
{z_3}\right) \nonumber \\
& & \nonumber \\
&Af (iv)& \frac{z_0}{z_3}\exp\left(\frac{z_2^2
-z_1z_3}{z_0z_3}\right) \nonumber
\end{eqnarray}

\begin{eqnarray}
&Af (v)& \frac{z_2}{z_3}\exp\left(\frac{z_0z_3
-z_1z_2}{z_3^2}\right) \nonumber \\
& & \nonumber \\
&Af (vi)& \frac{{\left(z_0z_3-z_1z_2\right)}^\kappa}
{z_2^{\kappa+1}z_3^{\kappa-1}}\nonumber \\
& & \nonumber \\
&Af (vii)& \frac{z_0z_3-z_1z_2}{z_3^2}\exp\left(-
\frac{z_2}{z_3}\right)\nonumber
\end{eqnarray}
 o\`u $\kappa$ et $\lambda_i \in \C$.
\vs
\end{thm}

En particulier, on obtient le th\'eor\`eme annonc\'e dans
l'introduction. Les symboles Ab et Af distinguent les cas ab\'eliens
et affines.\vs

La proposition \ref{degre} assure qu'un $\mathscr{L}$-feuilletage
$\mathscr{F}$ de degr\'e $2$ sur $\C \P(3)$ est
associ\'e \`a une sous-alg\`ebre $\mathscr{L} \subset
\chi(\C \P(3))$ de dimension $2$. D'apr\`es la
classification des alg\`ebres de Lie complexes de dimension $2$,
l'alg\`ebre $\mathscr{L}$ est soit ab\'elienne, soit isomorphe \`a
l'alg\`ebre de Lie du groupe des transformations affines $\{z
\mapsto az+b\}$. La d\'emonstration du th\'eor\`eme \ref{clas} se fait
en plusieurs \'etapes suivant la nature de l'alg\`ebre
$\mathscr{L}$. Elle proc\`ede par classification effective des
sous-alg\`ebres de matrices de dimension deux.\vs

Bien que cela n'apparaisse pas directement dans le texte nous
avons souvent utilis\'e le lemme suivant qui permet d'\'eliminer
certaines configurations :
\vs

\begin{lem}\label{elimination}
Soit $\mathscr{F}$ un
$\mathscr{L}$-feuilletage sur $\C \P(n)$.
Supposons qu'il existe un \'el\'ement $X \in \mathscr{L}$ dont la
matrice associ\'ee est de rang $1$. Alors $\deg \mathscr{F} < n-1$
et il existe une projection $\tau \hspace{1mm} \colon \hspace{1mm}
\C \P(n) \dashrightarrow \C \P
(n-1)$, un feuilletage $\mathscr{F}'$ sur $\C \P
(n-1)$ tels que $\mathscr{F}=\tau^*\mathscr{F}'$.
\vs
\end{lem}

\begin{proof}[D\'emonstration] Soit $X_1=X$, $X_2$, $\ldots$, 
$X_{n-1}$ des \'el\'ements de
$\mathscr{L}$ tels qu'en un point g\'en\'erique $m$, les champs
$R(m)$, $X_1(m)$, $\ldots$, $X_{n-1}(m)$ soient ind\'ependants. La
condition de rang entra\^ine que le champ
$X$ s'\'ecrit $$\ell X_0$$ o\`u $\ell$ est une forme lin\'eaire et
$X_0$ un champ constant.

On peut choisir des coordonn\'ees dans lesquelles le champ $X_0$
est de la forme {\large $\frac{\partial}{\partial z_0}$} ; on a alors
$$i_R i_{X_0}\ldots i_{X_{n-1}} \hspace{1mm} dz_0\wedge \ldots
\wedge dz_n=P \omega$$ o\`u $P$ est un polyn\^ome homog\`ene
\'eventuellement constant et $\omega$ une $1$-forme d\'efinissant
$\mathscr{F}$. On constate que $\deg \omega <n-1$ et que la
$1$-forme $\omega$ ne d\'epend pas de $z_0$.
\end{proof}

\vs

\subsection{L'alg\`ebre $\mathscr{L}$ est ab\'elienne.}

Si $A$ est une matrice carr\'ee on \'ecrit $A=A_S+A_N$ o\`u $A_S$
est semi-simple, $A_N$ nilpotente et $[A_S,A_N]=0$. Si $A$ et $B$
commutent, alors $[A_S,B_S]=[A_N,B_N]=[A_S,B_N]=[A_N,B_S]=0$ :
c'est la \guillemotleft \hspace{1mm} jordanisation simultan\'ee
\guillemotright.

Par Jordanisation simultan\'ee des matrices, on \'etablit la :
\vs

\begin{pro} Soit $\mathscr{F}$ un $\mathscr{L}$-feuilletage de degr\'e $2$ sur
$\C\P(3)$ associ\'e \`a une alg\`ebre $\mathscr{L}$ ab\'elienne.
On est \`a conjugaison pr\`es dans l'une des configurations suivantes :

\begin{itemize}

\item[(i)] l'alg\`ebre $\tilde{\mathscr{L}}$ est diagonale,

\item[(ii)] l'alg\`ebre
 $\tilde{\mathscr{L}}$ est engendr\'ee par
les champs $$z_1\frac{\partial}{\partial z_0}+(z_2+z_3)\frac
{\partial}{\partial z_2}+z_3\frac{\partial}{\partial z_3},
\hspace{2mm} (z_2+\kappa z_3)\frac{\partial}{\partial z_2}
+z_3\frac{\partial}{\partial z_3} \hspace{2mm} \mbox{et } R,$$

\item[(iii)] l'alg\`ebre $\tilde{\mathscr{L}}$ est engendr\'ee par
les champs $$\kappa z_0 \frac{\partial}{\partial
z_0}+z_2\frac{\partial}{\partial z_1} +z_3\frac{\partial}{\partial
z_2}, \hspace{2mm} z_0\frac{\partial}{\partial z_0}
+z_3\frac{\partial}{\partial z_1} \hspace{2mm} \mbox{et } R.$$

\item[(iv)] l'alg\`ebre $\tilde{\mathscr{L}}$ est engendr\'ee par
les champs
$$z_1
\frac{\partial}{\partial z_1}+z_3\frac{\partial}{\partial z_2},
\hspace{2mm} z_0 \frac{\partial}{\partial z_0}+\kappa z_1
\frac{\partial}{\partial z_1} \hspace{2mm} \mbox{et } R,$$

\item[(v)] l'alg\`ebre $\tilde{\mathscr{L}}$ est engendr\'ee par
les champs $$z_1\frac{\partial}{\partial z_0}
+z_2\frac{\partial}{\partial z_1}+z_3\frac{\partial}{\partial
z_2}, \hspace{2mm} (z_2+\kappa z_3)\frac{\partial}{\partial
z_0}+z_3 \frac{\partial}{\partial z_1} \hspace{2mm} \mbox{et }
R.$$
\end{itemize}
\end{pro}

\vs

\begin{proof}[D\'emonstration]  On d\'ecrit par jordanisation toutes les alg\`ebres de
matrices $\tilde{\mathscr{L}}=\C R \oplus \mathscr{L}'$ avec
$\mathscr{L}$ ab\'elienne ; si une telle alg\`ebre contient un
\'el\'ement de rang $1$ on ne la retient pas (\ref{elimination})
L'\'etude se fait \guillemotleft \hspace{1mm} \`a la main
\guillemotright \hspace{1mm} en examinant la nature du spectre
d'un \'el\'ement g\'en\'erique de $\tilde{\mathscr{L}}$. Tous
calculs faits on obtient les alg\`ebres suivantes g\'en\'er\'ees
par
$$R \sim \left(%
\begin{array}{cccc}
  1 & 0 & 0 & 0 \\
  0 & 1 & 0 & 0 \\
  0 & 0 & 1 & 0 \\
  0 & 0 & 0 & 1 \\
\end{array}%
\right)=Id$$ et les matrices de $X$ et $Y$  avec  \vspace{5mm}

\begin{itemize}
\item[(i)] $X=\left(%
\begin{array}{cccc}
  \alpha_0 & 0 & 0 & 0 \\
  0 & \alpha_1 & 0 & 0 \\
  0 & 0 & \alpha_2 & 0 \\
  0 & 0 & 0 & \alpha_3 \\
\end{array}%
\right) \hspace{6mm} Y=\left(%
\begin{array}{cccc}
  \beta_0 & 0 & 0 & 0 \\
  0 & \beta_1 & 0 & 0 \\
  0 & 0 & \beta_2 & 0 \\
  0 & 0 & 0 & \beta_3 \\
\end{array}%
\right)$ \vspace{5mm}

\item[(ii)] $X=\left(%
\begin{array}{cccc}
  1 & 0 & 0 & 0 \\
  0 & 1 & 1 & 0 \\
  0 & 0 & 1 & 0 \\
  0 & 0 & 0 & 0 \\
\end{array}%
\right) \hspace{6mm} Y=\left(%
\begin{array}{cccc}
  0 & 0 & 0 & 0 \\
  0 & 1 & \kappa & 0 \\
  0 & 0 & 1 & 0 \\
  0 & 0 & 0 & 0 \\
\end{array}%
\right) \hspace{3mm} \kappa \not =1$ \vspace{5mm}

\item[(iii)] $X=\left(%
\begin{array}{cccc}
  \kappa & 0 & 0 & 0 \\
  0 & 0 & 1 & 0 \\
  0 & 0 & 0 & 1 \\
  0 & 0 & 0 & 0 \\
\end{array}%
\right) \hspace{6mm} Y=\left(%
\begin{array}{cccc}
  1 & 0 & 0 & 0 \\
  0 & 0 & 1 & 0 \\
  0 & 0 & 0 & 0 \\
  0 & 0 & 0 & 0 \\
\end{array}%
\right) \hspace{3mm} \kappa \not = 1$\vspace{5mm}

\item[(iv)] $X=\left(%
\begin{array}{cccc}
  0 & 0 & 0 & 0 \\
  0 & 1 & 0 & 0 \\
  0 & 0 & 0 & 1 \\
  0 & 0 & 0 & 0 \\
\end{array}%
\right) \hspace{6mm} Y=\left(%
\begin{array}{cccc}
  1 & 0 & 0 & 0 \\
  0 & \kappa & 0 & 0 \\
  0 & 0 & 0 & 0 \\
  0 & 0 & 0 & 0 \\
\end{array}%
\right) \hspace{3mm} \kappa \not =0$\vspace{5mm}

\item[(v)] $X=\left(%
\begin{array}{cccc}
  0 & 1 & 0 & 0 \\
  0 & 0 & 1 & 0 \\
  0 & 0 & 0 & 1 \\
  0 & 0 & 0 & 0 \\
\end{array}%
\right) \hspace{6mm} Y=\left(%
\begin{array}{cccc}
  0 & 0 & 1 & \kappa \\
  0 & 0 & 0 & 1 \\
  0 & 0 & 0 & 0 \\
  0 & 0 & 0 & 0 \\
\end{array}%
\right)$ \vspace{5mm}
\end{itemize}

En prenant les champs de vecteurs associ\'es \`a ces cinq
configurations on obtient la liste annonc\'ee dans la
proposition.
\end{proof}

Nous allons maintenant pr\'esenter la liste des int\'egrales
premi\`eres des feuilletages associ\'es \`a ces cinq alg\`ebres.
\vs

\begin{itemize}
\item[(i)] Cas o\`u $\tilde{\mathscr{L}}$ est diagonale. Avec les notations
pr\'ec\'edentes on a $$X=\sum_{k=0}^3 \alpha_k z_k
 \frac{\partial}{\partial z_k}, \hspace{2mm} Y=\sum_{k=0}^3
 \beta_k z_k \frac{\partial}{\partial z_k}\mbox{ et } R$$ Le
 feuilletage $\mathscr{F}$ est d\'ecrit par $$i_Ri_Xi_Y
 \hspace{1mm} dz_0\wedge dz_1 \wedge dz_2 \wedge dz_3 =
z_0 z_1 z_2 z_3 \sum_{k=0}^3 \lambda_k \frac{dz_k}{z_k}$$ avec
$\displaystyle \sum_{k=0}^3 \lambda_k \alpha_k=\displaystyle
\sum_{k=0}^3 \lambda_k \beta_k=\displaystyle \sum_{k=0}^3
\lambda_k=0$. C'est le cas Ab (i) du th\'eor\`eme 2.1.
\vs
\vs

\item[(ii)] Traitons le cas o\`u $\tilde{\mathscr{L}}$ est
engendr\'ee par les trois champs
$$z_1\frac{\partial}{\partial z_0}+(z_2+z_3)\frac
{\partial}{\partial z_2}+z_3\frac{\partial}{\partial z_3},
\hspace{2mm} (z_2+\kappa z_3)\frac{\partial}{\partial z_2}
+z_3\frac{\partial}{\partial z_3} \hspace{2mm} \mbox{et } R,
\hspace{2mm} \kappa\not =1$$ correspondant \`a la deuxi\`eme
configuration de matrices. La
$1$-forme $\omega$ annul\'ee par ces trois champs s'\'ecrit
$$(1-\kappa)z_1z_3^2dz_0+((\kappa-1) z_0 -\kappa
z_1)z_3^2dz_1-z_1^2z_3dz_2+(z_2+\kappa z_3)z_1^2dz_3$$

\vs

Dans ce cas, on constate que si $\kappa=0$, alors
$$\frac{\omega}{z_1^2z_3^2} =d\left(\frac{z_0}{z_1}-
\frac{z_2}{z_3}\right)$$
qui correspond au premier cas de Ab (ii).

\vs

En fait pour tout $\kappa$, $z_1^2z_3^2$ est un facteur
int\'egrant de $\omega$. On v\'erifie que
$$-\frac{\omega}{z_1^2z_3^2}=(\kappa-1) d\left(\frac{z_0}
{z_1}\right)+d\left(\frac{z_2}{z_3}\right)+\kappa\left(
\frac{dz_1}{z_1}-\frac{dz_3}{z_3}\right), \hspace{3mm} 
\kappa \not =1$$ C'est le cas g\'en\'eral de Ab (ii) du
th\'eor\`eme 2.1. Une int\'egrale
premi\`ere du feuilletage d\'ecrit par $\omega$ est $$(\kappa -1)
\frac{z_0}{z_1}+\frac{z_2}{z_3}+\kappa \log\left(\frac{z_1}
{z_3}\right) \mbox{ ou son exponentielle }$$ autrement dit, 
\`a conjugaison
pr\`es $$\frac{z_1z_2-z_0z_3}{z_1z_3} \hspace{2mm} \mbox{ si }
\kappa=0$$ $$\frac{z_1}{z_3} \exp\left(\frac{z_0z_3+z_2z_1}{z_1
z_3}\right) \hspace{2mm} \mbox{ si } \kappa \not =0$$ 

\vs

On note que lorsque $\kappa=1$ le feuilletage correspondant
est de degr\'e $1$.
\vs

\item[(iii)] Cas o\`u $\tilde{\mathscr{L}}$ est engendr\'ee par
les trois champs
$$\kappa z_0\frac{\partial}{\partial z_0} + z_2\frac{\partial}{\partial z_1}
+z_3\frac{\partial}{\partial z_2}, \hspace{2mm}
z_0\frac{\partial}{\partial z_0} +z_3\frac{\partial}{\partial z_1}
\hspace{2mm} \mbox{et } R$$ correspondant \`a la troisi\`eme
configuration de matrices. La $1$-forme $\omega$ annul\'ee par
ces trois champs s'\'ecrit \`a multiplication pr\`es
$$-z_3^3dz_0+z_0z_3^2dz_1+(\kappa z_0z_3-z_0z_2)z_3dz_2$$ $$+(z_0
z_3^2-\kappa z_0z_2z_3-z_0z_1z_3+z_0z_2^2)dz_3$$ Ici encore on
trouve un facteur int\'egrant de $\omega$ \'evident : $z_0z_3^3$.
On a $$\frac{\omega}{z_0z_3^3}=\frac{dz_0}{z_0}-\frac{dz_3}{z_3}
+d\left(\frac{z_2^2-2\kappa z_2z_3-2z_1z_3}{2z_3^2}\right)$$ Une
int\'egrale premi\`ere de cette forme et du feuilletage
correspondant est donn\'ee par $$\frac{z_0}{z_3}\exp \left(
\frac{z_2^2-2z_1z_3-2\kappa z_2z_3}{2z_3^2}\right)$$ soit \`a conjugaison pr\`es
$$\frac{z_0}{z_3}\exp \left(\frac{z_2^2-z_1z_3}{z_3^2}\right)$$
C'est le cas Ab (iii).
\vs

\item[(iv)] Cas o\`u $\tilde{\mathscr{L}}$ est engendr\'ee par les
trois champs
$$z_1
\frac{\partial}{\partial z_1}+z_3\frac{\partial}{\partial z_2},
\hspace{2mm} z_0 \frac{\partial}{\partial z_0}+\kappa z_1
\frac{\partial}{\partial z_1} \hspace{2mm} \mbox{et } R.$$  La
$1$-forme $\omega$ annul\'ee par ces trois champs s'\'ecrit \`a
multiplication pr\`es
$$-\kappa z_1z_3^2dz_0+z_0z_3^2dz_1-z_0z_1z_3dz_2+z_0z_1((\kappa-1)
z_3+z_2)dz_3$$ En cherchant une famille d'applications lin\'eaires
laissant le feuilletage invariant, on construit une sym\'etrie et
on trouve que la $1$-forme
$\omega$ admet pour facteur int\'egrant $z_0z_1z_3^2$ ; ce que
l'on peut constater de visu :
$$\frac{\omega}{z_0z_1z_3^2}=-\kappa\frac{dz_0}{z_0}+\frac{dz_1}{z_1}+
(\kappa-1)\frac{dz_3}{z_3}-d\left(\frac{z_2}{z_3}\right)$$ Elle
d\'ecrit un $\mathscr{L}$-feuilletage de degr\'e $2$ sur
$\C \P(3)$ dont une int\'egrale premi\`ere est
$$\frac{z_1z_3^{\kappa-1}}{z_0^{\kappa}} \exp\left(
\frac{z_2}{z_3}\right)$$ Il s'agit du cas Ab (iv).
\vs

\item[(v)] Cas o\`u $\tilde{\mathscr{L}}$ est engendr\'ee par les
trois champs
$$z_1\frac{\partial}{\partial z_0}
+z_2\frac{\partial}{\partial z_1}+z_3\frac{\partial}{\partial
z_2}, \hspace{2mm} (z_2+\kappa z_3)\frac{\partial}{\partial
z_0}+z_3 \frac{\partial}{\partial z_1} \hspace{2mm} \mbox{et }
R$$ correspondant \`a la cinqui\`eme configuration de matrices. La
 $1$-forme $\omega$ annul\'ee par ces trois champs s'\'ecrit
$$-z_3^3dz_0+z_3^2(z_2+\kappa z_3)dz_1+z_3(z_1z_3-z_2^2-\kappa
z_2z_3)dz_2$$ $$+(z_0z_3^2-2z_1z_2z_3+z_2^3+\kappa z_2^2z_3-\kappa
z_1z_3^2)dz_3$$ 

\vs

Le polyn\^ome $z_3^4$ est un facteur int\'egrant de $\omega$ et
$$\frac{\omega}{z_3^4}=d\left(-\frac{z_0}{z_3}+\kappa \frac{z_1}{z_3}
-\kappa\frac{z_2^2}{2z_3^2}-\frac{z_2^3}{3z_3^3}+\frac{z_1z_2}{z_3^2}\right)$$
Cette forme d\'ecrit un $\mathscr{L}$-feuilletage de degr\'e $2$ sur
$\C \P(3)$ dont le lieu singulier est une droite
$$\mathscr{S}ing \hspace{1mm} \mathscr{F}=(z_2=z_3=0)$$ et
d'int\'egrale premi\`ere
$$\frac{-z_0z_3^2+\kappa z_1z_3^2+z_1z_2z_3-\frac{z_2^3}{3}-
\frac{\kappa z_2^2z_3}{2}}{z_3^3}$$
Celle-ci s'\'ecrit \`a conjugaison
et composition \`a gauche par une homographie pr\`es
$$\frac{z_0z_3^2+z_1z_2z_3+z_2^3}{z_3^3}$$ C'est le cas Ab (v).
Dans la carte affine
$z_3=1$, cette int\'egrale premi\`ere s'\'ecrit
$$C(z_0,z_1,z_2)=z_0+z_1z_2+z_2^3$$ ce qui nous donne un polyn\^ome de degr\'e
$3$ induisiant un $\mathscr{L}$-feuilletage dans $\C^3$. On constate que
pour tout $\rho \in \C$, on a
$$C(\rho^3z_0,\rho^2z_1,\rho z_2)=\rho^3 C(z_0,z_1,z_2),$$ en
particulier pour $\rho=e^t$, on a $$C(e^{3t}z_0,e^{2t}z_1,e^t
z_2)=e^{3t} C(z_0,z_1,z_2).$$ Autrement dit le champ lin\'eaire
$$3z_0\frac{\partial}{\partial z_0}+2z_1\frac{\partial}{\partial
z_1}+z_2\frac{\partial}{\partial z_2}$$ est une sym\'etrie du
feuilletage $\mathscr{F}$ : le feuilletage $\mathscr{F}$ admet \`a
la fois une int\'egrale premi\`ere rationnelle non triviale et une
sym\'etrie. Chaque cubique $z_0+z_1z_2+z_2^3=$ cte est une surface
r\'egl\'ee : les droites d'\'equations
$$\left \{
\begin{array}{c}
  z_2=\mu \hspace{10mm} \\
  z_0+\alpha z_1=\lambda \\
\end{array} \right.$$ sont contenues dans cette surface. Lorsqu'on coupe
$\C \P(5)$ par le 2-plan $\{z_5=z_2,\hspace{1mm}
z_1=~z_3\}$, on constate que le polyn\^ome de Gordan-N\oe ther
homog\'en\'eis\'e
$$\frac{z_0z_1^2+z_1z_2z_4+z_2^2z_5}{z_3^3}$$ co\"{\i}ncide avec
l'int\'egrale premi\`ere du $\mathscr{L}$-feuilletage
$\mathscr{F}$. Outre le fait que le lieu singulier soit une droite
exactement ce dernier fait est remarquable.
\vs
\end{itemize}

\textbf{Probl\`eme} : classifier les $\mathscr{L}$-feuilletages de
$\C\P(n)$ qui ont pour ensemble singulier un sous-espace
lin\'eaire de codimension $2$. Les feuilletages de degr\'e $0$ en
font partie ; pour $n=3$, y a-t-il d'autres exemples que le
pr\'ec\'edent ? \vs

\subsection{L'alg\`ebre $\mathscr{L}$ est isomorphe \`a l'alg\`ebre
du groupe des transformations affines.}

Pour classifier les feuilletages de degr\'e $2$ associ\'es \`a une
alg\`ebre de Lie isomorphe \`a l'alg\`ebre du groupe des
transformations affines, on va classifier les alg\`ebres de Lie
engendr\'ees par deux champs lin\'eaires $X$ et $Y$ de
$\C^4$ satisfaisant $[X,Y]=Y$. Comme l'alg\`ebre $<X,Y>$
est r\'esoluble, par triangulation on constate que $Y$ est
nilpotent. D'apr\`es le lemme \ref{elimination} le rang de la 
matrice associ\'ee au champ $Y$ est $2$ ou $3$, le cas o\`u $Y$
est de rang $1$ n'est donc pas retenu.
\vs

Commen\c{c}ons par supposer que le rang de la matrice associ\'ee
au champ $Y$ est $3$, notons $-A$ la matrice de $X$ et $B$ celle
de $Y$. Comme $B$ est nilpotente de rang $3$ elle s'\'ecrit \`a
conjugaison pr\`es $$\left(%
\begin{array}{cccc}
  0 & 1 & 0 & 0 \\
  0 & 0 & 1 & 0 \\
  0 & 0 & 0 & 1 \\
  0 & 0 & 0 & 0 \\
\end{array}%
\right)$$ On remarque que la matrice $$A_0=\left(%
\begin{array}{cccc}
  3 & 0 & 0 & 0 \\
  0 & 2 & 0 & 0 \\
  0 & 0 & 1 & 0 \\
  0 & 0 & 0 & 0 \\
\end{array}%
\right)$$ satisfait $A_0B-BA_0=B$ ; par suite la matrice $A-A_0$
commute \`a $B$. Du calcul direct du commutateur de $B$,
on d\'eduit que $A$ est du type $$\left(%
\begin{array}{cccc}
  \lambda+3 & \kappa_1 & \kappa_2 & \kappa_3 \\
  0 & \lambda+2 & \kappa_1 & \kappa_2 \\
  0 & 0 & \lambda+1 & \kappa_1 \\
  0 & 0 & 0 & \lambda \\
\end{array}%
\right)$$ En particulier $A$ est diagonalisable donc s'\'ecrit \`a
conjugaison pr\`es $$\left(%
\begin{array}{cccc}
  \lambda+3 & 0 & 0 & 0 \\
  0 & \lambda+2 & 0 & 0 \\
  0 & 0 & \lambda+1 & 0 \\
  0 & 0 & 0 & \lambda \\
\end{array}%
\right)$$ Une matrice $B$ satisfaisant $AB-BA=B$ est alors de la
forme
$$\left(%
\begin{array}{cccc}
  0 & \beta_1 & 0 & 0 \\
  0 & 0 & \beta_2 & 0 \\
  0 & 0 & 0 & \beta_3 \\
  0 & 0 & 0 & 0 \\
\end{array}%
\right)$$ Comme $B$ est suppos\'e \^etre de rang $3$, les
$\beta_i$ sont non nuls donc \'egaux \`a $1$ \`a conjugaison
pr\`es. Finalement le $\mathscr{L}$-feuilletage correspondant est
dans des coordonn\'ees homog\`enes bien choisies tangent aux
champs
$$R=\sum_{i=0}^3 z_i \frac{\partial}{\partial z_i}, \hspace{2mm}
X=3z_1 \frac{\partial}{\partial z_1}+2z_2\frac{\partial}{\partial
z_2}+ z_3\frac{\partial}{\partial z_3}$$
$$\mbox{ et } Y=z_1\frac{\partial}{\partial z_0}+z_2\frac{\partial}{\partial
z_1}+z_3\frac{\partial}{\partial z_2}$$ Il s'agit du feuilletage
exceptionnel qui admet pour int\'egrale premi\`ere
$$\frac{{\left(z_0z_3^2-z_1z_2z_3+\frac{z_2^3}{3}
\right)}^2}{{\left(z_1z_3-\frac{z_2^2}{2}\right)}^3}$$ C'est le cas Af (i).
\vs

Supposons maintenant que la matrice associ\'ee au champ $Y$ soit
de rang $2$. Il y a \`a conjugaison pr\`es deux types de matrices
$(4,4)$ nilpotentes de rang $2$

\vs

$$B_1=\left(%
\begin{array}{cccc}
  0 & 0 & 0 & 0 \\
  0 & 0 & 1 & 0 \\
  0 & 0 & 0 & 1 \\
  0 & 0 & 0 & 0 \\
\end{array}%
\right) \mbox{ et } B_2=\left(%
\begin{array}{cccc}
  0 & 1 & 0 & 0 \\
  0 & 0 & 0 & 0 \\
  0 & 0 & 0 & 1 \\
  0 & 0 & 0 & 0 \\
\end{array}%
\right) \sim \left(%
\begin{array}{cccc}
  0 & 0 & 1 & 0 \\
  0 & 0 & 0 & 1 \\
  0 & 0 & 0 & 0 \\
  0 & 0 & 0 & 0 \\
\end{array}%
\right)$$\vs

\begin{itemize}
\item[2.2.1.] Cas $B=B_1$.

L'\'egalit\'e $[A,B]=B$ entra\^ine que $A$ s'\'ecrit sous la forme $$\left(%
\begin{array}{cccc}
  \lambda_1 & 0 & 0 & \kappa_4 \\
  \kappa_5 & 2+\lambda & \kappa_7 & \kappa_8 \\
  0 & 0 & 1+\lambda & \kappa_7 \\
  0 & 0 & 0 & \lambda \\
\end{array}%
\right)$$ Puisqu'on s'int\'eresse \`a l'alg\`ebre engendr\'ee par
les champs $X$, $Y$ et $R$ on peut supposer que $\kappa_7=0$.
Remarquons qu'avec les notations habituelles un changement de base
$e_0 \rightarrow e_0+\alpha e_1$ n'alt\`ere pas la matrice $B$.
Par suite si $\lambda_1 \not =2+\lambda$, on peut supposer que $A$
est du type $$A_1=\left(%
\begin{array}{cccc}
  \lambda_1 & 0 & 0 & \kappa_4 \\
  0 & 2+\lambda & 0 & \kappa_8 \\
  0 & 0 & 1+\lambda & 0 \\
  0 & 0 & 0 & \lambda \\
\end{array}%
\right)$$ et sinon du type $$A_2=\left(%
\begin{array}{cccc}
  2+\lambda & 0 & 0 & \kappa_4 \\
  \kappa_5 & 2+\lambda & 0 & \kappa_8 \\
  0 & 0 & 1+\lambda & 0 \\
  0 & 0 & 0 & \lambda \\
\end{array}%
\right)$$

\begin{itemize}
\item[a. ]Si $\lambda_1 \not =2+\lambda$, on fait le changement de base $$e_4
\rightarrow e_4+ae_1+be_2=e'_4$$ (qui n'alt\`ere pas $B$) pour obtenir
$Ae'_4=\lambda e'_4$ soit encore
$$a(\lambda_1-\lambda)+\kappa_4=0$$ $$2b+\kappa_8=0$$
Ainsi si $\lambda \not = \lambda_1$, on peut par
conjugaison \'eliminer $\kappa_4$ et $\kappa_8$, sinon on peut
\'eliminer $\kappa_8$. \vs

\begin{itemize}
\item[a.1. ] Si $\lambda \not = \lambda_1$, on se ram\`ene \`a une
alg\`ebre du type $$R, \hspace{2mm} Y=z_2\frac{\partial} {\partial
z_1}+z_3\frac{\partial}{\partial z_2}, \hspace{2mm} X= \kappa
z_0\frac{\partial}{\partial z_0}+2z_1\frac{\partial} {\partial
z_1}+z_2\frac{\partial}{\partial z_2}$$

Les hyperplans $z_3=$ cte sont invariants par les deux champs $X$
et $Y$ ; ainsi ces deux champs engendrent le feuilletage restreint
\`a la carte affine $z_3=1$ ; ils s'\'ecrivent
$$\tilde{X}=X_{|z_3=1}=\kappa z_0
\frac{\partial}{\partial z_0}+2z_1\frac{\partial}{\partial z_1}+
z_2\frac{\partial}{\partial z_2}$$
$$\tilde{Y}=Y_{|z_3=1}=z_2\frac{\partial}{\partial z_1}+
\frac{\partial}{\partial z_2}$$ Le champ $\tilde{Y}$ a une
\guillemotleft \hspace{1mm} partie constante \guillemotright
\hspace{1mm}; nous allons essayer de le redresser par un automorphisme
polynomial. Le flot de $\tilde{Y}$ est
$$\varphi \hspace{1mm} \colon \hspace{1mm} (z;t) \mapsto
(z_0,z_1+z_2t+\frac{t^2}{2},z_2+t)$$ On d\'efinit maintenant le
diff\'eomorphisme $$H \hspace{1mm} \colon \hspace{1mm}
(z_0,z_1,z_2) \mapsto
\varphi(z_0,z_1,0;z_2)=(z_0,z_1+\frac{z_2^2}{2},z_2)$$ Par
construction $H$ conjugue {\large $\frac{\partial}{\partial z_2}$}
\`a $\tilde{Y}$. Le feuilletage $\mathscr{F}'$ associ\'e \`a
l'alg\`ebre engendr\'ee par les champs $H_*^{-1}\tilde{X}$ et
$H_*^{-1}\tilde{Y}$ est d\'ecrit par une $1$-forme $\Omega'$ qui,
en particulier, annule le champ $H_*^{-1}\tilde{Y}=${\large
$\frac{\partial}{\partial z_2}$} donc s'\'ecrit $$A'(z_0,z_1)
dz_0 +B'(z_0,z_1)dz_1$$ Le diff\'eomorphisme $H$ laisse le plan
$z_2=0$ invariant point par point. On constate que le champ
$$\tilde{X}_{|z_2=0}=\kappa z_0\frac{\partial}{\partial z_0}+2z_1
\frac{\partial}{\partial z_1}$$ est tangent au feuilletage et au
plan $z_2=0$ ; il d\'ecrit donc le feuilletage en restriction au plan
$z_2=0$. On d\'eduit une int\'egrale premi\`ere du champ
$\tilde{X}$ $$\frac{z_0^2}{z_1^\kappa}$$ et donc de $\mathscr{F}'$. Ce
qui donne apr\`es composition par $H^{-1}$ et homog\'en\'eisation l'int\'egrale
premi\`ere de $\mathscr{F}$ suivante :
$$\frac{z_0^2z_3^{2(\kappa-1)}}{{\left(z_1z_3-\frac{z_2^2}{2}
\right)}^\kappa}$$ C'est le cas Af (ii).
\vs

\item[a.2. ] Si $\lambda=\lambda_1$, l'alg\`ebre
$\tilde{\mathscr{L}}$ est engendr\'ee par
$$R, \hspace{2mm} z_2\frac{\partial}{\partial z_1}+z_3
\frac{\partial}{\partial z_2}, \hspace{2mm} \kappa_4z_3
\frac{\partial}{\partial z_0}+2z_1\frac{\partial}{\partial z_1}+
z_2\frac{\partial}{\partial z_2}$$

Si $\kappa_4=0$ alors les champs {\large $\frac{\partial}{\partial
z_1}$} et {\large $\frac{\partial}{\partial z_2}$} sont tangents
au feuilletage qui ne peut donc \^etre un feuilletage de degr\'e
$2$ sur $\C \P(3)$ ; par homoth\'etie on peut donc
se ramener \`a $\kappa_4=1$.

Par la m\^eme m\'ethode que pr\'ec\'edemment, en redressant $z_2
\frac{\partial}{\partial z_1}+z_3\frac{\partial}{\partial z_2}$
dans la carte affine $z_3=1$, on montre que le feuilletage
$\mathscr{F}$ admet pour int\'egrale premi\`ere
$$\frac{z_3^2}{2z_1z_3-z_2^2}\exp\left(\frac{2z_0}{z_3}\right)$$
ce qui donne le cas Af (iii) \`a conjugaison pr\`es.
\vs
\end{itemize}

\item[b. ]Si $\lambda_1=2+\lambda$, par un changement de base $$e_4 \rightarrow
e_4+ae_1+be_2$$ on peut se ramener \`a $\kappa_4=\kappa_8=0$. En
effet,
$$A_2(e_4+ae_1+b
e_2)=0$$ \'equivaut \`a $$a=-\frac{\kappa_4}{2}, \hspace{2mm}
 b=\frac{a\kappa_4\kappa_5-2\kappa_8}{4}$$
Par homoth\'etie on se ram\`ene \`a $\kappa_5=0$ ou $1$. \vs

Commen\c{c}ons par traiter le cas $\kappa_5=0$ ; alors
l'alg\`ebre $\tilde{\mathscr{L}}$ est engendr\'ee par les champs
$$R, \hspace{2mm} z_2\frac{\partial}{\partial
z_1}+z_3\frac{\partial} {\partial z_2} \mbox{ et
}2z_0\frac{\partial}{\partial z_0}+ 2z_1\frac{\partial}{\partial
z_1}+z_2\frac{\partial}{\partial z_2}$$ qui est un cas particulier
du cas a.1. ($\kappa=2$)
\vs

Finalement traitons le cas $\kappa_5=1$ ; l'alg\`ebre
$\tilde{\mathscr{L}}$ est du type $$R, \hspace{2mm}
z_2\frac{\partial}{\partial z_1}+z_3\frac{\partial} {\partial z_2}
\mbox{ et }2z_0\frac{\partial}{\partial z_0}+(z_0+2z_1)
\frac{\partial}{\partial z_1}+z_2\frac{\partial} {\partial z_2}
$$
Toujours par la m\^eme m\'ethode (qui consiste \`a redresser le
champ $z_2\frac{\partial}{\partial z_1}+z_3\frac{\partial}
{\partial z_2}$ dans la carte affine $z_3=1$) on montre que le
feuilletage $\mathscr{F}$ admet pour int\'egrale premi\`ere
$$\frac{z_0}{z_3}\exp\left(\frac{z_2^2-2z_1z_3}{z_0z_3}\right)$$
on obtient le cas Af (iv) \`a conjugaison pr\`es.
\vs
\end{itemize}

\item[2.2.2.] Cas $B=B_2$.

Ce cas est caract\'eris\'e par $\ker B=im \hspace{1mm}B$ et $B$
est un isomorphisme de $<e_2,e_3>$ dans $<e_0 ,e_1>$. En
particulier pour chaque base $f_0,f_1$ de $\ker B$, il existe une
base $f_2,f_3$ de $<e_2,e_3>$ telle que $Bf_2=f_0$ et $Bf_3=f_1$.
L'\'egalit\'e $AB-BA=B$ implique que $\ker B$ est invariant par
$A$. On peut donc supposer que l'on
a une base $e_0,e_1,e_2,e_3$ dans laquelle $A$ et $B$ sont du type
suivant (on jordanise la restriction de $A$
\`a $\ker B$) $$A=A_1=\left(%
\begin{array}{cccc}
  \lambda_1 & 0 & * & * \\
  0 & \lambda_2 & * & * \\
  0 & 0 & * & * \\
  0 & 0 & * & * \\
\end{array}%
\right) \mbox{ avec } \lambda_1 \not =\lambda_2 \mbox{ ou } A=A_2=\left(%
\begin{array}{cccc}
  \lambda & \varepsilon & * & * \\
  0 & \lambda & * & * \\
  0 & 0 & * & * \\
  0 & 0 & * & * \\
\end{array}%
\right)$$ $$\mbox{et } B=\left(%
\begin{array}{cccc}
  0 & 0 & 1 & 0 \\
  0 & 0 & 0 & 1 \\
  0 & 0 & 0 & 0 \\
  0 & 0 & 0 & 0 \\
\end{array}%
\right)$$ On remarque que par homoth\'etie on se ram\`ene \`a
$\varepsilon \in \{0,1\}$. \vs

\begin{itemize}
\item[c. ] $A$ est du type $A_1$.

On a
$$ABe_3-BAe_3=Be_3 \Leftrightarrow BAe_3=(\lambda_2-1)Be_3$$ Il en
r\'esulte que $$Ae_3=\kappa_3e_0+\kappa_5e_1+(\lambda_2-1)e_3$$ De
la m\^eme fa\c{c}on, on obtient $$Ae_2=\kappa_2e_0+\kappa_4e_1
+(\lambda_1-1)e_2$$ et $$A=\left(%
\begin{array}{cccc}
  \lambda_1 & 0 & \kappa_2 & \kappa_3 \\
  0 & \lambda_2 & \kappa_4 & \kappa_5 \\
  0 & 0 & \lambda_1-1 & 0 \\
  0 & 0 & 0 & \lambda_2-1 \\
\end{array}%
\right)$$ On note que les changements de base du type $$e_2
\rightarrow e_2 +\alpha e_0+\beta e_1=\underline{e_2}$$ $$e_3
\rightarrow e_3 + \gamma e_0+\delta e_1=\underline{e_3}$$
n'alt\`erent ni $B$, ni la forme de $A$. On a avec les notations
pr\'ec\'edentes $$A\underline{e_2}-(\lambda_1-1)\underline{e_2}=
(\kappa_2+\alpha)e_0+(\beta(\lambda_2-\lambda_1+1)+\kappa_4)e_1$$
$$A\underline{e_3}-(\lambda_2-1)\underline{e_3}=
(\kappa_3+\gamma(\lambda_1-\lambda_2+1))e_0+(\delta+\kappa_5)e_1$$
On constate qu'en choisissant $\delta=-\kappa_5$ et
$\alpha=-\kappa_2$, on peut supposer que $\kappa_2=\kappa_5=0$ ;
de plus si $\lambda_2-\lambda_1 \not = \pm 1$, on peut supposer
que $\kappa_3=\kappa_4=0$.

En enlevant $(\lambda_2-1)R$ \`a $X$ repr\'esent\'e par la matrice
$A$, on se ram\`ene au cas o\`u l'alg\`ebre $\tilde{\mathscr{L}}$
est engendr\'ee par $R$ et les trois champs $$
z_2\frac{\partial}{\partial z_0}+z_3\frac{\partial}{\partial
z_1},\hspace{2mm} ((\mu+1)z_0+
\kappa_3z_3)\frac{\partial}{\partial z_0}+(z_1+\kappa_4z_2)
\frac{\partial}{\partial z_1}+\mu z_2\frac{\partial}{\partial
z_2}$$ o\`u $\mu=\lambda_1-\lambda_2$ ($\not =0$ par hypoth\`ese).

On va suivre la m\'ethode utilis\'ee pr\'ec\'edemment mais nous la
d\'etaillerons car elle est ici un peu plus technique. Comme
l'hyperplan $z_3=0$ est invariant par tous ces champs on se place
dans la carte affine $z_3=1$ o\`u l'on a $$\tilde{Y}=Y_{|z_3
=1}=z_2\frac{\partial}{\partial z_0}+\frac{\partial}{\partial
z_1}$$ $$\tilde{X}=X_{|z_3=1}=((\mu+1)z_0+ \kappa_3)
\frac{\partial}{\partial z_0}+(z_1+\kappa_4z_2) \frac{\partial}
{\partial z_1}+\mu z_2\frac{\partial}{\partial z_2}$$ L'hyperplan
$z_1=0$ est transverse \`a $\tilde{Y}$ donc au feuilletage
$\mathscr{F}$. Notons $$\varphi(z_0,z_1,z_2;t)=(z_0+tz_2,z_1+t,
z_2)$$ le flot de $\tilde{Y}$, alors le diff\'eomorphisme
$$H(z_0,z_1,z_2)=\varphi(z_0,0,z_2;z_1)=(z_0+z_1z_2,z_1,z_2)$$
redresse le champ $\tilde{Y}$ sur {\large
$\frac{\partial}{\partial z_1}$} et trivialise le feuilletage.
Dans l'hyperplan $z_1=0$, le feuilletage $\mathscr{F}_{|z_1=0}$
est donn\'e par
$$(\tilde{X}-\kappa_4z_2 \tilde{Y})_{|z_1=0}=((\mu+1)z_0+\kappa_3
-\kappa_4z_2^2) \frac{\partial}{\partial z_0}+\mu
z_2\frac{\partial}{\partial z_2}$$ La $1$-forme
$$\overline{\omega}=\mu z_2dz_0-((\mu+1)z_0
+\kappa_3-\kappa_4z_2^2)dz_2$$ d\'ecrit donc $\mathscr{F}_{|z_1=0}$.
\vs

\begin{itemize}
\item[c.1. ]Si $\mu=-1$, le feuilletage restreint est d\'ecrit par la
 $1$-forme $$d\left(z_0-\frac{\kappa_4}{2}z_2^2\right)+\kappa_3
 \frac{dz_2}{z_2}$$ d'o\`u une int\'egrale premi\`ere, si $\kappa_3
  \not =0$, du type :
$$z_2\exp\left(\frac{z_0} {\kappa_3}-\frac{\kappa_4}{2\kappa_3}z_2
^2\right)$$ % soit \`a conjugaison pr\`es $z_2\exp(z_0-z_2^2)$
ce qui donne apr\`es composition avec $H^{-1}$ et
homog\'en\'eisation
$$\frac{z_2}{z_3}\exp\left(\frac{2z_0z_3-2z_1z_2-\kappa_4z_2^2}
{2\kappa_3z_3^2}\right)$$ qui correspond au cas Af (v) \`a
conjugaison pr\`es. Si $\kappa_3=0$, le feuilletage est de degr\'e
$1$.
\vs

\item[c.2. ] Si $\mu \not =\pm 1$, puisqu'on peut se ramener \`a
$\kappa_3=\kappa_4=0$, on a $$\overline{\omega} =\mu z_2dz_0-
(\mu+1)z_0dz_2$$ On en d\'eduit,
pour le feuilletage restreint \`a $z_1=0$, l'int\'egrale premi\`ere
suivante
$$\frac{z_0^\mu}{z_2^{\mu+1}}$$ Apr\`es composition par
$H^{-1}$ et homog\'en\'eisation, le feuilletage admet
pour int\'egrale premi\`ere
$$\frac{{\left(z_0z_3-z_1z_2\right)}^\mu}{z_2^{\mu+1}z_3^{\mu-1}}$$
qui donne le cas Af (vi).
\vs

\item[c.3. ] Si $\mu=1$, la $1$-forme $\overline{\omega}$ s'\'ecrit
$$z_2dz_0-(2z_0+\kappa_3-\kappa_4z_2^2)dz_2$$ Quitte \`a composer
$\overline{\omega}$ par $$\phi \hspace{1mm} \colon
\hspace{1mm}(z_0,z_2) \mapsto (z_0-\frac{\kappa_3}{2},z_2)$$ on a
$$\overline{\omega}=z_2dz_0-(2z_0-\kappa_4
z_2^2) dz_2$$ ou encore, \`a multiplication pr\`es : $$d\left(
\frac{z_0}{z_2^2}\right)+\kappa_4\frac{dz_2}{z_2}$$ d'o\`u 
l'int\'egrale premi\`ere ($\kappa_4 \not=0$ sinon le degr\'e chute)
$$z_2\exp\left( \frac{z_0}{\kappa_4z_2^2}\right)$$ Ce qui donne apr\`es
composition avec $\phi^{-1}$, $H^{-1}$ et homog\'en\'eisation
$$\frac{z_2}{z_3}
\exp\left(\frac{2z_0z_3-2z_1z_2+\kappa_3z_3^2}{2\kappa_4z_2^2}\right)$$
Quitte \`a permuter $z_2$ et $z_3$ puis $z_0$ et $z_1$ on retombe
sur le cas c.1.
\vs
\end{itemize}

\item[d. ] Supposons que $A$ est du type $A_2$.

Si $\varepsilon=0$, alors $\mu=0$ et en reprenant les calculs
pr\'ec\'edents, on constate que le degr\'e du feuilletage chute.
\vs

Si $\varepsilon=1$, alors on a $$ABe_2-BAe_2=Be_2
\Leftrightarrow BAe_2=(\lambda-1)Be_2$$ Il en r\'esulte que
$$Ae_2=\kappa_2e_0+\kappa_4e_1$$ De
la m\^eme fa\c{c}on, on obtient $$Ae_3=\kappa_3e_0+\kappa_5e_1
+e_2+(\lambda-1)e_3$$ et $$A=\left(%
\begin{array}{cccc}
  \lambda & 1 & \kappa_2 & \kappa_3 \\
  0 & \lambda & \kappa_4 & \kappa_5 \\
  0 & 0 & \lambda-1 & 1 \\
  0 & 0 & 0 & \lambda-1 \\
\end{array}%
\right) $$ On note que les changements de base du type $$e_2
\rightarrow e_2 +\alpha e_0+\beta e_1=\underline{e_2}$$ $$e_3
\rightarrow e_3 + \gamma e_0+\delta e_1=\underline{e_3}$$
n'alt\`erent ni $B$, ni la forme de $A$. On a $$A\underline{e_2}
-(\lambda-1)\underline{e_2}=(\kappa_2+\alpha+\beta)e_1+\kappa_4+
\beta$$ Donc quitte \`a prendre $\beta=-\kappa_4$ et
$\alpha=\kappa_4-\kappa_2$, on peut supposer que
$\kappa_2=\kappa_4=0$. On a alors $$A\underline{e_3}-(\lambda-1)
\underline{e_3}-\underline{e_2}=(\kappa_3+\delta+\gamma)e_1+
(\kappa_5+\delta)e_2$$ Puis quitte \`a prendre
$\delta=-\kappa_5$ et $\gamma=\kappa_5-\kappa_3$, on se ram\`ene
\`a $\kappa_3=\kappa_5=0$. Finalement \`a conjugaison pr\`es $A$
s'\'ecrit $$\left(%
\begin{array}{cccc}
  \lambda & 1 & 0 & 0 \\
  0 & \lambda & 0 & 0 \\
  0 & 0 & \lambda-1 & 1 \\
  0 & 0 & 0 & \lambda-1 \\
\end{array}%
\right)$$ et notre feuilletage est tangent aux champs de vecteurs
$$R, \hspace{2mm} Y=z_2\frac{\partial}{\partial z_0}+z_3
\frac{\partial}{\partial z_1} \mbox{ et } X=(z_0+z_1)
\frac{\partial}{\partial z_0}+z_1\frac{\partial}{\partial z_1}+
z_3\frac{\partial}{\partial z_2}$$ Comme l'hyperplan $z_3=0$ est
invariant par tous ces champs on se place dans la carte affine
$z_3=1$ o\`u %l'on a $$\tilde{Y}=Y_{|z_3=1}=z_2 \frac{\partial}
%{\partial z_0}+ \frac{\partial}{\partial z_1}$$
%$$\tilde{X}=X_{|z_3=1}=(z_0+z_1)\frac{\partial}{\partial z_0}
%+z_1\frac{\partial}{\partial z_1}+\frac{\partial}{\partial z_2}$$
la restriction de notre feuilletage \`a $z_1=0$ est donn\'e par le
champ $$z_0\frac{\partial}{\partial z_0}
+\frac{\partial}{\partial z_2}$$ autrement dit par l'int\'egrale
premi\`ere $$z_0\exp(-z_2)$$
Apr\`es composition par $H^{-1}$ et homog\'en\'eisation, on
constate que notre feuilletage admet pour int\'egrale premi\`ere
$$\frac{z_0z_3-z_1z_2}{z_3^2} \exp\left(-\frac{z_2}{z_3}\right)$$
c'est le cas Af (vii).
\vs
\end{itemize}
\end{itemize}

%%%%%%%%%%%%%%%%%%%%%%%%%%%%%%%%%%%%%%%%%%%%%%%%%%%%%%%%%%%%%%%%%%%%
%%%%%%%%%%%%%%%%%%%%%%%%%%%%%%%%%%%%%%%%%%%%%%%%%%%%%%%%%%%%%%%%%%%%

\chapter{$\mathscr{L}$-feuilletages quadratiques sur $\C
\P(n)$.}

\vs

L'\'etude des $\mathscr{L}$-feuilletages de codimension $1$ sur
$\C\P(n)$ pour $n\geq 4$ s'av\`ere d\'elicate, non
seulement pour des raisons de taille de calculs mais aussi pour
des raisons conceptuelles. Bien s\^ur nous savons qu'en degr\'es
$0$ et $1$ tous les feuilletages de codimension $1$ sont des
$\mathscr{L}$-feuilletages. Par contre d\'ej\`a en degr\'e $2$ on
se heurte au probl\` eme suivant : majorer la dimension de
$\mathscr{L}$. Nous avons vu qu'en degr\'e maximal $n-1$ la
dimension de $\mathscr{L}$ est pr\'ecis\'ement $n-1$ ce qui permet
tout du moins en petite dimension  d'esp\'erer une classification
en se r\'eferrant \`a la description des alg\`ebres de Lie. Nous
la pr\'esentons en dimension $4$. %La taille des calculs explose
%de la dimension $3$ \`a la dimension $4$ et nous avons utilis\'e
%pour certains d'entre eux le logiciel Maple. Chaque fois qu'un
%r\'esultat d\'epend de cette utilisation nous l'indexerons par une
%ast\'erisque : lemme$^*$ ou proposition$^*$. M\^eme si l'on ne
%peut en toute rigueur consid\'erer qu'il s'agit l\`a de
%r\'esultats au sens classique, on obtient tout de m\^eme au
%minimum une grande liste d'exemples.

\vs

\section{$\mathscr{L}$-feuilletages quadratiques sur $\C \P(n)$,
$n \geq 4$.}

\vs

Nous donnons quelques propri\'et\'es g\'en\'erales susceptibles
d'\^etre utiles pour classifier les $\mathscr{L}$-feuilletages de
degr\'e $2$. La plupart sont cons\'equences de la description
g\'en\'erale des feuilletages de degr\'e $2$ sur $\C
\P(n)$ faite dans \cite{[Ce-LN]} et que nous avons pr\'esent\'ee
dans le cas $n=3$ (chapitre $4$). Nous donnerons aussi
quelques exemples. \vs

\begin{pro} Soit $\mathscr{F}$ un feuilletage de
degr\'e $2$ sur $\C \P(n)$ d\'efini en
coordonn\'ees homog\`enes par la $1$-forme $\omega \in
\Omega_3^1(\C^{n+1})$. On est dans l'une des situations
suivantes :

1. $\mathscr{F}$ est exceptionnel, i.e. il existe $\tau
\hspace{1mm} \colon \hspace{1mm} \C \P(n)
\dashrightarrow \C \P(3)$
 lin\'eaire telle que $\mathscr{F}=\tau^*\mathscr{F}_\Gamma$

2. Il existe un feuilletage $\mathscr{F}'$ de degr\'e $2$ sur
$\C\mathbb {P}(2)$, une application lin\'eaire
$\tau\hspace{1mm}\colon\hspace{1mm} \C \P(n)
\dashrightarrow \C \P(2)$ tels que
$\mathscr{F}=\tau^*\mathscr{F}'$.

3. $\mathscr{F}$ a une int\'egrale premi\`ere du type
$\frac{Q_1}{Q_2}$ o\` u $Q_1$ et $Q_2$ sont des formes
quadratiques.

4. $\mathscr{F}$ poss\`ede un hyperplan invariant $\mathscr{H}$ et
$\omega$ un facteur int\'egrant de l'un des types suivants :

\begin{itemize}
\item[4.1. ] $L_1^4$

\item[4.2. ] $L_1^3L_2$

\item[4.3. ] $L_1^2L_2^2$

\item[4.4. ] $L_1^2L_2L_3$

\item[4.5. ] $L_1L_2L_3L_4$

\item[4.6. ] $L_1^2Q$

\item[4.7. ] $L_1L_2Q$

\item[4.8. ] $L_1C$

o\`u les $L_i$ d\'esignent des formes lin\'eaires distinctes, $Q$
une forme quadratique et $C$ une cubique.
\vs
\end{itemize}
\end{pro}

\begin{proof}[D\'emonstration]
C'est une application directe de \cite{[Ce-LN]}.
\end{proof}

\vs

\begin{rem}
A partir des facteurs int\'egrants on peut \'evidemment donner les
diff\'erents types possibles d'int\'egrales premi\` eres.
\end{rem}

\vs

On sait que dans le cas des feuilletages exceptionnels, il s'agit
de $\mathscr{L}$-feuilletages.

Par contre la proposition \ref{obstruction} asssure que dans le cas 2,
$\mathscr{F}$ n'est jamais un $\mathscr{L}$-feuilletage.

Dans le cas 3 o\`u $\mathscr{F}$ poss\`ede une int\'egrale
premi\` ere $\frac{Q_1}{Q_2}$, alors si $\mathscr{F}$ est un
$\mathscr{L}$-feuilletage, les formes quadratiques $Q_1$ et $Q_2$
ne peuvent \^etre g\'en\'eriques. En effet, si $Q_1$ et $Q_2$ sont
g\'en\'eriques, alors la surface $Q_1=Q_2=0$ est une surface de
Del Pezzo ; ceci implique par un argument d\'ej\`a rencontr\'e
dans la preuve du lemme \ref{courbeell} qu'il ne peut y avoir de champ de
vecteurs tangent \`a
 la fois aux niveaux de $Q_1$ et de $Q_2$.

\vs

\subsection{Exemples de type 4.2.}

\vs

Ils se traitent relativement ais\'ement : \vs

\begin{pro} Soit $\mathscr{F}$ un
$\mathscr{L}$-feuilletage de degr\'e $2$ sur $\C
\P(4)$ ayant une int\'egrale premi\`ere du type $z_1\exp
H$, $\deg H=2$, dans une carte affine ad-hoc. Alors \`a
conjugaison pr\`es $H$ est de l'un des types suivants
$$z_2+z_3z_4$$ $$z_2+z_3^2
$$ $$z_2+z_1z_3+z_4^2$$\vs
\end{pro}

\begin{proof}
On remarque qu'un champ affine $X$ annulant
$\frac{dz_1}{z_1}+dH$ est n\'ecessairement du type $$X=\lambda
z_1\frac{\partial}{\partial z_1} +Z$$ o\`u $$Z \in \hspace{1mm}
<\frac{\partial} {\partial z_2}, \hspace{1mm}
\frac{\partial}{\partial z_3}, \hspace{1mm}
\frac{\partial}{\partial z_4}>$$ D'autre part pour que
$\mathscr{F}$ soit un $\mathscr{L}$-feuilletage il est
n\'ecessaire que pour l'un des $X$, $\lambda$ soit non nul disons
\'egal \`a $1$. En \'ecrivant que $$1+(\lambda z_1\frac{\partial
H}{\partial z_1}+Z(H))=0$$ on constate que $H(0,z_2,z_3,z_4)$ est un
polyn\^ome quadratique qui est n\'ecessairement une submersion.
Par suite quitte \`a changer les coordonn\'ees on peut supposer
que $$H(0,z_2,z_3,z_4)= z_2+q(z_3,z_4)$$ o\`u $q$ est une forme
quadratique en $z_3$, $z_4$. On remarque aussi que $H$ n'a pas de
terme en $z_1^2$ si bien que $$H=z_2+z_1L(z_2,z_3,z_4) +
q(z_3,z_4)$$ On examine les cas o\`u $L\equiv 0$, $L=z_2$,
$L=z_2+z_3$, $L=z_3$ qui couvrent \`a conjugaison pr\`es tous les
autres.

Dans la premi\`ere \'eventualit\'e, toujours  \`a conjugaison
pr\`es, on peut supposer que $$H=z_2+z_3z_4$$ ou bien
$$H=z_2+z_3^2$$ Dans ces deux cas, un calcul \'el\'ementaire
montre qu'il s'agit bien de $\mathscr{L}$-feuilletages, le second
apparaissant comme un pull-back lin\'eaire d'un
$\mathscr{L}$-feuilletage sur $\C\P(3)$.

Consid\'erons le cas o\`u $$H=z_2+z_1z_2+z_1z_3+q(z_3,z_4)$$ Avec les
notations pr\'ec\'edentes si $X=z_1\frac{\partial}{\partial
z_1}+Z$, avec $Z=\displaystyle \sum_{i=2}^4 A_i \frac{\partial}{\partial
z_i}$, est tangent \`a $\mathscr{F}$, on obtient
$$1+z_1z_2+A_2+A_3\left(z_1+\frac{\partial q}{\partial z_3}\right)+A_4
\frac{\partial}{\partial z_4} \equiv 0$$ En posant
$q=az_3^2+bz_3z_4+cz_4^2$ on obtient
$$(1+z_1z_2)+A_2+A_3(z_1+2az_3+bz_4)+A_4(bz_3+2cz_4)\equiv 0 $$
Sur l'intersection des deux hyperplans $$z_1+2az_3+bz_4=0$$ et
$$bz_3+ 2cz_4=0 $$ on a $$\left(1-(2az_3+bz_4)z_2\right)+f \equiv
0$$ o\`u $f$ est affine donc $a=b=0$. Finalement \`a conjugaison pr\`es on a
$$H=z_2+z_1z_3+\kappa z_4^2$$ avec $\kappa \in \{0,1\}$. On constate que
$\mathscr{F}$ est un $\mathscr{L}$-feuilletage puisqu'on trouve parmi les champs
tangents \`a $\mathscr{F}$ :
$$z_1\frac{\partial}{\partial z_1}-\frac{\partial}{\partial
z_2}-z_3\frac{\partial}{\partial z_3}, \hspace{2mm} z_1\frac{\partial}{\partial
z_2}-\frac{\partial}{\partial z_3}, \hspace{2mm} 2\kappa
z_4\frac{\partial}{\partial z_2}-\frac{\partial}{\partial z_4}$$
\end{proof}

\vs

\subsection{Un cas sp\'ecial. }

\vs

Voici une m\'ethode qui pourrait s'av\'erer fructueuse en
particulier dans l'\'etude du cas quadratique. Consid\'erons
$\mathscr{F}$ un $\mathscr{L}$-feuilletage de degr\'e $2$ sur
$\C\P(4)$ ayant une int\'egrale premi\`ere de type
$\frac{Q_1}{Q_2}$ avec $Q_i$ formes quadratiques. Supposons
l'alg\` ebre $\mathscr{L}$ non nilpotente ; alors $\mathscr{L}'$
contient un champ $X$ semi-simple : $$X=\displaystyle \sum_{i=0}^k
\lambda_iz_i\frac{\partial}{\partial z_i} \hspace{8mm} \lambda_i
\not =0$$ On peut supposer $k \geq 1$ sinon $\mathscr{F}$ serait
un pull-back. Examinons le cas o\` u $rg \hspace{1mm} X =2$, i.e.
$k=1$ : $$X=\lambda_0z_0\frac{\partial}{\partial z_0}+\lambda_1
z_1 \frac{\partial}{\partial
z_1}\hspace{8mm}\lambda_0\lambda_1\not=0$$ Soit $Q$ une forme
quadratique annul\'ee par $X$ ; alors $Q$ s'\'ecrit $$\varepsilon
z_0z_1+q(z_2,z_3,z_4)$$ Remarquons que $X(Q_1)=X(Q_2)=0$. On se
ram\`ene \`a \'etudier les cas $$\frac{Q_1}{Q_2}=
\frac{q_1}{q_2}$$ et
$$\frac{Q_1}{Q_2}=\frac{z_0z_1+q_1(z_2,z_3,z_4)}
{q_2(z_2,z_3,z_4)}$$ Dans le premier cas $\mathscr{F}$ est un
pull-back. Dans le second cas, on va faire une discussion suivant
le rang de $q_2$.

Si $rg \hspace{1mm} q_2=1$, i.e. $q_2=z_2^2$ (\`a conjugiason
pr\`es) alors $\mathscr{F}$ est en fait de degr\'e un.

Supposons $q_2$ de rang $3$, autrement dit de la forme $z_2^2+z_3^2+
z_4^2$ ; on doit trouver un champ $Y$ annulant d\'enominateur et num\'erateur.
Un tel champ s'\'ecrit $$Y=A_0\frac{\partial}{\partial z_0}+A_1
\frac{\partial}{\partial z_1}+B_2\frac{\partial}{\partial z_2}+
B_3\frac{\partial}{\partial z_3}+B_4\frac{\partial}{\partial
z_4}=A+B$$ L'\'egalit\'e $A(z_0z_1)=0$ implique que $$A=\lambda\left(z_0
\frac{\partial}{\partial z_0}-z_1\frac{\partial}{\partial z_1}\right)$$ On a
aussi $$B(q_1)=B(q_2)=0$$ En particulier on peut voir $B$ dans
$so(3,\C)$. Mais tous les \'el\'ements non nuls de
$so(3,\C)$ sont conjugu\'es (\`a multiplication pr\`es).
On peut donc supposer que $$B=z_3\frac{\partial}{\partial
z_2}-z_2\frac{\partial}{\partial z_3}$$ et
$$q_1=a(z_2^2+z_3^2)+bz_4^2$$ Il nous faut pour des raisons de
dimension ponctuelle trouver un troisi\`eme champ $B'$ annulant
$q_1$ et $q_2$ du type pr\'ec\'edent. Par suite tout champ
annulant $q_2$ annule $q_1$ et $q_1=\alpha q_2$. Finalement
$\mathscr{F}$ poss\`ede une int\'egrale premi\`ere du type
$$\frac{z_0z_1}{z_2^2+z_3^2+z_4^2}$$ que nous avons
d'ailleurs rencontr\'ee dans le principe de construction d'exemple.

Terminons par le cas o\`u $rg \hspace{1mm} q_2=2$, i.e.
$q_2=z_2^2+z_3^2$. En suivant la d\'emarche pr\'ec\'edente on doit
trouver deux champs du type $$Y_i=z_2\frac{\partial}{\partial
z_3}-z_3\frac{\partial}{\partial z_3}+L_i\frac{\partial}{\partial
z_4}  \hspace{5mm} i=1,2$$ ind\'ependants c'est-\`a-dire, en
particulier, non $\C$-colin\'eaires et annulant $q_1$ ; ceci
n'est possible que si $\frac{\partial}{\partial z_4}$ annule
$q_1$ cas o\`u $\mathscr{F}$ est un pull-back. \vs

\subsection{Exemple de type 3. }

\vs

On consid\`ere un feuilletage de degr\'e $2$ sur
$\C\P(4)$ ayant pour int\'egrale premi\`ere
$$\frac{Q_1}{Q_2}=\frac{z_0^2+z_1^2+z_2^2+z_3^2+z_4^2}
{\lambda_0z_0^2+\lambda_1z_1^2+\lambda_2z_2^2+\lambda_3z_3^2+
\lambda_4z_4^2}$$ On peut supposer que $\lambda_0=0$,
$\lambda_4=1$ et qu'un autre $\lambda_i$ est non nul sinon le
feuilletage serait de degr\'e $1$. Les champs $X$ qui annulent
$Q_1$ sont du type $$\sum_{0\leq i<j \leq 4}
\mu_{i,j}\left(z_i\frac{\partial}{\partial z_j}-z_j
\frac{\partial} {\partial z_i}\right)$$ Si de plus $X$ annule 
$Q_2$, on a $$\sum_{0\leq i<j \leq 4}
\mu_{i,j}(\lambda_i-\lambda_j)z_iz_j=0$$ Si $\lambda_1$, $\lambda_2$
et $\lambda_3$ sont non nuls et diff\'erents de $1$, on aura
$$\mu_{0,j}=0 \mbox{ pour } j=0,\ldots,4$$ $$\mu_{i,4}=0 \mbox{
pour } i=1,2,3$$ Ainsi un champ $X$ qui annule $Q_1$ et $Q_2$ n'a
pas de composante sur $\frac{\partial}{\partial z_0}$ et sur
$\frac{\partial}{\partial z_4}$. La condition de dimension
ponctuelle pour l'alg\`ebre des champs annulant $Q_1$ et $Q_2$ ne
peut \^etre satisfaite. On peut donc supposer toujours sous
l'hypoth\`ese $\lambda_1\lambda_2 \not =0$ que $\lambda_3=1$.
Si $\lambda_1$ et $\lambda_2$ sont encore diff\'erents de $1$, on
aura toujours pour $X$ annulant $Q_1$ et $Q_2$ les \'egalit\'es
$$\mu_{0,j}=\mu_{1,4}=\mu_{2,4}=\mu_{1,3}=\mu_{2,3}=0$$ et $X$
sera du type $$\mu_{1,2}\left(z_1\frac{\partial}{\partial
z_2}-z_2\frac{\partial}{\partial z_1}\right)+\mu_{3,4}\left(z_3
\frac{\partial} {\partial z_4}-z_4\frac{\partial}{\partial
z_3}\right) $$ ce qui est insuffisant pour obtenir la dimension
ponctuelle souhait\'ee. Par suite on peut supposer $\lambda_2=1$.
On se ram\`ene donc \`a
$$\frac{Q_1}{Q_2}=\frac{z_0^2+z_1^2+z_2^2+z_3^2+z_4^2}{z_0^2+
\lambda z_1^2}$$ Comme $\lambda \not =0$, on constate que la seule
possibilit\'e est $\lambda=1$. Le feuilletage a donc une
int\'egrale premi\`ere de type $$\frac{z_2^2+z_3^2+
z_4^2}{z_0^2+z_1^2}$$ vue dans le principe de construction
d'exemples et dans le cas sp\'ecial pr\'ec\'edent. On d\'emontre
de fa\c{c}on analogue que les autres configurations de $\lambda$
conduisent au m\^eme mod\`ele. D'o\`u la : \vs

\begin{pro} Soit $\mathscr{F}$ un
$\mathscr{L}$-feuilletage de degr\'e $2$ sur
$\C\P(4)$ ayant une int\'egrale premi\`ere
$\frac{Q_1}{Q_2}$ o\`u les $Q_i$ sont des formes quadratiques
diagonales. Si l'\'el\'ement g\'en\'erique du pinceau $Q_1-\kappa
Q_2$ est de rang maximum, alors $\mathscr{F}$ a une int\'egrale
premi\`ere du type $$\frac{z_2^2+ z_3^2+ z_4^2}{z_0^2+z_1^2}$$
\vs
\end{pro}

\begin{rem}
Un pinceau de formes quadratiques diagonales $Q_1=\displaystyle
\sum_{i=0}^k z_i^2$, $Q_2=\displaystyle \sum_{i=0}^n
\lambda_iz_i^2$ tel que l'\'el\'ement g\'en\'erique ne soit pas de
rang maximal satisfait $$\lambda_i=0 \mbox{ pour } i>k$$ En
particulier le feuilletage associ\'e est un pull-back.
\end{rem}

\vs

\subsection{}

Lors de l'examen de certaines configurations de
$\mathscr{L}$-feuilletages de degr\'e trois sur $\C\P(4)$
associ\'es \`a des alg\`ebres r\'esolubles de dimension trois,
nous avons trouv\'e quelques d\'eg\'en\'erescences conduisant \`a
des $\mathscr{L}$-feuilletages de degr\'e deux poss\'edant des
int\'egrales premi\`eres du type $\frac{Q_1}{Q_2}$ avec $Q_i$
forme quadratique ; pr\'ecis\'ement soient
$$Q_1=z_0z_4-p_1z_3^2-p_2z_3z_4-\frac{1}{2}z_2^2+p_3z_4^2$$
$$Q_2=z_1z_4-q_1z_3^2-q_2z_3z_4-z_2z_3-q_3z_4^2$$

\vs

On constate que l'\'el\'ement g\'en\'erique du pinceau $\lambda_1
Q_1+\lambda_2Q_2$ est d\'eg\'en\'er\'e. \vs

\begin{pro}
Le feuilletage $\mathscr{F}$ de degr\'e $2$ sur $\C\P(4)$
donn\'e par les niveaux de $\frac{Q_1}{Q_2}$ est un $\mathscr{L}$-feuilletage.
L'alg\`ebre de Lie $\mathscr{L}$ associ\'ee est r\'esoluble de dimension trois ;
si $\tilde{X} \in \tilde{\mathscr{L}}$, le champ $\tilde{X}$ s'\'ecrit $$
\left(b_1z_0+(b_3-e_3q_2+q_2b_1-2q_1b_2)z_2+(2p_1b_2-p_2b_1+e_3p_2)z_3\right.$$ 
$$\left.+(2
p_3b_1 -2p_3e_3+p_2b_2)z_4\right)\frac{\partial}{\partial z_0}$$ $$
+\left(b_1z_1+b_2z_2+b_3z_3+(q_2b_2-2q_3b_1+2q_3e_3)z_4 \right)
\frac{\partial}{\partial z_1}$$ $$+\left(e_3z_2+(b_3-e_3q_2+q_2b_1-2q_1b_2)z_4
\right)\frac{\partial}{\partial z_2}$$
$$+\left(e_3z_3+b_2z_4\right)\frac{\partial}{\partial z_3}
+(2e_3-b_1)z_4\frac{\partial}{\partial z_4}$$ o\`u $b_1$, $b_2$,
$b_3$, $e_3 \in \C$. Ce $\mathscr{L}$-feuilletage ne v\'erifie pas
la condition de \textbf{r\'egularit\'e} le long de $z_4=0$.
\end{pro}

\begin{proof}[D\'emonstration] L'expression des champs $\tilde{X}$ a \'et\'e obtenue via Maple.
L'hyperplan $z_4=0$ n'est pas invariant par le feuilletage puisque
$\frac{Q_1}{Q_2}$ n'y est pas constant. Par contre il est invariant par tous les
champs $\tilde{X}$. On constate par un calcul \'el\'ementaire qu'ils sont
d\'ependants le long de cet hyperplan.
\end{proof}

\vs

\subsection{Ubiquit\'e de l'exemple de Gordan-N\oe ther.}

\vs

Revenons \` a l'exemple de Gordan-N\oe ther
$$P=z_1^2z_3+z_1z_2z_4+z_2^2z_5$$ dans $\C^5$ et
consid\'erons sa restriction \`a la carte $z_1=1$
$$P_1=z_3+z_2z_4+z_2^2z_5$$ On constate que les champs
$$z_2\frac{\partial}{\partial z_3}-\frac{\partial}{\partial z_4},
\hspace{2mm} z_2\frac{\partial}{\partial
z_4}-\frac{\partial}{\partial z_5} \mbox{ et }
z_2\frac{\partial}{\partial z_2}-z_4\frac{\partial}{\partial
z_4}-2z_5\frac{\partial}{\partial z_5}$$ sont tangents aux 
niveaux de $P_1$. Le feuilletage par les niveaux de $P_1$ 
est donc une fois encore un $\mathscr{L}$-feuilletage de 
$\C^4$, lequel s'\'etend bien s\^ur \`a $\C\P(4)$ ; il a 
pour int\'egrale premi\`ere $$\frac{z_1^2z_3+z_1z_2z_4+
z_2^2z_5}{z_1^3}$$ qui d\'efinit un feuilletage de degr\'e $2$. 
\vs

Maintenant on consid\`ere la restriction de $P$ \`a la carte
$z_3=1$ ; on obtient le polyn\^ome
$$P_2=z_1^2+z_1z_2z_4+z_2^2z_5$$ qui est annul\'e par les champs
$$z_2\frac{\partial}{\partial z_1}-2\frac{\partial}{\partial z_4}, \hspace{2mm}
z_2\frac{\partial}{\partial z_4}-z_1\frac{\partial}{\partial z_5},
\hspace{2mm} z_2\frac{\partial}{\partial z_2}-z_4\frac{\partial}{\partial z_4}
-2z_5\frac{\partial}{\partial z_5}$$ et produit encore un exemple de degr\'e $2$
sur $\C\P(4)$ d'int\'egrale premi\`ere $$\frac{z_1^2z_3+z_1z_2
z_4+z_2^2z_5}{z_3^3}$$ Ces deux exemples, tous deux de type 4.1., sont non
\'equivalents. \vs
\vs

\section{Exemples de type 4.1. }

\vs

\subsection{Cas sp\'ecial } (o\`u l'on retrouve l'ubiquit\'e de
Gordan-N\oe ther).

\vs

Nous allons nous int\'eresser \`a un cas sp\'ecial de $\mathscr{L}$-feuilletage
quadratique de type 4.1. sur $\C\P(4)$. Soit $\mathscr{F}$ un
$\mathscr{L}$-feuilletage sur $\C\P(4)$ ayant une int\'egrale
premi\`ere du type $\frac{C}{z_0^3}$ o\`u $C$ est un polyn\^ome cubique. La
restriction de $\mathscr{F}$ \`a la carte affine $z_0=1$ d\'efinit un
$\mathscr{L}$-feuilletage de degr\'e $2$ sur $\C^4$ ayant pour
int\'egrale premi\`ere $P(z_1,z_2,z_3,z_4)=C(1,z_1,z_2,z_3,z_4)$ . Nous
allons supposer que $P$ est un polyn\^ome homog\`ene.

\begin{thm}
Soit $P$ un polyn\^ome homog\`ene minimal de degr\'e $3$ sur
$\C^4$ d\'efinissant un $\mathscr{L}$-feuilletage. Alors \`a
conjugaison pr\`es $P$ est de l'un des types suivants :
\begin{eqnarray}
&1.& z_1z_2^2\nonumber \\
& & \nonumber \\
&2.& z_1z_2z_3\nonumber \\
& & \nonumber \\
&3.& z_1(z_2^2-z_3z_4)\nonumber \\
& & \nonumber \\
&4.& z_1^2z_3+z_1z_2(a_1z_1+z_2+z_3+z_4)+z_2^2z_4\nonumber \\
& & \nonumber \\
&5.& z_1^2z_3+z_1z_2(z_3+z_4)+z_2^2z_4\nonumber \\
& & \nonumber \\
&6.& z_1^2z_3+z_1z_2z_4+a_2z_2^3\nonumber
\end{eqnarray}
\end{thm}

\vs

\begin{rems}
1. Consid\'erons le polyn\^ome de
\textbf{Gordan-N\oe ther}
$$P=z_1^2z_3+z_1z_2z_5+z_2^2z_4$$ 
Si l'on restreint $P$ \`a $z_2=0$ on obtient \`a permutation
pr\`es des coordonn\'ees le premier cas 1. En restreignant \`a
$z_3=z_4=0$ on trouve 2 ; la restriction de $P$ \`a $z_4=0$ est
\`a conjugaison pr\`es le cas 3. En
restreignant $P$ \`a l'hyperplan $a_1z_1+z_2+z_3+z_4$, on obtient
le cas 4. En le restreignant \`a $z_5=z_3+z_4$, on retrouve le
cas 5. Enfin le cas 6. s'obtient en permutant $z_4$ et $z_5$
et en restreignant $P$ \`a $z_5=a_2z_2$.

2. Par conjugaison on peut supposer $a_2\in \{0,1\}$.
\end{rems}

\vs

\begin{proof}
D'apr\`es le th\'eor\`eme \ref{cormal}, l'origine n'est
pas une singularit\'e isol\'ee de $dP$ sinon le feuilletage
$\mathscr{F}$ serait de degr\'e $1$. Il s'ensuit que $dP$ s'annule
sur un nombre fini de
droites, ou sur une surface ou encore sur une hypersurface. Nous
allons traiter ces diff\'erentes \'eventualit\'es au cas par
cas.\vs

\begin{itemize}
\item[1. ]  Si $dP$ s'annule sur une hypersurface, $P$ n'est pas
r\'eduit et s'\'ecrit $z_1^2L$ o\`u $L$ d\'esigne une forme
lin\'eaire. Puisque $P$ est minimal, $L$ est non colin\'eaire \`a
$z_1$. Le $\mathscr{L}$-feuilletage $\mathscr{F}$ d\'ecrit par le
polyn\^ome $P$ est alors
 de degr\'e $1$ et $P$ est du type $z_1z_2^2$. \vs

\item[2. ] Consid\'erons le cas o\`u $dP$ s'annule au moins sur une droite par
exemple d'\'equation $z_2=z_3=z_4=0$. En particulier le polyn\^ome $P$ n'a de
mon\^ome ni en $z_1^3$, ni en $z_1^2z_i$ ; par suite il est de la forme $$z_1
q(z_2,z_3,z_4)+r(z_2,z_3,z_4)$$ o\`u $q$ est une forme quadratique et $r$ un
polyn\^ome cubique. Un champ de vecteurs $X$ annulant $P$ doit \^etre tangent au
lieu singulier de $P$ et en particulier \`a l'axe des $z_1$. Par suite, il
s'\'ecrit $$A_1(z_1,z_2,z_3,z_4)\frac{\partial}{\partial z_1}+\displaystyle
\sum_{i=2}^4A_i(z_2,z_3,z_4)\frac{\partial}{\partial z_i}$$ Remarquons que si
$q\equiv 0$, alors $\mathscr{F}$ est un pull-back ; ce cas a \'et\'e trait\'e
dans le corollaire de l'introduction.

Si l'on \'ecrit $X$ sous la forme $A_1\frac{\partial}{\partial z_1}+Y$, on a
$$\left\{ \begin{array}{c} A_1q+z_1Y(q)=0\\ Y(r)=0 \hspace{10mm}\end{array}
\right.$$ Si $q\not \equiv 0$, alors $A_1$ est divisible par $z_1$, autrement
dit $$X=\lambda_1z_1 \frac{\partial}{\partial z_1}+Y$$ Supposons que pour tout
champ $X$ annulant le polyn\^ome $P$ on ait $\lambda_1=0$. La dimension
ponctuelle g\'en\'erique de tels $X$ \'etant trois, le
polyn\^ome $P$ est annul\'e par $\frac{\partial}{\partial z_2}$,
$\frac{\partial}{\partial z_3}$ et
$\frac{\partial}{\partial z_4}$ ce qui est absurde. Ainsi on peut trouver des
champs annulant $P$ du type $$X_1=z_1\frac{\partial }{\partial z_1} +Y_1,
\hspace{2mm} X_2=Y_2 \mbox{ et } X_3=Y_3$$ g\'en\'eriquement ind\'ependants ; il
en est de m\^eme pour les champs $Y_2$ et $Y_3$. Les
\'egalit\'es $$\left\{\begin{array}{c}Y_2(q)=Y_2(r)=0\nonumber \\ Y_3(q)=Y_3(r)=
0\nonumber\end{array}\right.$$ impliquent que les composantes connexes des
niveaux de $q$
et $r$ sont les m\^emes.

Si la forme quadratique $q$ est de rang $2$ ou $3$, alors on a $r=0$ et $P$
s'\'ecrit \`a conjugaison pr\`es $$z_1(z_2^2-z_3z_4)$$ que nous avons d\'ej\`a
rencontr\'e ou $$z_1z_2z_3$$ Dans les deux cas on obtient un
$\mathscr{L}$-feuilletage.

Reste \`a consid\'erer le cas o\`u la forme quadratique $q$ est de rang
$1$ : le polyn\^ome $P$ est alors de la forme $$z_1z_2^2+\varepsilon
z_2^3=z_2^2(z_1+\varepsilon z_2)$$ qui est conjugu\'e \`a $z_1z_2^2$ et le
feuilletage obtenu est de degr\'e $1$.\vs

\item[3. ] Cas o\`u $dP$ s'annule sur une surface $S$. \vs

\begin{itemize}
\item[3.1. ] Supposons qu'une
composante de $S$ soit un plan par exemple d'\'equation $z_1=z_2=0$, alors $P$
est du type $$z_1^2L_1+z_1z_2L_2+z_2^2L_3$$ o\`u $L_i$ d\'esigne une forme
lin\'eaire. On a l'alternative suivante :

- le feuilletage $\mathscr{F}$ est un pull-back (et donc encore du type
$z_1z_2^2$)

- les formes lin\'eaires $z_1$, $z_2$, $L_1$, $L_2$ et $L_3$ forment un
syst\`eme de rang $4$. On se ram\`ene alors aux deux cas suivants $$z_1^2z_3+
z_1z_2L+z_2^2z_4$$  $$z_1^2z_3+z_1z_2z_4+z_2^2L$$ o\`u $L$ d\'esigne une forme
lin\'eaire. Il nous faut maintenant caract\'eriser les $L$ qui produisent un
$\mathscr{L}$-feuilletage. Comme le lieu singulier doit \^etre invariant par
tout champ $X$ annulant $P$, un tel champ, n\'ecessairement lin\'eaire, est de
la forme $$(b_1z_1+b_2z_2)\frac{\partial}{\partial z_1}+(c_1z_1+c_2z_2)
\frac{\partial}{\partial z_2}$$ $$+\left(\displaystyle\sum_{i=1}^4 \alpha_iz_i\right)
\frac{\partial}{\partial z_3}+\left(\displaystyle \sum_{i=1}^4\beta_iz_i\right)
\frac{\partial}{\partial z_4}$$ Posons $$L=\displaystyle\sum_{i=1}^4 a_iz_i$$

\begin{itemize}
\item[3.1.1. ] Commen\c{c}ons par consid\'erer le cas o\`u $P$ est de la forme
$$z_1^2z_3+z_1z_2L+z_2^2z_4$$
L'\'egalit\'e $X(P)=0$ se traduit par le syst\`eme suivant :
$$\left\{ \begin{array}{c} 2b_1+c_1a_3+\alpha_3=0 \hspace{42mm} \\
2a_1b_1+2a_2c_1+a_1c_2+\alpha_2+\alpha_1a_3+\beta_1a_4 \hspace{7mm}\\
b_1a_2+2a_1b_2+2a_2c_2+\alpha_2a_3+\beta_2a_4+\beta_1=0 \\
b_1a_4+2c_1+a_4c_2+\alpha_4a_3+\beta_4a_4=0 \hspace{14mm}\\
2b_2+b_1a_3+a_3c_2+\alpha_3a_3+\beta_3a_4=0 \hspace{14mm}\\
b_2a_2+\beta_2=0 \hspace{52mm}\\
b_2a_4+2c_2+\beta_4=0 \hspace{42mm}\\
b_2a_3+\beta_3=0 \hspace{52mm}\\
a_1c_1+\alpha_1=0 \hspace{52mm}\\
a_4c_1+\alpha_4=0 \hspace{52mm}
\end{array} \right.$$

\vs

C'est un syst\`eme de dix \'equations \`a douze inconnues. Pour que la
dimension ponctuelle des champs $X$ annulant $P$ soit g\'en\'eriquement trois,
il nous faut au moins trois solutions ind\'ependantes. Ceci impose des
conditions sur les $a_i$ qui r\'egissent le syst\`eme. En substituant on obtient
que ce syst\`eme a un espace de solutions de dimension au moins trois si et
seulement si $$(a_2a_4-1)(2a_1+a_1a_4^2-3a_2a_4)=(a_3a_4-1)=0$$ \vs

\begin{itemize}
\item[3.1.1.1. ]
Lorsque $a_2a_4=a_3a_4=1$,
quitte \`a conjuguer par une application lin\'eaire ad-hoc on se ram\`ene \`a
$a_2=a_3=a_4=1$, autrement dit $$P=z_1^2z_3+z_1z_2
(a_1z_1+z_2+z_3+z_4)+z_2^2z_4$$ On peut alors r\'esoudre le syst\`eme
pr\'ec\'edent et l'espace vectoriel des champs annulant $P$ est engendr\'e par
les trois champs
$$X_1=z_2\frac{\partial}{\partial z_3}-z_1\frac{\partial}{\partial z_4}$$
$$X_2=(z_1+z_2)\frac{\partial}{\partial z_1}-2z_3\frac{\partial}{\partial z_3}
-(2a_1z_1+z_2+z_3+z_4)\frac{\partial}{\partial z_4}$$ $$X_3=(z_1+z_2)
\frac{\partial}{\partial z_2}-(a_1z_1+z_3+z_4)\frac{\partial}{\partial z_3}-2
(z_1+z_4)\frac{\partial}{\partial z_4}$$ On
constate que ces trois champs engendrent une alg\`ebre de dimension ponctuelle
trois. \vs

\item[3.1.1.2. ] Passons au cas o\`u $a_3a_4=1$ et $2a_1+a_1a_4^2-3a_2a_4=0$.
Comme $a_3$ et $a_4$ sont non nuls, on se ram\`ene par une application
lin\'eaire de type $$(z_1,z_2,z_3,z_4)\mapsto (\rho_1 z_1,\rho_2z_2,\rho_3z_3,
\rho_4z_4)$$ ad-hoc aux deux cas suivants : $$P_1=z_1^2z_3+z_1z_2(z_3+z_4)+
z_2^2z_4$$ et $$P_2=z_1^2z_3+z_1z_2(z_1+z_2+z_3+z_4)+z_2^2z_4$$ Dans le premier
cas, les champs suivants annulent $P_1$ :
$$z_1\frac{\partial}{\partial z_1}+z_2\frac{\partial}{\partial
 z_2}-2z_3\frac{\partial}{\partial z_3}-2z_4\frac{\partial}{\partial z_4}$$
$$z_2\frac{\partial}{\partial z_1}-z_2\frac{\partial}{\partial z_2}+(z_4-z_3)
\frac{\partial}{\partial z_4}$$
$$z_2\frac{\partial}{\partial z_3}-z_1\frac{\partial}{\partial z_4}$$ Ainsi
$P_1$ d\'efinit bien un $\mathscr{L}$-feuilletage.

Le second cas est un cas particulier du cas 3.1.1.1. avec $a_1=1$.\vs
\end{itemize}

\item[3.1.2. ] Consid\'erons le cas o\`u $P$ est de la forme $$z_1^2z_3+z_1z_2z_4
+z_2^2(a_1z_1+a_2z_2+a_3z_3+a_4z_4)$$ Notons que si $a_4\not =0$, en posant
$Z_4=a_1z_1+a_2z_2+a_3z_3+a_4z_4$ on se ram\`ene au cas pr\'ec\'edent. Nous
supposerons donc que $a_4=0$. Quitte \`a changer $z_4$ en
$z_4-a_1z_2$ on a $a_1=0$ ; autrement dit $P$
s'\'ecrit $$z_1^2z_3+z_1z_2z_4+z_2^2(a_2z_2+a_3z_3)$$ On proc\`ede comme
pr\'ec\'edemment ; l'\'egalit\'e $X(P)=0$ conduit au syst\`eme $$\left\{
\begin{array}{c}
b_1+\beta_4+c_2=0 \hspace{9mm}\\
b_2+\alpha_4a_3=0 \hspace{13mm} \\
c_1+\alpha_4=0 \hspace{17mm} \\
2b_1+\alpha_3=0 \hspace{15mm}\\
2b_2+2a_3c_1+\beta_3=0 \hspace{2mm}\\
a_3(2c_2+\alpha_3)=0 \hspace{9mm} \\
3a_2c_1+\alpha_1a_3+\beta_2=0\\
3a_2c_2+\alpha_2a_3=0 \hspace{9mm}\\
\alpha_1=0 \hspace{25mm}\\
\alpha_2+\beta_1=0 \hspace{17mm}
\end{array}\right.$$
\end{itemize}

\begin{itemize}
\item[3.1.2.1. ] Supposons que $a_3=0$, i.e. que le polyn\^ome $P$ s'\'ecrive $$
P=z_1^2z_3+z_1z_2z_4+a_2z_2^3$$ On a l'alternative suivante : \vs

- $a_2 \not =0$, le polyn\^ome $P$ est alors annul\'e par les champs
$$z_1\frac{\partial}{\partial z_1}-2z_3\frac{\partial}{\partial z_3}
-z_4\frac{\partial}{\partial z_4}$$ $$z_1\frac{\partial}{\partial z_2}-
z_4\frac{\partial}{\partial z_3}-3a_2z_2\frac{\partial}{\partial z_4}$$
$$z_2\frac{\partial}{\partial z_3}-z_1\frac{\partial}{\partial z_4}$$ Il
d\'efinit donc bien un $\mathscr{L}$-feuilletage.\vs

- $a_2=0$, le polyn\^ome $P$ d\'ecrit un $\mathscr{L}$-feuilletage dont
l'alg\`ebre associ\'ee, de dimension $4$, est engendr\'ee par les champs
$$z_1\frac{\partial}{\partial z_1}-2z_3\frac{\partial}{\partial z_3}-z_4
\frac{\partial}{\partial z_4}$$ $$z_1\frac{\partial}{\partial z_2}-z_4
\frac{\partial}{\partial z_3}$$ $$z_2\frac{\partial}{\partial z_2}-z_4
\frac{\partial}{\partial z_4}$$ $$z_2\frac{\partial}{\partial z_3}-z_1
\frac{\partial}{\partial z_4}$$

On constate que la dimension de l'alg\`ebre passe de trois \`a quatre lorsque
$a_2$ s'annule.\vs

\item[3.1.2.2. ] Supposons maintenant que $a_3\not =0$. A conjugaison pr\`es par
des applications du type $$(z_1,z_2,z_3,z_4)\mapsto (\rho_1z_1,\rho_2z_2,\rho_3
z_3,\rho_4z_4)$$ on se ram\`ene \`a l'\'etude des deux cas
suivants $$z_1^2z_3+z_1z_2z_4+z_2^2z_3$$ $$z_1^2z_3+z_1z_2z_4+z_2^2(z_2+z_3)$$

Les champs qui annulent ces deux polyn\^omes forment un espace de dimension $2$
donc les feuilletages d\'ecrits par ces polyn\^omes ne sont pas des
$\mathscr{L}$-feuilletages.
\vs
\end{itemize}

\item[3.2. ] Pla\c{c}ons nous dans l'\'eventualit\'e o\`u la surface $S$
contient une composante non plane.
On peut alors supposer que trois axes de coordonn\'ees sont dans $S$, par
exemple les axes $z_1$, $z_2$ et $z_3$. Par un argument d\'ej\`a
rencontr\'e, le polyn\^ome $P$ est affine en $z_1$, $z_2$, $z_3$ et quadratique
en $z_4$ ; on constate que $P$ est du type $$\varepsilon z_1z_2z_3+z_4^2(a_1z_1+
a_2z_2+a_3z_3)$$ On remarque que $\varepsilon \not =0$ sinon le feuilletage
d\'ecrit par $P$ est de degr\'e $1$. D\'eterminons les valeurs des $a_i$ pour
lesquelles ces polyn\^omes ont une surface singuli\`ere. \vs

\begin{itemize}
\item[3.2.1. ] Si $a_1a_2a_3 \not =0$, on se ram\`ene par homoth\'etie \`a
$$P=z_1z_2z_3+z_4^2(z_1+z_2+z_3)$$ et un calcul \'el\'ementaire montre que
$codim \hspace{1mm} \mathscr{S}ing \hspace{1mm} P=3$.
\vs

\item[3.2.2.] Si deux des $a_i$ sont non nuls, on se ram\`ene \`a
$$P=z_1z_2z_3 +z_4^2(z_1+z_2)$$ et
comme pr\'ec\'edemment le lieu singulier de $P$ est de codimension
$3$.\vs

\item[3.2.3. ] Reste donc le cas o\`u un seul des $a_i$ est non nul, alors $P$
est de la forme : $$z_1(z_2z_3+z_4^2)$$ cas d\'ej\`a rencontr\'e.
\end{itemize}
\end{itemize}
\end{itemize}
\end{proof}

\subsection{Exemples de $\mathscr{L}$-feuilletages quadratiques
satisfaisant $\dim \mathscr{L}=3$.}

\vs

Lors de l'\'etude des $\mathscr{L}$-feuilletages de degr\'e $3$ sur
$\C\P(4)$ correspondant aux alg\`ebres $\mathscr{L}_\alpha$,
nous avons rencontr\'e une d\'eg\'en\'erescence qui produit un
feuilletage de degr\'e $2$ de type 4.1 ; l'alg\`ebre $\mathscr{L}$ est ici de
dimension trois (cf. chapitre 6). Voici cet exemple : \vs

L'alg\`ebre $\mathscr{L}$ est engendr\'ee par les champs
$$X=\left(-2 z_0+(\xi_5c- h)z_2+\frac{c\xi_{13}+\delta+2\xi_5 e}{2}z_4\right)
\frac{\partial}{\partial z_0}$$$$+(\xi_5z_0-3z_1+\xi_7z_2+\xi_8z_3
+(h\xi_{13}+2i-\xi_5\delta)z_4)\frac{\partial}
{\partial z_1}$$$$+(-z_2+\xi_{13}z_3)\frac{\partial}
{\partial z_2}+(\xi_5z_3-z_4)\frac{\partial}{\partial
z_4},$$ $$Z=(cz_2+\delta z_3+ez_4)\frac{\partial} {\partial
z_0}+(z_0+hz_2+iz_3)\frac{\partial}{\partial z_1}+
z_3\frac{\partial}{\partial z_4},$$ %qui correspond au cas ou
%$a=b=m=\ell=0$, $f\not =0$ et $j=0$ parce qu'on peut l'enlever.
$$Y=z_3\frac{\partial}{\partial z_0}+z_4\frac{\partial}{\partial z_1}$$
\vs

L'alg\`ebre engendr\'ee par $X$, $Y$ et $Z$ est r\'esoluble et de type
$\mathscr{L}_\alpha$ (avec $\alpha=2$), c'est-\`a-dire poss\`ede une
pr\'esentation de la forme $$\{[\tilde{X},\tilde{Y}]=\tilde{Y}, \hspace{2mm}
[\tilde{X},\tilde{Z}]=\alpha \tilde{Z},
\hspace{2mm}[\tilde{Y},\tilde{Z}]=0\},\hspace{2mm} \alpha \not = 0$$

On montre que le feuilletage associ\'e \`a
cette alg\`ebre admet pour int\'egrale premi\`ere
$$\frac{3{\left(z_2-\xi_{13}z_3
\right)}^3}{z_0z_3z_4-z_1z_3^2-\frac{e}{3}z_4^3-\frac{cz_2+
\delta z_3}{2}z_4^2+(hz_2+iz_3)z_3z_4+\kappa z_3^3-3\xi_5h
z_2z_3^2}$$ o\`u $$\kappa=-\xi_{13}(\xi_7-\xi_5h)
-2(\xi_5i-\xi_8)$$ c'est-\`a-dire \`a conjugaison pr\`es $$\frac{z_2^3}{z_0z_3
z_4+z_1z_3^2+z_4^2(az_4+bz_2)}$$ avec $a$, $b \in\{0,1\}$.
\vs

\begin{rems}
1. En se restreignant \`a la carte affine $z_2=1$, on constate que le
polyn\^ome $$P=z_0z_3z_4+z_1z_3^2+az_4^3+bz_4^2$$ d\'efinit un
$\mathscr{L}$-feuilletage sur $\C^4$.

2. Le d\'enominateur $z_0z_3z_4+z_1z_3^2+z_4^2(az_4+bz_2)$
est encore lorsque $b$ est non nul de type Gordan-N\oe ther. 
Pour $b=0$, le feuilletage est un pull-back.
\end{rems}

%Les valeurs propres de la matrice $$\left(%
%\begin{array}{cc}
% 2\varepsilon & 0 \\
%  \frac{1}{2}(c\xi_{13}+\varepsilon \delta) & \varepsilon \\
%\end{array}%
%\right)$$ sont $2\varepsilon$ et $\varepsilon$ donc
%$\tilde{\alpha}=\tilde{2}$.
%\vs
%\vs

%%%%%%%%%%%%%%%%%%%%%%%%%%%%%%%%%%%%%%%%%%%%%%%%%%%%%%%%%%%%%%%%%%%%%%%
%%%%%%%%%%%%%%%%%%%%%%%%%%%%%%%%%%%%%%%%%%%%%%%%%%%%%%%%%%%%%%%%%%%%%%%

\chapter{$\mathscr{L}$-feuilletages de degr\'e trois sur $\C
\P(4)$.}

Soit $\mathscr{F}$ un $\mathscr{L}$-feuilletage de degr\'e $3$ sur
$\C \P(4)$, la proposition \label{degre} assure que $\dim
\mathscr{L}=3$. Ainsi la classification des
$\mathscr{L}$-feuilletages de degr\'e $3$ sur $\C
\P(4)$ repose sur celle des alg\`ebres de Lie de dimension
$3$ (\cite{[Fu-Ha]}). Pour classifier les alg\`ebres de Lie $\mathscr{L}$
de dimension $3$, on consid\`ere l'application bilin\'eaire
$$[\hspace{1mm},\hspace{1mm}] \hspace{1mm} \colon \hspace{1mm}
\mathscr{L} \times \mathscr{L} \to \mathscr{L}$$ et plus
particuli\`erement son rang.

Si $\dim [\mathscr{L},\mathscr{L}]=0$, l'alg\`ebre
$\mathscr{L}$ est ab\'elienne.

Si $\dim[\mathscr{L},\mathscr{L}]=1$, 
l'alg\`ebre $\mathscr{L}$ est d\'ecrite par les pr\'esentations
$$\{[X,Y]=[X,Z]=0, \hspace{1mm} [Y,Z]=Y\}$$ ou $$\{[X,Y]=[X,Z]
=0, \hspace{1mm} [Y,Z]=X\}$$

Si $\dim [\mathscr{L},\mathscr{L}]=2$,
l'alg\`ebre $\mathscr{L}$ est du type $$\{[X,Y]=Y, \hspace{1mm}
[X,Z]=\alpha Z, \hspace{1mm} [Y,Z]=0\}$$ avec $\alpha \in \C^*$ ou bien
$$\{[X,Y]=Y, \hspace{1mm} [X,Z]=Y+Z, \hspace{1mm} [Y,Z]=0\}$$.

Enfin si $\dim[\mathscr{L},\mathscr{L}]=3$, alors $\mathscr{L}$ est isomorphe
\`a $s\ell(2,\C)$.

Nous allons \'etudier au cas par cas les $\mathscr{L}$-feuilletages associ\'es.
\vs

\section{L'alg\`ebre $\mathscr{L}$ est isomorphe \`a $s\ell(2,
\C)$.} Par simple connexit\'e de $\mbox{SL}(2,\C)$, la
repr\'esentation de $s\ell(2,\C)$ dans $\mathscr{M}(5,\C)$
donn\'ee par $\mathscr{L}'$ assure l'existence d'une action
lin\'eaire de $\mbox{SL}(2,\C)$ sur $\C^5$. On est donc ramen\'e
\`a classifier les actions lin\'eaires de
$\mbox{SL}(2,\C)$ sur $\C^5$ que nous allons d\'ecrire.

Commen\c{c}ons donc par supposer que l'action est irr\'eductible ;
c'est alors \`a conjugaison pr\`es l'action naturelle de
$\mbox{SL}(2,\C)$ sur les polyn\^omes homog\`enes de degr\'e
quatre en deux variables. Cette action descend \`a l'espace
projectif pour donner un $\mathscr{L}$-feuilletage de degr\'e $3$ dont
les feuilles sont les niveaux g\'en\'eriques de la fonction $j$ ; nous
l'avons d\'ecrite au chapitre $2$. A
conjugaison pr\`es ce $\mathscr{L}$-feuilletage est unique
puisqu'il y a une unique action irr\'eductible de $\mbox{SL}(2,\C)$
 sur $\C^5$.

Supposons maintenant que l'action lin\'eaire de
$\mbox{SL}(2,\C)$ sur $\C^5$ soit r\'eductible ;
elle se d\'ecompose en actions irr\'eductibles avec les
configurations suivantes des diff\'erents facteurs : $$\C
\oplus \C^4,
\hspace{1mm} \C^2 \oplus \C^3,\hspace{1mm}
\C \oplus \C \oplus \C^3, \hspace{1mm}
\C \oplus \C^2 \oplus \C^2,$$
$$\C \oplus \C \oplus \C \oplus
\C^2,\mbox{ et } \C \oplus \C \oplus
\C \oplus \C \oplus \C.$$ On peut exclure
les deux derni\`eres configurations : l'action sur chaque facteur
$\C$
\'etant triviale, elles ne donneront pas d'orbites de dimension
quatre. Examinons les autres possibilit\'es au cas par cas :
\vs

(i) action de $\mbox{SL}(2,\C)$ sur $\C \oplus
\C^4$.

Sur le facteur $\C$ l'action est triviale et sur
$\C^4$ c'est l'action naturelle de $\mbox{SL}(2,\C)$ sur
$V_4$ (cf chapitre 2) \`a conjugaison pr\`es. Elle admet 
pour invariants $z_1$ et $\Delta(z_2,\ldots,z_5)$ o\`u $$\Delta(z_2,
z_3,z_4,z_5)=z_3^2z_4^2-4z_2z_4^3-4z_3^3z_5+18z_2z_3z_4z_5$$ Par suite
$\mathscr{F}$ poss\`ede l'int\'egrale premi\`ere $$
\frac{\Delta(z_2,\ldots,z_5)}{z_1^4}$$

A conjugaison pr\`es, le $\mathscr{L}$-feuilletage $\mathscr{F}$
est unique. Il a d\'ej\`a \'et\'e
\'etudi\'e dans la liste d'exemples (chapitre $2$) mais
pr\'esent\'e de mani\`ere diff\'erente.\vs

(ii) action de $\mbox{SL}(2,\C)$ sur $\C^2 \oplus
\C^3$.

Elle est conjugu\'ee \`a l'action de $\mbox{SL}(2,\C)$ sur
$V_2\oplus V_3$. On note les \'el\'ements de $V_2\oplus V_3$ sous la
forme  $$(z_1x+z_2y,z_3x^2+z_4xy+z_5y^2)$$ A chacun de ces \'el\'ements
on associe
$$Q_1=(z_1x+z_2y)^2=z_1^2x^2+2z_1z_2xy+z_2^2y^2 \hspace{1mm},$$
$$Q_2=z_3x^2+z_4xy+z_5y^2$$ et $\tilde{Q}_i$ la matrice de la forme
quadratique $Q_i$. Soit $P \in \mbox{SL}(2,\C)$ ; on a
$$\det(^tP(\tilde{Q}_1-\kappa \tilde{Q}_2)P)=\det(\tilde{Q}_1-
\kappa \tilde{Q}_2)$$ autrement dit $$\det(\tilde{Q}_1-
\kappa \tilde{Q}_2)=\kappa^2\left(z_3z_5-\frac{z_4^2}{4}\right)
-\kappa(z_3z_2^2+z_5z_1^2-z_1z_2z_4)$$ est un invariant de
l'action de $\mbox{SL}(2,\C)$ sur $\C^2 \oplus
\C^3$. Ainsi les niveaux g\'en\'eriques de $$(z_3z_5-
\frac{z_4^2}{4},\hspace{1mm} z_3z_2^2+z_5z_1^2-z_1z_2z_4)$$ d\'efinissent le
feuilletage de codimension $2$ associ\'e \`a
$\mathscr{L}'$ et $$\frac{{(z_3z_2^2+z_5z_1^2-z_1z_2z_4)}^2}{{\left(z_3z_5- \frac{z_4^2}{4}\right)}^3}$$ est une int\'egrale premi\`ere du
$\mathscr{L}$-feuilletage de degr\'e $3$ associ\'e \`a $\mathscr{L}$.\vs

\begin{rem}
Le num\'erateur de cette int\'egrale premi\`ere est le 
polyn\^ome de Gordan-N\oe ther ; on retrouve le feuilletage
 exceptionnel en coupant par un $\C 
\P(3)$ ad hoc. Ce $\C\P(3)$ est non g\'en\'erique puisqu'il
 y a une chute de degr\'e.
\end{rem}

\vs

(iii) action de $\mbox{SL}(2,\C)$ sur $\C \oplus
\C^2 \oplus \C^2$.

On \'ecrit l'\'el\'ement $g\in\mbox{SL}(2,\C)$ sous la forme $$\left(%
\begin{array}{cc}
  a & b \\
  c & d \\
\end{array}%
\right)$$ Comme l'action naturelle de $\mbox{SL}(2,\C)$
sur $\C^2$ correspond \`a l'action sur les formes
lin\'eaires, on peut voir l'action 
sur $\C^2 \oplus \C^2$ \`a conjugaison pr\`es comme suit
$$g.(z_2x+z_3y,z_4x+z_5y)\simeq \left(%
\begin{array}{cc}
  z_2 & z_3 \\
  z_4 & z_5 \\
\end{array}%
\right)\left(%
\begin{array}{cc}
  a & b \\
  c & d \\
\end{array}%
\right)^{-1}$$ On constate ainsi que $$\det \left(%
\begin{array}{cc}
  z_2 & z_3 \\
  z_4 & z_5 \\
\end{array}%
\right)=z_2z_5-z_3z_4$$ est un invariant et que les niveaux
g\'en\'eriques de $$(z_1,z_2z_5-z_3z_4)$$ sont les feuilles du
feuilletage de codimension $2$ associ\'e \`a $\mathscr{L}'$ avec les
notations habituelles. Ainsi
$$\frac{z_2z_5-z_3z_4}{z_1^2}$$ est une int\'egrale premi\`ere du
$\mathscr{L}$-feuilletage associ\'e \`a $\mathscr{L}$ qui est donc de
degr\'e $1$ (car d\'efini dans la carte affine $z_1=1$ par une
forme quadratique) ; cette action n'induit donc pas un
$\mathscr{L}$-feuilletage de degr\'e $3$ sur $\C
\P(4)$.\vs

(iv) action de $\mbox{SL}(2,\C)$ sur $\C \oplus
\C \oplus \C^3$.

Cette action est triviale sur chaque facteur $\C$ et sur
le facteur $\C^3$ c'est l'action naturelle sur V$_3$.
Ainsi le feuilletage associ\'e \`a $\mathscr{L}'$ est donn\'e par les
niveaux de $(z_1,z_2)$ et une int\'egrale premi\`ere du
$\mathscr{L}$-feuilletage associ\'e \`a $\mathscr{L}$ est {\Large
$\frac{z_1}{z_2}$}. Par suite ce feuilletage est de degr\'e $0$.
\vs

Finalement seules les deux premi\`eres actions conduisent \`a un
$\mathscr{L}$-feuilletage de degr\'e $3$ sur $\C
\P(4)$ :  \vs

\begin{thm} Soit $\mathscr{F}$ un
$\mathscr{L}$-feuilletage de degr\'e $3$ sur
$\C\P(4)$. Si $\mathscr{L}$ est isomorphe \` a
$s\ell(2,\C)$, alors le feuilletage $\mathscr{F}$ admet
\`a conjugaison pr\`es l'une des int\'egrales premi\`eres
suivantes : $$\frac{z_3^2z_4^2-4z_2z_4^3-4
z_3^3z_5+18z_2z_3z_4z_5}{z_1^4}$$
$$\frac{{(z_3z_2^2+z_5z_1^2-z_1z_2
z_4)}^2}{{\left(z_3z_5-\frac{z_4^2}{4}\right)}^3}$$
\end{thm}

\section{L'alg\`ebre $\mathscr{L}$ est ab\'elienne}

En classifiant les triplets de matrices $(5,5)$ qui commutent on se ram\`ene
\`a \'etudier les sept cas qui suivent :
\begin{eqnarray}
&(i)& \tilde{\mathscr{L}} \mbox{est diagonale}, \nonumber \\
& & \nonumber \\
& & \mbox{ ou bien $\tilde{\mathscr{L}}$ est engendr\'ee par les champs du
type suivant } \nonumber \\
& & \nonumber \\
&(ii)&  \kappa z_1\frac{\partial}{\partial
z_0}+z_2\frac{\partial}{\partial z_2},\hspace{2mm} \xi
z_1\frac{\partial}{\partial z_0}+z_3\frac{
\partial}{\partial z_3}, \hspace{2mm} z_1\frac{\partial}
{\partial z_0}+z_4\frac{\partial}{\partial z_4} \hspace{2mm}
\mbox{et } R. \nonumber \\
& & \nonumber \\
&(iii)& \kappa z_1\frac{\partial}{\partial z_0}+z_2\frac{\partial}
{\partial z_2}+z_3\frac{\partial}{\partial z_3}, \hspace{2mm}
z_1\frac{\partial}{\partial z_0}+z_4\frac{\partial}{\partial z_4},
\hspace{2mm} z_1 \frac{\partial}{\partial z_0}+z_3
\frac{\partial}{\partial z_2} \hspace{2mm} \mbox{et } R \nonumber \\
& & \nonumber \\
&(iv)& \kappa z_1\frac{\partial}{\partial
z_0}+z_3\frac{\partial}{\partial z_2} +z_4\frac{\partial}{\partial
z_3}, \hspace{2mm} z_1 \frac{\partial}{\partial
z_0}+z_4\frac{\partial}{\partial z_2}, \hspace{2mm}(z_0+\xi
z_1)\frac{\partial}{\partial z_0}+z_1\frac{\partial}{\partial
z_1}\hspace{2mm} \mbox{et } R. \nonumber \\
& & \nonumber \\
&(v)& z_1\frac{\partial}{\partial
z_1}+z_4\frac{\partial} {\partial
z_2}, \hspace{2mm} \kappa z_1\frac{\partial}{\partial
z_1}+z_3\frac{\partial}{\partial z_2}+z_4\frac{\partial}{\partial
z_3}, \hspace{2mm} \xi z_0 \frac{\partial}{\partial z_0}+z_1
\frac{\partial}{\partial z_1} \hspace{2mm} \mbox{et } R. \nonumber \\
& & \nonumber \\
&(vi)&  \kappa z_0\frac{\partial}{\partial z_0}+(z_2+\xi
z_4)\frac{\partial}{\partial z_1}+z_3\frac{\partial}{\partial z_2}+
z_4\frac{\partial}{\partial z_3}, \hspace{2mm} \beta z_0
\frac{\partial}{\partial z_0}+(z_3+\delta
z_4)\frac{\partial}{\partial z_2}, \nonumber \\
& & z_0
\frac{\partial}{\partial z_0}+z_4\frac{\partial}{\partial z_1}
\hspace{2mm} \mbox{et } R.  \nonumber 
\end{eqnarray}

\begin{eqnarray}
&(vii)& (z_1+\kappa z_4)\frac{\partial}{\partial z_0}+z_2
\frac{\partial}{\partial z_1}+z_3\frac{\partial}{\partial z_2}+z_4
\frac{\partial}{\partial z_3}, \hspace{2mm} (z_2+\xi z_4)\frac{\partial}
{\partial z_0}+z_3\frac{\partial}{\partial z_2}+ z_4
\frac{\partial}{\partial z_3}, \nonumber \\
& & (z_3+\beta z_4)
\frac{\partial}{\partial z_0}+z_4\frac{\partial}{\partial z_1}
\hspace{2mm} \mbox{et } R.\nonumber
\end{eqnarray}

\vs

\begin{rem}
On a \'elimin\'e les alg\`ebres qui donnent des feuilletages
de degr\'e strictement inf\'erieur \`a $3$.
\end{rem}

Examinons ces possibilit\'es au cas par cas :
\vs

\begin{itemize}

\item[(i)] Commen\c{c}ons par traiter le cas o\`u
$\tilde{\mathscr{L}}$ est diagonale ; une base de
$\tilde{\mathscr{L}}$ est
$$X=\sum_{k=0}^4 \kappa_k z_k \frac{\partial}{\partial z_k},
\hspace{2mm} Y=\sum_{k=0}^4 \mu_k z_k \frac{\partial}{\partial
z_k}, \hspace{2mm}Z=\sum_{k=0}^4 \nu_k z_k
\frac{\partial}{\partial z_k} \mbox{ et }R $$ Le feuilletage
$\mathscr{F}$ est d\'ecrit par
$$z_0 z_1 z_2 z_3 z_4 \sum_{k=0}^4 \xi_k
\frac{dz_k}{z_k}$$ avec $\displaystyle \sum_{k=0}^4 \kappa_k
\xi_k=\displaystyle \sum_{k=0}^4 \mu_k \xi_k=\displaystyle
\sum_{k=0}^4 \nu_k \xi_k=\displaystyle \sum_{k=0}^4 \xi_k=0$.
\vs
\vs

%En particulier si $\xi_0=\xi_1=\xi_2=1$, $\xi_3=-1$ et $\xi_4=-2$,
%autrement dit si le $\mathscr{L}$-feuilletage est d\'ecrit par la
%$1$-forme
%$$\frac{dz_0}
%{z_0}+\frac{dz_1}{z_1}+\frac{dz_2}{z_2}-\frac{dz_3}{z_3}-\frac{dz_4}{z_4},$$
%alors il admet l'int\'egrale premi\`ere rationnelle suivante
%$$\frac{z_0z_1z_2}{z_3z_4^2}.$$
%\vs

\item[(ii)] Cas o\`u $\tilde{\mathscr{L}}$ est
engendr\'ee par les champs
%$$\left(%
%\begin{array}{ccccc}
%  \lambda & s & 0 & 0 & 0 \\
%  0 & \lambda & 0 & 0 & 0 \\
%  0 & 0 & \mu_2 & 0 & 0 \\
%  0 & 0 & 0 & \mu_3 & 0 \\
%  0 & 0 & 0 & 0 & \mu_4 \\
%\end{array}%
%\right)$$
% autrement dit par
%$$\left(%
%\begin{array}{ccccc}
%  0 & \alpha & 0 & 0 & 0 \\
%  0 & 0 & 0 & 0 & 0 \\
%  0 & 0 & 1 & 0 & 0 \\
%  0 & 0 & 0 & 0 & 0 \\
%  0 & 0 & 0 & 0 & 0 \\
%\end{array}%
%\right),
%\left(%
%\begin{array}{ccccc}
%  0 & \beta & 0 & 0 & 0 \\
%  0 & 0 & 0 & 0 & 0 \\
%  0 & 0 & 0 & 0 & 0 \\
%  0 & 0 & 0 & 1 & 0 \\
%  0 & 0 & 0 & 0 & 0 \\
%\end{array}%
%\right),
%\left(%
%\begin{array}{ccccc}
%  0 & 1 & 0 & 0 & 0 \\
%  0 & 0 & 0 & 0 & 0 \\
%  0 & 0 & 0 & 0 & 0 \\
%  0 & 0 & 0 & 0 & 0 \\
%  0 & 0 & 0 & 0 & 1 \\
%\end{array}%
%\right)$$ et R

\vs

$$\kappa z_1\frac{\partial}{\partial z_0}+z_2\frac{\partial}{\partial
z_2},\hspace{2mm} \xi z_1\frac{\partial}{\partial z_0}+z_3\frac{
\partial}{\partial z_3}, \hspace{2mm} z_1\frac{\partial}
{\partial z_0}+z_4\frac{\partial}{\partial z_4} \hspace{2mm}
\mbox{et } R.$$

\vs

La $1$-forme annul\'ee par ces quatre champs
s'\'ecrit \`a multiplication pr\`es dans la carte affine $z_1=1$
$$dz_0-\kappa \frac{dz_2}{z_2}-\xi \frac{dz_3}{z_3}-\frac{dz_4}{z_4} $$
 Le $\mathscr{L}$-feuilletage d\'ecrit par cette $1$-forme
a pour int\'egrale premi\`ere
$$\frac{z_1^{\kappa+\xi+1}}{z_2^\kappa z_3^\xi z_4}\exp\left(
\frac{z_0}{z_1}\right)$$
\vs

\item[(iii)] Cas o\`u $\tilde{\mathscr{L}}$ est engendr\'ee par les champs
%$$\left(%
%\begin{array}{ccccc}
%  \lambda & s & 0 & 0 & 0 \\
%  0 & \lambda & 0 & 0 & 0 \\
%  0 & 0 & \mu & t & 0 \\
%  0 & 0 & 0 & \mu & 0 \\
%  0 & 0 & 0 & 0 & \mu_4 \\
%\end{array}%
%\right)$$
% autrement dit par
%$$\left(%
%\begin{array}{ccccc}
%  0 & \alpha & 0 & 0 & 0 \\
%  0 & 0 & 0 & 0 & 0 \\
%  0 & 0 & 1 & 0 & 0 \\
%  0 & 0 & 0 & 1 & 0 \\
%  0 & 0 & 0 & 0 & 0 \\
%\end{array}%
%\right),
%\left(%
%\begin{array}{ccccc}
%  0 & 1 & 0 & 0 & 0 \\
%  0 & 0 & 0 & 0 & 0 \\
%  0 & 0 & 0 & 0 & 0 \\
%  0 & 0 & 0 & 0 & 0 \\
%  0 & 0 & 0 & 0 & 1 \\
%\end{array}%
%\right),
%\left(%
%\begin{array}{ccccc}
%  0 & \beta & 0 & 0 & 0 \\
%  0 & 0 & 0 & 0 & 0 \\
%  0 & 0 & 0 & 1 & 0 \\
 % 0 & 0 & 0 & 0 & 0 \\
%  0 & 0 & 0 & 0 & 0 \\
%\end{array}%
%\right)$$ et R

\vs

$$\kappa z_1\frac{\partial}{\partial z_0}+z_2\frac{\partial}
{\partial z_2}+z_3\frac{\partial}{\partial z_3}, \hspace{2mm}
z_1\frac{\partial}{\partial z_0}+z_4\frac{\partial}{\partial z_4},
\hspace{2mm} z_1 \frac{\partial}{\partial z_0}+z_3
\frac{\partial}{\partial z_2} \hspace{2mm} \mbox{et } R.$$ 

\vs

La $1$-forme annul\'ee par ces champs s'\'ecrit \`a multiplication
pr\`es dans la carte affine $z_1=1$ $$-dz_0+d(\frac{z_2}{z_3})
+\kappa \frac{dz_3}{z_3}+\frac{dz_4}{z_4}$$ Le feuilletage
associ\'e \`a $\mathscr{L}$ admet donc pour int\'egrale premi\`ere
$$\frac{z_3^\kappa z_4}{z_1^{1+\kappa}}\exp\left(\frac{z_1z_2-z_0z_3}
 {z_1z_3}\right)$$\vs

\item[(iv)] Cas o\`u $\tilde{\mathscr{L}}$ est engendr\'ee par les
champs

\vs

$$\kappa
z_1\frac{\partial}{\partial z_0}+z_3\frac{\partial}{\partial z_2}
+z_4\frac{\partial}{\partial z_3}, \hspace{2mm} z_1
\frac{\partial}{\partial z_0}+z_4\frac{\partial}{\partial z_2},
\hspace{2mm}(z_0+\xi z_1)\frac{\partial}{\partial
z_0}+z_1\frac{\partial}{\partial z_1}\hspace{2mm} \mbox{et } R.$$

\vs

La $1$-forme annul\'ee par ces quatre champs s'\'ecrit \`a
multiplication pr\`es dans la carte affine $z_4=1$
$$d\left(\frac{z_0}{z_1}\right)-\xi \frac{dz_1}{z_1}+d(-z_2+
\frac{z_3^2}{2}-\kappa z_3) $$ Le $\mathscr{L}$-feuilletage 
d\'ecrit par cette $1$-forme a pour int\'egrale premi\`ere
$$\frac{z_0z_4^2-z_1z_2z_4-\kappa z_1z_3z_4+\frac{z_1z_3^2}
{2}}{z_1z_4^2}$$ si  $\xi=0$,
$$\left(\frac{z_1}{z_4}\right)^\xi\exp\left(\frac{z_0z_4^2-z_1z_2z_4-
\kappa z_1z_3z_4+\frac{z_1z_3^2}{2}}{z_1z_4^2}\right)$$
sinon. A conjugaison pr\`es on obtient :
$$\frac{z_0z_4^2-z_1z_2z_4+z_1z_3^2}{z_1z_4^2}$$
si  $\xi=0$, $$\left(\frac{z_1}{z_4}\right)^\xi
\exp\left(\frac{z_0z_4^2-z_1z_2z_4+z_1z_3^2}{z_1z_4^2}\right)$$
sinon.\vs

\item[(v)] Cas o\`u l'alg\`ebre $\tilde{\mathscr{L}}$ est
engendr\'ee par les champs

\vs

$$z_1\frac{\partial}{\partial z_1}+z_4\frac{\partial}{\partial
z_2}, \hspace{2mm} \kappa z_1\frac{\partial}{\partial
z_1}+z_3\frac{\partial}{\partial z_2}+z_4\frac{\partial}{\partial
z_3}, \hspace{2mm} \xi z_0 \frac{\partial}{\partial z_0}+z_1
\frac{\partial}{\partial z_1} \hspace{2mm} \mbox{et } R.$$ 

\vs

La $1$-forme $\omega$ annulant ces quatre champs s'\'ecrit \`a
multiplication pr\`es dans la carte affine $z_4=1$
$$-\frac{dz_0}{z_0}+\xi\frac{dz_1}{z_1}-\xi dz_2+\xi(z_3-\kappa)dz_3$$ Le
feuilletage associ\'e \`a l'alg\`ebre $\mathscr{L}$ admet pour
int\'egrale premi\`ere $$\frac{z_1^\xi z_4^{1-\xi}}{z_0}
\exp\left(\frac{\xi (z_3^2-2\kappa z_3z_4-2z_2z_4)}
{2z_4^2}\right)$$ soit \`a conjugaison pr\`es $$
\frac{z_1^\xi z_4^{1-\xi}}{z_0}
\exp\left(\frac{z_3^2-z_2z_4}{z_4^2}\right)$$
\vs

\item[(vi)] Cas o\`u l'alg\`ebre $\tilde{\mathscr{L}}$ est
engendr\'ee par les champs
$$\kappa z_0\frac{\partial}{\partial z_0}+z_2
\frac{\partial}{\partial z_1}+z_3\frac{\partial}{\partial z_2}+
z_4\frac{\partial}{\partial z_3}, \hspace{2mm} \beta z_0
\frac{\partial} {\partial z_0}+(z_3+\delta
z_4)\frac{\partial}{\partial z_1}+z_4\frac{\partial }{\partial z_2},$$
 $$z_0\frac{\partial}{\partial z_0}+z_4\frac{\partial}{\partial z_1}
\hspace{2mm} \mbox{et } R.$$ La $1$-forme $\omega$ annul\'ee par
ces quatre champs s'\'ecrit \`a multiplication pr\`es dans la
carte affine $z_4=1$
$$\frac{dz_0}{z_0}+d\left(-z_1+z_2z_3+(\delta-\beta)z_2-\kappa
z_3+\frac{\beta-\delta}{2}z_3^2-\frac{z_3^3}{3}\right)$$ Le
feuilletage associ\'e admet pour int\'egrale premi\`ere
$$\frac{z_0}{z_4}\exp\left(\frac{-z_1z_4^2+z_2z_3z_4+(\delta-
\beta)z_2z_4^2-\kappa z_3z_4^2-\frac{z_3^3}{3}+\frac{\beta-\delta}
{2}z_3^2z_4}{z_4^3}\right)$$\vs

que l'on peut conjuguer \`a $$\frac{z_0}{z_4}\exp\left(\frac{z_1z_4^2+z_2z_3z_4
+z_3^3}{z_4^3}\right)$$\vs

\item[(vii)] Cas o\`u $\tilde{\mathscr{L}}$ est engendr\'ee par
les champs
$$(z_1+\kappa z_4)\frac{\partial}{\partial z_0}+z_2\frac{\partial}
{\partial z_1}+z_3\frac{\partial}{\partial z_2}+z_4\frac{\partial}
{\partial z_3}, \hspace{2mm} (z_2+\xi z_4)\frac{\partial}
{\partial z_0}+z_3\frac{\partial}{\partial z_1}+ z_4
\frac{\partial}{\partial z_2},$$ $$(z_3+\beta z_4)
\frac{\partial}{\partial z_0}+z_4\frac{\partial}{\partial z_1}
\hspace{2mm} \mbox{et } R.$$ La $1$-forme $\omega$ qui annule ces
champs s'\'ecrit \`a multiplication pr\`es dans la carte affine
$z_4=1$
$$d\left(z_0-z_1z_3-\beta z_1-\xi z_2+\beta
z_2z_3-\frac{z_2^2}{2}+z_2z_3^2-\kappa
z_3-\frac{z_3^4}{4}-\frac{\beta}{3}z_3^3+\frac{\xi}{2}z_3^2\right)$$
Le feuilletage d\'ecrit par cette $1$-forme admet pour int\'egrale
premi\`ere
$$\frac{(z_0-\beta z_1-\xi z_2-\kappa z_3)z_4^3+(\frac{\xi}{2}z_3
-z_1+\beta z_2)z_3z_4^2-\frac{z_2^2z_4^2}{2}-\frac{z_3^4}
{4}+(z_2-\frac{\beta}{3}z_3)z_3^2z_4}{z_4^4}$$ que l'on peut conjuguer \`a
$$\frac{z_0z_4^3+z_1z_3z_4^2-\frac{z_2^2z_4^2}{2}-\frac{z_3^4}
{4}+z_2z_3^2z_4}{z_4^4} $$
\vs
\end{itemize}

Le th\'eor\`eme suivant r\'esume l'analyse pr\'ec\'edente ; elle a \'et\'e faite
directement sans utiliser Maple.

\begin{thm} Soit $\mathscr{F}$ un
$\mathscr{L}$-feuilletage de degr\'e $3$ sur
$\C\P(4)$. Si $\mathscr{L}$ est ab\'elienne, alors
le feuilletage $\mathscr{F}$ admet \`a conjugaison pr\`es l'une
des int\'egrales premi\`eres suivantes :
\begin{eqnarray}
& & z_0^{\kappa_0}
z_1^{\kappa_1} z_2^{\kappa_2} z_3^{\kappa_3} z_4^{\kappa_4}
\hspace{8mm} \sum_{j=0}^4 \kappa_j=0 \nonumber \\
& &\nonumber \\
& &\frac{z_1^{\kappa+\xi+1}}{z_2^\kappa z_3^\xi z_4}\exp\left(
\frac{z_0}{z_1}\right) \nonumber \\
& &\nonumber \\
& &\frac{z_3^\kappa z_4}{z_1^{1+\kappa}}\exp\left(\frac{z_1z_2-z_0z_3}
 {z_1z_3}\right) \nonumber \\
& &\nonumber \\
& &\frac{z_0z_4^2-z_1z_2z_4+z_1z_3^2}{z_1z_4^2} \nonumber \\
& &\nonumber \\
& &\left(\frac{z_1}{z_4}\right)^\xi
\exp\left(\frac{z_0z_4^2-z_1z_2z_4+z_1z_3^2}{z_1z_4^2}\right)
\nonumber 
\end{eqnarray}

\begin{eqnarray}
%& &\nonumber \\
& & \frac{z_1^\xi z_4^{1-\xi}}{z_0}
\exp\left(\frac{z_3^2-z_2z_4}{z_4^2}\right)\nonumber \\
& &\nonumber \\
& & \frac{z_0}{z_4}\exp\left(\frac{z_1z_4^2+z_2z_3z_4-
z_3^3}{z_4^3}\right)\nonumber \\
& &\nonumber \\
& & \frac{z_0z_4^3+z_1z_3z_4^2-\frac{z_2^2z_4^2}{2}-\frac{z_3^4}
{4}+z_2z_3^2z_4}{z_4^4} \nonumber
\end{eqnarray}
Les $\kappa_i$, $\kappa$ et $\xi$ sont des nombres complexes.
\end{thm}

\section[{$\mathscr{L}$ a pour pr\'esentation $\{[X,Y]=Y, 
\hspace{1mm} [X,Z]=\alpha Z, \hspace{1mm} [Y,Z]=0\}, \hspace{1mm}
\alpha \not = 0$.}]{L'alg\`ebre $\mathscr{L}=\mathscr{L}_\alpha$ a pour 
pr\'esentation $\{[X,Y]=Y, \hspace{1mm} [X,Z]=\alpha Z, \hspace{1mm}
[Y,Z]=0\}, \hspace{1mm} \alpha \not = 0$.}

\vs

On remarque que si l'on change $X$ en
$\frac{X}{\alpha}$ on a $$[\frac{X}{\alpha},Y]=\frac{Y}{\alpha},
\hspace{2mm} [Y,Z]=0, \hspace{2mm} [\frac{X}{\alpha},Z]=Z$$ On
obtient avec $X'=\frac{X}{\alpha}$, $Y'=Z$ et $Z'=Y$
$$[X',Y']=Y', \hspace{2mm} [Y',Z']=0, \hspace{2mm}
[X',Z']=\frac{Z'}{\alpha}$$ Par suite l'alg\`ebre $\mathscr{L}$
est indiff\'eremment d\'etermin\'ee par $\alpha$ ou
$\frac{1}{\alpha}$ que l'on note $\tilde{\alpha}$ (i.e. on note
$\tilde{\alpha}=\tilde{2}$ pour $\alpha=2$ ou
$\alpha=\frac{1}{2}$).
\vs

L'alg\`ebre $\C Y \oplus \C Z$ est sous-alg\`ebre
ab\'elienne de $\mathscr{L}$ qui est r\'esoluble. La pr\'esentation entra\^{i}ne
que $Y$ et $Z$ sont nilpotents ; on peut donc choisir une base dans laquelle $Y$
s'\'ecrit

$$N_{2,1}=\left(%
\begin{array}{ccccc}
  0 & 0 & 0 & 1 & 0 \\
  0 & 0 & 0 & 0 & 1 \\
  0 & 0 & 0 & 0 & 0 \\
  0 & 0 & 0 & 0 & 0 \\
  0 & 0 & 0 & 0 & 0 \\
\end{array}%
\right), \hspace{2mm} N_{2,2}=\left(%
\begin{array}{ccccc}
  0 & 1 & 0 & 0 & 0 \\
  0 & 0 & 1 & 0 & 0 \\
  0 & 0 & 0 & 0 & 0 \\
  0 & 0 & 0 & 0 & 0 \\
  0 & 0 & 0 & 0 & 0 \\
\end{array}%
\right),$$ $$N_{3,1}=\left(%
\begin{array}{ccccc}
  0 & 0 & 1 & 0 & 0 \\
  0 & 0 & 0 & 1 & 0 \\
  0 & 0 & 0 & 0 & 1 \\
  0 & 0 & 0 & 0 & 0 \\
  0 & 0 & 0 & 0 & 0 \\
\end{array}%
\right),\hspace{2mm}
N_{3,2}=\left(%
\begin{array}{ccccc}
  0 & 1 & 0 & 0 & 0 \\
  0 & 0 & 1 & 0 & 0 \\
  0 & 0 & 0 & 1 & 0 \\
  0 & 0 & 0 & 0 & 0 \\
  0 & 0 & 0 & 0 & 0 \\
\end{array}%
\right),$$ $$N_4=\left(%
\begin{array}{ccccc}
  0 & 1 & 0 & 0 & 0 \\
  0 & 0 & 1 & 0 & 0 \\
  0 & 0 & 0 & 1 & 0 \\
  0 & 0 & 0 & 0 & 1 \\
  0 & 0 & 0 & 0 & 0 \\
\end{array}%
\right)$$ suivant son ordre de nilpotence et son rang (on rappelle
que si $Y$ est de rang $1$, le degr\'e du feuilletage chute et en
fait le feuilletage est un pull-back) ; on cherche alors la forme
de la matrice $Z$ : soit $\mathscr{C}(N_i)$ le commutateur de
$N_i$
$$\mathscr{C}(N_{i,j})=\{M \in \mathscr{M}(5,\C), \hspace{1mm}
MN_{i,j}-N_{i,j}M=0\}$$ Un calcul fastidieux mais sans
difficult\'e donne la description des alg\`ebres de Lie
$\mathscr{C}(N_{i,j})$ :
$$\mathscr{C}(N_{2,1})=\left\{ \left(%
\begin{array}{ccccc}
  a & b & c & \delta & e \\
  f & g & h & i & j \\
  0 & 0 & k & \ell & m \\
  0 & 0 & 0 & a & b \\
  0 & 0 & 0 & f & g \\
\end{array}%
\right)\right\}$$ $$\mathscr{C}(N_{2,2})=\left\{ \left(%
\begin{array}{ccccc}
  a & b & c & \delta & e \\
  0 & a & b & 0 & 0 \\
  0 & 0 & a & 0 & 0 \\
  0 & 0 & f & g & h \\
  0 & 0 & i & j & k \\
\end{array}%
\right)\right\}$$
 $$\mathscr{C}(N_{3,1})=\left\{ \left(%
\begin{array}{ccccc}
  a & c & e & g & i \\
  0 & b & \delta & f & h \\
  0 & 0 & a & c & e \\
  0 & 0 & 0 & b & \delta \\
  0 & 0 & 0 & 0 & a \\
\end{array}%
\right)\right\}$$ $$\mathscr{C}(N_{3,2})=\left\{ \left(%
\begin{array}{ccccc}
  a & b & c & \delta & e \\
  0 & a & b & c & 0 \\
  0 & 0 & a & b & 0 \\
  0 & 0 & 0 & a & 0 \\
  0 & 0 & 0 & f & g \\
\end{array}%
\right)\right\}$$
$$\mathscr{C}(N_4)=\left\{ \left(%
\begin{array}{ccccc}
  a & b & c & \delta & e \\
  0 & a & b & c & \delta \\
  0 & 0 & a & b & c \\
  0 & 0 & 0 & a & b \\
  0 & 0 & 0 & 0 & a \\
\end{array}%
\right)\right\}$$

Comme $Z$ est nilpotente, $Z$ appartient suivant l'ordre de
nilpotence et le rang de $Y$
 \`a l'un des espaces suivants
\begin{eqnarray}
&\mathscr{C}'(N_{2,1})&=\left\{\left(%
\begin{array}{ccccc}
  a & b & c & \delta & e \\
  f & -a & h & i & j \\
  0 & 0 & 0 & \ell & m \\
  0 & 0 & 0 & a & b \\
  0 & 0 & 0 & f & -a \\
\end{array}%
\right), \hspace{1mm} a^2+bf=0\right\} \nonumber \\
& & \nonumber \\
&\mathscr{C}'(N_{2,2})&=\left\{\left(%
\begin{array}{ccccc}
  0 & b & c & \delta & e \\
  0 & 0 & b & 0 & 0 \\
  0 & 0 & 0 & 0 & 0 \\
  0 & 0 & f & g & h \\
  0 & 0 & i & j & -g \\
\end{array}%
\right), \hspace{1mm} g^2+jh=0 \right\} \nonumber \\
& & \nonumber \\
&\mathscr{C}'(N_{3,1})&=\left\{\left(%
\begin{array}{ccccc}
  0 & c & e & g & i \\
  0 & 0 & \delta & f & h \\
  0 & 0 & 0 & c & e \\
  0 & 0 & 0 & 0 & \delta \\
  0 & 0 & 0 & 0 & 0 \\
\end{array}%
\right)\right\}\nonumber \\
& & \nonumber \\
&\mathscr{C}'(N_{3,2})&=\left\{ \left(%
\begin{array}{ccccc}
  0 & b & c & \delta & e \\
  0 & 0 & b & c & 0 \\
  0 & 0 & 0 & b & 0 \\
  0 & 0 & 0 & 0 & 0 \\
  0 & 0 & 0 & f & 0 \\
\end{array}%
\right)\right\}\nonumber \\
&  & \nonumber \\
&\mathscr{C}'(N_4)&=\left\{\left(%
\begin{array}{ccccc}
  0 & b & c & \delta & e \\
  0 & 0 & b & c & \delta \\
  0 & 0 & 0 & b & c \\
  0 & 0 & 0 & 0 & b \\
  0 & 0 & 0 & 0 & 0 \\
\end{array}%
\right) \right\}\nonumber
\vs
\end{eqnarray}

Soit $E \in \mathscr{M}(n,\C)$ un sous-espace vectoriel de
$\mathscr{M}(n,\C)$ ; on note $r(E) \in \mathbb{N}$ le
rang g\'en\'erique des \'el\'ements de $E$ : $$r(E)=max \{rg M,
\hspace{1mm} M \in E\}$$ et $E(M,N_{i,j})$ l'espace vectoriel
engendr\'e par $M$ et $N_{i,j}$. Soit $M \in
\mathscr{C}'(N_{i,j})$ tel que $\dim E(M,N_{i,j})=2$. On introduit
la sous-vari\'et\'e alg\'ebrique
$$F_i=\bigcup_j \{M \in \mathscr{C}'(N_{i,j}), \hspace{1mm}
r(E(M,N_{i,j}))=i\}$$ les ensembles
$$G_{i,j}=\{ M \in \mathscr{C}'(N_{i,j}), \hspace{1mm}
r(E(M,N_{i,j})) \leq i\}$$ et la propri\'et\'e $\mathscr{P}$ : $$
\mathscr{P}(a,b,c,\delta,e,f,g,i,j,k,\ell) \Leftrightarrow a=b=c=f=g=0
\mbox{ ou } \left\{
\begin{array}{ccccc}
 ak=-f\ell  \hspace{22mm} \\
  \delta\ell k = -\ell^2i+\ell jk+ek^2  \\
  ck=-g\ell  \hspace{23mm}\\
  bk^2=-f\ell^2   \hspace{19mm}\\
\end{array}\right.$$

Le lemme suivant donne la description des $F_i$ ; sa preuve
utilise uniquement des arguments standards de mineurs. \vs

\begin{lem}
On a
\begin{eqnarray}
&G_{2,1}&=\left\{\left(%
\begin{array}{ccccc}
 a & b & c & \delta & e \\
  f & -a & h & i & j \\
  0 & 0 & 0 & \ell & m \\
  0 & 0 & 0 & a & b \\
  0 & 0 & 0 & f & -a \\
\end{array}%
\right), \hspace{1mm} \mathscr{P}(a,b,c,\delta,e,f,h,i,j,\ell,m) \right\}
\nonumber \\
&  & \nonumber \\
&G_{2,2}&=\left\{\left(%
\begin{array}{ccccc}
  0 & b & c & \delta & e \\
  0 & 0 & b & 0 & 0 \\
  0 & 0 & 0 & 0 & 0 \\
  0 & 0 & f & 0 & 0 \\
  0 & 0 & i & 0 & 0 \\
\end{array}%
\right) \right\}\nonumber \\
&  & \nonumber \\
&G_{3,1}&=\left\{\left(%
\begin{array}{ccccc}
  0 & c & e & g & i \\
  0 & 0 & \delta & f & h \\
  0 & 0 & 0 & c & e \\
  0 & 0 & 0 & 0 & \delta \\
  0 & 0 & 0 & 0 & 0 \\
\end{array}%
\right), \hspace{1mm} c\delta=0\right\} \nonumber \\
&  & \nonumber \\
&G_{3,2}&=\left\{\left(%
\begin{array}{ccccc}
  0 & 0 & c & \delta & 0 \\
  0 & 0 & 0 & c & 0 \\
  0 & 0 & 0 & 0 & 0 \\
  0 & 0 & 0 & 0 & 0 \\
  0 & 0 & 0 & 0 & 0 \\
\end{array}%
\right)\right\} \nonumber \\
&  & \nonumber \\
&G_4&=\left\{\left(%
\begin{array}{ccccc}
  0 & b & c & \delta & e \\
  0 & 0 & b & c & \delta \\
  0 & 0 & 0 & b & c \\
  0 & 0 & 0 & 0 & b \\
  0 & 0 & 0 & 0 & 0 \\
\end{array}%
\right)\right\} \nonumber
\end{eqnarray}
\vs
\end{lem}

\begin{rem}
Si l'on note $G_i=\displaystyle \bigcup_j G_{i,j}$,
on obtient \`a partir de la description des $G_{ij}$ celle des $F_i$ :
$$F_2=G_2, \hspace{1mm} F_3=G_3 \setminus G_2, \hspace{1mm}
F_4=G_4 \setminus (G_2 \cup G_3)$$ 
\end{rem}

\vs

\begin{proof}[D\'emonstration] Commen\c{c}ons par d\'ecrire l'ensemble
 $G_{2,1}$ ; l'annulation pour tout r\'eel $t$ des
mineurs $(3,3)$ de la matrice $$\left(%
\begin{array}{ccccc}
  a & b & c & \delta+t & e \\
  f & -a & h & i & j+t \\
  0 & 0 & 0 & \ell & m \\
  0 & 0 & 0 & a & b \\
  0 & 0 & 0 & f & -a \\
\end{array}%
\right)$$ est r\'ealis\'ee si et seulement si l'une des conditions
suivantes est satisfaite :
$$\left\{
\begin{array}{l} a=f=b=c=h=0 \\
b\ell^2=-fm^2, \hspace{1mm} a\ell=-fm, \hspace{1mm} \delta m \ell
= -m^2i+mj\ell+e\ell^2, \hspace{1mm} c\ell=-hm
\end{array} \right.$$
Pour l'ensemble $G_{2,2}$, les mineurs $(3,3)$ de la matrice $$
\left(%
\begin{array}{ccccc}
  0 & b+t & c & \delta & e \\
  0 & 0 & b+t & 0 & 0 \\
  0 & 0 & 0 & 0 & 0 \\
  0 & 0 & f & g & h \\
  0 & 0 & i & j & -g \\
\end{array}%
\right)$$ sont nuls pour tout r\'eel $t$ si et seulement si $g=h=j=0$. \vs

On d\'ecrit maintenant $G_{3,1}$ ; la matrice
$$\left(%
\begin{array}{ccccc}
  0 & c & e+t & g & i \\
  0 & 0 & \delta & f+t & h \\
  0 & 0 & 0 & c & e+t \\
  0 & 0 & 0 & 0 & \delta \\
  0 & 0 & 0 & 0 & 0 \\
\end{array}%
\right)$$ est de rang inf\'erieur ou \'egal \`a $3$ si et seulement si
$c\delta=0$.

Passons
\`a la description de $G_{3,2}$ ; l'annulation pour tout r\'eel
$t$ des mineurs $(4,4)$
de la matrice $$\left(%
\begin{array}{ccccc}
  0 & b+t & c & \delta & e \\
  0 & 0 & b+t & c & 0 \\
  0 & 0 & 0 & b+t & 0 \\
  0 & 0 & 0 & 0 & 0 \\
  0 & 0 & 0 & f & 0 \\
\end{array}%
\right) $$ est r\'ealis\'ee si et seulement si $b=e=f=0$.

Pour $G_4$, il n'y a rien \`a faire.
\end{proof}

\vs

On aborde maintenant la classification des
$\mathscr{L}$-feuilletages associ\'es aux alg\`ebres
$\mathscr{L}_\alpha$.\vs

On se permettra d'effectuer certaines combinaisons lin\'eaires dites permises
qui vont changer $\mathscr{L}'$ tout en gardant la m\^eme structure
d'alg\`ebre : par exemple, les alg\`ebres $<X,Y,Z>$ et $<X+\kappa R,Y,Z>$ sont
distinctes mais isomorphes et tangentes au m\^eme feuilletage.

\subsection{Cas o\`u les
matrices $Y$, $Z \in F_4$} C'est le cas o\`u $Y$ s'\'ecrit $$\left(%
\begin{array}{ccccc}
  0 & 1 & 0 & 0 & 0 \\
  0 & 0 & 1 & 0 & 0 \\
  0 & 0 & 0 & 1 & 0 \\
  0 & 0 & 0 & 0 & 1 \\
  0 & 0 & 0 & 0 & 0 \\
\end{array}%
\right)$$ et o\`u $Z$ est de la forme $$\left(%
\begin{array}{ccccc}
  0 & b & c & \delta & e \\
  0 & 0 & b & c & \delta \\
  0 & 0 & 0 & b & c \\
  0 & 0 & 0 & 0 & b \\
  0 & 0 & 0 & 0 & 0 \\
\end{array}%
\right)$$ \vs

Pour des raisons techniques de calcul, on cherche $X \in
\chi(\C \P(4))$, $\lambda$, $\mu$, $\beta$,
$\varepsilon \in \C$ et $Z$ du type pr\'ec\'edent
satisfaisant $$[X,Y]=\lambda Y+\mu Z, \hspace{1mm} [X,Z]=\beta
Y+\varepsilon Z$$ La raison est qu'un \'el\'ement g\'en\'erique
de $[\mathscr{L}_\alpha, \mathscr{L}_\alpha]$ peut \^etre de rang
$4$ alors que les $Y$ et $Z$ de la pr\'esentation
pr\'ec\'edente peuvent ne pas l'\^etre. Les \'equations trait\'ees
sont quadratiques et il n'est pas \'etonnant qu'arrivent plusieurs
cas ; ils correspondent aux diff\'erentes composantes irr\'eductibles
des ensembles d\'ecrits par ces \'equations.  On n'a retenu que les
solutions donnant des feuilletages de degr\'e trois, certaines donnant
des feuilletages de degr\'e plus bas ; ces solutions \'elimin\'ees
interviennent dans la classification des feuilletages de degr\'e
inf\'erieur \`a trois.
\vs

Etudions les feuilletages associ\'es aux solutions donn\'ees par Maple. Nous
pr\'esentons les calculs typiques et mentionnons par quelles techniques se
traitent les autres cas.
\vs

\begin{itemize}
\item[3.1.1.] Les champs $X$ et $Z$ s'\'ecrivent respectivement
\`a combinaisons lin\'eaires permises pr\`es
$$X=(b z_1+cz_2+\delta z_3+e z_4) \frac{\partial}{\partial z_0}+
(\kappa z_1+bz_2+cz_3+\delta z_4)\frac{\partial}{\partial z_1}$$
$$+(2\kappa z_2+ bz_3+cz_4)\frac{\partial}{\partial z_2}+(3\kappa
z_3+bz_4)\frac{\partial}{\partial z_3}+4\kappa z_4\frac{\partial}
{\partial z_4}, $$ $$Z=z_2 \frac{\partial} {\partial z_0}+
z_3\frac{\partial}{\partial z_1}+z_4 \frac{\partial}{\partial z_2}
\hspace{2mm}$$ Dans ce cas $$\left(%
\begin{array}{cc}
  \lambda & \mu \\
  \beta & \varepsilon \\
\end{array}%
\right)=\left(%
\begin{array}{cc}
  \kappa & 0 \\
  -b\kappa & 2\kappa \\
\end{array}%
\right)$$
\vs

Ceci implique que $\kappa\not =0$ ; par suite on peut supposer que $\kappa=1$.
On note que $\tilde{\alpha}=\tilde{2}$.

Quitte \`a faire des combinaisons lin\'eaires permises, le champ
$X$ s'\'ecrit $$(-4z_0+\delta z_3+ez_4) \frac{\partial}{\partial
z_0}+(-3z_1+\delta z_4) \frac{\partial}{\partial
z_1}-2z_2\frac{\partial}{\partial z_2}-
z_3\frac{\partial}{\partial z_3}$$ Les hyperplans $z_4=$ cte sont
alors invariants par les trois champs qui engendrent donc le feuilletage
restreint \`a la carte affine $z_4=1$. Dans cette carte affine, on a :
$$\tilde{X}=X_{|z_4=1}=(-4z_0+\delta
z_3+e)\frac{\partial}{\partial z_0}+(-3z_1+\delta)
\frac{\partial}{\partial z_1}-2z_2\frac{\partial}{\partial z_2}-
z_3\frac{\partial}{\partial z_3}$$ $$\tilde{Y}=Y_{|z_4=1}=
z_1\frac{\partial}{\partial z_0}+z_2 \frac{\partial}{\partial
z_1}+z_3\frac{\partial}{\partial z_2}+\frac{\partial}{\partial
z_3}$$ $$\tilde{Z}=Z_{|z_4=1}=z_2 \frac{\partial} {\partial z_0}+
z_3\frac{\partial}{\partial z_1}+\frac{\partial}{\partial z_2}$$
Les champs $\tilde{Y}$ et $\tilde{Z}$ ont une \guillemotleft
\hspace{1mm} partie
constante \guillemotright \hspace{1mm}
donc ne s'annulent pas et commutent ; nous allons essayer de les
redresser par un automorphisme polynomial. Le flot de $\tilde{Y}$
est
$$\varphi(z;t)=
(z_0+z_1t+z_2\frac{t^2}{2}+z_3\frac{t^3}{6}+\frac{t^4}{24},z_1+
z_2t+z_3\frac{t^2}{2}+\frac{t^3}{6},z_2+z_3t+\frac{t^2}{2},z_3+
t)$$
et celui de $\tilde{Z}$ $$\psi \hspace{1mm} \colon \hspace{1mm}
(z;s) \mapsto (z_0+z_2s+\frac{s^2}{2},z_1+z_3s,z_2+s, z_3)$$ On
d\'efinit maintenant le diff\'eomorphisme
$$H(z_0,z_1,z_2,z_3) =
\varphi(\psi(z_0, z_1,0,0;z_2);z_3)$$ qui est donn\'e par
$$(z_0+\frac{z_2^2}{2}+z_1z_3+\frac{z_2z_3^2}{2}+
\frac{z_3^4}{24},z_1+z_2z_3+\frac{z_3^3}{6},z_2+\frac{z_3^2}{2},z_3)$$
Par construction $H$ conjugue $\frac{\partial}{\partial z_3}$ \`a
$\tilde{Y}$ et $\frac{\partial}{\partial z_2}$ \`a $\tilde{Z}$. 

On note $\tilde{\mathscr{F}}$ le feuilletage associ\'e \`a
l'alg\`ebre engendr\'ee par les champs $\tilde{X}$, $\tilde{Y}$ et
$\tilde{Z}$ ; soit $\tilde{\Omega}$ une $1$-forme qui d\'efinit
$\tilde{\mathscr{F}}$. Le feuilletage $\mathscr{F}'=H^*\tilde{\mathscr{F}}$
est d\'ecrit par une $1$-forme $\Omega'$ qui annule les champs
$H_*^{-1}\tilde{Y}=\frac{\partial}{\partial z_3}$ et
$H_*^{-1}\tilde{Z}=\frac{\partial}{\partial z_2}$. Elle s'\'ecrit donc
du fait de l'int\'egrabilit\'e
$$\Omega'= A'(z_0,z_1)
dz_0 +B'(z_0,z_1)dz_1$$ Le diff\'eomorphisme $H$ laisse le plan
$z_2=z_3=0$ invariant point par point, donc $\Omega'=\mbox{cte
}\iota^*\tilde{\Omega}$ o\`u $$\iota \hspace{1mm} \colon
\hspace{1mm} (z_0,z_1) \mapsto (z_0,z_1,0,0)$$ et $\Omega'$
d\'efinit $\tilde{\mathscr{F}}_{|z_2=z_3=0}$. On constate que le
champ $$\tilde{X}_{|z_2=z_3=0}=(-4z_0+e)\frac{\partial}{\partial
z_0}+(-3z_1+\delta)\frac{\partial}{\partial z_1}$$ %est tangent au
%feuilletage et au $2$-plan $z_2=z_3=0$ donc
d\'ecrit le feuilletage en restriction au $2$-plan $z_2=z_3=0$. On d\'eduit une
int\'egrale premi\`ere de $\tilde{\mathscr{F}}_{|z_2=z_3=0}$ :
$$\frac{(z_0-\frac{e}{4})^3}{(z_1-\frac{\delta}{3})^4}$$
qui donne apr\`es composition avec $H^{-1}$ et homog\'en\'eisation
$$\frac{{\left(z_0z_4^3-\frac{z_2^2z_4^2}{2}-z_1z_3z_4^2+
z_2z_3^2z_4-\frac{z_3^4}{4}-\frac{e}{4}z_4^4\right)}^3}{{\left(z_1z_4^2
-z_2z_3z_4+\frac{z_3^3}{3}-\frac{\delta}{3}z_4^3\right)}^4}$$ 
\`a conjugaison pr\`es s'\'ecrit $$\frac{{\left(z_0z_4^3-
\frac{z_2^2z_4^2}{2}-z_3z_4^2
+z_2z_3^2z_4-\frac{z_3^4}{4}\right)}^3}{{\left(z_1z_4^2
-z_2z_3z_4+\frac{z_3^3}{3}\right)}^4} $$ A posteriori, on
constate que ce cas s'obtient \`a partir de 3.1.3. en
faisant $\kappa=0$.
\vs

\item[3.1.2.] Les champs $X$ et $Z$ sont \`a combinaisons
lin\'eaires permises pr\`es (Maple)
$$X=(\xi_1 z_1+ \xi_2 z_2+\xi_3 z_3+\xi_4 z_4)
\frac{\partial}{\partial z_0}+4(b\eta+\varsigma)
z_4\frac{\partial}{\partial z_4}$$ $$+\left( (b\eta+\varsigma) z_1+
\xi_1 z_2+(\xi_2+\eta \delta) z_3+ z_4 \right)
\frac{\partial}{\partial z_1}$$
$$+\left( 2(b\eta+\varsigma)z_2+\xi_1z_3 + (\xi_2+2\eta
\delta)z_4\right) \frac{\partial}{\partial z_2}+\left( 3(b\eta +
\varsigma)z_3+\xi_1z_4\right) \frac{\partial}{\partial z_3},$$
$$Z=(b z_1+\delta z_3) \frac{\partial} {\partial z_0}+
(bz_2+\delta z_4)\frac{\partial}{\partial z_1}+bz_3
\frac{\partial}{\partial
z_2}+bz_4\frac{\partial}{\partial z_3}$$ et $$\left(%
\begin{array}{cc}
  \lambda & \mu \\
  \beta & \varepsilon \\
\end{array}%
\right)=\left(%
\begin{array}{cc}
  \varsigma & \eta \\
  -b(3b\eta+2\varsigma) & 4\eta b+3\varsigma \\
\end{array}%
\right)$$ \vs

Les valeurs propres de cette matrice sont $3(\varsigma+b\eta)$ et $\varsigma+
b\eta$ donc $\tilde{\alpha}=\tilde{3}$. On remarque que $\varsigma+b\eta$ est
non nul.

On note que si $\delta=0$, la dimension ponctuelle de
$\mathscr{L}$ est deux. Quitte \`a changer $Z$ en
$\frac{1}{\delta}(Z-bY)$, le champ $Z$ s'\'ecrit
$$z_3\frac{\partial}{\partial z_0}+z_4\frac{\partial}{\partial
z_1}$$ A combinaisons lin\'eaires et conjugaison pr\`es, on a
$$X=(-4(b\eta+\varsigma)z_0+\xi_4 z_4)
\frac{\partial}{\partial z_0}+(-3(b\eta+\varsigma) z_1+\eta
\delta z_3)\frac{\partial}{\partial z_1}$$
$$+(-2(b\eta+\varsigma)z_2+2\mu \delta z_4)\frac{\partial}{\partial
z_2}-(b\eta+\varsigma)z_3\frac{\partial}{\partial z_3}$$ On traite ce cas de la
m\^eme mani\`ere qu'en 3.1.1 et une int\'egrale
premi\`ere de $\mathscr{F}$ est \`a conjugaison pr\`es
$$\frac{{\left(z_2z_4-\frac{z_3^2}{2}-\kappa z_4^2\right)}^2}{\left
(2(b\eta+\varsigma)(z_1z_3z_4^2+\frac{z_2z_3^2z_4}{2}-
\frac{z_3^4}{8})\right)+z_0z_4^3}$$ avec
$$\kappa=\frac{2\eta\delta}{2(b\eta+\varsigma)}$$
\vs

On peut en fait conjuguer la fonction rationnelle pr\'ec\'edente
\`a $$\frac{{\left(z_2z_4-z_3^2\right)}^2}{z_1z_3z_4^2+z_2z_3^2z_4-
z_3^4+z_0z_4^3}$$

\vs

\item[3.1.3.] A combinaisons lin\'eaires pr\`es on a suivant
Maple :
$$X=(\xi_1 z_1+ \xi_2 z_2+\xi_3 z_3+\xi_4 z_4)
\frac{\partial}{\partial z_0}-\frac{8 \varsigma
c^2}{\delta}z_4\frac{\partial}{\partial z_4}$$ $$+ \left(
-\frac{2\varsigma c^2}{\delta} z_1+ (\xi_1+\varsigma c) z_2+(\xi_2+
\varsigma \delta)
z_3+ \frac{4\xi_3 c+5\varsigma \delta^2}{4c} z_4\right)
\frac{\partial}{\partial z_1}$$ $$+\left( -\frac{4\varsigma
c^2}{\delta}z_2+(\xi_1+2\varsigma c)z_3 + (\xi_2+2\varsigma
\delta)z_4\right)\frac{\partial}{\partial z_2}$$ $$+\left(
-\frac{6\varsigma c^2}{\delta} z_3+(\xi_1+3\varsigma
c)z_4\right)\frac{\partial}{\partial z_3},$$
$$Z=(bz_1+cz_2+\delta z_3+\frac{5\delta^2}{4c}z_4) \frac{\partial}
{\partial z_0}+ (bz_2+cz_3+\delta z_4)\frac{\partial}{\partial
z_1}$$ $$+(bz_3+cz_4) \frac{\partial}{\partial
z_2}+bz_4\frac{\partial}{\partial z_3}$$ et $$\left(%
\begin{array}{cc}
  \lambda & \mu \\
  \beta & \varepsilon \\
\end{array}%
\right)=\left(%
\begin{array}{cc}
  -\frac{\varsigma(\delta b+2c^2)}{\delta} & \varsigma \\
  -\frac{b\varsigma(\delta b-2c^2)}{\delta} & \frac{\varsigma(\delta
  b-4c^2)}{\delta} \\
\end{array}%
\right)$$ On remarque que $c$, $\varsigma$, $\delta\not =0$ et que $\tilde{\alpha}
= \tilde{2}$.\vs

Notons $$\kappa=\frac{\delta}{c},\hspace{2mm}
\xi'_4=\frac{\kappa}{\varsigma c}(\xi_4-5\xi_2\kappa^2) \mbox{ et }
\xi'_3=\frac{\kappa}{\varsigma c}(\xi_3-\xi_2\kappa)$$ Quitte \`a faire
des combinaisons lin\'eaires de $Z$ avec $Y$, le champ $Z$
s'\'ecrit
$$(z_2+\kappa z_3+\frac{5}{4}\kappa^2 z_4)
\frac{\partial} {\partial z_0}+ (z_3+\kappa z_4)\frac{\partial}
{\partial z_1}+z_4\frac{\partial}{\partial z_2}$$ De m\^eme
\`a combinaisons lin\'eaires permises pr\`es on a $$X=(8z_0+\xi'_3z_3+\xi'_4z_4)
\frac{\partial}{\partial z_0}+(6z_1 +
\kappa z_2+\kappa^2 z_3+(\xi'_3+\frac{5}{4}\kappa^3)z_4)\frac{
\partial}{\partial z_1}$$ $$+(4z_2+2\kappa z_3+2\kappa^2z_4) \frac
{\partial}{\partial z_2}+(2z_3+3\kappa z_4)\frac{\partial}
{\partial z_3}$$ On trouve que $\mathscr{F}$ admet comme
int\'egrale premi\`ere
$$\frac{{\left((z_1-\kappa z_2)z_4^2+(\frac{\kappa z_3^2}{2}-z_2z_3)z_4+
\frac{z_3^3}{3}
+(\frac{\xi'_3}{6}-\frac{\kappa^3}{8})z_4^3\right)}^4}
{{\left((2z_0+\frac{\kappa^2}{2} z_2-3\kappa z_1)z_4^3+((3 \kappa
z_2-2z_1)z_3-\frac{\kappa^2}{4}z_3^2-z_2^2)z_4^2+(2z_2z_3^2-\kappa
z_3^3)z_4-\frac{z_3^4}{2} -\eta z_4^4\right)}^3}$$ avec
$$\eta=\frac{3\kappa \xi'_3}{8}+\frac{11
\kappa^4}{32}+\frac{\xi'_4} {4}$$
soit \`a conjugaison pr\`es %$$\frac{{\left(z_1z_4^2+(\frac{\kappa
%z_3^2}{2}-z_2z_3)z_4+\frac{z_3^3}{3}
%+(\frac{\xi'_3}{6}-\frac{\kappa^3}{8})z_4^3\right)}^4}
%{{\left(z_0z_4^3+((\kappa z_2-2z_1)z_3-\frac{\kappa^2}{4}z_3^2-z_2^2)z_4^2+(2z_2
%z_3^2-\kappa z_3^3)z_4-\frac{z_3^4}{2} -\eta z_4^4\right)}^3} $$
$$\frac{{\left(z_1z_4^2-z_2z_3z_4+\frac{z_3^3}{3}\right)}^4}
{{\left(z_0z_4^3-(2z_1z_3+z_2^2)z_4^2-\frac{z_3^4}{2}\right)}^3}$$
\vs
\end{itemize}

\begin{rems} 
1. Ce feuilletage, ainsi que celui produit en 3.1.1, 
n'appartient pas aux composantes connues de l'espace 
des feuilletages de degr\'e $3$ sur $\C\P(4)$ ; c'est
un candidat \`a \^etre un feuilletage \guillemotleft
\hspace{1mm} \textbf{exceptionnel} \guillemotright.

2. Une d\'eg\'en\'erescence de ce cas correspond \`a $\delta=0$ (multiplier les
constantes par $\delta$ et annuler $\delta$). On obtient l'alg\`ebre engendr\'ee
par les champs $$R,\hspace{2mm} X=4z_0\frac{\partial}{\partial z_0}+3z_1
\frac{\partial}{\partial z_1}+2z_2\frac{\partial}{\partial z_2}+z_3
\frac{\partial}{\partial z_3}$$ $$Y=z_1\frac{\partial}{\partial z_0}+z_2
\frac{\partial}{\partial z_1}+z_3\frac{\partial}{\partial z_2}+z_4
\frac{\partial}{\partial z_3} ,\hspace{2mm}Z=z_2\frac{\partial}{\partial z_0}+
z_3\frac{\partial}{\partial z_1}+z_4\frac{\partial}{\partial z_2}$$ qui
correspond en fait au cas $\kappa=0$, i.e. 3.1.1.

La d\'eg\'en\'erescence $\varsigma=0$ produit $\alpha=0$ et n'est donc pas
\` a consid\'erer ici.   

3. On constate que lorsque le rang g\'en\'erique de $<Y,Z>$ est
maximal l'invariant $\tilde{\alpha}$ est discret
et vaut $\tilde{2}$ ou $\tilde{3}$.
\end{rems}

\vs

\subsection{Les matrices $Y$, $Z\in F_3$.}

Comme pr\'ec\'edemment on cherche $X$ dans $\chi(\C \P(4))$,
$\lambda$, $\mu$, $\beta$, $\varepsilon \in \C$
satisfaisant $$[X,Y]=\lambda Y+\mu Z \mbox{ et } [X,Z]=\beta 
Y+\varepsilon Z$$ Gr\^ace \`a Maple, on obtient les sept 
solutions suivantes que l'on \'etudie au cas par cas. On n'a
retenu que les solutions donnant des feuilletages de degr\'e
$3$.
\vs

\begin{itemize}
\item[3.2.1.] Les champs $X$, $Y$ et $Z$ s'\'ecrivent
respectivement \`a combinaisons lin\'eaires permises pr\`es
$$X=\left( -\frac{1}{2}\sigma z_0 +\xi_1z_1+\xi_2z_2+\xi_3z_3 + \xi_4
z_4\right) \frac{\partial}{\partial z_0}$$ $$+\left(
-\frac{1}{2}\frac{ i\sigma
+2k(g\zeta-\xi_1)}{g}z_2+\xi_8z_3+\xi_9z_4\right) \frac{\partial}
{\partial z_1}$$ $$+\left(\frac{1}{2}\sigma z_2+(\xi_1+\zeta g)z_3+
(\xi_2+\zeta i)z_4\right) \frac{\partial}{\partial z_2}$$
$$+\left(\sigma z_3-\frac{1}{2}\frac{i\sigma-2\xi_1k}{g}z_4\right)
\frac{\partial}{\partial z_3}+\frac{3}{2}\sigma z_4
\frac{\partial}{\partial z_4},$$ $$Y=z_2\frac{\partial}{\partial
z_0}+z_3 \frac{\partial}{\partial z_1}+z_4\frac{\partial}{\partial
z_2} ,$$
$$Z=(gz_3+iz_4)\frac{\partial}{\partial z_0}+kz_4\frac{\partial}
{\partial z_1}, \hspace{1mm} k \not =0$$
et $$\left(%
\begin{array}{cc}
  \lambda & \mu \\
  \beta & \varepsilon \\
\end{array}%
\right)=\left(%
\begin{array}{cc}
  \sigma & \zeta \\
  0 & \frac{3\sigma}{2} \\
\end{array}%
\right)$$ \vs
On note que $\tilde{\alpha}=\tilde{\frac{3}{2}}$.

Le feuilletage associ\'e \`a l'alg\`ebre
engendr\'ee par les champs $X$, $Y$, $Z$ et $R$ admet pour
int\'egrale premi\`ere
$$\frac{z_4\left(
(\frac{g}{k}(z_2-z_1)z_3-\frac{z_2^2}{2}-\kappa
z_3^2+\frac{i}{k}z_2z_3)z_4+(z_0-\frac{i}{k}z_1-\eta z_3)z_4^2-2Cz_3^3+
\nu z_4^3\right)}{\left(z_3+\frac{\xi'_1}{\sigma}z_4\right)^4}$$ o\`u
$$\xi'_1=\frac{i\sigma-2\xi_1
k}{g},\hspace{2mm} \xi'_2=\frac{(\xi_2-\xi_8+\zeta
i)i-\xi_9g}{k}+\xi_3,$$ $$\xi'_3=i\frac{\xi_1+2\zeta
g}{k}+\frac{(\xi_2-\xi_8)g}{k},\hspace{2mm}
\xi'_4=\frac{g(\xi_1+\zeta g)}{k} \mbox{ et }
\xi'_5=\xi_4-\frac{i\xi_9} {k}$$
$$A=
\frac{\xi'_2}{\sigma}-\frac{2\xi'_1\xi'_3}{\sigma^2}+\frac{3\xi^{'2}_1
\xi'_4}{\sigma^3},\hspace{2mm} B=\frac{\xi'_3}{\sigma}-
\frac{3\xi'_4\xi'_1}{\sigma^2}, \hspace{2mm} C=\frac{\xi'_4}
{\sigma},$$ $$\xi''_1=\frac{\xi'_5}{\sigma}-\frac{\xi'_2\xi'_1}
{\sigma^2}+ \frac{\xi^{'2}_1\xi'_3}{\sigma^3}-
\frac{\xi^{'3}_1\xi'_4}{\sigma^4},\hspace{2mm} \kappa=B+\frac{6C
\xi'_1}{\sigma}, $$
$$\eta=\frac{2A}{3}+\frac{2\xi'_1B}{\sigma}+\frac{6C
\xi^{'2}_1}{\sigma^2} \mbox{ et } \nu=-\frac{\xi''_1}{2}
-\frac{B\xi^{'2}_1}{\sigma^2}-2C\frac{\xi_1{'^3}}{\sigma^3}-
\frac{2A\xi'_1}{3\sigma}$$ ou encore \`a conjugaison pr\`es
%$$\frac{z_4\left( (z_1z_3-\frac{z_2^2}{2})z_4+z_0z_4^2-2Cz_3^3\right)}
%{\left(z_3+\frac{\xi'_1}{\sigma}z_4\right)^4} $$
$$\frac{z_4(z_0z_4^2+z_1z_3z_4+z_2^2z_4+az_3^3)}{z_3^4}$$ o\`u
$a\in \{0,1\}$.
\vs

\begin{rem}
Dans la carte affine $z_3=1$, le
feuilletage admet pour int\'egrale premi\`ere \`a conjugaison pr\`es 
$$z_4(z_0z_4^2+z_1z_4+z_2^2z_4+a)$$ Il produit donc un polyn\^ome de
degr\'e quatre sur $\C^4$ d\'efinissant un
$\mathscr{L}$-feuilletage.
\end{rem}

\vs

\item[3.2.2.] Les champs $X$, $Y$ et $Z$ sont \`a combinaisons
lin\'eaires permises pr\`es
$$X=\left( \frac{(e-f)(\xi_7+\kappa k)}{k}z_0+\xi_3z_3+\xi_4z_4
\right) \frac{\partial}{\partial z_0}$$
$$+\left(
\xi_7z_2+(\xi_8-\xi_2)z_3+\xi_9z_4\right)\frac{\partial} {\partial
z_1}+\left(\frac{(e-f)\xi_7}{k}z_2+\kappa iz_4\right)
\frac{\partial}{\partial z_2}$$ $$+\left(2\kappa(f-e)z_3+(\xi_7+\kappa
k)z_4\right) \frac{\partial}{\partial z_3}+\frac{(f-e)(\kappa k
-\xi_7)}{k} z_4 \frac{\partial}{\partial z_4},$$
$$Y=z_2\frac{\partial}{\partial z_0}+z_3 \frac{\partial}{\partial
z_1}+z_4\frac{\partial}{\partial z_2},$$
$$Z=iz_4 \frac{\partial}{\partial z_0}+((f-e)z_3+kz_4)\frac{\partial}
{\partial z_1} \hspace{4mm} i \not = 0, \hspace{1mm} k \not =0$$
et $$\left(%
\begin{array}{cc}
  \lambda & \mu \\
  \beta & \varepsilon \\
\end{array}%
\right)=\left(%
\begin{array}{cc}
  \kappa(f-2e) & \kappa \\
  -ef\kappa & \kappa(2f-e) \\
\end{array}%
\right)$$ \vs

Les valeurs propres de cette matrice sont $\kappa(f-e)$ et $2\kappa(f-e)$ donc
$\tilde{\alpha}=\tilde{2}$.
\vs

On montre que le feuilletage associ\'e \`a
l'alg\`ebre engendr\'ee par les champs $X$, $Y$, $Z$ et $R$ admet
pour int\'egrale premi\`ere
$$\frac{z_4^{3\xi'_5-\xi'_1}\left(z_3+\frac{\xi'_6}{\xi'_5}z_4\right)^{\xi'_1}}
{{\left(z_1z_4^2-z_2z_3z_4+(\frac{f-e}{2i}z_3-\frac{k}{2i}z_4)(z_2^2-2z_0z_4)
-(A+2B\frac{\xi'_6}{\xi'_5})z_3z_4^2-Bz_3^2z_4-\nu
z_4^3\right)}^{\xi'_5}}$$ avec
$$\xi'_1=\frac{(f-e)(\xi_7-k\kappa)}{k},\hspace{2mm}
\xi'_5=\frac{(\kappa k+\xi_7)(f-e)}{k},\hspace{2mm} \xi'_6=\xi_7+\kappa
k,$$ $$A=\frac{\xi'_2\xi'_5- 2\xi'_3\xi'_6}{\xi'_5(\xi'_5
-\xi'_1)}, \hspace{2mm} B=\frac{\xi'_3}{2\xi'_5-\xi'_1} \mbox{ et
}$$ $$\nu=-\frac{\xi'_2\xi'_6}{\xi'_1\xi'_5}+\frac{\xi'_3
\xi^{'2}_6}{\xi'_1\xi^{'2}_5}+\frac{\xi'_4}{\xi'_1}-A\frac{\xi'_6}
{\xi'_5}-B\frac{\xi^{'2}_6}{\xi^{'2}_5}$$

\vs

Cette int\'egrale premi\`ere est conjugu\'ee \`a : 

\vs

$$\frac{z_4^{3\xi'_5-\xi'_1}z_3^{\xi'_1}}
{{\left(z_1z_4^2+z_2z_3z_4+z_2^2z_3+z_0^2z_4\right)}^{\xi'_5}} $$
\vs

\item[3.2.3.] A combinaisons lin\'eaires permises pr\`es on a
$$X=((-\nu z_0+\xi_3z_3+\xi_4z_4)\frac{\partial}{\partial
z_0} +(\xi_8z_3+\xi_9z_4) \frac{\partial}{\partial z_1}$$
$$+ ((-\nu-\delta)z_2+\varsigma iz_4) \frac{\partial}{\partial
z_2} - 2\delta z_3 \frac{\partial}{\partial z_3} +
(-\nu-2\delta)z_4 \frac{\partial}{\partial z_4},$$
$$Y=z_2\frac{\partial}{\partial z_0}+z_3 \frac{\partial}{\partial
z_1}+z_4\frac{\partial}{\partial z_2},$$
$$Z=iz_4 \frac{\partial}{\partial z_0}+(f-e)z_3
\frac{\partial}{\partial z_1} \hspace{4mm} i \not = 0, \hspace{1mm}
f-e \not =0$$ o\`u $\delta=\varsigma(e-f)$
et $$\left(%
\begin{array}{cc}
  \lambda & \mu \\
  \beta & \varepsilon \\
\end{array}%
\right)=\left(%
\begin{array}{cc}
  \varsigma(-2e+f) & \varsigma \\
  -e\varsigma f & \varsigma(2f-e) \\
\end{array}%
\right)$$
On note que $\tilde{\alpha}=\tilde{2}$.\vs

Posons $$\eta=\xi_8-\xi_2-\varsigma i-\frac{\xi_4}{i}(f-e)$$ et
$$\kappa=\frac{\xi_3(e-f)}{i}$$\vs

Si $\nu+2\delta=0$, le feuilletage associ\'e \`a
l'alg\`ebre engendr\'ee
 par les champs $X$, $Y$, $Z$ et $R$ admet pour int\'egrale premi\`ere
$$\left(\frac{z_3}{z_4}\right)^{\xi_9} \exp\left(\frac{\eta z_3
z_4^2+\frac{\kappa}{2}z_3^2z_4-\nu(z_1z_4^2-z_2z_3z_4-\frac{f-e}{i}z_0
z_3z_4+\frac{(f-e)}{2i}z_2^2z_3)}{z_4^3}\right) $$ c'est-\`a-dire \`a
conjugaison pr\`es $$\left(\frac{z_3}{z_4}\right)^{\xi_9} \exp\left(\frac{
z_1z_4^2+z_0z_3z_4+z_2^2z_3}{z_4^3}\right)$$
\vs

Si $2\delta+\nu\not =0$ et $2\delta-\nu \not =0$, alors le feuilletage admet
pour int\'egrale premi\`ere
$$\frac{z_3^{2\delta+\nu}z_4^{2(\nu-\delta)}}{\left(z_1z_4^2
-z_2z_3z_4-\frac{f-e}{i}z_0z_3z_4+\frac{f-e}{2i}z_2^2z_3+ \frac{\eta}{2
\delta}z_3z_4^2-\frac{\kappa}{2\delta-\nu}z_3^2z_4+\frac{\xi_9}{2\delta+\nu}
z_4^3\right)^\nu}$$ i.e. \`a conjugaison pr\`es $$\frac{z_3^{2\delta+\nu}
z_4^{2(\nu-\delta)}}{\left(z_1z_4^2+z_0z_3z_4+z_2^2z_3\right)^\nu} $$ \vs

Finalement si $2\delta+\nu\not =0$ et $\nu=2\delta$, l'int\'egrale premi\`ere
s'\'ecrit $$\left(\frac{z_3}{z_4}\right)^\beta
\exp\left(\frac{2\delta(z_1z_4^2-z_2z_3z_4-\frac{f-e}{i}z_0z_3z_4+\frac{f-e}
{2i}z_2^2z_3)}{z_3^2z_4}\right)$$ ou encore \`a conjugaison pr\`es
$$\left(\frac{z_3}{z_4}\right)^{\frac{\beta}{2\delta}}\exp\left(\frac{
z_1z_4^2+z_0z_3z_4+z_2^2z_3}{z_3^2z_4}\right) $$
\vs

\item[3.2.4.] Les champs $X$, $Y$ et $Z$ s'\'ecrivent \`a
combinaisons lin\'eaires pr\`es
$$X=\left(2\varsigma(f-e)z_0-2\varsigma g z_1+\frac{\varsigma i(2e-
f)+g(\xi_7+3\varsigma k)+i\eta}{\delta} z_3+\xi_4z_4\right)
\frac{\partial}{\partial z_0}$$
$$+\left(-2\varsigma \delta z_0+\xi_7z_2+(\xi_8+\frac{3\varsigma i
\delta+ \xi_7(f-e)+\varsigma k
(2f-e)+k\eta-\xi_8\delta}{\delta})z_3+\xi_9z_4\right)
 \frac{\partial}{\partial z_1}$$ $$+((\eta+\varsigma(2f-e))z_2-
 \varsigma g z_3+\varsigma iz_4)\frac{\partial}{\partial z_2}$$
$$+(-\varsigma \delta z_2+(\eta+\varsigma f)z_3+(\xi_7+\varsigma
 k) z_4)\frac{\partial}{\partial z_3}+2(\varsigma f+\eta)z_4
 \frac{\partial}{\partial z_4},$$
$$Y=z_2\frac{\partial}{\partial z_0}+z_3 \frac{\partial}{\partial
z_1}+z_4\frac{\partial}{\partial z_2},$$
$$Z=(gz_3+iz_4)\frac{\partial}{\partial z_0}+(\delta z_2+(f-e)z_3+
 kz_4)\frac{\partial}{\partial z_1}+\delta z_4\frac{\partial}
 {\partial z_3}$$
On traite ce cas comme d'habitude en trivialisant deux champs tangents dans la
carte affine $z_4=1$. Suivant les valeurs du param\`etre, on trouve les
int\'egrales premi\`eres suivantes : \vs
$$\frac{(P-\gamma^+Q)
^{\vartheta^+}(P-\gamma^-Q) ^{\vartheta^-}}{z_4^2}$$ si $e^2+\delta g
\not =0$,
$$\left(\frac{P-\frac{e}{2\delta}Q}{z_4^2}\right)^{\varsigma
\delta}\exp\left(\frac{\zeta P}{P-\frac{e}{2\delta}Q
}\right) $$ sinon. Les polyn\^omes $P$ et $Q$ sont d\'efinis par :
$$P(z_0,z_1,z_2,z_3,z_4)=z_0z_4-
\frac{gz_3^2}{2\delta}-\frac{iz_3z_4} {\delta}-\frac{z_2^2}{2}+
\kappa_0z_4^2$$ $$Q(z_0,z_1,z_2,z_3,z_4)=z_1z_4-
\frac{f-e}{2\delta}z_3^2-\frac{k} {\delta} z_3z_4-z_2z_3
-\kappa_1z_4^2$$ et
 $$\kappa_0=\frac{(2\varsigma
%\delta-\xi_9\delta+k\xi_7+\varsigma k^2)\varsigma g+(\xi_4\delta-i\xi_7-
i\varsigma k)\eta}{2\delta(\varsigma^2g \delta-\varsigma
e\eta-\eta^2)}, \hspace{2mm} \zeta=-\frac{2\delta}{e}\left(
\frac{\varsigma e}{2}+\eta\right)$$
$$\kappa_1=\frac{(\delta\varsigma\xi_4+2f(\varsigma e+
\eta))\delta\varsigma-\xi_9\delta(\eta+e\varsigma)+(\xi_7+
\varsigma k)(e\varsigma
k-i\delta\varsigma+\eta k) }{2(\varsigma^2 g\delta-\varsigma
e\eta-\eta^2)\delta}$$ $$\gamma^\pm=\frac{e \pm \sqrt{(e^2+4\delta g)}}{2\delta} \mbox{ et
}\vartheta^\pm=\frac{1}{2}\left(1\pm\frac{2\eta+e\varsigma}{
\sqrt{\varsigma^2(e^2+4g\delta)}}\right) $$\vs Ce qui donne \`a conjugaison
pr\`es $$\frac{{\left((z_0-\gamma^+z_1)z_4+\frac{\gamma^+(f-e)-g}{2\delta}z_3^2
-\frac{z_2^2}{2}+\gamma^+z_2z_3\right)}^{\vartheta_+}{\left((z_0-\gamma^-z_1)z_4+
\frac{\gamma^-(f-e)-g}{2\delta}z_3^2-\frac{z_2^2}{2}+\gamma^-z_2z_3\right)}^{
\vartheta_-}}{z_4^2} $$ si $e^2+\delta g\not =0$, $$\frac{(z_0-z_1)z_4+e
z_2z_3-\frac{z_2^2}{2}+\upsilon z_3^2}{z_4^2}
\exp\left(\frac{\zeta\left(z_0z_4-2g\delta z_3^2-\frac{z_2^2}{2}\right)}
{(z_0-z_1)z_4+ez_2z_3-\frac{z_2^2}{2}+\upsilon z_3^2}\right)$$ avec $$\upsilon=e(f-e)-2g\delta$$ sinon. On peut \'eventuellement simplifier
 encore ces expressions ; mais il faut tenir compte de l'annulation ou non 
de certains param\`etres. Nous en laissons le soin au lecteur.
\vs

%Les valeurs propres de la matrice $$\left(%
%\begin{array}{cc}
%  \lambda & \mu \\
%  \mu(g\delta-ef) & \mu(f+e)+\lambda \\
%\end{array}%
%\right)$$ sont $\rho^\pm=\frac{\mu}{2}(e+f)+\lambda\pm
%\frac{1}{2}\sqrt{\mu^2((e-f)^2+4g\delta)}$ donc
%$\tilde{\alpha}=\tilde{\frac{\rho^+}{\rho^-}}$.
%\vs

\item[3.2.5.] Etudions le cas suivant pr\'edit par Maple :
$$Y=z_2 \frac{\partial}{\partial z_0}+z_3 \frac{\partial}{\partial
z_1}+z_4 \frac{\partial}{\partial z_2}$$
$$Z=(z_1+gz_3+iz_4)\frac{\partial} {\partial z_0}+ (fz_3+kz_4)
\frac{\partial}{\partial z_1}+z_3 \frac{\partial}{\partial z_2}$$
Notons,  et ceci se constate dans les cartes $z_4=1$ puis
$z_3=1$ que l'on peut par conjugaison supposer que $g=i=0$. Lorsqu'on
demande \`a Maple $X,\lambda,\mu,\varepsilon,\beta$ tels que $[X,Y]=
\lambda X+\mu Y$ et $[X,Z]=\beta X+\varepsilon Y$ il donne %$$X=
%\left(\xi_1z_1+(\xi_8-\xi_1f)z_2+\xi_3z_3+\xi_4z_4\right)
%\frac{\partial}{\partial z_0}$$$$+\left((\zeta+\delta f)z_1+\delta k
%z_2+\xi_8z_3+\xi_1kz_4\right)\frac{\partial}{\partial z_1}+\left(
%\delta z_1+\zeta z_2+\xi_1z_3+(\xi_8-\xi_1f)z_4\right)\frac{\partial}
%{\partial z_2}$$$$+\left(2(\zeta+\delta f)z_3+2\delta k z_4\right)
%\frac{\partial}{\partial z_3}+(2\delta z_3+2\delta z_4)\frac{\partial}
%{\partial z_4} $$
%On constate que quitte
%\`a retrancher $\xi_1Z+\xi_8Y$ on peut supposer que $\xi_1=\xi_8=0$,
%i.e.
 $$X=\left(\xi_3z_3+\xi_4z_4\right)
\frac{\partial}{\partial z_0}$$ $$+\left((\zeta+\delta f)z_1+\delta k
z_2+\xi_1kz_4\right)\frac{\partial}{\partial z_1}+\left(
\delta z_1+\zeta z_2\right)\frac{\partial}{\partial z_2}$$ $$
+\left(2(\zeta+\delta f)z_3+2\delta k z_4\right)
\frac{\partial}{\partial z_3}+(2\delta z_3+2\delta z_4)\frac{\partial}
{\partial z_4}$$
et $$\left(%
\begin{array}{cc}
  \lambda & \mu \\
  \beta & \varepsilon \\
\end{array}%
\right)=\left(%
\begin{array}{cc}
  \zeta & \delta \\
  \delta k & \delta f + \zeta \\
\end{array}%
\right)$$ matrice qui doit \^etre inversible. \vs

Les m\'ethodes que nous avons appliqu\'ees jusqu'ici
conduisent \` a des calculs relativement p\'enibles. Nous allons
essayer une autre m\'ethode ; supposons que l'on puisse trouver $X'$
tel que par exemple $[X,X']=0$ et tel que $X'$ ne soit pas tangent au
feuilletage. Alors $X'$ est n\'ecessairement une sym\'etrie ; ce qui
produit un facteur int\'egrant qui en principe permet de calculer
l'int\'egrale premi\`ere. Cherchons par exemple, et cela conviendra
dans presque tous les cas ($\delta \not =0$), $X'$ sous la forme
$$(\xi'_3z_3+\xi'_4z_4)
\frac{\partial}{\partial z_0}+z_1\frac{\partial}{\partial z_1}+
z_2\frac{\partial}{\partial z_2}+2z_3\frac{\partial}{\partial z_3}+
2z_4\frac{\partial}{\partial z_4} $$ ceci
correspondant aux choix $\mu'=0$ et $\lambda'=1$. Pour que $[X,X']=0$
il faut et il suffit que l'on ait $$\left\{\begin{array}{c}
(\zeta+\delta f)\xi'_3+\delta\xi'_4=\xi_3 \\ \delta k \xi'_3+\zeta
\xi'_4=\xi_4 \hspace{9mm} \end{array}\right.$$ qui poss\`ede une solution unique
 puisque la matrice $$\left(%
\begin{array}{cc}
  \zeta & \delta \\
  \delta k & \delta f + \zeta \\
\end{array}%
\right)$$ est inversible. Dans le cas o\`u $\delta \not =0$ on v\'erifie
que $X'$ n'est pas tangent au feuilletage ; par contre, lorsque
$\delta=0$, le champ $X'$ est colin\'eaire \`a $X$. La r\'esolution explicite
de $\xi'_3$ et $\xi'_4$ donne $$ \xi'_3=\frac{\zeta\xi_3-\delta\xi_4}
{\nu} \hspace{1mm} \mbox{ et } \hspace{1mm}  \xi'_4=\frac{-\mu k\xi_3+(\zeta+
\delta f)\xi_4}{\nu}$$ avec $$\nu=\zeta(\zeta+\delta f)-\delta^2k$$

On proc\` ede alors au calcul de $$\omega=i_Xi_Yi_Zi_R dz_0 \wedge dz_1 \wedge
dz_2\wedge dz_3 \wedge dz_4$$ $1$-forme qui d\'ecrit le feuilletage
 et de $$P=\omega(X')=i_{X'}i_Xi_Yi_Zdz_0 \wedge dz_1 \wedge dz_2
\wedge dz_3 \wedge dz_4$$ facteur int\'egrant de cette $1$-forme. 
Si $\kappa=\sqrt{\delta^2f^2+4\delta^2 k}$, on constate que 
les deux hyperplans $$(\delta f+\kappa)z_3+2\delta
kz_4=0$$ et $$(\delta f-\kappa)z_3+2\delta kz_4=0$$ sont invariants par
le feuilletage car invariants par $X$ et \'evidemment par $Y$ et $Z$.
 Le polyn\^ome $P$ s'\'ecrit donc \` a constante multiplicative pr\` es
$$\left((\delta f+
\kappa)z_3 +2\delta kz_4\right)\left((\delta f-\kappa)z_3+2\delta
kz_4\right)Q$$ o\`u $$Q=
-2\nu (k+fz_3-z_3^2)z_0-\nu(fz_3+k)z_2^2+2\nu z_1z_2z_3$$ $$-
\nu z_1^2-\sigma z_3^3
+(\varsigma+f\sigma)z_3^2
+(k\sigma-\varsigma f)z_3
-k\varsigma$$ avec $$\varsigma=\xi_4\zeta-\delta k\xi_3+
\delta\xi_4f \mbox{ et } \sigma=\delta\xi_4-\zeta\xi_3$$
Le polyn\^ome $Q$, obtenu \`a l'aide de Maple, est de degr\'e $3$ ;
g\'en\'eriquement  sur les valeurs des param\`etres il est
irr\'eductible. Ainsi g\'en\'eriquement le feuilletage associ\'e \`a
l'alg\`ebre engendr\'ee par les champs $X$, $Y$, $R$ et $Z$ admet une
int\'egrale premi\`ere de la forme $$\left((\delta f+\kappa)z_3
+2\delta  k
z_4\right)^{\lambda_0}\left((\delta f-\kappa)z_3+2\delta kz_4
\right)^{\lambda_1}\tilde{Q}^{\lambda_2}$$
o\`u $\tilde{Q}$ est l'homog\'en\'eis\'e de $Q$ ; on constate que 
$\tilde{Q}$ vaut : 

\vs

$$-2\nu (kz_4^2+fz_3z_4-z_3^2)z_0
-\nu(fz_3+kz_4)z_2^2+2\nu z_1z_2z_3$$ $$-\nu z_1^2z_4-\sigma z_3^3+(\varsigma+
f\sigma)z_3^2z_4+(k\sigma-\varsigma f)z_3z_4^2-k\varsigma z_4^3$$

\vs

Les $\lambda_i$ se calculent \`a partir de $\frac{\omega}{\omega(X')}$. 

Pour les valeurs non g\'en\'eriques des param\`etres o\`u
le polyn\^ome $P$ n'est pas r\'eduit, on trouve des int\'egrales
premi\`eres contenant des exponentielles.\vs

Le polyn\^ome $P$ s'\'ecrit $\delta^2 \tilde{P}$ donc s'annule
pour $\delta=0$. Par contre $\tilde{P}$ ne s'annule pas pour
$\delta=0$ et par suite est un facteur int\'egrant du feuilletage
correspondant \`a $\delta=0$. On traite ce cas comme d'habitude
suivant les valeurs des param\`etres ; il ne donne pas de
ph\'enom\`ene nouveau.\vs

\item[3.2.6.]  A combinaisons lin\'eaires permises pr\`es on a
$$X=(-3\varsigma z_0+\xi_2z_2+\xi_3z_3+\xi_4z_4)\frac{\partial}
{\partial z_0}+(-2\varsigma z_1+\gamma cz_2+(\xi_2+\gamma\delta)z_3)
\frac{\partial} {\partial z_1}$$
$$+(-\varsigma z_2+2\gamma cz_3)\frac{\partial} {\partial z_2}+
(\xi_{0}z_3+(\xi_{1}-3\varsigma)z_4)\frac{\partial}{\partial z_4},$$
$$Y=z_1\frac{\partial}{\partial z_0}+z_2 \frac{\partial}{\partial
z_1}+z_3\frac{\partial}{\partial z_2},$$
$$Z=(z_2+\delta z_3)\frac{\partial}{\partial z_0}+z_3
\frac{\partial}{\partial z_1}$$
et $$\left(%
\begin{array}{cc}
  \lambda & \mu \\
  \beta & \varepsilon \\
\end{array}%
\right)=\left(%
\begin{array}{cc}
  \varsigma & -\frac{\varsigma\delta}{2} \\
  0 & 2\varsigma \\
\end{array}%
\right)$$\vs
On note que $\varsigma \not = 0$ et $\tilde{\alpha}=2$. \vs

Si $\xi_1-3\varsigma\not =0$, le
feuilletage associ\'e \`a l'alg\`ebre engendr\'ee par les
champs $X$, $Y$, $Z$ et $R$ admet pour int\'egrale premi\`ere
$$\frac{(z_4+\kappa z_3)^{3\varsigma}{\left(\xi_1(z_0z_3^2+\nu z_3^3
+\delta(\frac{z_2^2}{2}-z_1z_3)z_3-z_1z_2z_3+\frac{1}{3}z_2^3)-
\xi_4(z_4+\kappa z_3)z_3^2\right)}^{\xi_1-3\varsigma}}
{z_3^{3\xi_1-6\varsigma}}$$ avec
$$\kappa=\frac{\xi_0}{\xi_1-3\varsigma} \mbox{ et } \nu=\frac{1}
{3\varsigma}\left(\xi_4\kappa+\delta(\xi_2+\mu\delta)-\xi_3\right)$$
A conjugaison pr\`es on obtient $$\frac{z_4^{3\varsigma}{\left(z_0z_3^2+
z_1z_2z_3+z_2^3\right)}^{\xi_1-3\varsigma}}{z_3^{3\xi_1-6
\varsigma}} $$
\vs

Si $\xi_1=0$, alors le feuilletage est d\'ecrit par $$\left(\frac{z_4}{z_3}\right)^{\frac{\xi_4}{3\varsigma}}\exp\left( \frac{z_0z_3^2+\delta(\frac{1}{2}z_2^2z_3-z_1z_3^2)-z_1z_2z_3+
\frac{1}{3}z_2^3}{z_3^2z_4}\right)$$ ou encore \`a
conjugaison pr\`es $$\left(\frac{z_4}{z_3}\right)^{\frac{\xi_4}
{3\varsigma}}\exp\left(\frac{z_0z_3^2+z_1z_2z_3+z_2^3}{z_3^2z_4}\right)$$

Pour $\xi_1-3\varsigma=0$, on trouve comme int\'egrale premi\`ere
$$\frac{\left(z_0z_3^2-\delta z_1z_3^2+\frac{\delta}{2} z_2^2z_3
-z_1z_2z_3+\frac{1}{3}z_2^3+\rho z_3^3-\frac{\xi_4}{3
\varsigma}z_3^2z_4\right)}{z_3^3}\exp\left(\frac{3\varsigma z_4}
{\xi_0 z_3}\right)$$
avec $$\rho=\frac{1}{3\varsigma}\left(\delta(\xi_2
+\mu\delta)-\xi_3+\frac{\xi_0\xi_4}{3\varsigma}\right)$$ A
conjugaison pr\`es elle est de la forme $$\frac{z_0z_3^2+z_1z_2
z_3+z_2^3}{z_3^3}\exp\left(\frac{z_4}{z_3}\right) $$
\vs
\end{itemize}

\subsection{Cas o\`u les matrices $Y$, $Z\in F_2$.}

Maple nous a donn\'e trois types de champs $Z$ que l'on va examiner. \vs
\begin{itemize}
\item[3.3.1. ] Le premier cas est le suivant
$$Y=\left(%
\begin{array}{ccccc}
  0 & 1 & 0 & 0 & 0 \\
  0 & 0 & 1 & 0 & 0 \\
  0 & 0 & 0 & 0 & 0 \\
  0 & 0 & 0 & 0 & 0 \\
  0 & 0 & 0 & 0 & 0 \\
\end{array}%
\right) \mbox{ et } Z=\left(%
\begin{array}{ccccc}
  0 & b & c & \delta & e \\
  0 & 0 & b & 0 & 0 \\
  0 & 0 & 0 & 0 & 0 \\
  0 & 0 & f & 0 & 0 \\
  0 & 0 & i & 0 & 0 \\
\end{array}%
\right)$$ 
\vs 

Les champs $Z_{|z_2=0}$ et $Y_{|z_2=0}$ sont
colin\'eaires donc le feuilletage associ\'e \`a l'alg\`ebre
$\mathscr{L}$ ne peut \^etre un $\mathscr{L}$-feuilletage de
degr\'e $3$ sur $\C \P(4)$. En fait, il conduit \`a un
$\mathscr{L}$-feuilletage quadratique.\vs

\item[3.3.2.] Traitons le cas o\`u $$Y=\left(%
\begin{array}{ccccc}
  0 & 0 & 0 & 1 & 0 \\
  0 & 0 & 0 & 0 & 1 \\
  0 & 0 & 0 & 0 & 0 \\
  0 & 0 & 0 & 0 & 0 \\
  0 & 0 & 0 & 0 & 0 \\
\end{array}%
\right) \mbox{ et } Z=\left(%
\begin{array}{ccccc}
  0 & 0 & 0 & \delta & e \\
  0 & 0 & 0 & i & j \\
  0 & 0 & 0 & \ell & m \\
  0 & 0 & 0 & 0 & 0 \\
  0 & 0 & 0 & 0 & 0 \\
\end{array}%
\right)$$
\vs

Si $(\ell,m)=(0,0)$, alors le feuilletage associ\'e \`a
l'alg\`ebre $\mathscr{L}$ vient de $\C \P(2)$ donc
n'est pas de degr\'e $3$.

Soit $$P=\left(%
\begin{array}{ccccc}
  1 & \varepsilon & 0 & 0 & 0 \\
  0 & 1 & \varepsilon & 0 & 0 \\
  0 & 0 & 1 & \varepsilon & 0 \\
  0 & 0 & 0 & 1 & \varepsilon \\
  0 & 0 & 0 & 0 & 1 \\
\end{array}%
\right)$$ %alors $$P^{-1}=\left(%
%\begin{array}{ccccc}
%  1 & -\varepsilon & \varepsilon^2 & -\varepsilon^3 & \varepsilon^4 \\
%  0 & 1 & -\varepsilon & \varepsilon^2 & -\varepsilon^3 \\
%  0 & 0 & 1 & -\varepsilon & \varepsilon^2 \\
%  0 & 0 & 0 & 1 & -\varepsilon \\
%  0 & 0 & 0 & 0 & 1 \\
%\end{array}%
%\right)$$ et
alors $$PZP^{-1}=\left(%
\begin{array}{ccccc}
  0 & 0 & 0 & \delta+\varepsilon i & \varepsilon(j-\delta-\varepsilon i)+e \\
  0 & 0 & 0 & i+\varepsilon \ell & \varepsilon(m-\varepsilon \ell-i)+j \\
  0 & 0 & 0 & \ell & m-\varepsilon \ell \\
  0 & 0 & 0 & 0 & 0 \\
  0 & 0 & 0 & 0 & 0 \\
\end{array}%
\right)$$
\vs

En particulier si $\ell \not =0$, on peut, en posant
$\varepsilon=\frac{m}{\ell}$, se ramener \`a $m=0$.

Si $\ell=0$, alors $m\not =0$ et $i\not =0$ sinon le degr\'e chute,
et 
\vs

$$PZP^{-1}=\left(%
\begin{array}{ccccc}
  0 & 0 & 0 & \delta+\varepsilon i & \varepsilon(j-\delta-\varepsilon i)+e \\
  0 & 0 & 0 & i & \varepsilon(m-i)+j \\
  0 & 0 & 0 & 0 & m \\
  0 & 0 & 0 & 0 & 0 \\
  0 & 0 & 0 & 0 & 0 \\
\end{array}%
\right)$$ 
\vs

Quitte \`a permuter $z_3$ et $z_4$, %$$PZP^{-1}=\left(%
%\begin{array}{ccccc}
%  0 & 0 & 0 & \varepsilon(j-\delta-\varepsilon i)+e & \delta+\varepsilon i  \\
%  0 & 0 & 0 & \varepsilon(m-i)+j & i \\
%  0 & 0 & 0 & m & 0 \\
%  0 & 0 & 0 & 0 & 0 \\
%  0 & 0 & 0 & 0 & 0 \\
%\end{array}%
%\right) \mbox{ et } Y=\left(%
%\begin{array}{ccccc}
%  0 & 0 & 0 & 0 & 1 \\
%  0 & 0 & 0 & 1 & 0 \\
%  0 & 0 & 0 & 0 & 0 \\
%  0 & 0 & 0 & 0 & 0 \\
%  0 & 0 & 0 & 0 & 0 \\
%\end{array}%
%\right)$$
puis $z_0$ et $z_1$, on obtient
\vs

$$PZP^{-1}=\left(%
\begin{array}{ccccc}
  0 & 0 & 0 & \varepsilon(m-i)+j & i\\
  0 & 0 & 0 & \varepsilon(j-\delta-\varepsilon i)+e & \delta+\varepsilon i \\
  0 & 0 & 0 & m & 0 \\
  0 & 0 & 0 & 0 & 0 \\
  0 & 0 & 0 & 0 & 0 \\
\end{array}%
\right)$$ et $$Y=\left(%
\begin{array}{ccccc}
  0 & 0 & 0 & 1 & 0 \\
  0 & 0 & 0 & 0 & 1 \\
  0 & 0 & 0 & 0 & 0 \\
  0 & 0 & 0 & 0 & 0 \\
  0 & 0 & 0 & 0 & 0 \\
\end{array}%
\right)$$
\vs

 i.e. on s'est ramen\'e au cas $\ell \not=0$ et $m=0$.
\vs

Une fois $Y$ et $$Z=(\delta z_3+ez_4)\frac{\partial}{\partial z_0}+
(iz_3+jz_4)\frac{\partial}{\partial
z_1}+z_3\frac{\partial}{\partial z_2}$$ fix\'es, Maple donne :
$$X=(\xi_0z_0-2\mu ez_1+(2\mu(\delta-j)^2+2\xi_5 e+4i\mu e+(\xi_0-\xi_6)
(\delta-j))z_2+\xi_3z_3+\xi_4z_4)\frac{\partial}{\partial z_0}$$
$$+(\xi_5z_0+\xi_6z_1+(2i(\xi_0-\xi_6)+(\xi_5-2i\mu)(j-\delta)
)z_2+\xi_8z_3+\xi_9z_4)\frac{\partial}
{\partial z_1}$$ $$+(-\mu z_0+(2\xi_0+3\mu(\delta-j)-\xi_6)z_2+\xi_{13}z_3 +
\xi_{14}z_4)\frac{\partial}{\partial z_2}$$
$$+((\xi_0+\lambda+\mu\delta)z_3-\mu ez_4)\frac{\partial}{\partial
z_3}$$ $$+((\xi_5+i\mu)z_3+(\xi_6+\lambda+\mu j)z_4)
\frac{\partial}{\partial z_4},$$
\vs

Nous pr\'esentons le calcul de l'int\'egrale premi\`ere, la proc\'edure
d'int\'egration \'etant sensiblement diff\'erente de celles d\'ej\`a rencontr\'ees.

Quitte \`a retrancher $((\xi_0+\lambda+\mu\delta)z_3-\mu ez_4)R$ \`a $X$,
le champ $X$ n'a pas de composante en $\frac{\partial}{\partial z_3}$.
Ainsi les trois champs engendrent le feuilletage
restreint \`a la carte affine $z_3=1$ ; dans cette carte affine,
les champs $Y$ et $Z$ s'\'ecrivent
$$\tilde{Y}=Y_{|z_3=1}=\frac{\partial}{\partial z_0}+z_4
\frac{\partial}{\partial z_1}$$
$$\tilde{Z}=Z_{|z_3=1}=(\delta+ez_4)\frac{\partial}{\partial z_0}
+(i+jz_4)\frac{\partial} {\partial z_1}+\frac{\partial}{\partial z_2}$$
Redressons les champs $\tilde{Y}$ et $\tilde{Z}$ par un automorphisme
polynomial.

Le flot de $\tilde{Y}$ est
$$\varphi \hspace{1mm} \colon \hspace{1mm} (z;t) \mapsto
(z_0+t,z_1+ z_4t,z_2,z_4)$$ et celui de
$\tilde{Z}$
$$\psi \hspace{1mm} \colon \hspace{1mm} (z;s) \mapsto
(z_0+(\delta+ez_4)s,(i+jz_4)s+z_1,s+z_2,z_4)$$
On d\'efinit comme d'habitude le diff\'eomorphisme
\begin{eqnarray}
H(z_0,z_1,z_2,z_4) &=&
\varphi(\psi(0,z_1,0,z_4;z_2);z_0) \nonumber \\
&=&
(z_0+(\delta+ez_4)z_2,z_1+(i+jz_4)z_2+z_0z_4,z_2,z_4) \nonumber
\end{eqnarray}
%et $$H^{-1}=(z_0-(\delta+ez_4)z_2,z_1-z_0z_4+ez_2z_4^2-iz_2+(\delta-j)z_2z_4,
%z_2,z_4)$$
Par construction $H$ conjugue $\frac{\partial}{\partial z_0}$ \`a
$\tilde{Y}$ et $\frac{\partial}{\partial z_2}$ \`a $\tilde{Z}$. On
note $\tilde{\mathscr{F}}$ le feuilletage associ\'e \`a
l'alg\`ebre engendr\'ee par les champs $\tilde{X}=X_{|z_3=1}$,
$\tilde{Y}$, $\tilde{Z}$ et $R_{|z_3=1}$. On constate que le
champ $\tilde{X}_{|z_0=z_2=0}$, quitte \`a le modifier avec des
combinaisons ici polynomiales de $R$, $Y$ et $Z$, s'\'ecrit
$$\left((\xi_6-\xi_0-\lambda-\mu\delta)z_1+3\mu ez_1z_4+(\xi_9-i\xi_{14}-\xi_3
-\xi_{13}(j-\delta))z_4\right.$$ $$\left.+(\xi_{13}e+(\delta-j)\xi_{14}-\xi_4)
z_4^2+e\xi_{14}z_4^3+\xi_8-i\xi_{13}\right)\frac{\partial}{\partial z_1}$$ $$+
\left((\xi_5+\mu i)+(\xi_6-\xi_0+\mu(j-\delta))z_4+\mu ez_4^2\right)
\frac{\partial}{\partial z_4}$$ Il est tangent
\`a $\tilde{\mathscr{F}}$ et au $2$-plan $z_0=z_2=0$ ; donc il d\'ecrit
$\tilde{\mathscr{F}}$ en restriction \`a ce $2$-plan. Ce feuilletage
est aussi d\'efini par $$\omega=\left(\kappa_0+\kappa_1 z_4+\kappa_2
z_4^2\right)dz_1$$ $$+\left(\kappa_3z_1-3\kappa_2z_1z_4+\kappa_5z_4+\kappa_6
z_4^2+\kappa_7z_4^3+\kappa_8\right)dz_4$$ avec $$\kappa_0=\xi_5+\mu
i,\hspace{1mm}\kappa_1=\xi_6-\xi_0+\mu(j-\delta)$$ $$\kappa_2=\mu e,
\hspace{1mm}\kappa_3=-\xi_6+\xi_0+\lambda+\mu\delta,$$ $$\kappa_5=
-\xi_9+i\xi_{14}+\xi_3+\xi_{13}(j-\delta),\hspace{1mm}\kappa_6=(j-\delta)
\xi_{14}-\xi_{13}e+\xi_4,$$ $$\kappa_7=-e\xi_{14},\hspace{1mm}\kappa_8=i
\xi_{13}-\xi_8$$
On regarde la forme $\omega$ comme une \'equation diff\'erentielle lin\'eaire
avec second membre dont on cherche une solution particuli\`ere
$$z_1=a_0+a_1z_4+a_2z_4^2+a_3z_4^3$$ En d\'eveloppant on obtient
$$\left(\kappa_0a_1+\kappa_3a_0+\kappa_8\right)+\left(2\kappa_0a_2+(\kappa_1
+\kappa_3)a_1-3\kappa_2a_0+\kappa_5\right)z_4$$ $$+\left(3\kappa_0a_3+(2\kappa_1
+\kappa_3)a_2-2\kappa_2a_1+\kappa_6\right)z_4^2+\left((3\kappa_1+\kappa_3)a_3
-\kappa_2a_2+\kappa_7\right)z_4^3=0$$ 

\vs

La r\'esolution de cette \'equation
conduit \`a $$a_0=\frac{1}{f}\left(3\kappa_0^3\kappa_7+3\kappa_0^2\kappa_2\kappa_5
-9\kappa_0\kappa_1\kappa_2\kappa_8+5\kappa_0\kappa_1\kappa_3\kappa_5-6\kappa_1^3
\kappa_8-5
\kappa_0\kappa_2\kappa_3\kappa_8\right.$$ $$\left.-\kappa_0^2\kappa_3\kappa_6+6
\kappa_0\kappa_1^2\kappa_5-3\kappa_0^2
\kappa_1\kappa_6-11\kappa_1^2\kappa_3\kappa_8-6\kappa_1\kappa_3^2\kappa_8+
\kappa_0\kappa_3^2\kappa_5-\kappa_3^3\kappa_8\right),$$
$$a_1=-\frac{1}{f}\left(3\kappa_2\kappa_3^2\kappa_8+3\kappa_0\kappa_2\kappa_3
\kappa_5+6\kappa_1^2\kappa_3
\kappa_5+15\kappa_1\kappa_2\kappa_3\kappa_8+5\kappa_1\kappa_3^2
\kappa_5\right.$$ $$\left.-\kappa_0\kappa_3^2\kappa_6-3\kappa_0\kappa_1\kappa_3
\kappa_6+\kappa_3^3
\kappa_5+3\kappa_0^2\kappa_3\kappa_7+18\kappa_1^2\kappa_2\kappa_8+9
\kappa_0\kappa_2^2\kappa_8\right),$$ $$a_2=-\frac{1}{f}\left(6\kappa_2^2\kappa_3
\kappa_8+18\kappa_1\kappa_2^2\kappa_8+3\kappa_0
\kappa_2\kappa_3\kappa_6+9\kappa_0\kappa_1\kappa_2\kappa_6+2\kappa_2
\kappa_3^2\kappa_5+\kappa_3^3\kappa_6\right.$$
$$\left.-3\kappa_0\kappa_3^2\kappa_7
+4\kappa_1
\kappa_3^2\kappa_6+3\kappa_1^2\kappa_3\kappa_6-3\kappa_0\kappa_1\kappa_3
\kappa_7+6\kappa_1\kappa_2\kappa_3\kappa_5-9\kappa_0^2
\kappa_2\kappa_7\right),$$
$$a_3=-\frac{1}{f}\left(2\kappa_1^2\kappa_3\kappa_7+\kappa_1\kappa_2\kappa_3
\kappa_6+6\kappa_0\kappa_1\kappa_2\kappa_7+3\kappa_1\kappa_3^2\kappa_7+3\kappa_0
\kappa_2^2\kappa_6\right.$$ $$\left.+\kappa_2\kappa_3^2\kappa_6+
\kappa_3^3\kappa_7+5\kappa_0
\kappa_2\kappa_3\kappa_7+6\kappa_2^3\kappa_8+2\kappa_2^2\kappa_3\kappa_5
\right)$$ o\`u $$f=11\kappa_1^2\kappa_3^2+6\kappa_1\kappa_3^3+\kappa_3^4+8
\kappa_0\kappa_2\kappa_3^2$$ $$+6\kappa_1^3\kappa_3+24\kappa_0\kappa_1\kappa_2
\kappa_3+18\kappa_0\kappa_1^2\kappa_2+9\kappa_0^2\kappa_2^2$$
Notons $p(z)=a_0+a_1z+a_2z^2+a_3z^3$ et posons $z_1=Z_1+p(z_4)$ ; la
nouvelle forme s'\'ecrit alors $$\left(\kappa_0+\kappa_1z_4+\kappa_2
z_4^2\right)dZ_1+Z_1\left(\kappa_3-3\kappa_2z_4\right)dz_4$$ soit encore
$$\frac{dZ_1}{Z_1}+\frac{\kappa_3-3\kappa_2z_4}{\kappa_0+\kappa_1z_4+\kappa_4
z_2^2}dz_4$$ qui s'int\`egre explicitement par d\'ecomposition en \'el\'ements
simples ; notons que $\kappa_0+\kappa_1z_4+\kappa_4
z_2^2\not \equiv 0$. Il y a plusieurs cas suivant que le polyn\^ome $\kappa_0+
\kappa_1z_4+\kappa_4z_2^2 $ a des racines distinctes ou non.

Commen\c{c}ons par consid\'erer le cas o\`u le polyn\^ome $\kappa_0+\kappa_1z_4
+\kappa_4z_2^2 $ a deux racines distinctes $$\gamma^\pm=\frac{-\kappa_1 \pm
\sqrt{\kappa_1^2-4\kappa_0\kappa_2}}{2\kappa_2}$$ Le feuilletage est alors
d\'ecrit par la $1$-forme $$\frac{dZ_1}{Z_1}+\frac{1}{\kappa_2}\left(
\frac{b}{z_4-\gamma^-}+\frac{c}{z_4-\gamma^+}\right)dz_4$$ avec
$$c=\frac{\kappa_2(\kappa_3-3\kappa_2\gamma^+)}{\sqrt{\kappa_1^2-4\kappa_0
\kappa_2}},\hspace{2mm} b=-3\kappa_2-c$$ On en d\'eduit l'int\'egrale premi\`ere
%$$Z_1(z_4-\gamma^-)^b(z_4-\gamma^+)^c$$ soit $$(z_1-p(z_4))(z_4-\gamma^-)^b(z_4-
%\gamma^+)^c$$
$$\frac{{P(z_0,z_1,z_2,z_3,z_4)\left(z_4-\gamma^-z_3
\right)}^b{\left(z_4-\gamma^+z_3\right)}^c}{z_3^{3+b+c}}$$ o\`u
$P$ s'\'ecrit :
$$(z_1-iz_2-a_1z_4)z_3^2+((\delta-j)z_2-z_0)z_3z_4+
(ez_2-a_2z_3)z_4^2-a_0z_3^3-a_3z_4^3$$ soit \`a conjugaison pr\`es
%$$\frac{{\left(z_1z_3^2+z_0z_3z_4+ez_2z_4^2-a_3z_4^3\right)\left(z_4-\gamma^-z_3
%\right)}^b{\left(z_4-\gamma^+z_3\right)}^c}{z_3^{3+b+c}} $$
$$\frac{{\left(z_1z_3^2+z_0z_3z_4+ez_2z_4^2-az_4^3\right)\left(z_4-\varrho^-z_3
\right)}^b{\left(z_4-\varrho^+z_3\right)}^c}{z_3^{3+b+c}} $$ avec $a$,
$e \in \{0,1\}$ et les $\varrho^+$ et $\varrho^-$ se calculent en 
fonction des diff\'erents param\`etres.\vs

Finalement consid\'erons le cas o\`u le polyn\^ome $\kappa_0+\kappa_1z_4
+\kappa_4z_2^2 $ a une racine double $\gamma=-\frac{\kappa_1}{2\kappa_2}$
Le feuilletage est alors d\'ecrit par la $1$-forme $$\frac{dZ_1}{Z_1}+
\eta\frac{dz_4}{(z_4-\gamma)^2}-3\frac{dz_4}{z_4-\gamma}$$ avec $$\eta=
\frac{\kappa_3-3\kappa_2\gamma}{\kappa_2}$$ On obtient
comme int\'egrale premi\`ere %$$Z_1(z_4-\gamma)^B\exp\left(
%-\frac{A}{z_4-\gamma}\right)$$ soit $$(z_1-a_0-a_1z_4-a_2z_4^2-a_3z_4^3)
%(z_4-\gamma)^B\exp\left(-\frac{A}{z_4-\gamma}\right)$$ et finalement
$$\frac{P(z_0,z_1,z_2,z_3,z_4)}{{\left(z_4-\gamma
z_3\right)}^3}\exp\left(\frac{\eta z_4}{\gamma z_3-z_4}\right) $$
ou encore \`a conjugaison pr\`es $$\frac{\left(z_1z_3^2+z_0z_3z_4+ez_2z_4^2-
az_4^3\right)}{{\left(z_4-z_3\right)}^3}\exp\left(\frac{\upsilon z_4}{
z_3-z_4}\right)$$ avec $a$, $e \in \{0,1\}$.
\vs
%Notons $\zeta=\xi_6-\xi_0$ ; les valeurs propres de la matrice $$\left(%
%\begin{array}{cc}
%  \lambda & \mu \\
%  \mu(j\delta-ie-2j^2)-\xi_5 e-j\zeta &
%  \zeta+\lambda+\mu(3j-\delta) \\
%\end{array}%
%\right)$$ sont $\rho^\pm=\lambda-\frac{1}{2}(\mu\delta-\zeta-3\mu
%j)\pm\frac{1}{2}\sqrt{(2\mu (j-\delta)+\zeta)\zeta -4\mu
%e(\xi_5-i\mu)+\mu^2(j-\delta)^2}$ et
%$\tilde{\alpha}=\tilde{\frac{\rho^+}{\rho^-}}$.
%\vs

\item[3.3.3.] Pour finir traitons le cas o\`u $$Y=\left(%
\begin{array}{ccccc}
  0 & 0 & 0 & 1 & 0 \\
  0 & 0 & 0 & 0 & 1 \\
  0 & 0 & 0 & 0 & 0 \\
  0 & 0 & 0 & 0 & 0 \\
  0 & 0 & 0 & 0 & 0 \\
\end{array}%
\right) \mbox{ et } Z=\left(%
\begin{array}{ccccc}
  a & b & c & \delta & e \\
  f & -a & h & i & j \\
  0 & 0 & 0 & \ell & m \\
  0 & 0 & 0 & a & b \\
  0 & 0 & 0 & f & -a \\
\end{array}%
\right)$$ avec  $\left\{\begin{array}{ccccc}
  a\ell=-fm  \hspace{25mm} \\
  \delta m\ell = -m^2i+mj\ell+e\ell^2  \\
  c\ell=-hm  \hspace{25mm}\\
  b\ell^2=-fm^2   \hspace{22mm}\\
\end{array}\right.$
\vs
\vs

Il y a deux cas.
\vs

3.3.3.1. Maple donne l'alg\`ebre $\mathscr{L}$ engendr\'ee par les
champs
$$X=(-3\varepsilon z_0+\xi_1z_1+\xi_2z_2+\frac{(\delta\xi_1-\xi_2m)
\ell^2+2m\varepsilon(\delta \ell+mi-\ell j)}{b\ell^2}z_3+
\xi_4z_4)\frac{\partial}{\partial z_0}$$ $$+ (-2\varepsilon
z_1+\frac{2\xi_1i+\xi_2\ell-\varepsilon \delta}{2b}z_3)
\frac{\partial}{\partial z_1}$$ $$+ (-2\varepsilon z_2+\frac{\ell
\xi_1+m\varepsilon}{b}z_3+ \xi_{14}z_4)\frac{\partial} {\partial
z_2}+(-\varepsilon z_3+\xi_1z_4)\frac{\partial}{\partial z_3}$$
$$Z=(bz_1+\delta z_3+\frac{m(\delta\ell-mi-\ell j)} {\ell^2}z_4)
\frac{\partial}{\partial z_0}+iz_3\frac{\partial} {\partial z_1}$$
$$+(\ell z_3+mz_4) \frac{\partial}{\partial z_2}+ bz_4
\frac{\partial}{\partial z_3},$$ et $$\left(%
\begin{array}{cc}
2\varepsilon & 0 \\
\frac{1}{2}(\varepsilon\delta-\xi_2\ell) & \varepsilon \\
\end{array}%
\right) $$
\vs
Les valeurs propres de cette matrice sont $2\varepsilon$ et $\varepsilon$ donc
$\tilde{\alpha}=\tilde{2}$ ; on remarque que $\varepsilon \not =0$.
\vs

Le feuilletage associ\'e \`a cette alg\`ebre admet pour
int\'egrale premi\`ere
$$\frac{{\left(z_2z_4-\frac{m}{b}z_3z_4-\frac{\ell}{2b}z_3^2+\kappa z_4^2\right)
}^3}{{\left(\xi_2z_2z_4^2-\varepsilon z_0z_4^2+\varepsilon z_1z_3z_4+\varsigma
z_3z_4^2+(\frac{\delta\varepsilon-\ell\xi_2}{2b})z_3^2z_4-\frac{i\varepsilon}
{3b}z_3^3+(\xi_2\kappa-\varepsilon\nu)z_4^3\right)}^2}$$

\vs

avec $$\kappa=-\frac{\xi_{14}}{2\varepsilon}+\frac{\xi_1m}{2\varepsilon b},
\hspace{2mm}\varsigma=\frac{m\varepsilon(\delta\ell-mi-j\ell)
-m\ell^2\xi_2}{\ell^2b}$$
$$\nu=-\frac{1}{3\varepsilon}\left(\frac{\xi_2\xi_{14}}{2\varepsilon}-
\frac{\xi_1\xi_2 m}{2b\varepsilon}+\xi_4-\frac{\xi_1m(\delta\ell-mi-\ell j)}{b
\ell^2}\right)$$

\vs

soit \`a conjugaison pr\`es $$\frac{{\left(z_2z_4-\zeta z_3^2\right)}^3}{{\left(z_0z_4^2+z_1z_3z_4-\upsilon
z_3^3\right)}^2}$$ avec $\upsilon$, $\zeta \in \{0,1\}$.\vs

3.3.3.2. L'alg\`ebre $\mathscr{L}$ est engendr\'ee par les
champs
$$X=(-4\varepsilon z_0+\xi_1z_1+\xi_2z_2+\frac{\delta\xi_1\ell-hm
\xi_{14}+3\varepsilon m(\delta-j)-\xi_2m\ell}{b\ell}z_3+\xi_4z_4)
\frac{\partial}{\partial z_0}$$$$+(-3\varepsilon z_0+\frac{h(\ell
\xi_1+m\varepsilon}{b\ell}z_2+\frac{\xi_2\ell+h\xi_{14}-2\varepsilon
\delta}{2b}z_3)\frac{\partial}{\partial z_1}$$$$+ (-2\varepsilon
z_2+\frac{\ell \xi_1+m\varepsilon}{b}z_3+\xi_{14}z_4)
\frac{\partial}{\partial z_2}+(-\varepsilon z_3+\xi_1z_4)
\frac{\partial}{\partial z_3}$$ $$Z=(bz_1-\frac{hm}{\ell}z_2+
\delta z_3+ \frac{m(\delta-j)}{\ell} z_4)\frac{\partial}{\partial
z_0}+hz_2\frac{\partial}{\partial z_1}$$$$+(\ell z_3+mz_4)
\frac{\partial}{\partial z_2}+bz_4 \frac{\partial}{\partial
z_3}$$ et $$\left(%
\begin{array}{cc}
3\varepsilon & 0 \\
\frac{(b\ell(h\xi_{14}-\xi_2\ell))+2b\ell \varepsilon\delta-2hm
(\varepsilon m+\ell\xi_1)}{2b\ell} & \varepsilon \\
\end{array}%
\right) $$
\vs

On note que $\tilde{\alpha}=\tilde{3}$, que $\varepsilon \not =0$ et que
$b\ell\not =0$ (sinon $\tilde{\alpha}=\tilde{0}$).\vs

Le feuilletage associ\'e \`a cette alg\`ebre admet pour
int\'egrale premi\`ere $$\frac{{\left(z_2z_4-\ell P+\eta z_4^2\right)}^2}{2
\varepsilon \left(\zeta z_3z_4^3-z_0z_4^3-hz_2z_4P+z_1z_3z_4^2+2\tau
\rho^2\ell z_3^2z_4^2+\vartheta z_3^3z_4+\frac{h\ell}{8b^2}z_3^4
\right)+\kappa z_4^2\left(z_2z_4-\ell P+\eta z_4^2\right)}$$ \vs
avec $$\nu=\frac{1}{4\varepsilon}(-\frac{
1}{2\varepsilon}(\xi_2+\frac{\xi_1hm}{b\ell})(\frac{\xi_1m}{b}-
\xi_{14})+\xi_4-\frac{\xi_1m(\delta-j)}{b\ell}),$$
$$\kappa=\xi_2+\frac{\xi_1hm}{b\ell}, \hspace{2mm} \eta=\frac{1}
{2\varepsilon}(\frac{\xi_1m}{b}-\xi_{14}), \hspace{2mm}
\varsigma=\frac{\ell\varepsilon h}{2b}$$
$$\tau=m^2h+b\ell\delta, \hspace{2mm} \rho=\frac{1}{2b\ell}, \hspace{2mm}
 \vartheta=\frac{mh}{2b^2},$$
$$\zeta=2m\rho(\delta-j) \mbox{
et } P=\rho(2mz_4+\ell z_3)z_3$$\vs
Ce qui conduit \`a conjugaison pr\`es \`a %$$\frac{{\left(z_2z_4-\rho\ell^2z_3^2
%\right)}^2}{z_0z_4^3+z_1z_3z_4^2-h\rho\ell z_2z_3^2z_4+\vartheta z_3^3z_4
%+\frac{h\ell}{8b^2}z_3^4}$$
$$\frac{{\left(z_2z_4+\Upsilon_1 z_3^2
\right)}^2}{z_0z_4^3+z_1z_3z_4^2+\Upsilon_2 z_2z_3^2z_4+\upsilon z_3^3z_4
+\Upsilon_3 z_3^4}$$ o\`u $\Upsilon_i \in \{0,1\}$.
\vs
\vs
\end{itemize}

\begin{thm} (Maple) Soit $\mathscr{F}$ un
$\mathscr{L}$-feuilletage de degr\'e $3$ sur
$\C\P(4)$. Si $\mathscr{L}$ admet pour
pr\'esentation $\{[X,Y]=Y,\hspace{1mm} [X,Z]=\alpha Z,\hspace{1mm}
[Y,Z]=0\}$, $\alpha \not =0$, alors le feuilletage $\mathscr{F}$ 
poss\`ede un int\'egrale premi\`ere de l'un des types :
\begin{eqnarray}
& & \frac{{\left(z_0z_4^3-\frac{z_2^2z_4^2}{2}-z_3z_4^2
+z_2z_3^2z_4-\frac{z_3^4}{4}\right)}^3}{{\left(z_1z_4^2
-z_2z_3z_4+\frac{z_3^3}{3}\right)}^4}\nonumber \\
& & \nonumber \\
& & \frac{{\left(z_1z_4^2-z_2z_3z_4+\frac{z_3^3}{3}\right)}^4}
{{\left(z_0z_4^3-(2z_1z_3+z_2^2)z_4^2-\frac{z_3^4}{2}\right)}^3}
\nonumber \\
& & \nonumber \\
& & \frac{z_4(z_0z_4^2+z_1z_3z_4+z_2^2z_4+az_3^3)}{z_3^4}\nonumber \\
& & \nonumber \\
& & \frac{z_4^{3\kappa-\lambda}z_3^{\lambda}}
{{\left(z_1z_4^2+z_2z_3z_4+z_2^2z_3+z_0^2z_4\right)}^{\kappa}}  \nonumber \\
& & \nonumber \\
& & \left(\frac{z_3}{z_4}\right)^{\kappa} \exp\left(\frac{
z_1z_4^2+z_0z_3z_4+z_2^2z_3}{z_4^3}\right)\nonumber  \\
& & \nonumber \\
& & \frac{z_3^{\kappa+\lambda}
z_4^{2\lambda-\kappa}}{\left(z_1z_4^2+z_0z_3z_4+z_2^2z_3\right)^\lambda}
\nonumber \\
& & \nonumber \\
& & \left(\frac{z_3}{z_4}\right)^{\kappa}\exp\left(\frac{
z_1z_4^2+z_0z_3z_4+z_2^2z_3}{z_3^2z_4}\right)
 \nonumber \\
& & \nonumber \\
& &
\frac{{\left((z_0-\gamma^+z_1)z_4+(\gamma^+c-d)z_3^2
-\frac{z_2^2}{2}+\gamma^+z_2z_3\right)}^{\vartheta_+}{\left((z_0-\gamma^-z_1)z_4+
(\gamma^-c-d)z_3^2-\frac{z_2^2}{2}+\gamma^-z_2z_3\right)}^{
\vartheta_-}}{z_4^2} \nonumber \\
& & \nonumber \\
& & \frac{(z_0-z_1)z_4+c
z_2z_3-\frac{z_2^2}{2}+\upsilon z_3^2}{z_4^2}
\exp\left(\frac{\zeta\left(z_0z_4-2g\delta z_3^2-\frac{z_2^2}{2}\right)}
{(z_0-z_1)z_4+cz_2z_3-\frac{z_2^2}{2}+\upsilon z_3^2}\right)\nonumber \\
& & \nonumber \\
& & z_3^{\lambda_3}z_4^{\lambda_4}Q^{\lambda_2} \mbox{, ou } Q
\mbox{ est une certaine cubique, et ses d\'eg\'en\'erescences (cf } 3.2.5.) \nonumber \\
& & \nonumber \\
& & \frac{z_4^\kappa{\left(z_0z_3^2+
z_1z_2z_3+z_2^3\right)}^{\lambda-\kappa}}{z_3^{3\lambda-2
\kappa}} \nonumber \\
& & \nonumber \\
& &\left(\frac{z_4}{z_3}\right)^{\kappa}\exp\left(\frac{z_0z_3^2
+z_1z_2z_3+z_2^3}{z_3^2z_4}\right)
\nonumber %\\
%& & \nonumber \\
\end{eqnarray}

\begin{eqnarray}
& &\frac{z_0z_3^2+z_1z_2z_3+z_2^3}{z_3^3}\exp
\left(\frac{z_4}{z_3}\right) \nonumber \\
& & \nonumber \\
& & \frac{\left(z_1z_3^2+z_0z_3z_4+ez_2z_4^2-
az_4^3\right)}{{\left(z_4-z_3\right)}^3}\exp\left(\frac{\upsilon z_4}{
z_3-z_4}\right) \nonumber\\
& & \nonumber \\
& & \frac{{\left(z_1z_3^2+z_0z_3z_4+ez_2z_4^2-az_4^3\right)\left(z_4-\varrho^-z_3
\right)}^\kappa{\left(z_4-\varrho^+z_3\right)}^\lambda}{z_3^{3+\kappa
+\lambda}}\nonumber \\
& & \nonumber\\
& & \frac{{\left(z_2z_4-az_3^2\right)}^3}{{\left(z_0z_4^2+
z_1z_3z_4-ez_3^3\right)}^2}\mbox{ o\`u } \xi, \upsilon \in \{0,1\}\nonumber\\
& & \nonumber\\
& & \frac{{\left(z_2z_4+\Upsilon_1 z_3^2
\right)}^2}{z_0z_4^3+z_1z_3z_4^2+\Upsilon_2 z_2z_3^2z_4+\upsilon z_3^3z_4
+\Upsilon_3 z_3^4}\nonumber
\end{eqnarray}

o\`u $a$, $e$, $\Upsilon_i \in  \{0,1\}$ et les autres param\`etres
sont des nombres complexes.
\end{thm}

\vs

\section[{$\mathscr{L}$ a pour pr\'esentation $\{[X,Y]=[X,Z]=0,
\hspace{1mm} [Y,Z]=Y\}$.}]{L'alg\`ebre $\mathscr{L}$ a pour
pr\'esentation $\{[X,Y]=[X,Z]=0, \hspace{1mm} [Y,Z]=Y\}$.}

Comme pr\'ec\'edemment on d\'efinit les combinaisons lin\'eaires
permises par les combinaisons lin\'eaires de $X$, $Y$ et $Z$ avec
$X$, $Y$, $Z$ et $R$ qui respectent la structure de l'alg\`ebre
$\mathscr{L}$. \vs

La r\'esolubilit\'e de $\mathscr{L}$ et l'\'egalit\'e $[Y,Z]=Y$
assurent que $Y$ est nilpotent. Les matrices $X$ et $Y$ commutent
d'o\`u l'existence d'une base dans laquelle on peut jordaniser $X$
et $Y$ ; on obtient modulo $R$ les configurations suivantes : \vs
$$X=\left(%
\begin{array}{ccccc}
  \nu_0 & \beta & 0 & 0 & 0 \\
  0 & \nu_0 & 0 & 0 & 0 \\
  0 & 0 & 0 & \delta & 0 \\
  0 & 0 & 0 & 0 & 0 \\
  0 & 0 & 0 & 0 & \nu_2 \\
\end{array}%
\right) \mbox{ et } Y=\left(%
\begin{array}{ccccc}
  0 & \varepsilon & 0 & 0 & 0 \\
  0 & 0 & 0 & 0 & 0 \\
  0 & 0 & 0 & \eta & 0 \\
  0 & 0 & 0 & 0 & 0 \\
  0 & 0 & 0 & 0 & 0 \\
\end{array}%
\right)$$

$$X=\left(%
\begin{array}{ccccc}
  \nu_0 & \beta & 0 & 0 & 0 \\
  0 & \nu_0 & 0 & 0 & 0 \\
  0 & 0 & 0 & \delta & \vartheta \\
  0 & 0 & 0 & 0 & \delta \\
  0 & 0 & 0 & 0 & 0 \\
\end{array}%
\right) \mbox{ et } Y=\left(%
\begin{array}{ccccc}
  0 & \varepsilon & 0 & 0 & 0 \\
  0 & 0 & 0 & 0 & 0 \\
  0 & 0 & 0 & \phi & \eta \\
  0 & 0 & 0 & 0 & \phi \\
  0 & 0 & 0 & 0 & 0 \\
\end{array}%
\right)$$

$$X=\left(%
\begin{array}{ccccc}
  \nu_0 & 0 & 0 & 0 & 0 \\
  0 & \nu_1 & 0 & 0 & 0 \\
  0 & 0 & 0 & \varepsilon & \delta \\
  0 & 0 & 0 & 0 & \varepsilon \\
  0 & 0 & 0 & 0 & 0 \\
\end{array}%
\right) \mbox{ et } Y=\left(%
\begin{array}{ccccc}
  0 & 0 & 0 & 0 & 0 \\
  0 & 0 & 0 & 0 & 0 \\
  0 & 0 & 0 & \vartheta & \beta \\
  0 & 0 & 0 & 0 & \vartheta \\
  0 & 0 & 0 & 0 & 0 \\
\end{array}%
\right)$$

$$X=\left(%
\begin{array}{ccccc}
  \nu & 0 & 0 & 0 & 0 \\
  0 & 0 & \eta & \vartheta & \phi \\
  0 & 0 & 0 & \eta & \vartheta \\
  0 & 0 & 0 & 0 & \eta \\
  0 & 0 & 0 & 0 & 0 \\
\end{array}%
\right) \mbox{ et } Y=\left(%
\begin{array}{ccccc}
  0 & 0 & 0 & 0 & 0 \\
  0 & 0 & \beta & \delta & \varepsilon \\
  0 & 0 & 0 & \beta & \delta \\
  0 & 0 & 0 & 0 & \beta \\
  0 & 0 & 0 & 0 & 0 \\
\end{array}%
\right)$$

$$X=\left(%
\begin{array}{ccccc}
  0 & \eta & \vartheta & \phi & \nu \\
  0 & 0 & \eta & \vartheta & \phi \\
  0 & 0 & 0 & \eta & \vartheta \\
  0 & 0 & 0 & 0 & \eta \\
  0 & 0 & 0 & 0 & 0 \\
\end{array}%
\right) \mbox{ et } Y=\left(%
\begin{array}{ccccc}
  0 & \lambda & \beta & \delta & \varepsilon \\
  0 & 0 & \lambda & \beta & \delta \\
  0 & 0 & 0 & \lambda & \beta \\
  0 & 0 & 0 & 0 & \lambda \\
  0 & 0 & 0 & 0 & 0 \\
\end{array}%
\right)$$
\vs

On cherche alors $Z$ tel que $[X,Z]=0$ et $[Y,Z]=Y$, on obtient
les sept configurations suivantes que l'on traite au cas par cas.
Certaines configurations s'av\`erent redondantes comme nous le
verrons ; comme cela ne se constate qu'\`a posteriori nous les
avons laiss\'es. \vs

\begin{itemize}
\item[4.1.] D'apr\`es Maple, les champs $X$, $Y$ et $Z$ s'\'ecrivent
$$X=z_0\frac{\partial}{\partial z_0}+z_1\frac{\partial} {\partial
z_1},\hspace{3mm} Y=\varepsilon z_1 \frac{\partial}{\partial
z_0}+z_3 \frac{\partial}{\partial z_2}, \hspace{2mm} \varepsilon
\not = 0$$
$$Z=(\xi_0 z_0+\xi_1 z_1)\frac{\partial}{\partial z_0}+
(\xi_0-1)z_1 \frac{\partial}{\partial
z_1}+\xi_2z_3\frac{\partial}{\partial z_3}$$
$$+((\xi_2+1)z_2+\xi_3z_3+\xi_4z_4)
\frac{\partial}{\partial z_2}+(\xi_5z_3+\xi_6z_4)
\frac{\partial}{\partial z_4}$$ Le feuilletage $\mathscr{F}$
associ\'e \`a cette alg\`ebre est d\'ecrit dans la carte affine
$z_3=1$ par la $1$-forme $\omega$ donn\'ee par
$$((\xi_2-\xi_6)z_4-\xi_5)(z_1dz_0-z_0dz_1)+\varepsilon(\xi_5+(\xi_6-\xi_2)
z_4)z_1^2dz_2$$ $$+(z_0+(\xi_1-\varepsilon \xi_3)z_1-
\varepsilon z_1z_2-\varepsilon \xi_4z_1z_4)z_1dz_4$$ Supposons
que $\lambda=\xi_6-\xi_2 \not =0$ ; on peut r\'esorber le
terme $\xi_5$ et, comme $\varepsilon \not =0$, par
une homoth\'etie on se ram\`ene \`a
$$\lambda(z_0dz_1-z_1dz_0)z_4+\lambda z_1^2z_4dz_2+ (z_0+\xi'_1
z_1-z_1z_2-\xi_4 z_1z_4)dz_4$$ On proc\`ede ici par \'eclatement
Posons $z_0=tz_1$, $\xi'_1=\xi_1
-\xi_3-\frac{\xi_4\xi_5}{\lambda}$ puis $T=t-z_2+\xi'_1$ ;
on obtient
$$z_1^2\left((T-\xi_4z_4)dz_4- \lambda z_4 dT\right)$$ qui a pour
facteur int\'egrant
$$z_1^2z_4\left((1-\lambda)T-
\xi_4z_4\right)$$
Si $\lambda \not =1$, le feuilletage admet pour int\'egrale
premi\`ere
$$\frac{z_1^\lambda z_3^{\lambda-1}z_4}{((1-\lambda)
(z_0z_3-z_1z_2+\xi'_1z_1z_3)-\xi_4z_1z_4)^\lambda}$$ et
$$\frac{z_4}{z_3}
\exp\left(\frac{z_0z_3-z_1z_2+\xi'_1z_1z_3}{z_1z_4}\right)$$ sinon. \vs

Ce qui donne \`a conjugaison pr\`es $$\frac{z_1^\lambda z_3^{\lambda-1}z_4}
{{\left(z_0z_3+z_1z_2\right)}^\lambda}$$ si $\lambda\not =1$ et $$\frac{z_4}{z_3
}\exp\left(\frac{z_0z_3+z_1z_2}{z_1z_4}\right) $$ sinon.\vs

Supposons que $\lambda=0$ ; on remarque que si $\xi_5=0$, le
degr\'e du feuilletage chute. La $1$-forme $\omega$ s'\'ecrit alors
$$-\xi_5(z_1dz_0-z_0dz_1)+\varepsilon\xi_5z_1^2dz_2$$ $$+
\left(z_0+(\xi_1
- \varepsilon\xi_3)z_1-\varepsilon z_1z_2-\varepsilon\xi_4z_1
z_4\right)z_1dz_4$$ Si $\xi_4 \not = 0$, posons $t=\frac{z_0}
{z_1}$, $T=-t+\varepsilon
z_2$ et $Z_4=z_4-\frac{\xi_1-\varepsilon \xi_3}{\varepsilon
\xi_4}$; quitte \`a diviser par $z_1^2$, on a
$$\omega=\xi_5dT-(T+\varepsilon\xi_4Z_4)dZ_4$$ i.e. $$\frac{dT}{dZ_4}
=\frac{T+\varepsilon\xi_4Z_4}{\xi_5}$$ Ainsi $T=C\exp z_4-
\varepsilon \xi_4(Z_4+\xi_5)$ et $$C=(T+\varepsilon\xi_4
(Z_4+\xi_5))\exp(-z_4)$$
est une int\'egrale premi\`ere du feuilletage. Finalement on obtient $$
\left(\frac{\varepsilon z_1z_2-z_0z_3+\varepsilon
\xi_4z_1z_4-\varepsilon \xi_4(\xi_5-\nu)z_1z_3}{z_1z_3}
\right)\exp\left(\frac{-z_4+\nu z_3}{z_3}\right)$$ o\`u
$$\nu=\frac{\xi_1-
\varepsilon\xi_3}{\varepsilon\xi_4}$$ soit \`a conjugaison pr\`es $$\left(
\frac{z_1z_2+z_0z_3}{z_1z_3}\right)\exp\left(\frac{z_4}{z_3}\right)$$
\vs

Pour $\xi_4=0$, on trouve l'int\'egrale premi\`ere
$$\frac{z_0z_3-\varepsilon z_1z_2-(\xi_1-\varepsilon\xi_3)z_1
z_3}{z_1z_3}\exp\left(-\frac{z_4}{\xi_5z_3}\right)$$ soit \`a conjugaison
pr\`es $$\frac{z_0z_3+z_1z_2}{z_1z_3}\exp\left(\frac{z_4}{z_3}\right) $$
On constate que les int\'egrales premi\`eres trouv\'ees dans les cas
$\xi_4\not =0$ et $\xi_4=0$ sont conjugu\'ees.
\vs
\vs

\item[4.2.] Le deuxi\`eme cas donn\'e par Maple est le suivant
$$X=(z_0+\beta z_1)\frac{\partial}{\partial z_0}+z_1
\frac{\partial}{\partial z_1}, \hspace{3mm}
Y=(z_3+(\xi_9-\xi_3)z_4) \frac{\partial}{\partial z_2}+z_4
\frac{\partial}{\partial z_3}$$
$$Z=\xi_1z_1\frac{\partial}{\partial z_0}+(z_2+\xi_3 z_3+
 \xi_4z_4)\frac{\partial}{\partial z_2}+\xi_3z_4 \frac{\partial}
{\partial z_3}-z_4\frac{\partial}{\partial z_4}$$
Posons $$Z_2=\frac{2z_2-z_3^2+2(\xi_3-
\xi_9)z_3+\xi_9
\xi_3+\xi_4 -\xi_9^2}{2}$$ alors la $1$-forme annulant
les champs $X$, $Y$, $Z$ et $R$ s'\'ecrit $$\frac{z_1dz_0-
(z_0+\beta z_1)dz_1}{z_1^2}- \frac{\xi_1}{2}\frac{dZ_2}{Z_2}$$ ce
qui conduit \`a l'int\'egrale premi\`ere $$\frac{z_1^{2\beta}
{\left(z_2z_4-\frac{z_3^2}{2}+(\xi_3-\xi_9)z_3z_4+(\xi_3
(\xi_3-\xi_9)+\xi_4)z_4^2\right)}^{\xi_1}}
{z_4^{2(\beta+\xi_1)}}\exp\left( \frac{2z_0}{z_1}\right)$$ soit
\`a conjugaison pr\`es $$\frac{z_1^{2
\beta}{\left(z_2z_4-\frac{z_3^2}{2}\right)}^{\xi_1}}
{z_4^{2(\beta+\xi_1)}}\exp\left( \frac{z_0}{z_1}\right)$$
\vs
\vs

\item[4.3.] Comme troisi\`eme cas on a
%$$X=z_0\frac{\partial}{\partial z_0}+z_1 \frac{\partial} {\partial
%z_1}, \hspace{3mm}
% Y=(z_3+(\xi_{19}-\xi_3)z_4)\frac{\partial}{\partial z_2}+z_4
%\frac{\partial}{\partial z_3}, \hspace{3mm} R \mbox{ et }$$
%$$Z=(\xi_0z_0+\xi_1z_1)\frac{\partial}{\partial
%z_0}+(\xi_5z_0+\xi_6z_1) \frac{\partial}{\partial
%z_1}+(2z_2+\xi_3z_3+\xi_{14}z_4)\frac{\partial}{\partial
%z_2}+(z_3+\xi_{19}z_4)\frac{\partial} {\partial z_3}$$
%La $1$-forme $\omega$ annulant ces champs
%s'\'ecrit dans la carte
%affine $z_4=1$ $$(2z_2-(2z_3+\xi_{19})(\xi_{19}-\xi_3)
%-z_3^2+\xi_{14}) (z_1dz_0-z_0dz_1)$$
%$$-((\xi_0-\xi_6)z_0z_1+\xi_1z_1^2 -\xi_5 z_0^2) (dz_2
%-(z_3+(\xi_{19} -\xi_3))dz_3)$$ Posons $Z_2=\frac{1}{2}(
%2z_2-2(\xi_{19}-\xi_3)z_3-z_3^2+\xi_{14} -\xi_{19}(\xi_{19}
%-\xi_3))$ alors $\omega$ s'\'ecrit
%$$\frac{z_1dz_0-z_0dz_1} {(\xi_0-\xi_6)z_0z_1+\xi_1z_1^2
%-\xi_5 z_0^2}+\frac{dZ_2}{Z_2}$$ c'est-\`a-dire \`a conjugaison
%pr\`es
%$$d(\frac{z_0}{z_1})+\frac{dZ_2}{Z_2} \mbox{ ou } \lambda(\frac{dz_0}{z_0}-
%\frac{dz_1}{z_1})+\frac{dZ_2}{Z_2}$$ Le feuilletage d\'ecrit par
%$\omega$ admet donc pour int\'egrale premi\`ere
%$$\frac{1}{2}\left( 2z_2z_4-2(\xi_{19}-\xi_3)z_3z_4-z_3^2+(\xi_{14}
%-\xi_{19}(\xi_{19} -\xi_3))z_4^2\right)\exp\left(\frac{z_0}
%{z_1}\right)
%$$ ou
%$$ \frac{z_0^\lambda\left(2z_2z_4-2(\xi_{19}-\xi_3)z_3z_4-z_3^2+(\xi_{14}
%-\xi_{19}(\xi_{19} -\xi_3))z_4^2\right)}{2z_1^\lambda z_4^2}$$

$$X=z_0\frac{\partial}{\partial z_0}+z_1 \frac{\partial} {\partial
z_1}, \hspace{3mm}
 Y=(z_3-\xi_3z_4)\frac{\partial}{\partial z_2}+z_4
\frac{\partial}{\partial z_3}, \hspace{3mm} R \mbox{ et }$$
$$Z=(\xi_0z_0+\xi_6z_1) \frac{\partial}{\partial z_1}+(2z_2+
\xi_3z_3)\frac{\partial}{\partial z_2}+z_3\frac{\partial}
{\partial z_3}$$

ou bien si $\xi_0=\xi_6$ $$X=z_0\frac{\partial}{\partial z_0}+z_1
\frac{\partial}{\partial z_1}, \hspace{3mm}
 Y=(z_3-\xi_3z_4)\frac{\partial}{\partial z_2}+z_4
\frac{\partial}{\partial z_3}, \hspace{3mm} R \mbox{ et }$$
$$Z=(\xi_0z_0+z_1) \frac{\partial}{\partial
z_0}+\xi_0z_1\frac{\partial}{\partial z_1}+(2z_2+
\xi_3z_3)\frac{\partial}{\partial z_2}+z_3\frac{\partial}
{\partial z_3}$$

\vs

Dans le premier cas, apr\`es \'etude dans la carte affine $z_4=1$, on
trouve l'int\'egrale premi\`ere
$$\left(\frac{z_0}{z_1}\right)^{\frac{2}{\xi_6-\xi_0}} \left(
\frac{z_2z_4-\frac{z_3^2}{2}+\xi_3z_3z_4}{z_4^2}\right)$$
qui est conjugu\'ee \`a
$$\left(\frac{z_0}{z_1}\right)^{\frac{2}{\xi_6-\xi_0}} \left(
\frac{z_2z_4-z_3^2}{z_4^2}\right)$$
et dans le second cas
$$\frac{z_2z_4-\frac{z_3^2}{2}+\xi_3z_3z_4}{z_4^2}
\exp\left(-\frac{2z_0}{z_1}\right)$$
qui se ram\`ene \`a
$$\frac{z_2z_4-z_3^2}{z_4^2}\exp\left(\frac{z_0}{z_1}\right)$$
C'est un cas particulier du 4.2 ($\beta=0$, $\xi_1=1$)\vs

\item[4.4.] Les champs $X$, $Y$ et $Z$ s'\'ecrivent
dans le quatri\`eme cas :
$$X=z_0\frac{\partial}{\partial z_0}+z_1\frac{\partial} {\partial
z_1}, \hspace{3mm}
Y=z_1\frac{\partial}{\partial z_0}+z_4 \frac{\partial}{\partial
z_2}$$
$$Z=-z_1\frac{\partial}{\partial z_1}+(z_2+\xi_3z_3)
\frac{\partial}{\partial z_2}+(\xi_8z_3+\xi_9z_4)
\frac{\partial}{\partial z_3}$$\vs

La $1$-forme $\omega$ annulant ces champs s'\'ecrit dans la carte
affine $z_4=1$
$$(\xi_8z_3+\xi_9)(z_1dz_0-z_0dz_1-z_1^2dz_2)
-z_1(z_0-z_1(z_2+\xi_3z_3z_1))dz_3$$ autrement dit
%quitte \`a diviser z_1^2(\xi_8z_3-\xi_9)
$$d\left(\frac{z_0}{z_1}-z_2\right)-\left(\frac{z_0}{z_1}-z_2-\xi_3z_3\right)
\frac{dz_3}{\xi_8z_3+\xi_9}$$ Posons $u=\frac{z_0}{z_1}-
z_2$ ; alors \`a multiplication pr\`es $\omega$ est du type $$(\xi_8z_3+
\xi_9)du+(u-\xi_3z_3)dz_3$$ Notons, si $\xi_8\not =0$ :
$$\kappa=\frac{\xi_9}{\xi_8}, \hspace{2mm} Z_3=z_3+\kappa
\mbox{ et } U=u-\xi_3\kappa$$ Alors le feuilletage est d\'ecrit par
$$\xi_8Z_3dU+(U-\xi_3Z_3)dZ_3$$ et admet pour int\'egrale
premi\`ere $$Z_3{\left((\xi_8)+1)U-\xi_3Z_3\right)}^{\xi_8}$$ si
$\xi_8\not =-1$ et $$Z_3^{\xi_3}\exp\left(\frac{U}{Z_3}
\right)$$ sinon.\vs

On obtient dans les coordonn\'ees initiales l'int\'egrale
premi\`ere $$\frac{(z_3+\kappa
z_4){\left((\xi_8+1)(z_0z_4-z_1z_2+\xi_3\kappa z_1z_4)-
\xi_3(z_3z_4+\kappa z_4^2)\right)}^{\xi_8}}{z_1^{\xi_8}
z_4^{1+\xi_8}}$$ qui est conjugu\'ee \`a $$
\frac{z_3{\left(z_0z_4-z_1z_2\right)}^{\xi_8}}{z_1^{\xi_8}
z_4^{1+\xi_8}}$$ si $\xi_8\not =1$, et
$$\left(\frac{z_3+\kappa z_4}
{z_4}\right)\exp\left(\frac{z_0z_4-z_1z_2-\xi_3\kappa z_1z_4}
{z_1(z_3+\kappa z_4)}\right)$$ que l'on conjugue \`a $$
\left(\frac{z_3}{z_4}\right)\exp\left(\frac{z_0z_4-z_1z_2}
{z_1z_3}\right)$$ sinon. On remarque que cette derni\`ere
int\'egrale premi\`ere se ram\`ene \`a celle obtenue en 4.1.
pour $\lambda\not =1$. \vs

Lorsque $\xi_8=0$, la $1$-forme $\omega$ s'\'ecrit $$
\xi_9du+(u-\xi_3z_3)dz_3$$ et a pour int\'egrale
premi\`ere $$(u-\xi_3z_3+\xi_3\xi_9)\exp\left(\frac{z_3}
{\xi_9}\right)$$ ce qui conduit \`a
$$\frac{z_0z_4-z_1z_2-\xi_3z_1z_3+\xi_3\xi_9z_1z_4}{z_1
z_4}\exp\left(\frac{z_3}{\xi_9z_4}\right)$$ que l'on transforme
en : $$\frac{z_0z_4-z_1z_2}{z_1
z_4}\exp\left(\frac{z_3}{z_4}\right)$$ Cette int\'egrale premi\`ere est
encore conjugu\'ee \`a celle obtenue en 4.1. pour $\lambda=0$.
\vs

\item[4.5.] Gr\^ace \`a Maple, on obtient pour le cinqui\`eme
cas
$$X=z_0\frac{\partial}{\partial z_0}+z_1\frac{\partial}{\partial
z_1}, \hspace{3mm} Y= z_1 \frac{\partial}{\partial z_0}+z_3
\frac{\partial}{\partial z_2}+z_4\frac{\partial}{\partial z_3}$$
$$Z=\xi_1z_1\frac{\partial}{\partial z_0}-z_1\frac{\partial}
{\partial z_1}+2z_2\frac{\partial}{\partial z_2}+z_3
\frac{\partial}{\partial z_3}$$\vs
La $1$-forme qui annule ces champs s'\'ecrit dans
la carte affine $z_4=1$ $$(2z_2-z_3^2)(z_1dz_0-z_0dz_1-z_1^2dz_3)+z_1(z_0-\xi_1
z_1-z_1z_3)(dz_2-z_3dz_3)$$ soit \`a coefficient
multiplicatif
pr\`es $$2\frac{d\left(\frac{z_0}{z_1}-\xi_1-z_3\right)}
{\frac{z_0}{z_1}-\xi_1-z_3}-\frac{d(2z_2-z_3^2)}{2z_2-z_3^2}$$
On constate que
$$\frac{{\left(z_0z_4-\xi_1z_1z_4-z_1z_3\right)}^2}{z_1^2
\left(2z_2z_4-z_3^2\right)}$$ est une int\'egrale premi\`ere
rationnelle du feuilletage de $\C\P(4)$ associ\'e, ce qui
donne \`a
conjugaison pr\`es $$\frac{{\left(z_0z_4-z_1z_3\right)}^2
}{z_1^2\left(z_2z_4-z_3^2\right)} $$
\vs

\item[4.6.] D'apr\`es Maple l'avant derni\`ere possibilit\'e
est la suivante : $$X=z_0\frac{\partial}{\partial z_0}+\nu z_1
\frac{\partial}{\partial z_1}, \hspace{3mm}
Y=(z_3+\kappa z_4) \frac{\partial}{\partial z_2}+ z_4
\frac{\partial}{\partial z_3}, \hspace{3mm} R \mbox{ et}$$
$$Z=(\xi_6-\xi_1-\nu(\xi_0-\xi_1))z_1\frac{\partial}{\partial
z_1}$$
$$+(2z_2-\kappa z_3+(\xi_4-\xi_9\kappa)
z_4)\frac{\partial}{\partial z_2}+z_3\frac{\partial}{\partial
z_3}$$ o\`u $$\kappa=\xi_9-\xi_3$$\vs

La $1$-forme annulant les champs $X$, $Y$, $Z$ et $R$
s'\'ecrit dans la carte affine $z_4=1$
$$(\beta-2\kappa z_3-z_3^2+2z_2)(\nu z_1dz_0- z_0dz_1)+
\gamma z_0 z_1(dz_2-\kappa-z_3)dz_3$$ o\`u
$$\beta=\xi_4 -\xi_9\kappa \mbox{ et }
\gamma=\xi_6+(\nu-1)\xi_1-\nu \xi_0$$

Cette $1$-forme s'\'ecrit
aussi %apr\`es division par $z_0z_1\beta-2(\xi_9-\xi_3)z_3-z_3^2+2z_2$
$$\nu \frac{dz_0}{z_0}-\frac{dz_1}{z_1}+\frac{\gamma}{2} \frac{d(\beta
+2z_2-2\kappa z_3-z_3^2)}{\beta +2z_2-2\kappa z_3-z_3^2}$$ Le
$\mathscr{L}$-feuilletage d\'ecrit par cette forme admet donc pour
int\'egrale premi\`ere
$$\frac{z_0^{2\nu} {(\beta z_4^2+2z_2z_4-2\kappa z_3z_4-
z_3^2)}^\gamma}{z_4^{2(\gamma+\nu-1)}z_1^2}$$ que l'on conjugue \`a  $$\frac{z_0^{2\nu}{(z_2z_4-z_3^2)}^\gamma}{z_4^{2(\gamma+\nu-1)}z_1^2}$$
\vs

\item[4.7.] Enfin la derni\`ere possibilit\'e est celle o\`u
les champs $X$, $Y$ et $Z$ s'\'ecrivent
$$X=z_0\frac{\partial}{\partial z_0}+ z_4\frac{\partial}{\partial
z_1}, \hspace{3mm} Y=( z_3+\varepsilon
z_4)\frac{\partial}{\partial z_1}+ z_4 \frac{\partial}{\partial
z_2}, \hspace{3mm} R \mbox{ et}$$
$$Z=((\xi_9+\varepsilon)z_2+(\xi_8-\xi_9)z_3+(\xi_9-\xi_0+\xi_6-
\varepsilon \xi_9)z_4)\frac{\partial}{\partial z_1}$$ $$+(z_2+
\xi_3z_3)\frac{\partial}{\partial z_2}+(-z_3+\xi_9z_4)
\frac{\partial}{\partial z_3}$$\vs

La $1$-forme annulant les champs $X$, $Y$, $Z$ et $R$ s'\'ecrit
 dans la carte affine $z_4=1$ $$(z_3-\xi_9)(dz_0-z_0dz_1+z_0(\varepsilon+ z_3)
)dz_2$$ $$+\left(z_0(z_2(z_3-\xi_9)+z_3((z_3 +\varepsilon)
\xi_3-\xi_8+\xi_9)-\xi_{9}+\varepsilon
\xi_9+\xi_0\xi_6)\right)dz_3$$
qui admet comme int\'egrale premi\`ere
$$\frac{z_0(z_3-\xi_{18}z_4)^\zeta}{z_4^{1+\zeta}} \exp
\left(\frac{-z_1z_4+\varepsilon z_2z_4+z_2z_3+ \frac{\xi_3}{2}z_3^2
+\upsilon z_3z_4}{z_4^2}\right)$$ o\`u
$$\zeta=((\xi_3+1)(\xi_9+\varepsilon)-\xi_8)\xi_9-\xi_9
+\xi_0\xi_6$$ et $$\upsilon=((\varepsilon +\xi_9)\xi_3-\xi_8+\xi_9)$$
soit \`a conjugaison pr\`es $$\frac{z_0(z_3-
z_4)^\zeta}{z_4^{1+\zeta}} \exp\left(\frac{z_1z_4+z_2z_3}{z_4^2}
\right) $$
\end{itemize}

\begin{thm} (Maple) Soit $\mathscr{F}$ un
$\mathscr{L}$-feuilletage de degr\'e $3$ sur
$\C\P(4)$. Si $\mathscr{L}$ admet pour
pr\'esentation $\{[X,Y]=0,\hspace{1mm} [X,Z]=0,\hspace{1mm}
[Y,Z]=Y\}$, alors le feuilletage $\mathscr{F}$ poss\`ede \`a
conjugaison pr\`es l'une des int\'egrales premi\`eres suivantes :
\begin{eqnarray}
& & \frac{z_1^\lambda z_3^{\lambda-1}z_4}
{{\left(z_0z_3+z_1z_2\right)}^\lambda}\nonumber \\
& & \frac{z_4}{z_3
}\exp\left(\frac{z_0z_3+z_1z_2}{z_1z_4}\right) \nonumber \\
& & \nonumber \\
& & \frac{z_0z_3+z_1z_2}{z_1z_3}\exp\left(\frac{z_4}{z_3}\right) \nonumber \\
& & \nonumber \\
& & \frac{z_1^{2
\lambda}{\left(z_2z_4-\frac{z_3^2}{2}\right)}^{\kappa}}
{z_4^{2(\lambda+\kappa)}}\exp\left( \frac{z_0}{z_1}\right)\nonumber \\
& & \nonumber \\
& & \left(\frac{z_0}{z_1}\right)^{\kappa} \left(
\frac{z_2z_4-z_3^2}{z_4^2}\right)\nonumber \\
& & \nonumber \\
& & \frac{z_3{\left(z_0z_4-z_1z_2\right)}^\kappa}{z_1^{\kappa}
z_4^{\kappa+1}}\nonumber \\
& & \nonumber \\
& & \nonumber \\
& & \frac{{\left(z_0z_4-z_1z_3\right)}^2
}{z_1^2\left(z_2z_4-z_3^2\right)}\nonumber \\
& & \nonumber \\
& & \frac{z_0^{2\kappa}{(z_2z_4-
z_3^2)}^\lambda}{z_4^{2(\lambda+\kappa-1)}z_1^2}\nonumber \\
& & \nonumber \\
& & \frac{z_0(z_3-z_4)^\kappa}{z_4^{1+\kappa}}
\exp\left(\frac{z_1z_4+z_2z_3 }{z_4^2}\right) \nonumber
\end{eqnarray}
o\`u $\kappa$ et $\lambda$ d\'esignent des nombres complexes.
\vs
\end{thm}

\section[{$\mathscr{L}$ a pour pr\'esentation
 $\{[X,Y]=[X,Z]=0, \hspace{1mm} [Y,Z]=X\}$}.]
{L'alg\`ebre $\mathscr{L}$ a pour pr\'esentation
 $\{[X,Y]=[X,Z]=0, \hspace{1mm} [Y,Z]=X\}$.}

On se propose de v\'erifier le : \vs

\begin{thm} Il n'existe pas de
$\mathscr{L}$-feuilletage de degr\'e $3$ sur $\C
\P(4)$ associ\'e \`a l'alg\`ebre $\{[X,Y]=0,\hspace{2mm}
[X,Z]=0,\hspace{2mm} [Y,Z]=X\}$.
\vs
\end{thm}

La d\'emonstration se fait en deux \'etapes, chacune correspondant \`a l'un
des lemmes qui suivent. \vs

L'\'egalit\'e $[Y,Z]=X$ implique que $X$ est nilpotente ; par
jordanisation $X$ s'\'ecrit
$$\left(%
\begin{array}{ccccc}
  0 & \varepsilon_1 & 0 & 0 & 0 \\
  0 & 0 & \varepsilon_2 & 0 & 0 \\
  0 & 0 & 0 & \varepsilon_3 & 0 \\
  0 & 0 & 0 & 0 & \varepsilon_4 \\
  0 & 0 & 0 & 0 & 0 \\
\end{array}%
\right)$$ o\`u $\varepsilon_i \in \{0,1\}$. Notons $$X_4=\left(%
\begin{array}{ccccc}
  0 & 1 & 0 & 0 & 0 \\
  0 & 0 & 0 & 0 & 0 \\
  0 & 0 & 0 & 1 & 0 \\
  0 & 0 & 0 & 0 & 0 \\
  0 & 0 & 0 & 0 & 0 \\
\end{array}%
\right)$$ On a le : \vs

\begin{lem} Si $X$ n'est pas conjugu\'e \`a $X_4$, alors l'alg\`ebre
$\mathscr{L}$ n'est pas associ\'ee \`a un feuilletage de degr\'e
$3$ sur $\C \P(4)$.\vs
\end{lem}

\begin{proof}[D\'emonstration] Si $X$ est de rang $4$, comme $Y$ et $Z$ commutent avec $X$,
les
matrices $Y$ et $Z$ sont de la forme $$\left(%
\begin{array}{ccccc}
 \nu & \kappa & \beta & \gamma & \delta \\
  0 & \nu & \kappa & \beta & \gamma \\
  0 & 0 & \nu & \kappa & \beta \\
  0 & 0 & 0 & \nu & \kappa \\
  0 & 0 & 0 & 0 & \nu \\
\end{array}%
\right)$$ Mais alors $[Y,Z]=0$, ce qui est exclu. \vs

Supposons que $X$ soit de rang $3$ ; alors $X$ est du type $$X_1=\left(%
\begin{array}{ccccc}
  0 & 1 & 0 & 0 & 0 \\
  0 & 0 & 1 & 0 & 0 \\
  0 & 0 & 0 & 1 & 0 \\
  0 & 0 & 0 & 0 & 0 \\
  0 & 0 & 0 & 0 & 0 \\
\end{array}%
\right) \mbox{ ou } X_2=\left(%
\begin{array}{ccccc}
  0 & 1 & 0 & 0 & 0 \\
  0 & 0 & 1 & 0 & 0 \\
  0 & 0 & 0 & 0 & 0 \\
  0 & 0 & 0 & 0 & 1 \\
  0 & 0 & 0 & 0 & 0 \\
\end{array}%
\right)$$

Commen\c{c}ons par supposer que $X$ est du type $X_1$ ; $Y$ et $Z$
commutent
\`a $X$ donc s'\'ecrivent quitte \`a soustraire un multiple de $X$ $$\left(%
\begin{array}{ccccc}
  \kappa & 0 & \beta & * & * \\
  0 & \kappa & 0 & \beta & 0 \\
  0 & 0 & \kappa & 0 & 0 \\
  0 & 0 & 0 & \kappa & 0 \\
  0 & 0 & 0 & * & * \\
\end{array}%
\right)$$ Mais alors $[Y,Z]$ est de la forme $$\left(%
\begin{array}{ccccc}
  0 & 0 & 0 & * & * \\
  0 & 0 & 0 & 0 & 0 \\
  0 & 0 & 0 & 0 & 0 \\
  0 & 0 & 0 & 0 & 0 \\
  0 & 0 & 0 & * & 0 \\
\end{array}%
\right)$$ et ne peut donc pas \^etre \'egal \`a $X$.

Finalement supposons que $X$ soit du type $X_2$ ; $Y$ et $Z$
commutent \`a $X$ donc s'\'ecrivent respectivement quitte \`a
soustraire un multiple de $X$
$$\left(%
\begin{array}{ccccc}
0 & 0 & \xi_1 & \xi_2 & \xi_3 \\
0 & 0 & 0 & 0 & \xi_2 \\
0 & 0 & 0 & 0 & 0 \\
0 & \xi_4 & \xi_5 & \xi_6 & \xi_7 \\
0& 0 & \xi_4 & 0 & \xi_6 \\
\end{array}%
\right) \mbox{ et }
\left(%
\begin{array}{ccccc}
0 & 0 & \kappa_1 & \kappa_2 & \kappa_3 \\
0 & 0 & 0 & 0 & \kappa_2 \\
0 & 0 & 0 & 0 & 0 \\
0 & \kappa_4 & \kappa_5 & \kappa_6 & \kappa_7 \\
0& 0 & \kappa_4 & 0 & \kappa_6 \\
\end{array}%
\right)$$ Alors la condition $[Y,Z]=X$ ne peut \^etre satisfaite.\vs

Si $X$ est de rang $2$, alors $X$ est de la forme $X_4$ ou $$X_3=\left(%
\begin{array}{ccccc}
  0 & 1 & 0 & 0 & 0 \\
  0 & 0 & 1 & 0 & 0 \\
  0 & 0 & 0 & 0 & 0 \\
  0 & 0 & 0 & 0 & 0 \\
  0 & 0 & 0 & 0 & 0 \\
\end{array}%
\right)$$

Si $X$ est du type $X_3$ ; comme $Y$
et $Z$ commutent
avec $X$, ils s'\'ecrivent quitte \`a soustraire un
multiple de $X$ $$\left(%
\begin{array}{ccccc}
  \kappa & 0 & * & * & * \\
  0 & \kappa & 0 & 0 & 0 \\
  0 & 0 & \kappa & 0 & 0 \\
  0 & 0 & * & * & * \\
  0 & 0 & * & * & * \\
\end{array}%
\right)$$ Mais alors $[Y,Z]$ est de la forme $$\left(%
\begin{array}{ccccc}
  0 & 0 & * & * & * \\
  0 & 0 & 0 & 0 & 0 \\
  0 & 0 & 0 & 0 & 0 \\
  0 & 0 & * & * & * \\
  0 & 0 & * & * & * \\
\end{array}%
\right)$$ et ne peut donc pas \^etre \'egal \`a $X$.
\vs

Si $X$ est de rang $1$, alors le feuilletage associ\'ee \`a
$\mathscr{L}$ est en fait un pull-back.
\end{proof}

Reste donc \`a traiter le cas o\`u $X=X_4$ ; les matrices $X$ et
$Y$ commutent donc sont simultan\'ement jordanisables. Il nous 
reste parmi les configurations rencontr\'ees en $4.$ la suivante :
$$ X=\left(%
\begin{array}{ccccc}
  0 & 1 & 0 & 0 & 0 \\
  0 & 0 & 0 & 0 & 0 \\
  0 & 0 & 0 & 1 & 0 \\
  0 & 0 & 0 & 0 & 0 \\
  0 & 0 & 0 & 0 & 0 \\
\end{array}%
\right)\mbox{ et }  Y=\left(%
\begin{array}{ccccc}
  \nu_0 & \beta & 0 & 0 & 0 \\
  0 & \nu_0 & 0 & 0 & 0 \\
  0 & 0 & 0 & \delta & 0 \\
  0 & 0 & 0 & 0 & 0 \\
  0 & 0 & 0 & 0 & \nu_2 \\
\end{array}%
\right)$$
\vs

\begin{lem} Cette configuration ne g\'en\`ere pas
un $\mathscr{L}$-feuilletage de
degr\'e $3$ sur $\C \P(4)$. \vs
\end{lem}

\begin{proof}[D\'emonstration] A combinaisons lin\'eaires permises pr\`es $Y$ et $Z$
s'\'ecrivent $$\left(%
\begin{array}{ccccc}
\nu_0 & 0 & 0 & 0 & 0 \\
0 & \nu_0 & 0 & 0 & 0 \\
0 & 0 & 0 & \delta & 0 \\
0 & 0 & 0 & 0 & 0 \\
0 & 0 & 0 & 0 & \nu_2 \\
\end{array}%
\right)\mbox{ et }
\left(%
\begin{array}{ccccc}
0 & 0 & \zeta_2 & \zeta_3 & \zeta_4 \\
\zeta_5 & \zeta_6 & \zeta_7 & \zeta_8 & \zeta_9 \\
\zeta_{10} & \zeta_{11} & \zeta_{12} & \zeta_{13} & \zeta_{14} \\
\zeta_{15} & \zeta_{16} & \zeta_{17} & \zeta_{18} & \zeta_{19} \\
\zeta_{20} & \zeta_{21} & \zeta_{22} & \zeta_{23} & \zeta_{24} \\
\end{array}%
\right) $$

L'\'egalit\'e $[Y,Z]-X=0$ s'\'ecrit alors $$
\left(%
\begin{array}{ccccc}
0 & -1 & -\nu_0\zeta_2 & \zeta_2\delta-\nu_0\zeta_3 &
\zeta_4(\nu_2-\nu_0)
\\
0 & 0 & -\nu_0\zeta_7 & \zeta_7\delta-\nu_0\zeta_8 & \zeta_9
(\nu_2-\nu_0) \\
\zeta_{10}\nu_0-\delta\zeta_{15} & \zeta_{11}\nu_0-\delta
\zeta_{16} & -\zeta_{17}\delta & \delta(\zeta_{12}-
\zeta_{18})-1 & \zeta_{14}\nu_2-\delta\zeta_{19} \\
\zeta_{15}\nu_0 & \zeta_{16}\nu_0 & 0 & \zeta_{17}\delta & \zeta_{19}\nu_2 \\
\zeta_{20}(\nu_0-\nu_2) & \zeta_{21}(\nu_0-\nu_2) &
-\nu_2\zeta_{22} &
\zeta_{22}\delta-\nu_2\zeta_{23} & 0\\
\end{array}%
\right)$$ ce qui est absurde (premi\`ere ligne, 2\`eme colonne).
\end{proof}

\section[{$\mathscr{L}$ a pour pr\'esentation 
$\{[X,Y]=Y, \hspace{1mm} [X,Z]=Y+Z, \hspace{1mm}
[Y,Z]=0\}$.}]{L'alg\`ebre $\mathscr{L}$ a pour pr\'esentation 
$\{[X,Y]=Y, \hspace{1mm} [X,Z]=Y+Z, \hspace{1mm}
[Y,Z]=0\}$.}

De m\^eme on a le :

\begin{thm} Il n'existe pas de
$\mathscr{L}$-feuilletage de degr\'e $3$ sur $\C
\P(4)$ associ\'e \`a l'alg\`ebre d\'ecrite par $\{[X,Y]=Y,\hspace{1mm}
[X,Z]=Y+Z,\hspace{1mm} [Y,Z]=0\}$.
\vs
\end{thm}

La d\'emonstration se fait en deux \'etapes, chacune correspondant \`a l'un
des lemmes qui suivent. \vs

La r\'esolubilit\'e de $\mathscr{L}$ et l'\'egalit\'e $[X,Y]=Y$
implique que $Y$ est nilpotente, alors par jordanisation $Y$
s'\'ecrit
$$\left(%
\begin{array}{ccccc}
  0 & \varepsilon_1 & 0 & 0 & 0 \\
  0 & 0 & \varepsilon_2 & 0 & 0 \\
  0 & 0 & 0 & \varepsilon_3 & 0 \\
  0 & 0 & 0 & 0 & \varepsilon_4 \\
  0 & 0 & 0 & 0 & 0 \\
\end{array}%
\right)$$ o\`u $\varepsilon_i \in \{0,1\}$. Notons $$Y_4=\left(%
\begin{array}{ccccc}
  0 & 1 & 0 & 0 & 0 \\
  0 & 0 & 0 & 0 & 0 \\
  0 & 0 & 0 & 1 & 0 \\
  0 & 0 & 0 & 0 & 0 \\
  0 & 0 & 0 & 0 & 0 \\
\end{array}%
\right)$$ on a \vs

\begin{lem} Si $Y$ n'est pas conjugu\'e \`a $Y_4$, alors l'alg\`ebre
$\mathscr{L}$ n'est pas associ\'ee \`a un feuilletage de degr\'e
$3$ sur $\C \P(4)$.\vs
\end{lem}

\begin{proof}[D\'emonstration] Supposons que $Y$ soit de rang $4$. Comme $Z$ commute \`a
$Y$, la matrice $Z$ est quitte \`a soustraire $Y$ de la forme
$$\left(%
\begin{array}{ccccc}
  0 & 0 & \beta & \gamma & \delta \\
  0 & 0 & 0 & \beta & \gamma \\
  0 & 0 & 0 & 0 & \beta \\
  0 & 0 & 0 & 0 & 0 \\
  0 & 0 & 0 & 0 & 0 \\
\end{array}%
\right)$$ et l'\'egalit\'e $[X,Y]=Y$ implique que $$X=\left(%
\begin{array}{ccccc}
  \kappa' & \beta' & \gamma' & \delta' & \varepsilon' \\
  0 & \kappa'+1 & \beta' & \gamma' & \delta' \\
  0 & 0 & \kappa'+2 & \beta' & \gamma' \\
  0 & 0 & 0 & \kappa'+3 & \beta' \\
  0 & 0 & 0 & 0 & \kappa'+4 \\
\end{array}%
\right)$$ Ainsi $[X,Z]-Z$ s'\'ecrit $$\left(%
\begin{array}{ccccc}
  0 & 0 & * & * & * \\
  0 & 0 & 0 & * & * \\
  0 & 0 & 0 & 0 & * \\
  0 & 0 & 0 & 0 & 0 \\
  0 & 0 & 0 & 0 & 0 \\
\end{array}%
\right)$$ et ne peut dont pas \^etre \'egal \`a $Y$. \vs

Supposons que $Y$ soit de rang $3$, alors $Y$ est du type $$Y_1=\left(%
\begin{array}{ccccc}
  0 & 1 & 0 & 0 & 0 \\
  0 & 0 & 1 & 0 & 0 \\
  0 & 0 & 0 & 1 & 0 \\
  0 & 0 & 0 & 0 & 0 \\
  0 & 0 & 0 & 0 & 0 \\
\end{array}%
\right) \mbox{ ou } Y_2=\left(%
\begin{array}{ccccc}
  0 & 1 & 0 & 0 & 0 \\
  0 & 0 & 1 & 0 & 0 \\
  0 & 0 & 0 & 0 & 0 \\
  0 & 0 & 0 & 0 & 1 \\
  0 & 0 & 0 & 0 & 0 \\
\end{array}%
\right)$$

Si $Y$ est du type $Y_1$ ;
comme $[Z,Y]=0$, on obtient $$Z=\left(%
\begin{array}{ccccc}
  \kappa & \beta & \gamma & * & *\\
  0 & \kappa & \beta & \gamma & 0 \\
  0 & 0 & \kappa & \beta & 0 \\
  0 & 0 & 0 & \kappa & 0 \\
  0 & 0 & 0 & * & * \\
\end{array}%
\right)$$ L'\'egalit\'e $[X,Y]=Y$ implique que $X$ s'\'ecrit $$\left(%
\begin{array}{ccccc}
  \kappa' & \beta' & \gamma' & * & * \\
  0 & \kappa'+1 & \beta' & \gamma' & 0 \\
  0 & 0 & \kappa'+2 & \beta' & 0 \\
  0 & 0 & 0 & \kappa'+3 & 0 \\
  0 & 0 & 0 & * & * \\
\end{array}%
\right)$$ Alors $[X,Z]-Z$ est de la forme $$\left(%
\begin{array}{ccccc}
  -\kappa & 0 & \gamma & * & * \\
  0 & -\kappa & 0 & \gamma & 0 \\
  0 & 0 & -\kappa & 0 & 0 \\
  0 & 0 & 0 & -\kappa  & 0 \\
  0 & 0 & 0 & * & * \\
\end{array}%
\right)$$ donc diff\'erent de $Y$.

Si $Y$ est du type $Y_2$, alors
comme $[Z,Y]=0$ et $Z$ est nilpotent, on obtient $$Z=\left(%
\begin{array}{ccccc}
0 & 0 & \xi_0 & \beta & \xi_1\\
0 & 0 & 0 & 0 & \beta \\
0 & 0 & 0 & 0 & 0 \\
0 & \kappa & \xi_2 & 0 & \xi_3\\
0 & 0 & \kappa & 0 & 0 \\
\end{array}%
\right)$$\vs

L'\'egalit\'e $[X,Y]=Y$ implique que $X$ s'\'ecrit modulo $R$
$$\left(%
\begin{array}{ccccc}
0 & b & c & \delta & e\\
0 & 1 & b & 0 & \delta\\
0 & 0 & 2 & 0 & 0 \\
0 & f & g & a & h \\
0 & 0 & f & 0 & a+1 \\
\end{array}%
\right) $$
Alors $[X,Z]-Z-Y$ s'\'ecrit $$\left(%
\begin{array}{ccccc}
0 & \beta f-\kappa\delta-1 & \xi_0+\beta g+\xi_1f-\delta\xi_2 -e
\kappa & \beta (a-1) & \beta (h-b)+\xi_1a-\delta\xi_3 \\
0 & 0 & \beta f-\kappa\delta-1 & 0 & \beta(a-1)\\
0 & 0 & 0 & 0 & 0 \\
0 & -a\kappa & \kappa (b-h)+\xi_2(1-a)+\xi_3f & 0 &
\kappa\delta-\beta f-1\\
0 & 0 & -a\kappa & 0 & 0 \\
\end{array}%
\right)$$

Cette matrice doit \^etre identiquement nulle. En particulier
on a $$\beta f-\kappa\delta-1=\kappa\delta-\beta f-1=0$$ ce qui
est absurde.
\vs

Supposons maintenant que $Y$ soit de rang $2$ ; alors $Y$ est du
type $Y_4$ ou : $$Y_3=\left(%
\begin{array}{ccccc}
  0 & 1 & 0 & 0 & 0 \\
  0 & 0 & 1 & 0 & 0 \\
  0 & 0 & 0 & 0 & 0 \\
  0 & 0 & 0 & 0 & 0 \\
  0 & 0 & 0 & 0 & 0 \\
\end{array}%
\right)$$

Si $Y$ est de la forme $Y_3$, alors
l'\'egalit\'e $[Y,Z]=0$ conduit
\`a $$Z=\left(%
\begin{array}{ccccc}
  \kappa & \beta & * & * & * \\
  0 & \kappa & \beta & 0 & 0 \\
  0 & 0 & \kappa & 0 & 0 \\
  0 & 0 & * & * & * \\
  0 & 0 & * & * & * \\
\end{array}%
\right)$$ Ensuite l'\'egalit\'e $[X,Y]=Y$ implique
que $X$ est de la forme $$\left(%
\begin{array}{ccccc}
  \kappa' & \beta' & * & * & * \\
  0 & \kappa'+1 & \beta' & 0 & 0 \\
  0 & 0 & \kappa'+2 & 0 & 0 \\
  0 & 0 & * & * & * \\
  0 & 0 & * & * & * \\
\end{array}%
\right)$$ En r\'esulte que $[X,Z]-Z$ s'\'ecrit $$\left(%
\begin{array}{ccccc}
  -\kappa & 0 & * & * & * \\
  0 & -\kappa & 0 & 0 & 0 \\
  0 & 0 & -\kappa & 0 & 0 \\
  0 & 0 & * & * & * \\
  0 & 0 & * & * & * \\
\end{array}%
\right)$$ et ne peut donc \^etre \'egal \`a $Y$.
\vs

Si $X$ est de rang $1$, alors le feuilletage associ\'ee \`a
$\mathscr{L}$ est en fait un pull-back.
\end{proof}

Reste donc \`a traiter le cas o\`u $Y=Y_4$ ; les matrices $Y$ et
$Z$ commutent, d'o\`u l'existence d'une base dans laquelle on peut
jordaniser $Y$ et $Z$. Par ailleurs la r\'esolubilit\'e de
$\mathscr{L}$ et l'\'egalit\'e $[X,Z]=Y+Z$ implique que $Z$ est
nilpotente ; donc parmi les configurations du 4. il nous reste
seulement la configuration
 $$Y=\left(%
\begin{array}{ccccc}
  0 & 1 & 0 & 0 & 0 \\
  0 & 0 & 0 & 0 & 0 \\
  0 & 0 & 0 & 1 & 0 \\
  0 & 0 & 0 & 0 & 0 \\
  0 & 0 & 0 & 0 & 0 \\
\end{array}%
\right) \mbox{ et } Z=\left(%
\begin{array}{ccccc}
  0 & \beta & 0 & 0 & 0 \\
  0 & 0 & 0 & 0 & 0 \\
  0 & 0 & 0 & \delta & 0 \\
  0 & 0 & 0 & 0 & 0 \\
  0 & 0 & 0 & 0 & 0 \\
\end{array}%
\right)$$
\vs

\begin{lem} Cette configuration ne g\'en\`ere pas
un $\mathscr{L}$-feuilletage de
degr\'e $3$ sur $\C \P(4)$. \vs
\end{lem}

\begin{proof}[D\'emonstration] L'\'egalit\'e $[X,Y]=Y$ implique que $X$ s'\'ecrit \`a
combinaisons lin\'eaires permises pr\`es
$$\left(%
\begin{array}{ccccc}
0 & 0 & \zeta_1 & \zeta_2 & \zeta_3 \\
0 & 1 & 0 & \zeta_4 & 0 \\
\zeta_5 & \zeta_6 & \zeta_7 & \zeta_8 & \zeta_9 \\
0 & \zeta_{10} & 0 & \zeta_7+1 & 0 \\
0 & \zeta_{11} & 0 & \zeta_{12} & \zeta_{13}\\
\end{array}%
\right)$$ Mais alors on a $$[X,Z]-Y-Z=\left(%
\begin{array}{ccccc}
0 & -1 & 0 & \beta \zeta_4-\zeta_1\delta & 0 \\
0 & 0 & 0 & 0 & 0\\
0 & \delta \zeta_{10}-\zeta_5\beta & 0 & -1 & 0 \\
0 & 0 & 0 & 0 & 0 \\
0 & 0 & 0 & 0 & 0 \\
\end{array}%
\right)$$ d'o\`u $-1=0$ ce qui est exclu.
\end{proof}

%%%%%%%%%%%%%%%%%%%%%%%%%%%%%%%%%%%%%%%%%%%%%%%%%%%%%%%%%%%%%%%%%%%%%%%%%%
%%%%%%%%%%%%%%%%%%%%%%%%%%%%%%%%%%%%%%%%%%%%%%%%%%%%%%%%%%%%%%%%%%%%%%%%%%

\chapter[Codimension quelconque.]{$\mathscr{L}$-feuilletages
de degr\'e $0$ et de codimension quelconque sur $\C\P(n)$.}

\vs

Par d\'efinition, un $\mathscr{L}$-feuilletage de degr\'e $0$ et de codimension
$p$ sur $\C\P(n)$ est associ\'e au champ de $p$-plans d\'efini
par une $p$-forme homog\`ene $\Omega$ de degr\'e $1$ sur $\C^{n+1}$
annulant le champ radial $R$. Plus pr\'ecis\'ement une telle $1$-forme s'\'ecrit
$$\Omega=\displaystyle \sum_{I=(i_1,\ldots,i_p)} L_I dz_{i_1}\wedge \ldots
\wedge dz_{i_p}$$ o\`u $L_I$ d\'esigne une forme lin\'eaire. On demande que
$i_R\Omega \not =0$. L'int\'egrabilit\'e se traduit comme suit :
en tout point $m \in \C^{n+1}$ g\'en\'erique, i.e. $\Omega(m) \not =0$,
la $1$-forme $\Omega$ s'\'ecrit localement $$\Omega=\omega_1\wedge\ldots\wedge
\omega_p$$ o\`u les $\omega_i$ sont des $1$-formes locales telles que le
syst\`eme de Pfaff engendr\'e par les $\omega_i$ soit int\'egrable : $$\omega_1
\wedge\ldots\wedge\omega_p\wedge d\omega_i=0 \hspace{6mm} \forall \hspace{1mm}
i \in\{1,\ldots,p\}$$

\vs

\begin{thm} Un feuilletage $\mathscr{F}$ de
codimension $p$ et de degr\'e $0$ sur $\C\P(n)$
est un $\mathscr{L}$-feuilletage.
\end{thm}

\vs

\begin{proof}[D\'emonstration]
Elle r\'esulte de la description des formes $\Omega$ qui suit.

Si $p=n-1$, la condition $i_R\Omega=0$ implique l'existence d'un champ $X_0$ tel
que $$\Omega=i_{X_0}i_R dz_0\wedge\ldots\wedge dz_n$$ Comme $\deg \Omega=1$, on
peut prendre pour $X_0$ un champ constant, par exemple
$\frac{\partial}{\partial z_0}$. Visiblement, $\mathscr{F}$ est un
$\mathscr{L}$-feuilletage et l'alg\`ebre $\tilde{\mathscr{L}}$ est engendr\'ee
par les champs $$z_j\frac{\partial}{\partial z_0}, \hspace{1mm} j=0, \ldots,
n+1 \hspace{1mm} \mbox{ et } R$$ Les $\frac{z_i}{z_j}$, $1\leq i\leq j$, sont
des int\'egrales premi\`eres de $\mathscr{F}$. Dans la carte affine $z_0=1$, le
feuilletage est associ\'e au champ radial de $\C^n$ : $\displaystyle
\sum_{i=1}^n z_i\frac{\partial}{\partial z_i}$.

On suppose que $p<n-1$. Soit $m \in \C^{n+1}$ un point g\'en\'erique, la
$1$-forme $\Omega$ s'\'ecrit $$\omega_0\wedge \ldots\wedge\omega_{p-1}$$ o\`u
$\omega_i$ d\'esigne une $1$-forme locale en $m$. Le th\'eor\`eme de Frobenius
classique assure l'existence de submersions $f_1$, $\ldots$, $f_p$
ind\'ependantes telles que $$\Omega\wedge df_i=0$$ On peut supposer que $f_i=z_i
+\mbox{ termes de degr\'e sup\'erieur}$ pour tout $i=0,\ldots,p-1$. On a
l'\'egalit\'e $d\Omega\wedge df_i=0$ qui, comme la $1$-forme $d\Omega$ est \`a
coefficients constants, conduit pour $i=0,\ldots,p-1$ \`a 
$$d\Omega\wedge dz_i=0$$ En r\'esulte que $$\Omega=dz_0\wedge \ldots\wedge
dz_{p-1}\wedge \eta$$ o\`u $\eta$ est une $1$-forme \`a coefficients constants :
il existe une forme lin\'eaire $L$ tell que $\eta=dL$. L'identit\'e d'Euler 
$$L_R\Omega=
i_Rd\Omega=(p+1)\Omega$$ assure que $$dz_0\wedge\ldots\wedge dz_{p-1}\wedge dL \not
= 0$$ On peut donc supposer que $L=z_p$ ; on obtient alors
\begin{eqnarray}
\Omega&=&i_R(dz_0\wedge\ldots\wedge dz_p) \nonumber \\
&=& z_0dz_1\wedge\ldots\wedge dz_p-z_1dz_0\wedge dz_2\wedge\ldots\wedge dz_p+
\ldots\nonumber
\end{eqnarray}
On remarque que les champs de vecteurs $$z_j\frac{\partial}{\partial z_\ell},
\hspace{1mm} \ell \geq p+1 \mbox{ et } R_{ij}=z_i\frac{\partial}{\partial z_i} +
z_j\frac{\partial}{\partial z_j}, \hspace{1mm} i< j \leq p$$ annulent
$\Omega$. Ceci d\'emontre que $\mathscr{F}$ est un $\mathscr{L}$-feuilletage ;
il poss\`ede comme int\'egrales premi\`eres $$\frac{z_i}{z_j} \hspace{6mm} i
< j \leq p$$ Les feuilles r\'eguli\`eres sont d'adh\'erence des $p$-plans.
\end{proof}

Rappelons que lors de l'\'etude de l'exemple de Gordan-N\oe ther nous avons
rencontr\'e des $\mathscr{L}$-feuilletages de codimensions diff\'erentes de $1$
; les feuilletages de degr\'e $1$ sur $\C \P(n)$ de codimension
$p$, avec $1<p<n-1$, ne sont actuellement pas classifi\'es. Voici quelques
exemples de feuilletages de codimension $2$ obtenus par intersection de
feuilletages de codimension $1$ par exemple dans $\C\P(4)$. On
g\'en\'eralisera facilement.

Comme on l'a vu l'intersection de deux pinceaux d'hyperplans $\frac{z_0}{z_1}=$
cte et $\frac{z_0}{z_2}=$ cte d\'efinit un feuilletage de degr\'e $0$ sur
$\C\P(4)$. Ce feuilletage peut \^etre vu comme
d\'eg\'en\'erescence de l'intersection de deux pinceaux d'hyperplans
g\'en\'eriques $\frac{z_0}{z_1}=$ cte et $\frac{z_2}{z_3}=$ cte qui lui aussi
est un $\mathcal{L}$-feuilletage de degr\'e $1$ cette fois.

Consid\'erons maintenant un feuilletage logarithmique
$\mathcal{F}_\omega$, de degr\'e $1$, donn\'e par
$$\omega_1=z_0z_1z_2\displaystyle \sum_{i=0}^2 \lambda_i
\frac{dz_i}{z_i}$$ o\`u $\lambda_0+\lambda_1+\lambda_2=0$. Nous
allons l'intersecter avec un pinceau d'hyperplans ; le degr\'e du
feuilletage associ\'e sera $1$ ou $2$ suivant que la position du
pinceau d'hyperplans est g\'en\'erique ou non par rapport \`a
$\omega_1$. Par exemple le feuilletage de codimension $2$ obtenu
par intersection de $\mathscr{F}_{\omega_1}$ et de
$\frac{z_3}{z_4}=$ cte d\'efinit un feuilletage de degr\'e $2$ sur
$\C\P(4)$. Il poss\`ede les deux int\'egrales
premi\`eres $z_0^{\lambda_0}z_1^{\lambda_1}z_2^{\lambda_2}$ et
$\frac{z_2}{z_3}$. On peut faire d\'eg\'en\'erer ce pinceau sur le
suivant $\frac{z_0}{z_3}=$ cte et obtenir ainsi un feuilletage de
degr\'e $1$. Ces deux feuilletages sont des
$\mathscr{L}$-feuilletages. Si l'on consid\`ere un autre type de
d\'eg\'en\'erescence du pinceau, par exemple le pinceau
$\frac{z_0+z_1}{z_3}=$ cte, on obtient un feuilletage de degr\'e
$2$ et de codimension $2$ qui n'est pas un
$\mathscr{L}$-feuilletage.

\bibliographystyle{plain}
\bibliography{biblio}
\nocite{*}

\end{document}